\documentclass[upref]{amsart}
\usepackage{latexsym,amsmath,amstext,mathrsfs,amssymb,amsthm,amscd,amsopn,verbatim,graphics}

\usepackage{times}
\numberwithin{equation}{section}

\DeclareMathOperator{\id}{id}

\DeclareMathOperator{\conj}{c}

\DeclareMathOperator{\tange}{T}

\DeclareMathOperator{\Eta}{H}
\DeclareMathOperator{\Mu}{M}

\DeclareMathOperator{\coh}{H}

\DeclareMathOperator{\pr}{pr}

\newcommand{\unint}{\mathrm{I}}

\newtheorem{Thm}{Theorem}[section]
\newtheorem{Prop}[Thm]{Proposition}
\newtheorem{Lem}[Thm]{Lemma}
\newtheorem{Cor}[Thm]{Corollary}

\theoremstyle{remark}
\newtheorem{Rem}[Thm]{Remark}
\newtheorem*{Ack}{Acknowledgment}

\theoremstyle{definition}

\newtheorem{Def}[Thm]{Definition}
\newtheorem{Exa}[Thm]{Example}

\newtheorem*{Reg*}{Regularization procedure}

\newcommand{\bbR}{{\mathbb{R}}}

\newcommand{\calH}{\mathcal{H}}
\newcommand{\calS}{\mathcal{S}}

\newcommand{\calG}{\mathcal{G}}
\newcommand{\calU}{\mathcal{U}}

\newcommand{\calE}{\mathcal{E}}
\newcommand{\calF}{\mathcal{F}}

\newcommand{\tildep}{\widetilde{p}}

\newcommand\qq{\rm}

\newcommand\cme[1]{{\qq Commun.\ Math.\ Helv.\ \bf #1}}

\newcommand\Kth[1]{{\qq $K$-Theory \bf #1}}

\begin{document} 

\title[Principal bundles with groupoid\ldots]{Principal bundles with groupoid structure: local vs.\ global theory and nonabelian {\v C}ech cohomology}

\author[C.~A.~Rossi]{Carlo~A.~Rossi}
\address{Dept.\ of Mathematics---Technion---32000, Haifa---Israel}  
\email{crossi@techunix.technion.ac.il}

\begin{abstract}
The aim of this paper is to review and discuss in detail local aspects of principal bundles with groupoid structure.
Many results, in particular from the second and third section, are already known to some extents, but, due to the lack of a ``unified'' point of view on the subject, I decided nonetheless to (re)define all the main concepts and write all proofs; however, some results are reformulated in a more elegant way, using the division map and the generalized conjugation of a Lie groupoid.
In the same framework, I discuss later generalized groupoids and Morita equivalences from a local point of view; in particular, I found a (so far as I know) new characterization of generalized morphisms coming from nonabelian {\v C}ech cohomology, which allows one to view generalized morphisms as a generalization of classical descent data.
I found also a factorization formula for the division map, which is the crucial point in the local formulation of Morita equivalences.
\end{abstract}

\maketitle

\tableofcontents

\section{Introduction}\label{sec-intro}
Principal bundles are one of the main objects of study in many areas of mathematics and physics, e.g.\ differential geometry, algebraic topology, Topological Quantum Field Theories, gauge theory.
In particular, the study of principal bundles encodes also the study of many structures on them, like e.g.\ connections, basic differential forms, curvature of connections, etc.
In a previous paper~\cite{CR}, for the purpose of studying the properties of the so-called generalized Wilson loop observable for $BF$-theories (roughly speaking, higher-dimensional analoga of $3$-dimensional Chern--Simons theory), which is a formal series of differential forms on the space of loops on a manifold $M$ which mimics the shape of the path-ordered exponential giving an explicit representative of trace of the holonomy in some representation of the structure group $G$ of a bundle $P$ over $M$, I pursued in detail the idea of viewing holonomy w.r.t.\ a connection as a gauge transformation of some bundle over the space of loops in $M$; more generally, I introduced the concept of generalized gauge transformations and interpreted the parallel transport w.r.t.\ a connection as a generalized gauge transformation between two particular bundles one the space of curves in $M$.
Moreover, the flatness of the connection w.r.t.\ which one considers parallel transport has some consequences: namely, the parallel transport is a horizontal section w.r.t.\ a covariant derivative coming from pull-back of the previous connection.
From a more geometrical point of view, this is equivalent to the well-known fact that the holonomy map w.r.t.\ a flat connection, restricted to loops with a fixed base point, factors to a map from the fundamental group of $M$ w.r.t.\ the chosen base point to $G$.

Later on, I took notice that generalized gauge transformations may be introduced also in the more general framework of principal bundles with structure groupoid: roughly speaking, they can be viewed as manifolds, projecting down to base spaces, acted on from the right (or from the left) by a Lie groupoid, so that the action is free and transitive on each fiber.
This result is an easy consequence of the existence of a division map also for principal bundles with structure groupoid; I refer to~\cite{CR1} for more details on the subject.

At that point, I thought it would be interesting to pursue analoga of the interpretation of parallel transport w.r.t.\ connections also for principal bundles with structure groupoid.
Although there is a huge amount of literature about principal bundles with structure groupoid, the only explicit reference to connections on them I found in the forthcoming book of M{\oe}rdijk and Mr{\v c}un~\cite{Moer2} and in~\cite{Moer4}, again by M{\oe}rdijk and Mr{\v c}un; let me just point out that they discuss connections w.r.t.\ a foliated structure of the base space of the bundles they consider.
They also discuss briefly the notion of flat connections w.r.t.\ a foliated structure.
This was one of my starting point towards an attempt to find a convenient notion of connections and flatness of connections on principal bundles with structure groupoid.

On the other hand, I began also to pursue properties of principal bundles with groupoid structure without any reference to other structures; this I just did for a better understanding of the subject.
In the mathematical literature, I found a huge amount of work on this field, e.g.~\cite{Con},~\cite{H} and~\cite{H1},~\cite{HS},~\cite{Moer1} and~\cite{Mrcun}, to cite the main ones.
However, each one of the cited authors had his personal point of view about the way of dealing with principal bundles with groupoid structure; e.g.\ Connes~\cite{Con}, Hilsum and Skandalis~\cite{HS} and H{\"a}fliger~\cite{H1} prefer to think in local terms, introducing the analogon for groupoids of the nonabelian {\v C}ech first cohomology group of a manifold $M$ (the base of the bundle) with values in the structure groupoid $\calG$, generalizing ideas of Grothendieck, and viewing thus isomorphism classes of principal bundles with structure groupoid as cohomology classes in this framework.
To be more precise, in~\cite{H}, the author introduced nonabelian {\v C}ech cohomology (at degree $1$) of a topological space with values in sheaf of (topological) groupoids; the first nonabelian {\v C}ech cohomology group of $M$ with values in a sheaf of topological groupoids canonically associated to a topological groupoid is in one-to-one correspondence with the set of isomorphism classes of topological principal bundles with groupoid structure over $M$.
On the other hand, M{\oe}rdijk~\cite{Moer1}, Mr{\v c}un~\cite{Mrcun} and M{\oe}rdijk and Mr{\v c}un~\cite{Moer2} prefer to view principal bundles as global objects, pointing out to nonabelian cohomology theory from a slightly different perspective, using the division map.

Principal bundles with groupoid structure are interesting mathematical objects by themselves, since they encode a differential-geometric analogon of the algebraic notion of {\em bimodules} between rings or algebras, namely generalized morphisms between Lie groupoids.
Generalized morphisms between two Lie groupoids $\calG$ and $\calH$ are bibundles, i.e.\ (roughly speaking) principal bundles with structure groupoid $\calH$ over the manifold of objects of $\calG$, such that $\calG$ acts from the left on the bundle in a compatible way. 
Generalized morphisms play a central r{\^o}le in many areas of modern research: a well-known fact (which I will reprove later) is that $G$-equivariant principal $H$-bundles over a manifold $M$, acted on from the left by $G$, are in one-to-one correspondence with generalized morphisms from the action groupoid $G\ltimes M$ to $H$, viewed as a trivial Lie groupoid.
Connes~\cite{Con} and Hilsum and Skandalis~\cite{HS} present a nonabelian {\v C}ech cohomological version of generalized morphisms.
Last, but not least, let me spend two words on the notion of Morita equivalences between Lie groupoids: this notion mimics the notion of Morita equivalences in the category of bimodules.
In fact, a Morita equivalence between Lie groupoids is an isomorphism in the category of generalized morphisms; Morita equivalences between Lie groupoids are central objects of study in the theory of Lie groupoids.

Motivated by the results of~\cite{CR1}, I tried to pursue a ``universal'' point of view.
The key tool, due to my previous work~\cite{CR} in ordinary principal bundles, is the use of the mathematical object that MacKenzie~\cite{McK} calls ``division map'': in~\cite{CR1} I analyzed carefully the properties of the division map, and I will use it in this paper to pursue all local properties of principal bundles with groupoid structure. 
The division map of MacKenzie, in the framework of ordinary principal bundles, is, in the more general context of principal bundles with structure groupoid, what M{\oe}rdijk~\cite{Moer1} calls a ``cocycle over $M$ with values in $\calG$'' (or also a division map); a cocycle with values in $\calG$ depends on a surjective submersion from the total space $P$ to the base manifold $M$, and always comes in pair with a map from the base manifold $M$ to the manifold of objects of the structure groupoid $\calG$, which I call the {\em momentum} of the principal bundle. 
The main feature of a cocycle over $M$ with values in $\calG$ is that it contains all the informations one needs to characterize an action of $\calG$ on $P$ to be free and transitive on each fiber of the surjective submersion.
It is also what I use to link the global point of view to the local one, hence drawing a bridge between the two ``cohomology'' philosophies of~\cite{Con},~\cite{H} and~\cite{H1},~\cite{HS} and\cite{Moer1},~\cite{Mrcun}.
  
Let me now explain the structure of the paper.
In Section~\ref{sec-divmap}, I review the definition of principal bundles with groupoid structure in a ``global'' way, following closely~\cite{Moer2}.
I also review the definition and the main properties of the division map, for whose detailed exposition I refer to~\cite{Moer3} and~\cite{CR1}.

In Section~\ref{sec-localdata}, I review and discuss in detail so-called {\em local trivializing data} over a smooth manifold with values in a Lie groupoid (see also~\cite{Con} and~\cite{H}): I show that there is a one-to-one correspondence between local trivializing data and principal bundles with groupoid structure, once open coverings of the base space are chosen. 
Afterwards, I construct examples of principal bundles with groupoid structure for some known groupoids, using local trivializing data: among other things, I classify completely principal bundles with action groupoids as structure groupoids.
Moreover, I review the definition and the local characterization of morphisms between principal bundles with structure groupoid (always by choosing open coverings of the base space).
To get rid of the choice of open coverings, one tries to reformulate the theory of local trivializing data and local morphisms between them in the framework of nonabelian {\v C}ech cohomology: this is what I do at the end in a more elegant way than in the classical references, using the generalized conjugation of Lie groupoids.
The computations that I do in the framework of nonabelian {\v C}ech cohomology are crucial for what comes in the next Section.

In Section~\ref{sec-hilskan}, I first review the notion of {\em generalized morphism}, following closely~\cite{Mrcun} and~\cite{Moer2}, using a ``global'' point of view, focusing in particular on the additional properties the division map of generalized morphisms has to satisfy.
Later on, I rewrite in local terms, using the arguments of Section~\ref{sec-localdata}, the notion of generalized morphism, introducing the notion of {\em local generalized morphism}, and I then show that local generalized morphisms and generalized morphisms are in bijective correspondence.
Using the notion of local generalized morphism, I classify completely generalized morphisms between action groupoids.
A local construction of generalized morphisms was already pursued in~\cite{Con} and~\cite{HS}, but the local point of view that I take differs slightly from their, in the sense that I consider not only what in~\cite{HS} is called a ``cocycle on a groupoid $\calG$ with values in a groupoid $\calH$'', but with an additional equation, which is related to nonabelian {\v C}ech cohomology in a (so far as I know) new way.
This additional equation is the key point in the characterization of generalized morphisms that I adopted: to mention one fact, in the classical references, the components of a local generalized morphism are labelled as cocycles.
They are, more precisely, not cocycles, but coboundaries in nonabelian {\v C}ech cohomology between two cocycles obtained by two distinct pull-back procedures.
At a global level, this means that there is a morphism between principal bundles, which satisfies also an additional ``cocycle'' condition similar in spirit to the one appearing in classical Descent Theory.
Classical Descent Theory, to be precise, is proved to be a generalized morphism in this setting; therefore, the theory of generalized morphisms may be viewed as a generalization of classical Descent Theory.
Using this characterization of generalized morphisms, I came, by what I should call ``Serendipity'', to the result that I wanted initially to pursue in~\cite{CR}, namely a new characterization of flat bundles, motivated by the fact that flat bundles give rise to horizontal gauge transformations: in fact, isomorphism classes of flat $G$-bundles over $M$ correspond uniquely to equivalence classes of generalized morphisms from the fundamental groupoid of $M$ to $G$.
This I will not pursue here in detail, deserving to it a forthcoming paper; in fact, the announced result could be also generalized to principal bundles with structure groupoid, giving a purely algebraic definition of flat connections, which generalizes the well-known correspondence between flat bundles and conjugacy classes of of representations of the fundamental group in the structure group of the bundles.
Moreover, recently, a complete analysis of connections on so-called principal $G$-bundles over the Lie groupoid $\Gamma$ was pursued in~\cite{L-GTX}, using methods from the theory of simplicial manifolds.
Since connections on principal bundles with structure groupoid may be viewed as a generalized morphism from the quasi-groupoid of curves in $M$ to the structure group of the bundle, connections on principal $G$-bundles over a Lie groupoid $\Gamma$, for a general Lie groupoid $\Gamma$, can be characterized as two generalized morphisms from distinct (quasi) Lie groupoids with the same manifold of objects to the same groupoid $G$, which satisfy an additional compatibility condition, which should be expected to correspond to the vertical differential in the {\v C}ech--De Rham bicomplex used in~\cite{L-GTX}, which measures the condition for a quasi-connection to be a true connection.
This I plan to pursue also somewhere later.

In Section~\ref{sec-genmorcomp}, I first review the composition of generalized morphisms from a global point of view, following again~\cite{Mrcun} and~\cite{Moer2}; in particular, I focus on the computation of the division map of the composition of two generalized morphisms.
I then proceed by analyzing the concept of composition of local generalized morphisms.
This needs a ``refinement'' trick, due to the (tautological!) local nature of local generalized morphisms; once this point is clear, the composition of local generalized morphisms can be easily defined, and moreover, it is shown that composition of local generalized morphisms is equivalent to composition of generalized morphisms.
Let me just notice that, by the arguments of the final subsection of Section~\ref{sec-hilskan}, the composition of generalized morphisms is expected to correspond to an operation in nonabelian {\v C}ech cohomology; I do not intend to pursue this topic here, however I plan to deserve to the abstract cohomological aspects of the theory of generalized morphisms a forthcoming work.

In Section~\ref{sec-morita}, I first review the concept of {\em Morita equivalence}, again using as main references~\cite{Mrcun} and~\cite{Moer2}: I review the concept of {\em canonical division maps} of a Morita equivalence and the concept of the inverse $P^{-1}$ of a Morita equivalence $P$, showing that it is also a Morita equivalence and concentrating on the computation of its canonical division maps and their relationship with those of the Morita equivalence $P$.
In particular, using arguments of~\cite{CR1}, I find a {\em Factorization formula} for the (at first sight complicated) division map of the composition of a Morita equivalence $P$ with its inverse $P^{-1}$; this is the starting point from which I derive a notion of {\em local Morita equivalence}, which I prove to be completely equivalent to the ``global'' one.
Also Morita equivalences, having local counterparts, whose properties are strictly related to the composition of generalized morphisms, are expected to receive an interpretation in nonabelian {\v C}ech cohomology, recalling arguments from Section~\ref{sec-genmorcomp}; as for the final topic of the preceding Section, I will treat this topic also separately.

\begin{Ack}
I thank A.\ S.\ Cattaneo and G.\ Felder for reading the manuscript carefully and for many discussions; I also acknowledge the pleasant atmosphere at the Department of Mathematics of the Technion, where this work was accomplished.
\end{Ack}

\section{Principal bundles with groupoid structure: the division map}\label{sec-divmap}
I devote the first section to a reminder of the definition of principal bundles with groupoid structure borrowed from~\cite{Moer2}; in particular, I going to recall the definition of the {\em division map} for principal bundles.
I borrow the name ``division map'' from MacKenzie~\cite{McK}; M{\oe}rdijk~\cite{Moer1} calls the division map a {\em cocycle with values in a groupoid $\calG$} (although he considers it jointly with a smooth map, which is for me the momentum), and he already states some of its properties.
In~\cite{CR1}, I analyzed its properties carefully, in particular emphasizing its nature as a bundle map and its equivariance properties, which will prove to be a basic ingredient of later computations. 
The division map plays a fundamental r{\^o}le in the local description of principal bundles, and its main properties build the groundstone leading the local description of principal bundles; this r{\^o}le was already known to M{\oe}rdijk~\cite{Moer1}.
Let me only skip the introductory part to the theory of Lie groupoids which will be used throughout the paper, referring to~\cite{CR1} for the main conventions and notations.

I recall now the definition of principal bundles with groupoid structure.
\begin{Def}\label{def-princgroupoid}
A {\em principal bundle $P$ with groupoid structure $\calG$ over the manifold $M$} is a $4$-tuple $\left(P,\pi,\varepsilon,M\right)$, where $i)$ $P$ and $M$ are smooth manifolds and $ii)$ the pair $(P,\varepsilon)$ defines a structure of right $\calG$-space on $P$.
(Notice that usually, the right $\calG$-action is simply denoted as right multiplication; if otherwise a particular notation is needed, I will use the notation $\Psi$ or $\Psi_P$ for the right-action map.)

Moreover, the following requirements must hold:
\begin{itemize}
\item[i)] the map $\pi$ is a surjective submersion from $P$ to $M$; 
\item[ii)] the map $\pi$ is $\calG$-invariant, i.e.\ the following diagram commutes
\[
\begin{CD} 
P\times_{\varepsilon}\calG  @>\Psi>> P \\
@V\pr_1 VV  @VV\pi V\\
P @>\pi>> M
\end{CD}\quad ;
\] 
\item[iii)] the map $\left(\pr_1,\Psi\right)$ defined via
\begin{align*}
\left(\pr_1,\Psi\right)\colon P\times_\varepsilon \calG&\to P\times_M P,\\ 
(p,g)&\mapsto (p,pg),
\end{align*}
is a diffeomorphism; by $P\times_M P$ is meant
\[
P\times_M P\colon=\left\{(p,q)\in P\times P\colon \pi(p)=\pi(q)\right\}.
\]
\end{itemize} 
The map $\varepsilon$ is sometimes called the {\em (right) momentum of the bundle $P$}. 
\end{Def}

Let me notice that the first two requirements in Definition~\ref{def-princgroupoid} are the same as for ordinary principal bundles with structure group; the third one is most peculiar, but it may be viewed as another way of saying that the groupoid $\calG$ operates {\em freely and transitively on each fiber of $P$}.
Namely: assume first that the identity holds
\[
pg=p,\quad p\in P,g\in\calG\quad\text{such that $t_{\calG}(g)=\varepsilon(p)$}.
\]
It follows that the (a priori) distinct pairs $(p,g)$ and $(p,\iota_{\calG}(\varepsilon(p)))$, both in $P\times_\varepsilon \calG$, are mapped by the diffeomorphism $(\pr_1,\Psi)$ to the same image, namely $(p,p)$; hence, 
\[
g=\iota_{\calG}(\varepsilon(p)).
\] 
If, on the other hand, one takes any two points $p$ and $q$ of $P$, lying in the same fibre of $\pi$, which means that 
\[
\pi(p)=\pi(q),
\]
hence the pair $(p,q)$ belongs to $P\times_MP$, since $\left(\pr_1,\Psi\right)$ is a diffeomorphism, one has immediately that
\[
q=pg,\quad g\in\calG\quad \varepsilon(q)=s_{\calG}(g),\quad \varepsilon(p)=t_{\calG}(g),
\]
whence also $g\in \calG_{\varepsilon(q),\varepsilon(p)}$.

\subsection{First examples of principal bundles: unit bundle, pull-back bundle and trivial bundle}\label{ssec-unittriv}
In this subsection I define three particularly important principal bundles, which will play also a fundamental r{\^o}le in the local description of general principal bundles, namely the {\em unit bundle}, which is the basic groundstone for the subsequent theory, the natural notion of {\em pull-back bundle of a principal bundle $P$}, from which, using the unit bundle, from which I can define the notion of {\em trivial bundle}.
Notice that the notion of trivial bundle is not uniquely determined as in the case of ordinary principal bundles, but, in fact, there can be more than one trivial bundle over the same base space.

\begin{Def}\label{def-unitbun}
The {\em unit bundle of the Lie groupoid $\calG$} consists of the $4$-tuple $\left(\calG,t_\calG,s_\calG,X_\calG\right)$ (thus, it is a bundle over the manifolds of points of $\calG$), and the right $\calG$-action on itself is given by right multiplication; it is usually denoted by $\mathcal{U}_{\calG}$.
\end{Def}

Let me briefly sketch the proof of the fact that $\calU_\calG$, for any Lie groupoid $\calG$, is really a principal bundle in the sense of Definition~\ref{def-princgroupoid}.
That the projection $\pi=t_\calG$ is a surjective surjection follows immediately from the definition of Lie groupoid; similarly, the axioms of a Lie groupoid imply immediately that it is $\calG$-invariant in the above sense.
I come now to the last part of the proof: namely, let me consider the map $(\pr_1,\Psi)$ on $\calU_\calG$
\[
\calU_\calG\times_{s_{\calG}}\calG\ni (g_1,g_2)\mapsto (g_1,g_1g_2)\in \calU_\calG\times_{X_{\calG}}\calU_\calG.
\]
The inverse map thereof is simply given by
\[
\calU_\calG\times_{X_{\calG}}\calU_\calG\ni (g_1,g_2)\mapsto (g_1,g_1^{-1}g_2)\in\calU_\calG\times_{s_{\calG}}\calG.
\]
Notice that the previous map makes sense: in fact, if the pair $(g_1,g_2)$ belongs to $\calU_\calG\times_{X_{\calG}}\calU_\calG$, this means that
\[
t_\calG(g_1)=s_{\calG}(g_1^{-1})=t_{\calG}(g_2)\quad\text{and}\quad s_{\calG}(g_1)=t_{\calG}(g_1^{-1}).
\]
The axiom of a Lie groupoid imply that the above map is smooth, hence it is a diffeomorphism.
The map introduced above can be really thought of as a division map; in fact, this is the context where MacKenzie derived its name from. 

I define now the pull-back bundle of a general principal bundle in the sense of Definition~\ref{def-princgroupoid}.
\begin{Def}\label{def-pullbun}
If the $4$-tuple $\left(P,\pi,\varepsilon,N\right)$ is a principal $\calG$-bundle over $N$ and $M\overset{f}\to N$ is a smooth map from the manifold $M$ to the manifold $N$, the {\em pull-back bundle $f^*P$ of $P$ w.r.t.\ $f$} is defined by the $4$-tuple $\left(f^*P,\pr_1,\varepsilon\circ\pr_2,M\right)$, where the space $f^*P$ is  
\[
f^*P\colon=\left\{(m,p)\in M\times P\colon f(m)=\pi(p)\right\},
\]
and $\pr_i$, $i=1,2$, denotes projection onto the $i$-th term of $f^*P$.
\end{Def}

\begin{Lem}\label{lem-pullbackbun}
The $4$\-tuple $\left(f^*P,\pr_1,\varepsilon\circ\pr_2,M\right)$ is a principal bundle in the sense of Definition~\ref{def-princgroupoid}.
\end{Lem} 
\begin{proof}
The bundle projection $\pr_1$ is clearly a surjective submersion.
Since the right $\calG$ action, which is defined along the map
\[
(m,p)\overset{\varepsilon\circ\pr_2}\to \varepsilon(p),
\]
takes the explicit form
\[
f^*P\times_{\varepsilon\circ\pr_2}\calG\ni (m,p;g)\mapsto(m,pg),
\]
the bundle projection is also clearly $\calG$-invariant.

If two points $(m_1,p)$ and $(m_2,q)$ of $f^*P$ belong to the same fibre, it follows
\[
m_1=m_2\Rightarrow f(m_1)=\pi(p)=f(m_2)=\pi(q)\Leftrightarrow q=pg_{p,q},
\] 
for some element $g_{p,q}\in\calG$, since $P$ is a principal bundle.
Thus, the map 
\[
f^*P\times_{\varepsilon\circ\pr_2}\calG\ni (m,p;g)\to (m,p;m,pg)
\]
is a diffeomorphism, where the smooth inverse is given explicitly by
\[
f^*P\times_{M} f^*P\ni(m,p;m,q)\to \left(m,p;g_{p,q}\right), 
\]
and the claim follows.
\end{proof}

Now, I give the definition of trivial bundles over a base manifold $M$.
\begin{Def}\label{def-trivbun}
Given a groupoid $\calG$ and a smooth map $\alpha$ from a manifold $M$ to the manifold of objects $X_{\calG}$ of $\calG$, I consider the pull-back bundle $\alpha^*\calU_{\calG}$ of the unit bundle of $\calG$.

By its very definition, the total space of this bundle has the form
\[
\alpha^*\calU_{\calG}=\left\{(m,g)\in M\times \calG\colon \alpha(m)=t_{\calG}(g)\right\}.
\]
The bundle $\alpha^*\calU_{\calG}$ is called the {\em trivial $\calG$-bundle over $M$ w.r.t.\ $\alpha$}. 
\end{Def}
Lemma~\ref{lem-pullbackbun} implies that $\alpha^*\calU_\calG$ is in fact a principal bundle.

\begin{Exa}\label{exa-trivial1}
Recall that any Lie group $G$ may be viewed more generally as a Lie groupoid, where the manifold of objects is simply a point $*$, and target, source and identity map are defined accordingly.

There is only one map $\alpha$ from a manifold $M$ to the point $*$, mapping all $M$ onto $*$. 
Hence, ``the'' trivial bundle $\alpha^*\calU_G$ takes the form:
\begin{align*}
\alpha^*\calU_G&=\left\{(m,g)\in M\times G\colon \alpha(m)=*=t_G(g)\right\}=\\
&=M\times G,
\end{align*}
which coincides with the usual definition of trivial principal bundle over $M$.
\end{Exa}

\begin{Rem}\label{rem-trivial}
Notice that, while there is only one trivial principal $G$-bundle over a manifold $M$, with $G$ a group, there can be in principle many {\em distinct} trivial $\calG$-bundles over the same base. 
\end{Rem}

\begin{Exa}\label{exa-trivial2}
One can consider the manifold $M$ to be a point $*$; then, the map $\alpha$ simply sends the point $*$ to some $x\in X_\calG$ and this is clearly a smooth map.
The associated trivial bundle $\alpha^*\calU_\calG$ is the subset of $\calG$ of ``arrows'' arriving at $x$: namely, the base space can be immediately identified with the point $x$ and the total space is by definition
\begin{align*}
\alpha^*\calU_\calG&=\left\{(*,g)\in \left\{*\right\}\times\calG\colon t_{\calG}(g)=\alpha(*)=x\right\}\cong\\ 
&\cong \calG_{\bullet,x}.
\end{align*}
Hence, if consider e.g.\ the action groupoid $G\ltimes M$, for a smooth manifold acted on from the left by a Lie group $G$, then the trivial bundle over a point $*$, mapped to the point $m$ in $M$, is simply the $G$-orbit in $M$ through the point $m$.
\end{Exa}

\begin{Rem}
Observe that the ``momentum map'' $\varepsilon$, along which the right action of $\calG$ on $P$ is defined, is a surjective submersion in the case of a trivial bundle, as it is the composition of two surjective submersions.
\end{Rem}

I introduce now the following
\begin{Def}\label{def-isombun}
Given two principal bundles $\left(P,\pi,\varepsilon,M\right)$ and $\left(\widetilde{P},\widetilde{\pi},\widetilde{\varepsilon},M\right)$ over the same base manifold $M$ and with the same structure groupoid $\calG$, a {\em morphism of principal bundles from $P$ to $\widetilde{P}$} is a smooth map $\tau$ from $P$ to $\widetilde{P}$ enjoying the two requirements:
\begin{itemize}
\item[i)] $\tau$ is fibre-preserving, i.e.\ the following identity must hold:
\[
\widetilde{\pi}\circ \tau=\pi.
\] 
\item[ii)] $\tau$ is $\calG$-equivariant, i.e.\
\[
\widetilde{\varepsilon}\circ\tau=\varepsilon,\quad \tau(pg)=\tau(p)g,\quad \forall (p,g)\in P\times_{\varepsilon}\calG.
\] 
\end{itemize}
\end{Def}
\begin{Rem}
Notice that the first identity in $ii)$ of Definition~\ref{def-isombun} implies that both terms in the second identity are well-defined.
\end{Rem}

In the terminology introduced in~\cite{CR1}, a morphism between two principal bundles $P$ and $\widetilde{P}$ over the same base manifold and with the same structure groupoid is a {\em fibre-preserving twisted equivariant map between the $\calG$-spaces $P$ and $\widetilde{P}$}, when coupled with the identity map of $M$.

I recall from~\cite{CR1} that any morphism of principal bundles over the same base and with the same structure groupoid is invertible.
Hence, two principal bundles $P$ and $\widetilde{P}$ over the same base space and with the same structure groupoid, for which there exists a morphism in the sense of Definition~\ref{def-isombun}, are said to be {\em isomorphic}.

\subsection{The division map: definition and memento of main properties}\label{ssec-defdivmap}
In this subsection I discuss the {\em division map} of a general principal bundle $P$; for the name of the map, canonically associated to $P$, I have followed the convention adopted by MacKenzie~\cite{McK} for ordinary principal bundles.
First of all, I need a preliminary Lemma, whose proof may be found in~\cite{CR1}.
\begin{Lem}\label{lem-fibprod}
The $4$-tuple $\left(P\odot\widetilde{P},\overline{\pi},\varepsilon\times\widetilde{\varepsilon},M\right)$, where the manifold $P\odot\widetilde{P}$ is defined by
\[
P\odot\widetilde{P}\colon=\left\{(p,\tildep)\in P\times \widetilde{P}\colon \pi(p)=\widetilde{\pi}(\tildep)\right\},
\]
and the projection $\overline{\pi}$ is 
\[
\overline{\pi}(p,\tildep)=\pi(p)=\widetilde{\pi}(\tildep),
\]
defines a principal $\calG^2$\-bundle over $M$, which is called the fibred product bundle of $P$ and $\widetilde{P}$.
\end{Lem}

\begin{Rem}
It is customary to denote the fibred product bundle of $P$ and $\widetilde{P}$ by $P\times_M\widetilde{P}$, but I prefer to use the previous notation, which reminds somehow of the Whitney sum notation, whose analogon in the framework of principal bundles is exactly the fibred product operation.
\end{Rem}

I recall the definition of the generalized conjugation of a Lie groupoid $\calG$, referring to~\cite{CR1} for a more detailed description.

The generalized conjugation of a Lie groupoid $\calG$ consists of an action of the product groupoid $\calG^2$ of $\calG$ with itself on $\calG$; as such (see again~\cite{Moer2} or~\cite{CR1} for more details), it consists of a $3$-tuple $\left(\calG^2,J_{\conj},\Psi_{\conj}\right)$, with $J_{\conj}$ the momentum of the action and $\Psi_{\conj}$ the explicit action map.

The momentum $J_{\conj}$ of the generalized conjugation is simply 
\[
J_{{\conj}}\!\left(g\right)\colon=\left(t_\calG(g),s_\calG(g)\right),\quad \forall g\in\calG.
\]
Thus, the manifold $\calG^2\times_{J_{\conj}}\calG$, where the action makes sense, takes the form
\[
\calG^2\times_{J_{\conj}}\calG=\left\{\left(g_1,g_2;g_3\right)\in\calG^3\colon \begin{cases}
s_\calG(g_1)&=t_\calG(g_3)\\
s_\calG(g_2)&=s_\calG(g_3)
\end{cases}\right\}.
\]
Define then the action map $\Psi_{\conj}$ of the generalized conjugation from $\calG^2\times_{J_{\conj}}\calG$ to $\calG$ as
\begin{equation}\label{eq-genconjgroup}
\Psi_{\conj}\!\left(g_1,g_2;g_3\right)\colon= g_{1}g_3 g_2^{-1}.
\end{equation}

\begin{Prop}\label{prop-genconj}
The triple $\left(\calG^2,J_{\conj},\Psi_{\conj}\right)$ defines a left $\calG^2$-action on $\calG$, which I call the {\em generalized conjugation of $\calG$}. 
\end{Prop}
See~\cite{CR1} for the proof.

\begin{Rem}\label{rem-genconj}
Let me notice that there is a similar, still distinct, left $\calG^2$-action on $\calG$; in fact, one can consider the map momentum map $\overline{J}_{\conj}$ from $\calG$ to $X_\calG\times X_\calG$ given by
\[
\overline{J}_{\conj}(g)\colon=\left(s_{\calG}(g),t_{\calG}(g)\right),
\]
whence
\[
\calG^2\times_{\overline{J}_{\conj}}\calG=\left\{(g_1,g_2;g_3)\in\calG^3\colon \begin{cases}
s_\calG(g_1)&=s_{\calG}(g_3)\\
s_{\calG}(g_2)&=t_\calG(g_3)
\end{cases}
\right\},
\]
and the action map $\overline{\Psi}_{\conj}$ from $\calG^2\times_{\overline{J}_{\conj}}\calG$ to $\calG$ via
\[
\overline{\Psi}_{\conj}\!\left(g_1,g_2;g_3\right)\colon=g_2g_3g_1^{-1}.
\]
It is not difficult to verify that the triple $\left(\calG,\overline{J}_{\conj},\overline{\Psi}_{\conj}\right)$ defines also a left $\calG^2$-action on $\calG$.
\end{Rem}

\begin{Rem}\label{rem-conjright}
The maps $J_{\conj}$ and $\overline{J}_{\conj}$ define also right $\calG^2$-actions on $\calG$, the {\em right generalized conjugations}: namely, on the set $\calG\times_{J_{\conj}}\calG^2$, resp.\ $\calG\times_{\overline{J}_{\conj}}\calG^2$, define the map $\Psi_{\conj}^R$, resp.\ $\overline{\Psi}_{\conj}^R$, by the formula
\begin{align*}
(g_3;g_1,g_2)&\overset{\Psi_{\conj}^R}\mapsto g_1^{-1}g_3g_2,\quad\text{resp.}\\
(g_3;g_1,g_2)&\overset{\overline{\Psi}_{\conj}^R}\mapsto g_2^{-1}g_3g_1.
\end{align*}
\end{Rem}

Define now the {\em division map} of a general principal bundle $P$.
\begin{Def}\label{def-divmap}
Given a principal bundle $\left(P,\pi,\varepsilon,M\right)$ with structure groupoid $\calG$, the division map $\phi_P$ of $P$ is defined by the requirement
\begin{equation}\label{eq-divmap}
q=p \phi_P(p,q),\quad \pi(p)=\pi(q).
\end{equation}
\end{Def}
First of all, notice that the division map, because of Equation (\ref{eq-divmap}), is defined on the fibred product bundle $P\odot P$, and that it is, in fact, the second component of the smooth inverse of the canonical map $\left(\pr_1,\Psi\right)$ from $P\times_\varepsilon \calG$ to $P\odot P$.
Namely, the inverse of the map $\left(\pr_1,\Psi\right)$ can be factorized as follows
\[
\left(\pr_1,\Psi\right)^{-1}(p,q)=\left(\Phi_{P,1}(p,q),\Phi_{P,2}(p,q)\right),
\]
where $\Phi_{P,1}(p,q)$ belongs to $P$ and $\Phi_{P,2}(p,q)$ belongs to $\calG$, for any pair $(p,q)$ in the fibred product $P\odot P$.
From the very definition of inverse map, it follows easily
\begin{align*}
\left(\pr_1,\Psi\right)\!\left(\Phi_{P,1}(p,q),\Phi_{P,2}(p,q)\right)&=(\Phi_{P,1}(p,q),\Phi_{P,1}(p,q)\Phi_{P,2}(p,q))=\\
&=(p,q),
\end{align*}
whence it follows that
\[
\Phi_{P,1}=\pr_1,\quad \Phi_{P,2}=\phi_P.
\]

Let me now list the main properties of the division map $\phi_P$
\begin{Prop}\label{prop-prodivmap}
The map $\phi_P$ from $P\odot P$ to $\calG$ has the following properties:
\begin{itemize}
\item[i)] for any point $(p,q)$ of $P\odot P$, one has
\[
\phi_P(p,q)\in \calG_{\varepsilon(q),\varepsilon(p)}.
\] 
\item[ii)] On the diagonal submanifold of the total space of $P\odot P$, one has
\[
\phi_P(p,p)=\iota_{\calG}(\varepsilon(p)),\quad \forall p\in P.
\] 
\item[iii)] for any pair $(p,q)\in P\odot P$, the following equation holds
\[
\phi_P(p,q)=\phi_P(q,p)^{-1};
\] 
notice that the previous equation makes sense, since $(p,q)\in P\odot P$ implies that $(q,p)\in P\odot P$ also.
\item[iv)] The triple $\left(\phi_P,\id_{\calG^2},\id_{X_{\calG}^2}\right)$ is an equivariant map from the right $\calG^2$-space $P\odot P$ to the right $\calG^2$-space $\calG$ endowed with the right generalized conjugation defined by $\left(J_{\conj}^R,\Psi_{\conj}^R\right)$.  
\end{itemize}
\end{Prop}
See once again to~\cite{CR1} for the proof of Proposition~\ref{prop-prodivmap}

\begin{Exa}
Let me consider the unit bundle $\calU_\calG$ associated to a general Lie groupoid $\calG$, see Definition~\ref{def-unitbun}.
It is then easy to see that the division map $\phi_{\calU_{\calG}}=\phi_{\calG}$ of the unit bundle, defined on the space of pairs $(g_1,g_2)$ ending at the same point (which corresponds clearly to the fibred product of the unit bundle with itself), is the ``true'' division map
\[
(g_1,g_2)\mapsto g_1^{-1}g_2,\quad t_{\calG}(g_1)=t_{\calG}(g_2).
\]
\end{Exa}

\begin{Exa}
Given a principal bundle $P$ over a manifold $N$, and a smooth map $f$ from a manifold $M$ to $N$, it is not difficult to prove that the division map $\phi_{f^*P}$ of the pull-back bundle $f^*P$ is simply
\[
\phi_{f^*P}\!\left((m,p),(m,q)\right)=\phi_P(p,q),\quad f(m)=\pi(p)=\pi(q).
\]
Hence, the division map of the trivial bundle $\alpha^*\calU_{\calG}$ over $M$ associated to the smooth map $\alpha$ is simply
\[
\phi_\alpha\!\left((m,g_1),(m,g_2)\right)=g_1^{-1}g_2,\quad t_{\calG}(g_1)=t_{\calG}(g_2).
\]
\end{Exa}

\begin{Rem}
Given a principal bundle $\left(P,\pi,\varepsilon,M\right)$ with structure groupoid $\calG$, the pair $\left(\varepsilon,\phi_P\right)$ is called by M{\oe}rdijk~\cite{Moer1} a ``cocycle on $M$ with values in $\calG$; in the next section, I will explain in which sense this denomination has to be understood.
\end{Rem}

\section{Local data of principal bundles with groupoid structure}\label{sec-localdata}
In the previous section I introduced the main notion of principal bundle with groupoid structure and I discussed the first examples of principal bundles, namely the unit bundle and the trivial bundles associated to smooth maps; further, I introduced the division map of a general principal bundle $P$ and listed its main properties.

The main properties characterizing a principal bundle $P$ are encoded in the $\calG$-invariant surjective submersion $\pi$ and the division map, whose existence makes the action of $\calG$ on $P$ free and transitive on every fibre.

Now, I want to give a ``constructive'' definition of principal bundles, dealing with local data, such as in the ordinary case; in other words, I trivialize locally a general principal bundle, with the help of the submersivity of the projection $\pi$ and the division map.
Later, I will generalize the notion of local trivializing data, and I will also prove that local trivializing data provide an equivalent way of defining principal bundles.
Finally, with the help of local trivializing data, I will construct some examples of principal bundles for some particular Lie groupoids.

\subsection{Local sections of $\pi$ and trivializations of $P$}\label{ssec-localsec}
Recall that the projection $\pi$ of a principal bundle $P$ is a surjective submersion, i.e.\ the tangent map at any point of $P$ is a surjective linear map between the corresponding tangent space; hence, the Implicit Function Theorem implies the following useful 
\begin{Lem}\label{lem-loctriv}
Any principal $\calG$-bundle $\left(P,\pi,\varepsilon,M\right)$ is locally isomorphic to a trivial bundle, i.e.\ for any point $m$ of $M$ there is an open neighbourhood $U$, such that the restriction of $P$ to $U$ is isomorphic to a trivial bundle over $U$.
\end{Lem} 
\begin{proof}[See also~\cite{Moer2}]
Consider a general point $m$ of $M$ and choose a local section $\sigma_U$ of $\pi$ (it is possible by the Implicit Function Theorem, since $\pi$ is a surjective submersion) over an open neighbourhood $U=U_m$, and consider the (smooth) composite map 
\[
\varepsilon_U\colon=\varepsilon\circ\sigma_U.
\]
Consider the map
\[
\varepsilon_U^*\calU_{\calG}\ni(m,g)\overset{\varphi_U}\mapsto \sigma(m)g\in P;
\]
by the very definition of the map $\varepsilon_U$ and Definition~\ref{def-pullbun}, the map $\varphi_U$ is well-defined and smooth.

Since the restriction of $P$ to $U$ and $\varepsilon_U^*\calU_\calG$ are both principal bundles over the same manifold $U$ and with the same structure groupoid, to prove that both bundles are isomorphic via $\varphi_U$, it suffices to prove that $\varphi_U$ is $\calG$-equivariant and fibre\-preserving.

Let me prove first that it is fibre-preserving.
Namely, for any pair $(m,g)$ in $\varepsilon_U^*\calU_\calG$, one get
\begin{align*}
\left(\pi\circ \varphi_U\right)\!(m,g)&=\pi(\sigma_U(m)g)=\\
&=\pi(\sigma_U(m))=\\
&=m=\\
&=\pr_1(m,g),
\end{align*}
where I used the $\calG$-invariance of the projection $\pi$.

Now let me prove $\calG$-equivariance.
First, one has to show that $\varphi_U$ respects the momenta of the actions of $\calG$:
\begin{align*}
\left(\varepsilon\circ \varphi_U\right)\!(m,g)&=\varepsilon(\sigma_U(m)g)=\\
&=s_{\calG}(g),\quad \forall (m,g)\in \varepsilon_U^*\calU_\calG,
\end{align*}
by the very definition of the momentum for the action of $\calG$ on the trivial bundle $\varepsilon_U^*\calU_\calG$, and by the very properties of the momentum.
Furthermore, 
\begin{align*}
\varphi_U\!\left((m,g_1)g_2\right)&=\varphi_U(m,g_1g_2)=\\
&=\sigma_U(m)g_1g_2=\\
&=(\sigma_U(m)g_1)g_2=\\
&=\varphi_U(m,g_1)g_2,\quad \forall (m,g_1)\in\varepsilon_U^*P,s_{\calG}(g_1)=t_{\calG}(g_2).
\end{align*}
\end{proof} 

The isomorphism $\varphi_U$, associated to the local section $\sigma_U$ over $U$, is usually called a {\em local trivialization of $P$}.
Notice that the trivialization depends only on a choice of a local section of $\pi$, whereas the inverse of the trivialization depends additionally on the division map of $P$.
The question that arises naturally now is:

\fbox{\parbox{12cm}{\bf Given two trivializing open sets $U$, $V$, intersecting nontrivially, and corresponding sections of $\pi$, $\sigma_U$ and $\sigma_V$ respectively, giving rise to trivializations $\varphi_U$, resp.\ $\varphi_V$, is there an explicit relationship between the trivializations?}}

This is the question I want to answer in what follows, but, before entering into the details, I have to fix some notations
.
Given a local section $\sigma_U$ of $\pi$ over some open subset $U$ of $M$, denote by $\varepsilon_U$ the composite map
\[
\varepsilon_U\colon=\varepsilon\circ\sigma_U,
\]
which is a smooth map from $U$ to the manifold of objects $X_\calG$ of $\calG$, called the local momentum of $P$ w.r.t.\ the open set $U$; consequently, one can consider the trivial bundle $\varepsilon_U^*\calU_\calG$ associated to $\varepsilon_U$.
Given a principal bundle $P$ over $M$ and some open subset $U\subset M$, I will use the following short-hand notation for the restriction of $P$ to $U$:
\[
P_U\colon=\pi^{-1}\!(U).
\]
I need first the following technical
\begin{Lem}\label{lem-transmap}
If $U$, $V$ are two open subsets of $M$, intersecting nontrivially, over which there are trivializations $\varphi_U$ and $\varphi_V$, associated to sections $\sigma_U$ and $\sigma_V$ respectively, then 
\[
(\varepsilon_U^*\calU_\calG)_{U\cap V}\cong (\varepsilon_V^*\calU_\calG)_{U\cap V}.
\]
\end{Lem}
\begin{proof}
Given two trivializations $\varphi_U$ and $\varphi_V$ on the open sets $U$ and $V$ respectively, a morphism from $(\varepsilon_U^*\calU_\calG)_{U\cap V}$ to $(\varepsilon_V^*\calU_\calG)_{U\cap V}$ can be simply defined via
\[
\varphi_{VU}\colon=\varphi_V^{-1}\circ\varphi_U,
\]
where, of course, I consider the restrictions of the respective trivializations, so that the morphism is well-defined.

It remains to check that $\phi_{VU}$ is fibre-preserving and $\calG$-equivariant, which, on the other hand, are both consequences of Lemma~\ref{lem-loctriv}; hence, the claim follows.
\end{proof}

Now, I want to find an explicit expression for the isomorphism in Lemma~\ref{lem-transmap}.
Recalling that $\varphi_U$, resp.\ $\varphi_V$, is defined via a section $\sigma_U$ of $\pi$ over $U$, resp.\ $\sigma_V$ over $V$, one finds by a direct computation:
\begin{align*}
\varphi_{VU}(m,g)&=\varphi_V^{-1}\!\left(\sigma_U(m)g\right)=\\
&=\left(\pi(\sigma_U(m)g),\phi_P(\sigma_V(\pi(\sigma_U(m)g)),\sigma_U(m)g)\right)=\\
&=\left(m,\phi_P(\sigma_V(m),\sigma_U(m))g\right);
\end{align*}
in the previous computations, I made use of the $\calG$-invariance of the projection $\pi$, of the fact that $\sigma_U$ is a section of $\pi$ and of the equivariance properties of the division map of Proposition~\ref{prop-prodivmap}.

From now on, denote the map
\[
U\cap V\ni m\mapsto \phi_P(\sigma_V(m),\sigma_U(m))\in \calG
\]
simply by $\Phi_{VU}$, for any two local sections $\sigma_U$, $\sigma_V$ of $\pi$ over $U$, $V$.
I now want to analyse in detail the properties of the map $\Phi_{VU}$, which is called the {\em transition map from the trivialization $\varphi_U$ to the trivialization $\varphi_V$} or shortly the {\em trivialization from $U$ to $V$}.
\begin{Prop}\label{prop-protrans}
Given two open subsets $U$ and $V$ of $M$, intersecting nontrivially, and associated sections $\sigma_U$ and $\sigma_V$ respectively, the map $\Phi_{VU}$ enjoys the following properties:
\begin{itemize}
\item[i)] the following identities hold:
\[
t_{\calG}\circ\Phi_{VU}=\varepsilon_V,\quad s_{\calG}\circ\Phi_{VU}=\varepsilon_U,\quad \Phi_{UU}=\iota_{\calG}\circ \varepsilon_U.
\] 
\item[ii)] The following identity holds:
\[
\Phi_{UV}^{-1}=\Phi_{VU}.
\]
\item[iii)] For any three open subsets $U$,$V$ and $W$ of $M$, such that their triple intersection $U\cap V\cap W$ is nontrivial, the following identity holds:
\[
\Phi_{WU}(m)=\Phi_{WV}(m)\Phi_{VU}(m),\quad \forall m\in U\cap V\cap W.
\]
\end{itemize}
\end{Prop}
\begin{proof}
The proof of the first two statements follows directly from the definition of the maps $\varepsilon_U$ and $\varepsilon_V$ and from Proposition~\ref{prop-prodivmap}.

The third statement can be proved as follows: since $\Phi_{WU}$, $\Phi_{WV}$ and $\Phi_{VU}$ are all defined via sections of $\pi$, it follows from Equation (\ref{eq-divmap}) that
\begin{align*}
\Phi_{WU}(m)&=\phi_P(\sigma_W(m),\sigma_U(m))=\\
&=\phi_P(\sigma_V(m)\phi_P(\sigma_V(m),\sigma_W(m)),\sigma_U(m))=\\
&=\phi_P(\sigma_V(m),\sigma_W(m))^{-1}\phi_P(\sigma_V(m),\sigma_U(m))=\\
&=\phi_P(\sigma_W(m),\sigma_V(m))\phi_P(\sigma_V(m),\sigma_U(m))=\\
&=\Phi_{WV}(m)\Phi_{VU}(m),
\end{align*}
where I used also Proposition~\ref{prop-prodivmap}.
Notice that Property $i)$ implies that one can actually multiply $\Phi_{WV}(m)$ and $\Phi_{VU}(m)$.
\end{proof}

\subsection{Local trivializing data for principal bundles}\label{ssec-trivdata}
Motivated by the results of the last subsection, in particular Proposition~\ref{prop-protrans} regarding the properties of the transitions maps $\Phi_{VU}$, I want now to generalize the notion of transitions maps, and so I am led to the notion of {\em local trivializing data}.
Let me denote by $\mathfrak{U}=\left\{U_\alpha\right\}_\alpha$ an open cover of $M$; borrowing the notations from the algebro-geometric framework, I denote multiple intersections by
\[
U_{\alpha_1}\cap \cdots \cap U_{\alpha_p}=\colon U_{1\cdots p}\quad\text{or}\quad U_{\alpha_1}\cap \cdots \cap U_{\alpha_p}=\colon U_{\alpha_1\cdots \alpha_p}.
\]

\begin{Def}\label{def-trivdata}
Given a smooth manifold $M$ and a Lie groupoid $\calG$, local trivializing data over $M$ with values in $\calG$ (or shortly, local trivializing data, when the manifold $M$ and the Lie groupoid $\calG$ are clear from the context) consist of a $3$-tuple $\left(\mathfrak{U},\varepsilon_\alpha,\Phi_{\alpha\beta}\right)$, where $i)$ $\mathfrak{U}=\left\{U_\alpha\right\}_\alpha$ is an open cover of $M$, $ii)$ the $\varepsilon_\alpha$'s are smooth maps from $U_\alpha$ to the manifold of objects $X_\calG$ of the groupoid $\calG$, called {\em the local momenta of the data}, and $iii)$, for any two open sets $U_\alpha$ and $U_\beta$ of the cover $\mathfrak{U}$, intersecting nontrivially, smooth maps $\Phi_{\alpha\beta}$ from $U_\alpha\cap U_\beta$ to $\calG$, called also {\em transition maps or cocycle}, enjoying the properties:
\begin{itemize}
\item[a)] the following identities must hold:
\[
t_\calG\circ \Phi_{\alpha\beta}=\varepsilon_\alpha,\quad s_{\calG}\circ \Phi_{\alpha\beta}=\varepsilon_\beta,\quad \Phi_{\alpha\alpha}=\iota_{\calG}\circ\varepsilon_{\alpha};
\] 
\item[b)] for any three open subsets $U_\alpha$, $U_\beta$ and $U_\gamma$ of the cover $\mathfrak{U}$, such that their triple intersection $U_{\alpha\beta\gamma}$ is nontrivial, the following identity ({\em nonabelian cocycle identity}) must hold
\[
\Phi_{\alpha\gamma}(m)=\Phi_{\alpha\beta}(m)\Phi_{\beta\gamma}(m),\quad \forall m\in U_{\alpha\beta\gamma}.
\] 
\end{itemize}
\end{Def}
Notice that property $a)$ implies that the identity in $b)$ is well-defined.

The condition on $b)$ is called {\em cocycle condition}, because it is reminiscent of the ordinary cocycle condition for Lie groups; later, I will discuss in detail the cohomology theory behind it.

\begin{Rem}
The notion of local trivializing data can be found already in~\cite{Con} and~\cite{H}, although there is no explicit mentioning of the local momenta and the relationship between the cocycle $\Phi_{\alpha\beta}$ and the local momenta; in fact, the local momenta, in the cohomological framework, which I am also going to discuss in subsection~\ref{ssec-cohomint}, can be hidden as ``unit $0$-cochains associated to a $1$-cocycle with values in a sheaf of groupoids''.
I prefer to consider them explicitly, since I want to express the combination between $1$-cocycles as in~\cite{Con} and~\cite{H}, emphasizing the nonabelian {\v C}ech cohomological aspects, and the approach of M{\oe}rdijk~\cite{Moer1}, emphasizing the presence of a momentum. 
\end{Rem}

\begin{Exa}
Given a principal bundle $P$ over the manifold $M$ with structure groupoid $\calG$, Lemmata~\ref{lem-loctriv} and~\ref{lem-transmap} provide an example of local trivializing data, namely, after having chosen a countable open cover $\mathfrak{U}=\left\{U_\alpha\right\}_\alpha$ of $M$ and associated local sections $\sigma_\alpha$ of $\pi$ over $U_\alpha$, one has the local trivializing data
\[
\left(\mathfrak{U},\varepsilon_\alpha\colon=\varepsilon\circ \sigma_\alpha,\Phi_{\alpha\beta}\colon=\phi_P(\sigma_\alpha,\sigma_\beta)\right).
\]
Lemma~\ref{lem-transmap} is therefore the bridge to understand why M{\oe}rdijk called the division map, together with the momentum $\varepsilon$, a cocycle.
\end{Exa}

First of all, I need a technical 
\begin{Lem}\label{lem-trivdata}
Given local trivializing data $\left(\mathfrak{U},\varepsilon_\alpha,\Phi_{\alpha\beta}\right)$, and two open sets $U_\alpha$ and $U_\beta$ with nontrivial intersection $U_\alpha\cap U_\beta$, there is an isomorphism between the trivial bundles $(\varepsilon_\alpha^*\calU_\calG)_{U_{\alpha\beta}}$ and $(\varepsilon_\beta^*\calU_\calG)_{U_{\alpha\beta}}$.
\end{Lem}
\begin{proof}
Consider the following map from $(\varepsilon_\beta^*\calU_\calG)_{U_{\alpha\beta}}$ to $(\varepsilon_\alpha^*\calU_\calG)_{U_{\alpha\beta}}$:
\[
\varepsilon_\beta^*\calU_\calG\ni (m,g)\overset{\varphi_{\alpha\beta}}\mapsto (m,\Phi_{\alpha\beta}(m)g).
\]
One has to prove that $i)$ $\varphi_{\alpha\beta}$ is well-defined and that it maps really $(\varepsilon_\beta^*\calU_\calG)_{U_{\alpha\beta}}$ onto $(\varepsilon_\alpha^*\calU_\calG)_{U_{\alpha\beta}}$, that $ii)$ it is fibre-preserving and $iii)$ that it is $\calG$-equivariant. 

To prove $i)$, recall that the pair $(m,g)$ belongs to $(\varepsilon_\beta^*\calU_\calG)_{U_{\alpha\beta}}$ if and only if 
\[
t_{\calG}(g)=\varepsilon_\beta(m);
\]
thus, by Property $a)$ of Definition~\ref{def-trivdata}, the map $\phi_{\alpha\beta}$ is well-defined, and moreover
\begin{align*}
t_{\calG}\!\left(\Phi_{\alpha\beta}(m)g\right)&=t_{\calG}\!\left(\Phi_{\alpha\beta}(m)\right)=\\
&=\varepsilon_\alpha(m),
\end{align*}
whence it follows that $\varphi_{\alpha\beta}$ maps $(\varepsilon_\beta^*\calU_\calG)_{U_{\alpha\beta}}$ onto $(\varepsilon_\alpha^*\calU_\calG)_{U_{\alpha\beta}}$.

The proof of $ii)$ is trivial, hence it remains only to show the $\calG$-equivariance.
Recalling the definition of momentum map for trivial bundles, one gets
\begin{align*}
\varepsilon_{\varepsilon_\alpha^*\calU_\calG}\!\left(\varphi_{\alpha\beta}(m,g)\right)&=\varepsilon_{\varepsilon_\alpha^*\calU_\calG}\!\left((m,\Phi_{\alpha\beta}(m)g)\right)=\\
&=s_{\calG}\!\left(\Phi_{\alpha\beta}(m)g\right)=\\
&=s_{\calG}(g)=\\
&=\varepsilon_{\varepsilon_\beta^*\calU_\calG}(m,g). 
\end{align*}
Furthermore, recalling the definition of the right $\calG$-action on trivial bundles, $\calG$-equivariance follows immediately.

Since fibre-preserving, $\calG$-equivariant morphism between principal bundles over the same base space and with the same structure groupoid are invertible, the claim follows immediately.
\end{proof}

Lemma~\ref{lem-trivdata} is the groundstone of the ``constructive'' definition of principal bundles, which is the analogon of the local construction of ordinary principal bundles.
In fact, assume one is given local trivializing data $\left(\mathfrak{U},\varepsilon_\alpha,\Phi_{\alpha\beta}\right)$ over the manifold $M$ with values in the Lie groupoid $\calG$ in the sense of Definition~\ref{def-trivdata}; then consider the disjoint union of all trivial bundles $\varepsilon_\alpha^*\calU_\calG$
\[
Q_{\mathfrak{U}}\colon=\coprod_{\alpha}\varepsilon_\alpha^*\calU_\calG,
\]
i.e.\ the set consisting of all $3$-tuples of the form
\[
\left(\alpha,m,g\right),\quad\text{$\alpha$ is an index for the open cover $\mathfrak{U}$},\quad (m,g)\in\varepsilon_\alpha^*\calU_\calG.
\]
Introduce then the following equivalence relation on $Q_{\mathfrak{U}}$:
\begin{equation}\label{eq-glueing}
\left(\alpha,m_1,g_1\right)\sim(\beta,m_2,g_2)\Leftrightarrow\begin{cases}
&U_{\alpha\beta}\neq \emptyset,\\
&m_1=m_2\in U_{\alpha\beta},\\
&g_1=\Phi_{\alpha\beta}(m_1)g_2.
\end{cases}
\end{equation}
The equivalence relation makes sense in spite of Lemma~\ref{lem-trivdata}: in fact, Relation (\ref{eq-glueing}) means simply that, whenever one restricts the disjoint union $Q_{\mathfrak{U}}$ to double intersections of open sets of the cover $\mathfrak{U}$, one has an isomorphism between them.
Moreover, that the relation (\ref{eq-glueing}) is really an equivalence relation, it follows directly from Definition~\ref{def-trivdata}.

Consider now the quotient of $Q_{\mathfrak{U}}$ by the equivalence relation (\ref{eq-glueing}), which I denote by $P_{\mathfrak{U}}$:
\begin{equation}\label{eq-gluedprincbun}
P_{\mathfrak{U}}\colon=Q_{\mathfrak{U}}/ \sim .
\end{equation}
(I use the index $\mathfrak{U}$ so as to make clear the dependence on the chose cover of $M$.)

Define now two maps from $P_{\mathfrak{U}}$ to the manifold of objects $X_{\calG}$ of $\calG$ and to $M$ respectively:
\begin{equation}\label{eq-projmomentum}
\begin{cases}
\varepsilon_{\mathfrak{U}}\!\left(\left[\alpha,m_\alpha,g_\alpha\right]\right)&\colon=s_{\calG}(g_\alpha),\\
\pi_{\mathfrak{U}}\!\left(\left[\alpha,m_{\alpha},g_\alpha\right]\right)&\colon=m_\alpha,
\end{cases}
\end{equation}
where I used the notation $\left[\alpha,m_\alpha,g_\alpha\right]$ with square brackets for the equivalence class of the $3$-tuple $(\alpha,m_\alpha,g_\alpha)$ in $Q_{\mathfrak{U}}$.

First of all, one has to show that both maps are well-defined.
In fact, choosing some other representative $\left[\beta,m_\beta,g_\beta\right]$ for the class $\left[\alpha,m_\alpha,g_\alpha\right]$, then one would have obtained, by the definition of the equivalence relation (\ref{eq-glueing}):
\[
U_{\alpha\beta}\neq \emptyset,\quad m_\alpha=m_\beta,\quad g_\alpha=\Phi_{\alpha\beta}(m_\beta)g_\beta,
\]
whence it follows immediately that $\pi$ is well-defined.
On the other hand, by its very definition, $\varepsilon_\mathfrak{U}$ satisfies
\begin{align*}
\varepsilon_{\mathfrak{U}}\!\left(\left[\beta,m_\beta,g_\beta\right]\right)&=s_{\calG}(g_\beta)=\\
&=s_{\calG}\!\left(\Phi_{\beta\alpha}(m_\alpha)g_\alpha\right)=\\
&=s_{\calG}(g_\alpha)=\\
&=\varepsilon_{\mathfrak{U}}\!\left(\left[\alpha,m_\alpha,g_\alpha\right]\right).
\end{align*}
Hence, both maps are well-defined.
On the other hand, it is clear that $\pi_{\mathfrak{U}}$ is surjective: since $\mathfrak{U}$ is an open cover of $M$, for any point $m$ in $M$ one can choose an element $\alpha$, such that $m\in U_\alpha$.
Then, it is easy to verify that $\pi$ maps the equivalence class
\[
[\alpha,m,\iota_{\calG}(\varepsilon_\alpha(m))]
\]
onto $m$.

\begin{Thm}\label{thm-equivtrivdata}
Given local trivializing data $\left(\mathfrak{U},\varepsilon_\alpha,\Phi_{\alpha\beta}\right)$ over the manifold $M$ with values in the Lie groupoid $\calG$ in the sense of Definition~\ref{def-trivdata}, the $4$\-tuple
\[
\left(P_{\mathfrak{U}},\pi_{\mathfrak{U}},\varepsilon_{\mathfrak{U}},M\right),
\]
where the set $P_{\mathfrak{U}}$ is defined by Equation (\ref{eq-gluedprincbun}) and the maps $\pi_{\mathfrak{U}}$ and $\varepsilon_{\mathfrak{U}}$ are defined in Equation (\ref{eq-projmomentum}), is a principal bundle over $M$ with structure groupoid $\calG$ in the sense of Definition~\ref{def-princgroupoid}.
\end{Thm}
\begin{proof}
In order to show that the $4$-tuple $\left(P_{\mathfrak{U}},\pi_{\mathfrak{U}},\varepsilon_{\mathfrak{U}},M\right)$ defines a principal bundle in the sense of Definition~\ref{def-princgroupoid}, one has to show: $i)$ that $P$ is a smooth manifold, $ii)$ that the projection $\pi_{\mathfrak{U}}$ is a surjective, $\calG$-invariant submersion and $iii)$ that $\calG$ operates on $P_{\mathfrak{U}}$ freely and transitively on each fibre.

Let me first show $i)$. 
Since $P_{\mathfrak{U}}$ is obtained as a quotient of $Q_{\mathfrak{U}}$, notice first that $Q_{\mathfrak{U}}$ is a smooth manifold, as it is the disjoint union of smooth manifolds; moreover, the equivalence relation (\ref{eq-glueing}), by which one takes the quotient of $Q_{\mathfrak{U}}$, is defined by means of smooth maps.
Notice also that the natural mapping $Q_{\mathfrak{U}}\to P_{\mathfrak{U}}$ maps any trivial bundle $\varepsilon_\alpha^*\calU_\calG$ bijectively to $\pi_{\mathfrak{U}}^{-1}\!\left(U_\alpha\right)$: in fact, one can identify both sets via
\begin{equation}\label{eq-diffstr}
\varepsilon_\alpha^*\calU_\calG\ni(m_\alpha,g_\alpha)\mapsto \left[\alpha,m_\alpha,g_\alpha\right]\in \pi_{\mathfrak{U}}^{-1}(U_\alpha).
\end{equation}
Introduce at this point a differentiable structure on $P_{\mathfrak{U}}$ by requiring $a)$ the sets $\pi_{\mathfrak{U}}^{-1}(U_\alpha)$ to be open submanifolds of $P_{\mathfrak{U}}$ and that the maps (\ref{eq-diffstr}) are diffeomorphisms from the respective trivial bundles $\varepsilon_{\alpha}^*\calU_{\calG}$ to $\pi_{\mathfrak{U}}^{-1}(U_\alpha)$.
(Notice that any point $p\in P_{\mathfrak{U}}$ lies in some set $\pi_{\mathfrak{U}}^{-1}(U_\alpha)$, since $\pi_{\mathfrak{U}}$ is surjective and $\mathfrak{U}$ is an open cover of $M$.)
Hence, $P_{\mathfrak{U}}$ receives a smooth structure, which also makes the projection $\pi_{\mathfrak{U}}$ a smooth surjective submersion; in fact, this last fact follows directly from the equivalence relation (\ref{eq-glueing}).

It remains only to show that $\calG$ operates freely and transitively on each fibre via the momentum $\varepsilon_{\mathfrak{U}}$.
Define the right $\calG$-action $\Psi_{\mathfrak{U}}$ on $P_{\mathfrak{U}}$ by
\[
\left([\alpha,m_\alpha,g_\alpha];g\right)\overset{\Psi_{\mathfrak{U}}}\mapsto [\alpha,m_\alpha,g_\alpha g],\quad s_{\calG}(g_\alpha)=t_\calG(g).
\]
The action of $\calG$is well-defined, because
\[
[\alpha,m_\alpha,g_\alpha]=[\beta,m_\beta,g_\beta]\Leftrightarrow \begin{cases}
&U_{\alpha\beta}\neq\emptyset ,\\
&m_\alpha=m_\beta,\\
&g_\alpha=\Phi_{\alpha\beta}(m_\alpha)g_\beta,
\end{cases}
\]
whence it follows
\begin{align*}
\Psi_{\mathfrak{U}}\!\left([\beta,m_\beta,g_\beta];g\right)&=[\beta,m_\beta,g_\beta g]=\\
&=[\alpha,m_\alpha,\Phi_{\alpha\beta}(m_\alpha)(g_\beta g)]=\\
&=[\alpha,m_\alpha,g_\alpha g]=\\
&=\Psi_{\mathfrak{U}}\!\left([\alpha,m_\alpha,g_\alpha],g\right).
\end{align*}

Now, I want to show that the action $\Psi_{\mathfrak{U}}$ is free and transitive on every fibre.
First, assume that there is an element $g$ of $\calG$, such that, for some element $[\alpha,m_\alpha,g_\alpha]$ of $P_{\mathfrak{U}}$, the identity holds
\[
[\alpha,m_\alpha,g_\alpha]=[\alpha,m_\alpha,g_\alpha g]\Rightarrow g_\alpha =g_\alpha g\Rightarrow g=\iota_{\calG}(s_{\calG}(g_\alpha))=\iota_{\calG}\!\left([\varepsilon_{\mathfrak{U}}\!\left([\alpha,m_\alpha,g_\alpha]\right)\right),
\]
whence it follows that the action is free.

Assume now to have two points of $P_{\mathfrak{U}}$, say $[\alpha,m_\alpha,g_\alpha]$ and $[\beta,\widetilde{m}_\beta,\widetilde{g}_\beta]$, such that
\[
\pi_{\mathfrak{U}}\!\left([\alpha,m_\alpha,g_\alpha]\right)=m_\alpha=\widetilde{m}_\beta=\pi_{\mathfrak{U}}\!\left([\beta,\widetilde{m}_\beta,\widetilde{g}_\beta]\right).
\]
Hence, $U_{\alpha\beta}\neq\emptyset$, and therefore
\[
\widetilde{g}_\beta=\Phi_{\beta\alpha}(m_\alpha)\widetilde{g}_\alpha,
\]
for some $\widetilde{g}_\alpha$ such that the pair $(m_\alpha,\widetilde{g}_\alpha)$ belongs to $\varepsilon_\alpha^*\calU_\calG$; it follows
\[
[\beta,\widetilde{m}_\beta,\widetilde{g}_\beta]=[\alpha,m_\alpha,\widetilde{g}_\alpha].
\]
Therefore, the element $g_\alpha^{-1}\widetilde{g}_\alpha$ (which is well-defined by the properties of both factors) has the property of relating $[\alpha,m_\alpha,g_\alpha]$ and $[\beta,\widetilde{m}_\beta,\widetilde{g}_\beta]$ by right multiplication:
\begin{align*}
[\alpha,m_\alpha,g_\alpha]g_\alpha^{-1}\widetilde{g}_\alpha&=[\alpha,m_\alpha,g_\alpha g_\alpha^{-1}\widetilde{g}_\alpha ]=\\
&=[\alpha,m_\alpha,\widetilde{g}_\alpha ]=\\
&=[\alpha,m_\alpha,\Phi_{\beta\alpha}(m_\alpha)\widetilde{g}_\alpha\widetilde{g}_\alpha ]=\\
&=[\beta,\widetilde{m}_\beta,\widetilde{g}_\beta].
\end{align*}
Notice finally that, by the very definition of smooth structure on $P_{\mathfrak{U}}$, it follows that the momentum $\varepsilon_{\mathfrak{U}}$ of the $\calG$-action on $P_{\mathfrak{U}}$ is a smooth map and the right-action map $\Psi_{\mathfrak{U}}$ is also smooth.
\end{proof}

\begin{Cor}\label{cor-trivtrans}
Given local trivializing data $\left(\mathfrak{U},\varepsilon_\alpha,\Phi_{\alpha\beta}\right)$ over the manifold $M$ with values in the Lie groupoid $\calG$ in the sense of Definition~\ref{def-trivdata}, the principal bundle $P_{\mathfrak{U}}$, whose existence is guaranteed by Theorem~\ref{thm-equivtrivdata}, has over the open cover $\mathfrak{U}$ the functions $\Phi_{\alpha\beta}$ as transition functions.
\end{Cor}
\begin{proof}
By the very definition of smooth structure on $P_{\mathfrak{U}}$, the local trivialization $\varphi_\alpha$ over the open set $U_\alpha$ takes the form
\[
\varepsilon_\alpha^*\calU_\calG\ni (m_\alpha,g_\alpha)\overset{\varphi_\alpha^{-1}}\mapsto[\alpha,m_\alpha,g_\alpha]\in \pi_{\mathfrak{U}}^{-1}(U_\alpha).
\]
Therefore, taking two open sets $U_\alpha$ and $U_\beta$ intersecting nontrivially, one gets, by Lemma~\ref{lem-transmap}, the following isomorphism $\varphi_{\alpha\beta}$ between the restrictions $\left(\varepsilon_\alpha^*\calU_\calG\right)_{U_{\alpha\beta}}$ and $\left(\varepsilon_\beta^*\calU_\calG\right)_{U_{\alpha\beta}}$:
\begin{align*}
\varphi_{\alpha\beta}(m,g)&=\left(\varphi_\alpha\circ\varphi_\beta^{-1}\right)\!(m,g)=\\
&=\varphi_\alpha\!\left([\beta,m,g]\right)=\\
&=\varphi_\alpha\left([\alpha,m,\Phi_{\alpha\beta}(m)g]\right)=\\
&=(m,\Phi_{\alpha\beta}(m)g),\quad (m,g)\in \varepsilon_\beta^*\calU_\calG. 
\end{align*}
Hence, the claim follows.
\end{proof}

I have therefore proved the following fact:

\fbox{\parbox{12cm}{\bf There is an equivalence between principal bundles with groupoid structure in the sense of Definition~\ref{def-princgroupoid}, which can be trivialized over the open covering $\mathfrak{U}$, and local trivializing data in the sense of Definition~\ref{def-trivdata} w.r.t.\ open covering $\mathfrak{U}$.}}

\subsection{Examples of principal bundles with groupoid structure}\label{ssec-examples}
In this subsection I want to discuss some examples of principal bundles with groupoid structure, for particular Lie groupoids.
Let me start with the easiest example of {\em transitive Lie groupoid}, namely the product groupoid of a manifold $X$.

\subsubsection{Principal bundles with the product groupoid $X\times X$ of $X$ as groupoid structure}\label{sssec-productgroupoid}
The product groupoid $X\times X$ associated to a smooth manifold $X$ is defined as follows: 
\begin{itemize}
\item[i)] The product manifold $X\times X$ is the manifold of arrows of the product groupoid; 
\item[ii)] the manifold $X$ is the manifold of objects of the product groupoid; 
\item[iii)] the source map, resp.\ the target map, is defined as projection $\pr_1$ onto the second factor, resp.\ $\pr_2$ onto the first factor; the identity map is the diagonal map $\Delta_X$; 
\item[iv)] the product is simply defined as 
\[
(x,y)(y,z)\colon=(x,z).
\] 
\end{itemize}
It is easy to prove that $X\times X$ is a Lie groupoid.
I want now to describe precisely local trivializing data on $M$ with values in $X\times X$.
Following Definition~\ref{def-trivdata}, one needs $i)$ an open cover $\mathfrak{U}$ of $M$, $ii)$ smooth maps $\varepsilon_{\alpha}:U_{\alpha}\to X$ and $iii)$ smooth maps $\Phi_{\alpha\beta}:U_{\alpha\beta}\to X\times X$, for any nonempty intersection of any two open sets in $\mathfrak{U}$, satisfying additional identities, when composed with source and target map, and cocycle identities. 
The identities relating the ``cocycles'' $\Phi_{\alpha\beta}$ to the maps $\varepsilon_\alpha$ imply that there is only one possible cocycle, for any pair of open sets $U_\alpha$ and $U_\beta$, intersecting nontrivially: namely, such a map $\Phi_{\alpha\beta}$ takes the form
\[
\Phi_{\alpha\beta}(m)=\left(\Phi_{\alpha\beta}^1(m),\Phi_{\alpha\beta}^2(m)\right);
\]
recalling the definition of source and target map for the product groupoid, it follows immediately:
\[
\Phi_{\alpha\beta}^1=\pr_1\circ\Phi_{\alpha\beta}=\varepsilon_\alpha,
\]
and similarly
\[
\Phi_{\alpha\beta}^2=\varepsilon_\beta.
\]
It follows immediately that the maps $\Phi_{\alpha\beta}=(\varepsilon_\alpha,\varepsilon_\beta)$ satisfy the cocycle condition.

The trivial bundle associated to $\varepsilon_\alpha$ takes the form
\begin{align*}
\varepsilon_\alpha^*\calU_{X\times X}&=\left\{(m,(x,y))\in U_\alpha\times X\times X\colon \varepsilon_\alpha(m)=x\right\}\cong\\
&\cong U_\alpha\times X.
\end{align*}
The isomorphisms $\varphi_{\alpha\beta}$ take the form, using the above natural isomorphism between $\varepsilon_\alpha^*\calU_{X\times X}$ and the product manifold $U_\alpha\times X$, for any index $\alpha$:
\begin{align*}
U_{\alpha\beta}\times X\ni(m,x)&\mapsto(m,(\varepsilon_\beta(m),x))\overset{\varphi_{\alpha\beta}}\mapsto\\
&\overset{\varphi_{\alpha\beta}}\mapsto(m,(\varepsilon_\alpha(m),x))\mapsto\\
&\mapsto (m,x), 
\end{align*}
hence, one identifies in an obvious way the trivial bundles restricted to the intersections.
Thus, the equivalence relation (\ref{eq-glueing}) is the trivial relation, induced by the identity; therefore, the principal bundle $P_{\mathfrak{U}}$ equals the disjoint union $Q_{\mathfrak{U}}$, which turns out to be simply
\[
Q_{\mathfrak{U}}=\coprod_\alpha (U_\alpha\times X)=M\times X,
\]
since $\mathfrak{U}$ is an open cover of $M$.

\fbox{\parbox{12cm}{\bf Thus, for any open covering $\mathfrak{U}$ of $M$, there is only one principal bundle over $M$ with structure groupoid the product groupoid $X\times X$, namely $(M\times X,\pr_1,\pr_2,M)$.}}

\subsubsection{Principal bundles with a Lie group $G$ as structure groupoid}\label{sssec-liegroup}
It is already known that any Lie group $G$ can be made to a Lie groupoid, by putting $G$ itself as the manifold of arrows and a point $*$ as the manifold of objects; the source map, the target map and the identity map are defined respectively in a trivial way by putting
\[
s_G(g)=t_G(g)=*,\quad \forall g\in G;\quad \iota_G(*)=e,
\]
whereas the product is the usual product in $G$.

I want now to describe explicitly local trivializing data over $M$ with values in $G$, viewed as a Lie groupoid.
First of all, one needs, by Definition~\ref{def-trivdata}, an open cover $\mathfrak{U}$ of $M$ and maps $\varepsilon_\alpha$ from $U_\alpha$ to the manifold of objects of $G$; since this is simply a point, then all $\varepsilon_\alpha$ map the respective open set $U_\alpha$ onto the point $*$, and, as it was already remarked in Example~\ref{exa-trivial2} in Subsection~\ref{ssec-unittriv}, the trivial bundles $\varepsilon_\alpha^*\calU_G$ take the form $U_\alpha\times G$, which is {\em the} trivial bundle with structure group $G$ over $U_\alpha$.
Consider now the cocycles $\Phi_{\alpha\beta}$; they take their values in $G$, and satisfy additional identities relating them to the maps $\varepsilon_\alpha$, which are all equal, contracting the open sets to the point $*$.
Therefore, the identities relating the cocycles to the maps $\varepsilon_\alpha$ reduce simply to the equation
\[
\Phi_{\alpha\alpha}=e.
\]  
The cocycle condition becomes
\[
\Phi_{\alpha\beta}(m)\Phi_{\beta\gamma}(m)=\Phi_{\alpha\gamma}(m)\in G,\quad \forall m\in U_{\alpha\beta\gamma}\neq\emptyset,
\]
which is the ordinary cocycle condition for {\v C}ech cochains on $M$ w.r.t.\ the open cover $\mathfrak{U}$ with values in $G$.

Therefore, the principal bundle over $M$ with structure groupoid $G$ associated to the local trivializing data $\left(\mathfrak{U},\varepsilon_\alpha,\Phi_{\alpha\beta}\right)$ takes the form
\[
P_{\mathfrak{U}}=\coprod_{\alpha}(U_\alpha\times G)/ \sim,
\]
where the equivalence relation $\sim$ is simply
\[
U_{\alpha\beta}\times G\ni(\alpha,m_\alpha,g_\alpha)\sim(\beta,m_\beta,g_\beta)\in U_{\alpha\beta}\times G\Leftrightarrow \begin{cases}
&U_{\alpha\beta}\neq \emptyset,\\
&m_\alpha=m_\beta,\\
&g_\alpha=\Phi_{\alpha\beta}(m_\alpha)g_\beta. 
\end{cases}
\]
Hence, by the theory of ordinary principal bundles, it follows

\fbox{\parbox{12cm}{\bf Principal bundles over $M$ with groupoid structure $G$, for $G$ a Lie group, trivialized over the open covering $\mathfrak{U}$, are ordinary principal bundles over $M$ with structure group $G$, trivialized over the same open covering.}} 

\subsubsection{Principal bundles with the action groupoid $G\ltimes X$ as structure groupoid}\label{sssec-actiongroupoid}
Now comes a more interesting example of Lie groupoid, namely the {\em action groupoid} associated to a Lie group $G$ and a manifold $X$, on which $G$ acts from the left; let me recall briefly its definition.
The action groupoid $G\ltimes X$ associated to $G$ and $X$ is completely defined by the following requirements:
\begin{itemize}
\item[i)] the product manifold $G\times X$ is the manifold of arrows; 
\item[ii)] the manifold $X$ is the manifold of objects; 
\item[iii)] source map, target map and identity map are defined respectively via
\[
s_{G\ltimes X}(g,x)\colon=x,\quad t_{G\ltimes X}(g,x)\colon=gx,\quad \iota_{G\ltimes X}(x)\colon=(e,x);
\] 
\item[iv)] the product is defined via
\[
(g_1,g_2x)(g_2,x)\colon=(g_1g_2,x),\quad x\in X,\quad g_1,g_2\in G.
\] 
\end{itemize} 
Now consider local trivializing data $\left(\mathfrak{U},\varepsilon_\alpha,\Phi_{\alpha\beta}\right)$ over $M$ with values in $G\ltimes X$.
First of all, there is a family of maps $\varepsilon_\alpha$ from the open set $U_\alpha$ of the cover $\mathfrak{U}$ to $X$, and {\v C}ech cochains $\Phi_{\alpha\beta}$ over nontrivial intersections $U_{\alpha\beta}$ of any two open sets $U_\alpha$ and $U_\beta$.
Notice first that every cochain $\Phi_{\alpha\beta}$ can be written as
\[
\Phi_{\alpha\beta}(m)=(\Phi_{\alpha\beta}^G(m),\Phi_{\alpha\beta}^X(m)),
\]
where 
\[
\Phi_{\alpha\beta}^G\colon U_{\alpha\beta}\to G,\quad \Phi_{\alpha\beta}^X\colon U_{\alpha\beta}\to X.
\]
It follows from Condition $a)$ of Definition~\ref{def-trivdata} that
\[
s_{G\ltimes X}\circ \Phi_{\alpha\beta}=\Phi_{\alpha\beta}^X=\varepsilon_\beta,\quad t_{G\ltimes X}\circ \Phi_{\alpha\beta}=\Phi_{\alpha\beta}\varepsilon_\beta=\varepsilon_\alpha.
\]
Hence, Condition $a)$ of Definition~\ref{def-trivdata} establishes a relationship between the maps $\varepsilon_\alpha$ and the $G$-component of the cochains $\Phi_{\alpha\beta}$:
\begin{equation}\label{eq-actioncochain1}
\varepsilon_\alpha(m)=\Phi_{\alpha\beta}^G(m)\varepsilon_{\beta}(m),\quad \forall m\in U_{\alpha\beta}\neq \emptyset.
\end{equation}

Examine now the cocycle condition for the cochains $\Phi_{\alpha\beta}$: take any three open subsets $U_\alpha$, $U_\beta$ and $U_\gamma$ of the cover $\mathfrak{U}$, such that their triple intersection $U_{\alpha\beta\gamma}$ is nontrivial, then:
\begin{multline*}
\Phi_{\alpha\gamma}(m)=\Phi_{\alpha\beta}(m)\Phi_{\beta\gamma}(m)\Longleftrightarrow \\
\Longleftrightarrow\left(\Phi_{\alpha\gamma}^G(m),\varepsilon_\gamma(m)\right)=\left(\Phi_{\alpha\beta}^G(m),\Phi_{\beta\gamma}^G(m)\varepsilon_\gamma(m)\right)\left(\Phi_{\beta\gamma}^G(m),\varepsilon_\gamma(m)\right),\quad m\in U_{\alpha\beta\gamma},
\end{multline*}
whence the cocycle condition for the $G$-valued cochains $\Phi_{\alpha\beta}^G$:
\begin{equation}\label{eq-actioncochain2}
\Phi_{\alpha\gamma}^G(m)=\Phi_{\alpha\beta}^G(m)\Phi_{\beta\gamma}^G(m),\quad \forall m\in U_{\alpha\beta\gamma}.
\end{equation}
Let me describe now the trivial bundle $\varepsilon_\alpha^*\calU_{G\ltimes X}$: by definition, a pair 
\[
(m,(g,x)),\quad m\in U_\alpha,\quad (g,x)\in G\ltimes M,
\]
belongs to the trivial bundle $\varepsilon_\alpha^*\calU_{G\ltimes X}$ if and only if
\[
\varepsilon_\alpha(m)=gx\Rightarrow x=g^{-1}\varepsilon_\alpha(m).
\] 
Hence, the trivial bundle $\varepsilon_\alpha^*\calU_{G\ltimes X}$ can be diffeomorphically identified with the trivial $G$-bundle $U_\alpha\times G$ via the map
\[
\varepsilon_\alpha^*\calU_{G\ltimes X}\ni(m,(g,g^{-1}\varepsilon_\alpha(m)))\overset{\theta_\alpha}\mapsto (m,g)\in U_\alpha\times G,
\]
with inverse given by
\[
U_\alpha\times G\ni (m,g)\overset{\theta_\alpha^{-1}}\mapsto (m,(g,g^{-1}\varepsilon_\alpha(m)))\in \varepsilon_\alpha^*\calU_{G\ltimes X}.
\]
Furthermore, the isomorphisms $\varphi_{\alpha\beta}$, induced by the cocycles $\Phi_{\alpha\beta}$, between the restrictions of the trivial bundles to the (nontrivial) intersections of the open sets they are defined over take the following form, when identifying the trivial bundles $\varepsilon_\alpha^*\calU_{G\ltimes X}$ with the trivial $G$-bundles $U_\alpha\times G$:
\begin{align*}
U_{\alpha\beta}\times G\ni (m,g)&\overset{\theta_{\beta}^{-1}}\mapsto (m,(g,g^{-1}\varepsilon_\beta(m)))\in \left(\varepsilon_\beta^*\calU_{G\ltimes X}\right)\overset{\varphi_{\alpha\beta}}\mapsto\\
&\mapsto \left(m,\Phi_{\alpha\beta}(m)\!\left((g,g^{-1}\varepsilon_\beta(m)\right)\right)=\\
&=\left(m,\left(\Phi_{\alpha\beta}^G(m),\varepsilon_\beta(m)\right)\!\left((g,g^{-1}\varepsilon_\beta(m)\right)\right)=\\
&=\left(m,\left(\Phi_{\alpha\beta}^G(m)g,g^{-1}\Phi_{\beta\alpha}(m)\varepsilon_\alpha(m)\right)\right)=\\
&=\left(m,\left(\Phi_{\alpha\beta}^G(m)g,\left(\Phi_{\alpha\beta}^G(m)g\right)^{-1}\varepsilon_{\alpha}(m)\right)\right)\overset{\theta_\alpha}\mapsto\\
&\mapsto \left(m,\Phi_{\alpha\beta}^G(m)g\right)\in U_{\alpha\beta}\times G;
\end{align*}
notice that I made use of Equation (\ref{eq-actioncochain1}).

\fbox{\parbox{12cm}{\bf Hence, the isomorphisms $\varphi_{\alpha\beta}$ reduce to the isomorphisms between trivial $G$-bundles induced by the $G$-valued cocycles $\Phi_{\alpha\beta}^G$, and thus they give rise to a principal $G$-bundle $P_{\mathfrak{U}}^G$ over $M$.}}

Since the manifold $X$ is acted on from the left by $G$, one can form the associated fibre-bundle $P_{\mathfrak{U}}^G\times_G X$ with typical fibre $X$, and Equation (\ref{eq-actioncochain1}) implies that the functions $\varepsilon_\alpha$ glue together to form a global section of $P_{\mathfrak{U}}^G\times_G X$.
To prove this, recall the construction of associated bundles: given a principal $G$-bundle $P$ over $M$ and a manifold $X$, acted on from the left by $G$, one can form the quotient space of the product manifold $P\times X$ by the right $G$-action
\[
\left((p,x),g\right)\mapsto (pg,g^{-1}x);
\]
the projection $\pi_X$ is given by
\[
\pi_X\!\left([p,x]\right)\colon=\pi(p),
\]
$\pi$ being the projection of $P$ onto $M$.
Given trivializations $\varphi_\alpha$ over open subsets $U_\alpha$, belonging to some open cover $\mathfrak{U}$ of $M$, trivializations $\phi_{\alpha}^X$ of the associated bundle $P_{\mathfrak{U}}^G\times_G X$ are then 
\[
\pi_X^{-1}(U_\alpha)\ni [p,x]\overset{\varphi_\alpha^X}\mapsto\left(\pi(p),\left(\pr_G\circ\varphi_\alpha\right)(p)x\right)\in U_\alpha\times X, 
\]
where $\pr_G$ denotes the projection from $U_\alpha\times G$ onto $G$; it is not difficult to verify that the map is well-defined and that it is invertible, with inverse explicitly given by
\[
U_\alpha\times X\ni (m,x)\mapsto [\varphi_\alpha^{-1}(m,e),x]\in \pi_X^{-1}(U_\alpha).
\]
It is also easy to verify that the transition maps w.r.t.\ the above trivializations are simply given by
\[
U_{\alpha\beta}\times X\ni (m,x)\overset{\varphi_{\alpha\beta}}\mapsto\left(m,\Phi_{\alpha\beta}(m)x\right)\in U_{\alpha\beta}\times X,
\]
$\Phi_{\alpha\beta}$ being the transition functions of the bundle $P$ associated to the trivializations $\varphi_\alpha$.

Given now a section $\eta$ of the associated bundle $P_{\mathfrak{U}}^G\times_G X$, since $\pi_X\circ\eta=\id_M$, consider the composite map
\[
\eta_\alpha\colon=\varphi_\alpha^X\circ \eta\colon U_\alpha\to \pi_X^{-1}(U_\alpha)\to U_\alpha\times X, 
\]
which takes the form
\[
\eta_\alpha(m)=(m,\varepsilon_\alpha(m)),
\]
for $\varepsilon_\alpha$ a smooth map from $U_\alpha$ to the fibre $X$.
Now, since 
\[
\varphi_\alpha^X=\varphi_\alpha^X\circ\left(\varphi_\beta^X\right)^{-1}\circ\varphi_\beta^X
\]
on the intersection $U_{\alpha\beta}$, it follows
\begin{align*}
\eta_\alpha&=\varphi_\alpha^X\circ \eta=\\
&=\left(\varphi_\alpha^X\circ\left(\varphi_\beta^X\right)^{-1}\right)\circ\varphi_\beta^X\circ\eta=\\
&=\varphi_{\alpha\beta}^X\circ \eta_\beta,
\end{align*}
and the last identity takes the form
\[
\varepsilon_\alpha(m)=\Phi_{\alpha\beta}(m)\varepsilon_\beta(m),\quad \forall m\in U_{\alpha\beta},
\]
which is exactly Equation (\ref{eq-actioncochain1}).

On the other hand, given a principal $G$-bundle $P$ with trivializations $\varphi_\alpha$ over an open cover $\mathfrak{U}$, a family of functions $\varepsilon_\alpha$ associated to $\mathfrak{U}$, with a behaviour as in Equation (\ref{eq-actioncochain1}), gives rise to a global section $\eta$ in the following way: any point $m$ in $M$ belongs to some $U_\alpha$, because $\mathfrak{U}$ is an open cover of $M$, then
\[
\eta(m)\colon=\left[\varphi_\alpha^{-1}(m,e),\varepsilon_\alpha(m)\right]\in \pi_X^{-1}(U_\alpha).
\]
The section is well-defined, because, picking up another trivialization $\varphi_\beta$ and function $\varepsilon_\beta$, then 
\begin{align*}
\left[\varphi_\beta^{-1}(m,e),\varepsilon_\beta(m)\right]&=\left[\varphi_\alpha^{-1}(m,\Phi_{\alpha\beta}(m)),\varepsilon_\beta(m)\right]=\\
&=\left[\varphi_{\alpha}^{-1}(m,e),\Phi_{\alpha\beta}(m)\varepsilon_\beta(m)\right]=\\
&=\left[\varphi_{\alpha}^{-1}(m,e),\varepsilon_\alpha(m)\right],\quad\forall m\in U_{\alpha\beta}.
\end{align*}
Thus, I have proved the following important fact:

\fbox{\parbox{12cm}{\bf A principal bundle $P$ over a manifold $M$ with structure groupoid the action groupoid $G\ltimes X$, for a Lie group $G$ and a manifold $X$ acted on from the left by $G$, trivializable over the open covering $\mathfrak{U}$, is equivalent to a principal $G$-bundle $\widehat{P}$ over $M$, trivializable over the same covering, and a global section $\eta$ of the associated bundle $\widehat{P}\times_G X$; such a bundle I will call an $X$-pointed principal $G$-bundle over $M$.}}  

Let me consider some particular examples.
Consider an $m$-dimensional manifold $M$, and consider the obvious representation of $G=GL(m)$, the general linear group of $\bbR^m$, on $\bbR^m$; this induces in turn representations on the dual of $\bbR^m$, on tensor powers of $\bbR^m$ and/or its dual, on exterior powers of $\bbR^m$ and/or its dual, etc\ldots
Then, it is easy to verify the following equivalences:
\begin{itemize}
\item[i)] $X=\bbR^m$; then, the choice of a global vector field on $M$ corresponds to a principal $G\ltimes X$-bundle over $M$.
So, in particular, paracompact manifolds admit always such bundles.
Specifically, the principal $GL(m)$-bundle is the frame bundle $GL(M)$ of $M$ and the associated bundle is $\tange\!M$, the tangent bundle of $M$.
\item[ii)] $X=\bigwedge^m(\bbR^m)^*$; then, the choice of an orientation of $M$ corresponds to a principal $G\ltimes X$-bundle over $M$. 
This causes, by classical arguments, a reduction of the structure group $GL(m)$ to $SL(m)$, and the principal bundle in question corresponds to the {\em bundle of oriented frames $SL(M)$} over $M$, and the associated bundle is the line bundle $\bigwedge^m\tange^*M$.
\item[iii)] $X=\mathsf{Sym}^2_{>0}(m)$, the space of positive-definite, symmetric bilinear forms on $\bbR^m$; then, the choice of a Riemannian structure on $M$ corresponds to a principal $G\ltimes X$-bundle over $M$.
Again, paracompact bundles admit always such bundles.
Again by classical results, the principal bundle in question is a $GL(m)$-bundle, which admits a reduction to the orthonormal group $O(m)$; thus, the bundle is simply $O(M)$, the {\em bundle of orthonormal frames} over $M$, and the associated bundle is a subbundle of $\bigotimes^2\tange^*M$.
\end{itemize}
From the previous examples, one sees that in fact principal bundles over action groupoids encode many informations about the geometry of the base manifold.

\subsubsection{Some explicit constructions of principal bundles with structure groupoid $\calG(P)$}\label{sssec-gaugegroup}
Recall first the main features of the gauge groupoid $\calG(P)$ associated to an ordinary principal $G$-bundle $P\overset{\pi}\to X$, for a smooth manifold $X$ and for a Lie group $G$ (notice that the definition of gauge groupoid makes also sense for a principal bundle with {\em structure groupoid}:
\begin{itemize}
\item[i)] The manifold of arrows of $\calG(P)$ is set to be the quotient manifold of the product $P\times P$ w.r.t.\ the diagonal action of $G$.
\item[ii)] The manifold of objects of $\calG(P)$ is set to be the base manifold $X$ of the bundle $P$. 
\item[iii)] The source map, resp.\ the target map, of $\calG(P)$ is set to be
\[
s_{\calG(P)}\!\left(\left[p_1,p_2\right]\right)\colon=\pi(p_2),\quad\text{resp.}\quad t_{\calG(P)}\!\left(\left[p_1,p_2\right]\right)\colon=\pi(p_1),
\] 
where $\pi$ denotes the projection of the bundle $P$.
Notice that, by their very construction, both maps are well-defined; moreover, since $\pi$ is a surjective submersion, both maps $s_{\calG(P)}$ and $t_{\calG(P)}$ are also surjective submersions.

For the construction of the unit map of $\calG(P)$, consider an open cover $\mathfrak{U}$ of $X$, such that $P$ is trivializable over any open set of $\mathfrak{U}$; thus, smooth sections $\sigma_\alpha$ of $P$ over any open set $U_\alpha$ can be constructed, and one can set
\[
\iota_{\calG(P)}(x)\colon=\left[\sigma_\alpha(x),\sigma_\alpha(x)\right],\quad x\in U_\alpha.
\]
It is not difficult to verify that the unit map is well-defined, i.e.\ it does not depend on the choice of the section.
\item[iv)] The product in $\calG(P)$ is constructed by means of the division map of $P$ as follows: consider two composable elements of $\calG(P)$, say $[p_1,p_2]$ and $[q_1,q_2]$, in the sense that
\[
s_{\calG(P)}\!\left([p_1,p_2]\right)=\pi(p_2)=\pi(q_1)=t_{\calG(P)}\!\left([q_1,q_2]\right).
\] 
Hence, by the very definition of the division map, it follows:
\[
p_2=q_1\phi_P(q_1,p_2),\quad\text{or}\quad q_1=p_2\phi_P(p_2,q_1).
\]
Thus, it makes sense to set
\[
[p_1,p_2]\ [q_1,q_2]\colon=\left[p_1\phi_P(p_2,q_1),q_2\right].
\]
\end{itemize}
It is not difficult to verify that the $5$-tuple $\left(\calG(P),X,s_{\calG(P)},t_{\calG(P)},\iota_{\calG(P)}\right)$ defines a Lie groupoid, which is called the {\em gauge groupoid of the principal bundle $P$}.

Before going into the details of the construction of principal bundles with groupoid structure $\calG(P)$, notice that the manifold of arrows of $\calG(P)$, endowed with the natural projection $t_{\calG(P)}$ onto $X$, may be given another interpretation in terms of associated bundles.
In fact, observe first that the total space $P$ of the (ordinary) principal $G$-bundle may be viewed as a left $G$-space, where $G$ acts smoothly; namely, set the left $G$-action to be
\[
(g,p)\mapsto \rho_P^L(g)p\colon=pg^{-1}.
\]
Hence, one can consider the product manifold $P\times P$ with the following right $G$-action:
\[
\left((p_1,p_2),g\right)\mapsto \left(p_1g,g^{-1}p_2\right)=\left(p_1g,p_2g\right).
\]
Taking the quotient of $P\times P$ by the diagonal action of $G$ and considering the map from the quotient space to $X$ given by 
\[
[p_1,p_2]\mapsto \pi(p_1)
\]
shows that the manifold of arrows of the gauge groupoid is the total space of the associated bundle $P\times_G P$, defined by considering $P$ as a left $G$-space as pointed out above; the target map is then simply the projection from the total space onto the base space $X$.

For later purposes, let me compute the transition functions of the associated bundle $P\times_GP$; thus, consider an open trivializing cover $\mathfrak{U}$ of $P$, with trivializations $\varphi_\alpha^P$ and associated transition functions $\Phi_{\alpha\beta}^P$.
In order to avoid cumbersome notations, denote the transition maps of $P\times_GP$ associated to the transition maps $\Phi_{\alpha\beta}^P$, resp.\ the trivializations of $P\times_G P$ associated to the trivializations $\varphi_\alpha^P$, by $\Phi_{\alpha\beta}^{\calG(P)}$, resp.\ by $\varphi_\alpha^{\calG(P)}$; it is well-known that both transitions maps are related to each other by the formula
\[
\Phi_{\alpha\beta}^{\calG(P)}=\rho_{P}^L\!\left(\Phi_{\alpha\beta}^P\right),
\]
whence it follows that the composite maps $\varphi_\alpha^{\calG(P)}\circ \left(\varphi_\beta^{\calG(P)}\right)^{-1}$ take the explicit form
\[
U_{\alpha\beta}\times P\ni(x,p)\overset{\varphi_\alpha^{\calG(P)}\circ \left(\varphi_\beta^{\calG(P)}\right)^{-1}}\longmapsto \left(x,\rho_{P}^L\!\left(\Phi_{\alpha\beta}^P\right)\!(p)\right)\colon=\left(x,p\Phi_{\beta\alpha}^P(x)\right).
\]
Last, recall (\cite{BGV} for details) that the set $\Gamma\!\left(X,P\times_GP\right)$ of sections of the associated bundle $P\times_GP$ is in one-to-one correspondence with the set $C^{\infty}(P,P)^G$ of $G$-equivariant maps from $P$ to $P$, where the $G$-equivariance has to be w.r.t.\ the right $G$-action of $G$ on $P$ and the left action $\rho_P^L$; this is equivalent to the following requirement
\[
\forall \tau\in C^{\infty}(P,P)^G,\quad \forall p\in P,\quad \forall g\in G,\quad \tau(pg)=\rho_P^L\!\left(g^{-1}\right)\tau(p)=\tau(p)g.
\]
Moreover, observe that the composite map $\pi\circ\tau$ from $P$ to $X$ is $G$-invariant and is obviously smooth.

Translating this in the language of groupoids, the set $\Gamma\!\left(X,P\times_GP\right)$ is in one-to-one correspondence with the {\em bisections of the gauge groupoid $\calG(P)$}.

\paragraph{{\bf Two principal bundles over $X$}}\

I construct in the following two principal bundles over $X$ with structure groupoid $\calG(P)$, associated to a particular choice for the local momenta $\varepsilon_\alpha$.
Consider an open cover $\mathfrak{U}$ of $X$, in such a way that the $G$-bundle $P$ is locally trivializable over any open set $U_\alpha$ in the cover $\mathfrak{U}$.
I borrow the notations for trivializations and transition maps of $P$ w.r.t.\ open cover $\mathfrak{U}$ from above; notice that the (local) section of $P$ canonically associated to $\varphi_\alpha^P$ will be denoted by $\sigma_\alpha^P$.
Since I make an extensive use of the division map, it is better to write the explicit relationship between the local sections $\sigma_\alpha^P$ and transition maps $\Phi_{\alpha\beta}^P$:
\[
\Phi_{\alpha\beta}^P(x)=\phi_P\!\left(\sigma_\alpha^P(x),\sigma_\beta^P(x)\right),\quad x\in U_{\alpha\beta}.
\] 

A natural choice for the local momenta $\varepsilon_\alpha$ is simply:
\[
\varepsilon_\alpha\colon=\iota_\alpha\colon U_\alpha\hookrightarrow X, 
\]
i.e.\ the natural inclusions of the open sets $U_\alpha$ into $X$.

Consider a general cocycle $\Phi_{\alpha\beta}$ with values in $\calG(P)$, associated to the local momenta $\varepsilon_\alpha$.
First of all, any component of the cocycle is a smooth map from $U_{\alpha\beta}\subset X$ to the total space of the associated bundle $P\times_GP$; moreover, it has to satisfy
\[
t_{\calG(P)}\circ \Phi_{\alpha\beta}=\iota_\alpha=\id,
\] 
whence it follows that $\Phi_{\alpha\beta}$ represents a smooth section of $P\times_G P$ over $U_{\alpha\beta}$.
By the arguments at the beginning of this subsection, $\Phi_{\alpha\beta}$ is determined by a $G$-equivariant map $\tau_{\alpha\beta}$ from $P$, restricted to $U_{\alpha\beta}$, to $P$. 
Moreover, since 
\[
s_{\calG(P)}\circ\Phi_{\alpha\beta}=\id,
\]
it follows that $\tau_{\alpha\beta}$ preserves the fibres of $P$, hence it is a smooth $G$-equivariant map from $P$ restricted to $U_{\alpha\beta}$ to itself, which is moreover known to be invertible.

Hence, a general cocycle $\Phi_{\alpha\beta}$ takes the form
\[
\Phi_{\alpha\beta}(x)=\left[\sigma_\alpha^P(x),\tau_{\beta\alpha}\!\left(\sigma_\alpha^P(x)\right)\right],
\]
where $\tau_{\beta\alpha}$ lies in the gauge group of $P$ restricted to $U_{\alpha\beta}$.

Now let me check the cocycle condition:
\begin{align*}
\Phi_{\alpha\beta}(x)\Phi_{\beta\gamma}(x)&=\left[\sigma_\alpha^P(x),\tau_{\beta\alpha}\!\left(\sigma_\alpha^P(x)\right)\right]\left[\sigma_\beta^P(x),\tau_{\gamma\beta}\!\left(\sigma_\beta^P(x)\right)\right]=\\
&=\left[\sigma_\alpha^P(x)\phi_P\!\left(\tau_{\beta\alpha}\!\left(\sigma_\alpha^P(x)\right),\sigma_\beta^P(x)\right),\tau_{\gamma\beta}\!\left(\sigma_\beta^P(x)\right)\right]=\\
&=\left[\sigma_\alpha^P(x),\tau_{\gamma\beta}\!\left(\sigma_\beta^P(x)\right)\phi_P\!\left(\sigma_\beta^P(x),\tau_{\beta\alpha}\!\left(\sigma_\alpha^P(x)\right)\right)\right]=\\
&=\left[\sigma_\alpha^P(x),\tau_{\gamma\beta}\!\left(\tau_{\beta\alpha}\!\left(\sigma_\alpha^P(x)\right)\right)\right]\overset{!}=\\
&\overset{!}=\left[\sigma_\alpha^P(x),\tau_{\gamma\alpha}\!\left(\sigma_\alpha^P(x)\right)\right].
\end{align*}
Hence, it follows that the local gauge transformations $\tau_{\alpha\beta}$ satisfy the cocycle condition
\[
\tau_{\gamma\beta}\circ\tau_{\beta\alpha}=\tau_{\gamma\alpha}.
\] 
Since gauge transformations on $P$ are in bijective correspondence with smooth $G$-equivariant maps from $P$ to $G$ ($G$ viewed as a $G$-space via conjugation), one has $G$-equivariant maps $\Phi_{\alpha\beta}^\tau$ from $P$ restricted to $U_{\alpha\beta}$ to $G$, which satisfy the cocycle condition
\[
\Phi_{\alpha\gamma}^\tau(p)=\Phi_{\beta\gamma}^\tau(p)\Phi_{\alpha\beta}^\tau(p),\quad \forall p\in P\vert_{U_{\alpha\beta\gamma}}.
\]
There are many possible candidates: consider e.g.\ the two following cocycles:
\begin{itemize}
\item[i)] 
\begin{equation}\label{eq-gaugecoc1}
\tau_{\beta\alpha}(p)\colon=\left(\left(\varphi_\beta^P\right)^{-1}\circ\varphi_\alpha^P\right)\!(p),\quad \forall p\in P\vert_{U_{\alpha\beta}}, 
\end{equation}
or
\item[ii)] 
\begin{equation}\label{eq-gaugecoc2}
\tau_{\alpha\beta}(p)=\tau(p),\quad\forall\alpha,\beta,\quad \forall p\in P.
\end{equation}
\end{itemize}
Notice that the second case corresponds to a global gauge transformation, i.e.\ a $G$-equivariant morphism from the bundle $P$ to itself.
In general, observe that the trivial bundles $\varepsilon_\alpha^*\calU_{\calG(P)}$ may be identified with the product manifolds $U_\alpha\times P$.
Namely, by its very definition, a general element of $\varepsilon_\alpha^*\calU_{\calG(P)}$ takes the form
\[
\left(x;\left[p_1,p_2\right]\right),\quad U_\alpha\ni x=\pi(p_1),
\]
whence
\[
p_1=\sigma_\alpha(x)\phi_P(\sigma_\alpha(x),p_1)\Rightarrow\left(x;\left[p_1,p_2\right]\right)=\left(x;\left[\sigma_\alpha(x),p_2\phi_P(p_1,\sigma_\alpha(x))\right]\right). 
\]
Therefore, it makes sense to define the (smooth) map from the trivial bundle $\varepsilon_\alpha^*\calU_{\calG(P)}$ to $U_\alpha\times P$:
\begin{equation}\label{eq-difftriv}
\varepsilon_\alpha^*\calU_{\calG(P)}\ni\left(x;\left[p_1,p_2\right]\right)\mapsto (x,p_2\phi_P(p_1,\sigma_\alpha(x)))\in U_\alpha\times P, 
\end{equation}
which is invertible, with (smooth) inverse given by
\[
U_\alpha\times P\ni (x,p)\mapsto \left(x;\left[\sigma_\alpha(x),p\right]\right).
\]
Consider now any two open subsets $U_\alpha$ and $U_\beta$ with nontrivial intersection $U_{\alpha\beta}$; the diffeomorphism $\varphi_{\alpha\beta}$, induced by a general cocycle $\Phi_{\alpha\beta}^\tau$, takes the following form on the trivial $G$-bundle $U_{a\beta}\times P$:
\begin{align*}
U_{\alpha\beta}\times P\ni (x,p)&\mapsto (x,\left[\sigma_\beta(x),p\right])\mapsto \\
&\mapsto (x,\Phi_{\alpha\beta}(x)[\sigma_\beta(x),p])=\\
&=(x,\left[\sigma_\alpha(x),\tau_{\beta\alpha}(\sigma_\alpha(x))\right][\sigma_\beta(x),p])=\\
&=(x,\left[\sigma_\alpha(x),\sigma_\alpha(x)\Phi_{\alpha\beta}^\tau(\sigma_\alpha(x))\right][\sigma_\beta(x),p])=\\
&=(x,\left[\sigma_\alpha(x)\phi_P\!\left(\sigma_\alpha(x)\Phi_{\alpha\beta}^\tau(\sigma_\alpha(x)),\sigma_\beta(x)\right),p\right])=\\
&=(x,\left[\sigma_\alpha(x),p\phi_P\!\left(\sigma_\beta(x),\sigma_\alpha(x)\Phi_{\alpha\beta}^\tau(\sigma_\alpha(x))\right)\right])=\\
&=(x,\left[\sigma_\alpha(x),p\Phi_{\beta\alpha}^P(x)\Phi_{\alpha\beta}^\tau(\sigma_\alpha(x))\right])\mapsto \\
&\mapsto\left(x,p\Phi_{\beta\alpha}^P(x)\Phi_{\alpha\beta}^\tau(\sigma_\alpha(x))\right)\in U_{\alpha\beta}\times P,
\end{align*}
It is not difficult to verify that the maps 
\[
U_{\alpha\beta}\ni x\mapsto \Phi_{\alpha\beta}^{P,\tau}(x)\colon=\left(\Phi_{\beta\alpha}^P(x)\Phi_{\alpha\beta}^\tau(\sigma_\alpha(x))\right)^{-1}\in G
\]
satisfy the cocycle identity, hence they define a smooth (ordinary) principal $G$-bundle over $X$, which I denote by $P_\tau$ (so as to make explicit the dependence on $P$ (through the transition maps $\Phi_{\alpha\beta}^P$ and the family of gauge transformations $\tau_{\beta\alpha}$), by the usual gluing procedure, and hence the total space of the principal $\calG(P)$-bundle may be identified with the associated bundle $P_\tau\times_G P$; the (right) momentum of the principal $\calG(P)$-bundle $P_\tau\times_G P$ is simply 
\[
[p_\tau,p]\mapsto \pi(p)\in X,
\]
and the right $\calG(P)$-action on $P_\tau\times_G P$ takes the explicit form
\[
P_\tau\times_GP\ni[p_\tau,p],[p_1,p_2]\in \calG(P),\quad \pi(p)=\pi(p_1)\mapsto\left[p_\tau,p_2\phi_P(p_1,p)\right].
\]
The preceding formula may be checked directly by twisting the action of $\calG(P)$ on the trivial bundles $\varepsilon_\alpha^*\calU_{\calG(P)}$ by the diffeomorphisms (\ref{eq-difftriv}).

If one considers now e.g.\ the cocycle (\ref{eq-gaugecoc1}), one sees that any $\Phi_{\alpha\beta}^\tau\circ\sigma_\alpha$ equals $\Phi_{\alpha\beta}^P$, hence the cocycle $\Phi_{\alpha\beta}^{P,\tau}$ reduces to the identity element $e$ of $G$, independently from the choice of indices $\alpha$, $\beta$.
The bundle $P_\tau$ is thus simply the trivial principal $G$-bundle $X\times G$, and accordingly the associated principal $\calG(P)$-bundle $P_\tau\times_GP$ is given by the $4$-tuple $\left(X\times P,\pr_1,\pi\circ\pr_2,X\right)$.
Consider now the cocycle (\ref{eq-gaugecoc2}) for $\tau=\id_P$, then it follows easily
\[
\Phi_{\alpha\beta}^{P,\tau}=\Phi_{\alpha\beta}^P,
\]
whence it follows that $P_\tau=P$, and it is immediate to verify that the associated principal $\calG(P)$-bundle is simply $P\times_GP$, which, by previous arguments, equals the unit bundle $\calU_{\calG(P)}$ itself.
For a general global gauge transformation $\tau$ (which obviously induces a family of local gauge transformations obeying the cocycle condition), one gets a bundle $P_\tau$ in the same isomorphism class of $P$, hence, the $\calG(P)$-principal bundle may be also canonically identified with the unit bundle $\calU_{\calG(P)}$.

\paragraph{{\bf Two principal bundles over $P$}}\

Consider the base manifold of the $\calG(P)$-bundle to be $M=P$, the total space of the principal bundle defining the gauge groupoid $\calG(P)$.
Consider an open cover $\mathfrak{U}$ of $X$, the base space of $P$, which is assumed for simplicity to be a trivializing cover; I borrow the notations for trivializations, transition maps and sections of $P$ w.r.t.\ $\mathfrak{U}$ from the preceding paragraph.
Then, take the open cover of $P$ with elements given by 
\[
V_\alpha\colon=\pi^{-1}\!\left(U_\alpha\right).
\]
It makes thus sense to consider the local momenta $\varepsilon_\alpha$
\begin{align*}
\varepsilon_\alpha \colon V_\alpha&\to X,\\
p&\mapsto \pi(p)\in U_\alpha\quad\forall \alpha.  
\end{align*}
Therefore, with this choice of cover and local momenta, a cocycle $\Phi_{\alpha\beta}$ on $P$ with values in $\calG(P)$ is is represented by maps from any nontrivial intersection $V_{\alpha\beta}$ (which is contained in the restriction to $U_{\alpha\beta}$ of $P$) to the total space of the associated bundle $P\times_GP$, satisfying
\[
t_{\calG(P)}\circ\Phi_{\alpha\beta}=\pi,\quad s_{\calG(P)}\circ\Phi_{\alpha\beta}=\pi.
\]
Notice that any representative of a general cocycle $\Phi_{\alpha\beta}$ may be viewed as a local section over $V_{\alpha\beta}$ of the pull-back w.r.t.\ $t_{\calG(P)}$ of the associated bundle $P\times_GP$ over $X$.

Therefore, a natural choice for a local section of $t_{\calG(P)}^*\!\left(P\times_GP\right)$ over $V_{\alpha\beta}$ would be the section induced by the cocycle represented by 
\begin{equation}\label{eq-pullgaugecoc}
V_{\alpha\beta}\ni p\mapsto \Phi_{\alpha\beta}^X(\pi(p))\in P\times_GP, 
\end{equation}
where $\Phi_{\alpha\beta}^X$ is a general cocycle for $\calG(P)$ over $X$.
In the preceding paragraph, given a trivialization of $P$ over $X$, the general shape of such cocycles was computed:
\[
\Phi_{\alpha\beta}^X(x)=\left[\sigma_\alpha^P(x),\tau_{\beta\alpha}\!\left(\sigma_\alpha^P(x)\right)\right],
\]
for a family of local gauge transformations satisfying the cocycle relation
\[
\tau_{\beta\gamma}\circ\tau_{\alpha\beta}=\tau_{\alpha\gamma}.
\]
As in the preceding paragraph, it is better to work with $G$-equivariant maps $\Phi_{\beta\alpha}^\tau$ with values in $G$ canonically associated to the gauge transformations $\tau_{\alpha\beta}$; as seen before, these also satisfy some cocycle condition.
Consider now e.g.\ the two cocycles (\ref{eq-gaugecoc1}) and (\ref{eq-gaugecoc2}) of the preceding paragraph, the induced cocycles (\ref{eq-pullgaugecoc}) take the respective shapes:
\begin{align}
\label{eq-pullgaugecoc1} \Phi_{\alpha\beta}^X(\pi(p))&=\left[\sigma_\alpha(\pi(p)),\sigma_\beta(\pi(p))\right]=\left[\sigma_\alpha(\pi(p)),\sigma_\alpha(\pi(p))\Phi_{\alpha\beta}^P(\pi(p))\right],\quad p\in V_{\alpha\beta},\\ 
\label{eq-pullgaugecoc2} \Phi_{\alpha\beta}^X(\pi(p))&=\left[\sigma_\alpha(\pi(p)),\sigma_\alpha(\pi(p))\right]=\iota_{\calG(P)}(\pi(p)),\quad p\in V_{\alpha\beta}.
\end{align}
Keeping these results in mind, one finds the following expression for a possible cocycle on $P$ with values in $\calG(P)$:
\begin{align*}
V_{\alpha\beta}\ni p\mapsto \left[\sigma_\alpha^P(\pi(p)),\sigma_\alpha^P(\pi(p))\Phi_{\alpha\beta}^\tau(\sigma_\alpha^P(\pi(p)))\right].
\end{align*}
Notice, additionally, that there are also diffeomorphisms between the trivial bundles $\varepsilon_\alpha^*\calU_{\calG(P)}$ and $V_\alpha\times P$: namely, a general element of a trivial bundle takes the explicit form 
\begin{equation*}
\begin{aligned}
&\varepsilon_\alpha^*\calU_{\calG(P)}\ni\left(p;\left[p_1,p_2\right]\right),\quad \pi(p)=\pi(p_1)\Rightarrow p_1=p\phi_P(p,p_1),\ \text{whence}\\
&\Rightarrow\left(p;\left[p,p_2\phi_P(p_1,p)\right]\right)\overset{\xi_\alpha}\mapsto \left(p,p_2\phi_P(p_1,p)\right)\in V_\alpha\times P,
\end{aligned}
\end{equation*}
its inverse being simply
\[
V_\alpha\times P\ni (p,q)\overset{\xi_\alpha^{-1}}\mapsto (p;\left[p,q\right])\in \varepsilon_\alpha^*\calU_{\calG(P)}.
\]
With this identification, let me compute the transition maps between the ordinary trivial bundles $V_{\alpha\beta}\times P$ associated to the cocycle $\Phi_{\alpha\beta}$ as above:
\begin{align*}
V_{\alpha\beta}\times P\ni (p,q)&\overset{\xi_\beta}\mapsto(p;[p,q])\in\varepsilon_\beta^*\calU_{\calG(P)}\mapsto\\
&\mapsto\left(p;\Phi_{\alpha\beta}(p)\left[p,q\right]\right)=\\
&=\left(p;\left[\sigma_\alpha^P(\pi(p)),\sigma_\alpha^P(\pi(p))\Phi_{\alpha\beta}^\tau(\sigma_\alpha(\pi(p)))\right]\left[p,q\right]\right)=\\
&=\left(p;\left[\sigma_\alpha(\pi(p))\phi_P\!\left(\sigma_\alpha^P(\pi(p))\Phi_{\alpha\beta}^\tau(\sigma_\alpha(\pi(p))),p\right),q\right]\right)=\\ 
&=\left(p;\left[\sigma_\alpha(\pi(p)),q\phi_P\!\left(p,\sigma_\alpha^P(\pi(p)))\right)\Phi_{\alpha\beta}^\tau(\sigma_\alpha(\pi(p)))\right]\right)=\\
&=\left(p;\left[p,q\phi_P\!\left(p,\sigma_\alpha^P(\pi(p))\right)\Phi_{\alpha\beta}^\tau(\sigma_\alpha(\pi(p)))\phi_P\!\left(\sigma_\alpha(\pi(p)),p\right)\right]\right)=\\
&=\left(p;\left[p,q\Phi_{\alpha\beta}^\tau(p)\right]\right)\overset{\xi_{\alpha}^{-1}}\mapsto\\
&\mapsto\left(p,q\Phi_{\alpha\beta}^\tau(p)\right)=\\
&=\left(p,\rho_P^L\!\left(\Phi_{\alpha\beta}^\tau(p)^{-1}\right)q\right)\in V_{\alpha\beta}\times P. 
\end{align*}
By the previous arguments, it follows easily that the maps
\[
V_{\alpha\beta}\ni p\mapsto \Phi_{\alpha\beta}^\tau(p)^{-1}\in G,
\]
define a cocycle in $G$ over $P$ w.r.t.\ the open cover $\mathfrak{V}$.

Thus, by the usual gluing procedure, a smooth principal bundle $P_\tau$ over $P$ arises, and consequently the principal $\calG(P)$-bundle over $P$ defined by the family of gauge transformations $\tau_{\alpha\beta}$ is simply the associated bundle $P_\tau\times_GP$, with momentum map given by
\[
[p_\tau,p]\mapsto \pi(p),
\]
and right $\calG(P)$-action given by
\begin{align*}
&P_\tau\times_G P\ni [p_\tau,p],\quad [p_1,p_2]\in \calG(P),\quad \pi(p)=\pi(p_1)\Rightarrow\\
&[p_\tau,p][p_1,p_2]\colon=\left[p_\tau,p_2\phi_P(p_1,p)\right]\in P_\tau\times_GP.
\end{align*}

In particular, consider the cocycle (\ref{eq-pullgaugecoc1}); then, the transition maps of the bundle $P_\tau$ are simply the pull-backs of the transition maps of the bundle $P$, meaning that $P_\tau$ is in this case the tautological bundle $\pi^*P$ over $P$, and therefore the associated bundle $P_\tau\times_GP$ is simply the pull-back of $P\times_GP$ w.r.t.\ $\pi$ and the principal $\calG(P)$-bundle is the pull-back of the unit bundle w.r.t.\ the target map of $\calG(P)$ itself.

On the other hand, consider the cocycle (\ref{eq-pullgaugecoc2}); then one sees immediately that the transition maps of the bundle $P_\tau$ equal the identity maps, therefore the principal $\calG(P)$-bundle specified by is in this case is given by the $4$-tuple $\left(P\times P, pr_1,\pi\circ\pr_2,P\right)$, with right $\calG(P)$-action given by
\begin{align*}
&P\times P\ni (p_1,p_2),\quad \left[\widetilde{p}_1,\widetilde{p}_2\right]\in\calG(P),\quad \pi(p_2)=\pi(\widetilde{p}_1)\mapsto\\
&\mapsto(p_1,p_2)\left[\widetilde{p}_1,\widetilde{p}_2\right]\colon=\left(p_1,\widetilde{p}_2\phi_P\!\left(\widetilde{p}_1,p_2\right)\right). 
\end{align*}

\subsection{Morphisms of principal bundles}\label{ssec-morprincbun}
In this subsection I discuss morphisms of principal bundles over the same base manifold $M$ and with the same structure groupoid from a local point of view; in fact, the global point of view was already discussed extensively in~\cite{CR1}, constructing a generalization of gauge transformations.
I borrow the main notations and definitions from~\cite{CR1}: so, given two principal bundles $P$ and $Q$ over $M$ with structure groupoid $\calG$, a morphism $\Sigma$ between them is a (right) $\calG$-equivariant map, preserving projections and momenta.

I consider an explicit morphism $\Sigma$ from $P$ to $Q$; in~\cite{CR1}, it was shown that $\Sigma$ is automatically invertible by a global argument, and I denote by $K_\Sigma$ the unique map from $P\odot Q$ to $\calG$, equivariant with respect to right $\calG^2$-action on the fibred product $P\odot Q$ and the (right) generalized conjugation on $\calG$; see~\cite{CR1} for a complete description of such maps, called {\em generalized gauge transformations}, and their relationship with (iso)morphisms between principal bundles over the same base space and with the same structure groupoid.
Consider an open covering $\mathfrak{U}$ of $M$, such that there are sections $\sigma_{\alpha}^1$, resp.\ $\sigma_{\alpha}^2$, of $P$, resp.\ $Q$, over $U_\alpha$; as usual, denote by $\varepsilon_{\alpha}^1$, resp.\ $\varepsilon_{\alpha}^2$, the local momenta $\varepsilon_1\circ\sigma_{\alpha}^1$, resp. $\varepsilon_2\circ\sigma_{\alpha}^2$, and by $\Phi_{\alpha\beta}^1$, resp.\ $\Phi_{\alpha\beta}^2$, the transition functions w.r.t.\ the trivializations of $P$ w.r.t.\ the local sections $\sigma_{\alpha}^1$, resp.\ of $Q$ w.r.t.\ the local sections $\sigma_{\alpha}^2$.
Define the following maps
\begin{align*}
\Sigma_\alpha\colon U_\alpha&\to \calG,\\
x&\mapsto K_\Sigma\!\left(\sigma_{\alpha}^1(x),\sigma_{\alpha}^2(x)\right),
\end{align*}
where $K_\Sigma$ is the generalized gauge transformation associated to $\Sigma$.
\begin{Prop}\label{prop-locisomor}
The local maps $\Sigma_\alpha$ enjoy the following properties:
\begin{itemize}
\item[i)] (Compatibility with local momenta)
\[
t_{\calG}\circ \Sigma_\alpha=\varepsilon_{\alpha}^2,\quad s_{\calG}\circ \Sigma_\alpha=\varepsilon_{\alpha}^1.
\]
\item[ii)] (Coboundary relation)
\[
\Sigma_\beta(x)=\Phi_{\beta\alpha}^2(x)\Sigma_\alpha(x)\Phi_{\alpha\beta}^1(x),\quad \forall x\in U_{\alpha\beta},
\] 
provided $U_\alpha$ and $U_\beta$ intersect nontrivially.
\end{itemize}
\end{Prop}
\begin{proof}
I refer to~\cite{CR1} for the properties of generalized gauge transformations.
Then, the compatibility with local momenta follows and the coboundary relation follow immediately from the (global) compatibility of generalized gauge transformations with (global) momenta and from the $\calG^2$-equivariance of generalized gauge transformations, recalling that
\[
\sigma_{\beta}^1(x)=\sigma_{\alpha}^1(x)\Phi_{\alpha\beta}^1(x),\quad \sigma_{\beta}^2(x)=\sigma_{\alpha}^2(x)\Phi_{\alpha\beta}^2(x),\quad\forall x\in U_{\alpha\beta}.
\]
\end{proof}
The very definition of $\Sigma_\alpha$ implies that 
\[
\Sigma_\alpha=(\varphi_\alpha^2)^{-1}\circ \Sigma\circ \varphi_\alpha^1,
\]
where $\varphi_\alpha^1$, resp.\ $\varphi_\alpha^2$, is the trivialization of $P$ w.r.t.\ the local section $\sigma_{\alpha}^1$, resp.\ of $Q$ w.r.t.\ the local section $\sigma_{\alpha}^2$.

More abstractly, it makes sense to introduce the following
\begin{Def}\label{def-nonabcobound}
Let $\left(\mathfrak{U},\varepsilon_{\alpha}^1,\Phi_{\alpha\beta}^1\right)$, resp.\ $\left(\mathfrak{U},\varepsilon_{\alpha}^2,\Phi_{\alpha\beta}^2\right)$ be two local trivializing data over $M$ with values in the groupoid $\calG$, with the same open covering $\mathfrak{U}$.

A {\em local morphism} $\Sigma$ from $\left(\mathfrak{U},\varepsilon_{\alpha}^1,\Phi_{\alpha\beta}^1\right)$ to $\left(\mathfrak{U},\varepsilon_{\alpha}^2,\Phi_{\alpha\beta}^2\right)$ is a family of maps $\Sigma_\alpha$ from $U_\alpha$ to $\calG$, such that the following requirements hold:
\begin{itemize}
\item[a)] $\Sigma_\alpha$ puts the local moments $\varepsilon_{\alpha}^1$ and $\varepsilon_{\alpha}^2$ in relationship as follows:
\begin{equation*}
t_{\calG}\circ \Sigma_\alpha=\varepsilon_{\alpha}^2,\quad s_{\calG}\circ \Sigma_\alpha=\varepsilon_{\alpha}^1.
\end{equation*}
\item[b)] For any two nontrivially intersecting open sets $U_\alpha$ and $U_\beta$, the {\em nonabelian {\v C}ech cocycles} $\Phi_{\alpha\beta}^1$ and $\Phi_{\alpha\beta}^2$ are cohomologous w.r.t.\ $\Sigma$:
\begin{equation*}
\Sigma_\beta(x)=\Phi_{\beta\alpha}^2(x)\Sigma_\alpha(x)\Phi_{\alpha\beta}^1(x),\quad \forall x\in U_{\alpha\beta}.
\end{equation*}
\end{itemize} 
\end{Def}
By the previous computations, and since by Theorem~\ref{thm-equivtrivdata} the local trivializing data $\left(\mathfrak{U},\varepsilon_{\alpha}^1,\Phi_{\alpha\beta}^1\right)$ and $\left(\mathfrak{U},\varepsilon_{\alpha}^2,\Phi_{\alpha\beta}^2\right)$ give rise to principal $\calG$-bundles $P=P_{\mathfrak{U}}$ and $Q=Q_{\mathfrak{U}}$ respectively, it is natural to argue that the the morphism $\Sigma$ should give rise to a morphism from $P$ to $Q$, which I denote also by $\Sigma$, as follows by 
\begin{Thm}\label{thm-loctoglobmor}
Given two local trivializing data over $M$ with values in the groupoid $\calG$ $\left(\mathfrak{U},\varepsilon_{\alpha}^1,\Phi_{\alpha\beta}^1\right)$, resp.\ $\left(\mathfrak{U},\varepsilon_{\alpha}^2,\Phi_{\alpha\beta}^2\right)$, with the same open covering $\mathfrak{U}$, and a morphism $\Sigma$ between them in the sense of Definition~\ref{def-nonabcobound}, there is a morphism $\Sigma$ from $P$ to $Q$, the principal $\calG$-bundles associated to $\left(\mathfrak{U},\varepsilon_{\alpha}^1,\Phi_{\alpha\beta}^1\right)$ and $\left(\mathfrak{U},\varepsilon_{\alpha}^2,\Phi_{\alpha\beta}^2\right)$ respectively.
\end{Thm}  
\begin{proof}
Recall the construction of $P$, resp.\ $Q$, from the proof of Theorem~\ref{thm-equivtrivdata}:
\[
P=\coprod_{\alpha}\varepsilon_{\alpha}^{1*}\calU_{\calG}/ \sim,\quad Q=\coprod_{\alpha}\varepsilon_{\alpha}^{2*}\calU_{\calG}/ \sim,
\]
where the equivalence relation is induced by $\Phi_{\alpha\beta}^1$, resp.\ $\Phi_{\alpha\beta}^2$.
Given a family $\Sigma_\alpha$, define the morphism $\Sigma$ by the following rule:
\[
\Sigma\!\left([\alpha,x,g]\right)\colon=\left[\alpha,x,\Sigma_\alpha(x)g\right],\quad x\in U_\alpha,\quad g\in\calG,\quad \varepsilon_{\alpha}^1(x)=t_{\calG}(g),
\]
where the right-hand side of the previous equation has to be understood in $Q$.

First of all, one has to check that $\Sigma$ is well-defined, i.e.\ that it does not depend on the choice of the representatives and that, for $(\alpha,x,g)$ in $\varepsilon_{\alpha}^{1*}\calU_{\calG}$, the triple $(\alpha,x,\Sigma_\alpha(x)g)$ (which is well-defined since $s_{\calG}\circ \Sigma_\alpha=\varepsilon_{\alpha}^1$) belongs to $\varepsilon_{\alpha}^{2*}\calU_{\calG}$.
The second statement follows immediately from the relation $t_{\calG}\circ\Sigma_\alpha=\varepsilon_{\alpha}^2$. 
As for the first statement, let
\[
[\alpha,x,g]=[\beta,y,g_1];
\]
then $x=y$ in $U_{\alpha\beta}$ and $g_1=\Phi_{\beta\alpha}^1(x)^g$, whence it follows
\begin{align*}
\left[\beta,y,\Sigma_\beta(y)g_1\right]&=\left[\beta,x,\Sigma_\beta(x)\Phi_{\beta\alpha}^1(x)g\right]=\\
&=\left[\beta,x,\Phi_{\beta\alpha}^2(x)\Sigma_\alpha(x)g\right]=\\
&=\left[\alpha,x,\Sigma_\alpha(x)g\right].
\end{align*}
That $\Sigma$ preserves projections and momenta and is $\calG$-equivariant, it follows immediately by the definition of the projections, momenta and right $\calG$-actions for $P$ and $Q$.
\end{proof}
Therefore, the combined results of previous computations and Theorem~\ref{thm-loctoglobmor} can be resumed as follows:

\fbox{\parbox{12cm}{\bf There is an equivalence between morphisms of principal bundles over the same base space $M$ and with the same structure groupoid $\calG$ in the sense of~\cite{Moer2},~\cite{CR1}, trivializable over the same open covering $\mathfrak{U}$, and local morphisms between local trivializing data over the same manifold $M$ with values in the groupoid $\calG$ in the sense of Definition~\ref{def-nonabcobound}.}}

\subsubsection{Examples of morphisms of principal bundles}\label{sssec-examorprincbun}
\paragraph{\bf Principal bundles with a Lie group $G$ as structure groupoid}\label{par-liegroup}
Given a Lie group $G$, viewed as a (trivial )Lie groupoid over a point $\left\{*\right\}$ and a manifold $M$, it was proved in Subsubsection~\ref{sssec-liegroup} that a principal bundle $P$ with structure groupoid $G$ over $M$ is the same as a principal bundle with structure group $G$ over $M$.
Now, let me consider local trivializing data $\left(\mathfrak{U},\varepsilon_\alpha^1,\Phi_{\alpha\beta}^1\right)$ and $\left(\mathfrak{U},\varepsilon_\alpha^2,\Phi_{\alpha\beta}^2\right)$ over $M$ with values in $G$; then, clearly all local momenta $\varepsilon_\alpha^1$ and $\varepsilon_\alpha^2$ are trivial.
Consider now a local morphism $\Sigma$ from $\left(\mathfrak{U},\varepsilon_\alpha^1,\Phi_{\alpha\beta}^1\right)$ to $\left(\mathfrak{U},\varepsilon_\alpha^2,\Phi_{\alpha\beta}^2\right)$ in the sense of Definition~\ref{def-nonabcobound}; the compatibility  between the maps $\Sigma_\alpha$ and the local momenta are hollow, since the local momenta and the target and source map of $G$ are trivial.
It remains to check the significance of the nonabelian cohomological condition, namely:
\[
\Sigma_\beta(x)=\Phi_{\beta\alpha}^2(x)\Sigma_\alpha(x)\Phi_{\alpha\beta}^2(x),\quad \forall x\in U_{\alpha\beta}.
\] 
By classical arguments of the theory of usual principal $G$-bundles, it follows that $\Sigma$ defines a morphism between usual $G$-bundles, whence

\fbox{\parbox{12cm}{\bf Local morphisms $\Sigma$ between local trivializing data over $M$ with values in the Lie group $G$, viewed as a trivial Lie groupoid, correspond to (iso)morphisms between usual principal $G$-bundles.}}

\paragraph{\bf Principal bundles with the action groupoid $G\ltimes X$ as structure groupoid}\label{par-actgroup}
Consider a Lie group $G$ acting from the left on the manifold $X$, and the associated Lie groupoid $G\ltimes X$, the action groupoid; consider furthermore a manifold $M$.
From the results of Subsubsection~\ref{sssec-actiongroupoid}, one knows already that a principal bundle $P$ over $M$ with structure groupoid $G\ltimes X$ correspond to an $X$-pointed principal $G$-bundles $P$ over $M$; the $X$-point of $P$ is a global section of the associated bundle $P\times_G X$.
Given now two local trivializing data $\left(\mathfrak{U},\varepsilon_\alpha^1,\Phi_{\alpha\beta}^1\right)$ and $\left(\mathfrak{U},\varepsilon_\alpha^2,\Phi_{\alpha\beta}^2\right)$ over $M$ with values in the action groupoid $G\ltimes X$, I consider a local morphism $\Sigma$ between them.
Decomposing the maps $\Phi_{\alpha\beta}^1$, $\Phi_{\alpha\beta}^2$ and $\Sigma_\alpha$ as 
\[
\Phi_{\alpha\beta}^i(m)\colon=\left(\Phi_{\alpha\beta}^{G,i}(m),\Phi_{\alpha\beta}^{X,i}(m)\right),\quad \Sigma_\alpha(m)=\left(\Sigma_\alpha^G(m),\Sigma_\alpha^X(m)\right),
\] 
with maps
\begin{align*}
\Phi_{\alpha\beta}^{G,i}\colon U_{\alpha\beta}&\to G,\quad \Phi_{\alpha\beta}^{X,i}\colon U_{\alpha\beta}\to X,\\
\Sigma_\alpha^G\colon U_\alpha&\to G,\quad \Sigma_\alpha^X\colon U_\alpha\to X,
\end{align*}
where, again by the computations done in Subsubsection~\ref{sssec-actiongroupoid}, $\Phi_{\alpha\beta}^{X,i}=\varepsilon_{\beta}^i$, then the compatibility condition with local momenta in Definition~\ref{def-nonabcobound} implies immediately
\begin{align*}
\left(s_{G\ltimes X}\circ\Sigma_\alpha\right)(m)&=s_{G\ltimes X}\!\left(\Sigma_\alpha^G(m),\Sigma_\alpha^X(m)\right)=\\
&=\Sigma_\alpha^X(m)=\\
&=\varepsilon_\alpha^1(m),\quad\forall m\in U_\alpha,
\end{align*}
and
\begin{align*}
\left(t_{G\ltimes X}\circ\Sigma_\alpha\right)(m)&=t_{G\ltimes X}\!\left(\Sigma_\alpha^G(m),\Sigma_\alpha^X(m)\right)=\\
&=\Sigma_\alpha^G(m)\varepsilon_\alpha^1(m)=\\
&=\varepsilon_\alpha^2(m),\quad\forall m\in U_\alpha.
\end{align*}
On the other hand, the nonabelian cohomology condition in Definition~\ref{def-nonabcobound} can be rewritten as follows:
\begin{align*}
\Sigma_\beta(m)&=\left(\Sigma_\beta^G(m),\varepsilon_\beta^1(m)\right)=\\
&=\left(\Phi_{\beta\alpha}^{G,2}(m),\varepsilon^2_\alpha(m)\right)\left(\Sigma_\alpha^G(m),\varepsilon_\alpha^1(m)\right)\left(\Phi_{\alpha\beta}^{G,1}(m),\varepsilon_\beta^1(m)\right)=\\
&=\left(\Phi_{\beta\alpha}^{G,2}(m)\Sigma_\alpha^G(m)\Phi_{\alpha\beta}^{G,1}(m),\varepsilon_\alpha^1(m)\right),\quad \forall m\in U_{\alpha\beta},
\end{align*}
where I used $\Sigma_\alpha^G(m)\varepsilon_\alpha^1(m)=\varepsilon_\alpha^2(m)$.
Hence, one gets the nonabelian cohomological condition with values in the Lie group $G$:
\[
\Sigma_\beta^G(m)=\Phi_{\beta\alpha}^{G,2}(m)\Sigma_\alpha^G(m)\Phi_{\alpha\beta}^{G,1}(m),
\]
which corresponds to an isomorphism $\Sigma$ of principal $G$-bundles.
On the other hand, the identity
\[
\Sigma_\alpha^G(m)\varepsilon_\alpha^1(m)=\varepsilon_\alpha^2(m),\quad \forall m\in U_{\alpha\beta},
\]
has a significance: in fact, the local momenta $\varepsilon_\alpha^i$ are the local realizations of global sections $\eta_i$ of the principal $G$-bundles $P_i$ associated to the local trivializing data $\left(\mathfrak{U},\varepsilon_\alpha^i,\Phi_{\alpha\beta}^i\right)$, the so-called $X$-points of $P_i$, and the previous identity states that the isomorphism $\Sigma$, which induces also an isomorphism between the associated bundles $P_i\times_G X$ by the rule
\[
[m,x]\overset{\Sigma}\mapsto\left[m,\Sigma_\alpha^G(m)x\right],\quad x\in U_\alpha, 
\]
(the isomorphism $\Sigma$ is well-defined, because of the local construction of $P_i\times_G X$ and of the cocycle condition), maps the global section $\eta_1$ to $\eta_2$; this is evident from the local construction of the associated bundles.

\fbox{\parbox{12cm}{\bf Local morphisms $\Sigma$ between local trivializing data over $M$ with values in the action groupoid $G\ltimes X$ correspond to (iso)morphisms between $X$-pointed principal $G$-bundles; since an $X$-point of such a bundle $P$ corresponds to a global section of the associated bundle $P\times_G X$, (iso)morphisms of $X$-pointed bundles must preserve $X$-points.}}

The last condition translates, for the examples considered at the end of Subsubsection~\ref{sssec-actiongroupoid}, to the fact that the isomorphism of the tangent bundle, of the top exterior power of the cotangent bundle and of (a subspace of) the symmetric, covariant power of degree $2$ of the cotangent bundle of $M$ preserve global vector fields, orientation forms and Riemannian metrics respectively.

\subsection{Cohomological interpretation of principal bundles with structure groupoid}\label{ssec-cohomint}
In this Subsection, I want to point out and discuss a cohomological interpretation of principal bundles with structure groupoids.
In fact, it is well-known that isomorphism classes of principal $G$-bundles over a manifold $M$, $G$ being a Lie group, are in one-to-one correspondence with the first {\em nonabelian {\v C}ech cohomology group $\coh^1\!\left(M,G\right)$ of $M$ with values in $G$}, developed and discussed first by Grothendieck; I refer to~\cite{Gr1},~\cite{Gr2} and~\cite{Bry} for a complete discussion of the subject.
Moreover, H{\"a}fliger~\cite{H} introduced a natural generalization of the first nonabelian {\v C}ech cohomology group for the case of a topological groupoid; to be more precise, he discussed the first nonabelian {\v C}ech cohomology group $\coh^1\!\left(M,\calS_\calG\right)$, for $\calS$ being a sheaf of topological groupoids.
Since Lie groupoids have a smooth structure, I will consider here a sheaf of groupoids $\calS$ over a manifold $M$ and will define a slightly different notion of $1$-cochains and $1$-cocycles on it from the one of H{\"a}fliger; in fact, the two definitions are completely equivalent, although H{\"a}fliger does not mention explicitly the presence of local momenta, which, in my setting, makes definitions and computations more elegant.
In fact, for any open covering $\mathfrak{U}$ of $M$ and for a sheaf of Lie groupoids $\calS$ over $M$, $1$-cochains over $M$ with values in $\calS$ are shown to receive a natural left action of the groupoid of $0$-cochains over $M$ with values in $\calS$, which descends to an action on $1$-cocycles: the latter permits to define in a very explicit way the notion of cohomologous $1$-cocycles as orbits of a groupoid space, and, via refinement arguments, to define nonabelian cohomology classes of degree $1$ over $M$ with values in the sheaf of Lie groupoids $\calS$.
Finally, a Lie groupoid $\calG$ gives rise in a natural way to a sheaf $\calS_\calG$ of Lie groupoids over $M$, and local trivializing data in the sense of Definition~\ref{def-trivdata} modulo local morphisms in the sense of Definition~\ref{def-nonabcobound} correspond uniquely to nonabelian {\v C}ech cohomology classes in $\coh^1\!\left(M,\calS_\calG\right)$, cancelling moreover the dependence on open coverings of $M$, proving that the first nonabelian {\v C}ech cohomology group $\coh^1\!\left(M,\calS_\calG\right)$ classifies isomorphism classes of principal bundles over $M$ with structure groupoid $\calG$.

\begin{Def}\label{def-preshgroup}
Consider a smooth manifold $M$; then, a {\em sheaf $\calS$ of Lie groupoids over $M$} is a collection of Lie groupoids $\calS(U)$, for $U$ any open subset of $M$, such that, whenever there are $U$, $V$ open subsets of $M$ such that $U\subseteq V$, there is a morphism $\rho_{V,U}$ of Lie groupoids in the sense of Definition~\ref{def-groupmorph} from $\calS(V)$ to $\calS(U)$.
Moreover, the following compatibility condition must hold: whenever $U$, $V$ and $W$ are three open subsets of $M$, such that $U\subseteq V\subseteq W$, then the following identity must hold:
\[
\rho_{W,U}=\rho_{V,U}\circ \rho_{W,V},\quad \rho_{U,U}=\id_{\calS(U)}.
\] 
Additionally, the following gluing condition has to be satisfied: assume $U$ is an open subset of $M$, and $\mathfrak{U}$ is an open cover of $U$, and assume there are local sections $s_\alpha\in \calS(U_\alpha)$, such that
\[
\rho_{U_\alpha,U_{\alpha\beta}}(s_\alpha)=\rho_{U_\beta,U_{\alpha\beta}}(s_\beta),\quad U_{\alpha\beta}\neq \emptyset;
\]
then, there is a {\bf unique} section $s$ in $\calS(U)$, such that 
\[
s_\alpha=\rho_{U,U_\alpha}(s).
\]
\end{Def}
\begin{Rem}
It follows immediately from Definition~\ref{def-preshgroup} that a sheaf of Lie groupoids can be viewed as a functor associating to any open subset $U$ of $M$ two smooth manifolds, $\calS(U)$ and $X_{\calS(U)}$, the manifold of arrows and of objects respectively, plus (relative) target, source, unit and inversion maps.
Moreover, the restriction maps $\rho_{U,V}$, for $U\subseteq V$ open subsets of $M$, consist of two maps, namely the restriction $\rho_{U,V}$ on the manifolds of arrows and the restriction $r_{U,V}$ on the manifolds of objects, both compatible via (relative) target, source and unit maps.
It follows easily that the assignments $U\to X_{\calG(U)}$ and $U\to\calS(U)$ define two distinct sheaves of smooth manifolds, with restriction maps $\rho_{U,V}$ and $r_{U,V}$ respectively, for $U\subseteq V$ open subsets of $M$, called the (sheaf of arrows and of objects respectively.
In accordance to the usual terminology, I chose to denote a the manifold of arrows of a Lie groupoid $\calG$ by the same symbol, deserving a distinct notation for the manifold of objects; hopefully, this will not cause confusion.
\end{Rem}
\begin{Rem}
In~\cite{Moer3}, M{\oe}rdijk gives a different notion of sheaf of Lie groupoids.
Namely, a sheaf of Lie groupoids $\calS$ over $M$ is a Lie groupoid $\calS$, with manifold of objects $\calS_0$, such that there are two {\'E}tale maps (i.e.\ local diffeomorphisms) $p$ and $p_0$ from $\calS$, resp.\ $\calS_0$, to $M$, which are compatible with the structure of Lie groupoid of $\calS$.
This notion of sheaf of Lie groupoids is equivalent to the one I propose; still, for computational reasons, I will stick to Definition~\ref{def-preshgroup}.
\end{Rem}

\begin{Exa}
\begin{itemize}
\item[i)] A sheaf of Lie groups is a sheaf of Lie groupoids, with trivial sheaf of objects. 
\item[ii)] Let $M$ be a smooth manifold and $\calG$ be a Lie groupoid and consider, for $U$ an open subset of $M$, the sets 
\begin{align*}
C^\infty\!\left(U,X_\calG\right)&\colon=\left\{f\colon U\to X_{\calG}\colon \text{$f$ smooth}\right\},\\
C^\infty\!\left(U,\calG\right)&\colon=\left\{F\colon U\to \calG\colon \text{$F$ smooth}\right\}
\end{align*}
and consider the following maps 
\begin{align*}
s_U(F)&\colon=s_{\calG}\circ F,\\
t_U(F)&\colon=t_{\calG}\circ F,\\
j_U(F)&\colon=j_{\calG}\circ F,\quad\forall F\in C^\infty\!\left(U,\calG\right),\\
\iota_U(f)&\colon=\iota_\calG\circ f,\quad \forall f\in C^\infty\!\left(U,X_\calG\right), 
\end{align*}
and restriction maps $\rho_{U,V}$ and $r_{U,V}$ given simply by
\begin{align*}
\rho_{U,V}(F)&\colon=F\arrowvert_U,\quad F\in C^\infty\!\left(V,\calG\right),\\
r_{U,V}(f)&\colon=f\arrowvert_U,\quad f\in C^\infty\!\left(V,X_\calG\right). 
\end{align*}
It is easy to verify that the assignments $U\to C^\infty\!\left(U,\calG\right)$ and $U\to C^\infty\!\left(U,X_\calG\right)$ define a sheaf of (infinite-dimensional) Lie groupoids over $M$ associated to $\calG$, which I denote by $\calS_\calG$, and which I call the {\em canonical sheaf of Lie groupoids associated to $\calG$}.
Alternatively, the canonical sheaf of Lie groupoids associated to $\calG$ may be defined as the sheafification of the presheaf of germs of smooth maps from $M$ to $\calG$ (this way of constructing the canonical sheaf associated to $\calG$ makes it evident the relationship to the Definition in~\cite{Moer3}). 
\end{itemize}
\end{Exa}

Consider now an open covering $\mathfrak{U}$ of $M$ and a sheaf $\calS$ of Lie groupoids over $M$.
\begin{Def}\label{def-0coch}
A {\em $0$-cochain $\underline{\Sigma}$ over $M$ with values in $\calS$} (w.r.t.\ the open covering $\mathfrak{U}$) is given by a collection of sections $\Sigma_\alpha$ in the manifolds of arrows $\calS(U_\alpha)$.
A $0$-cochain $\underline{\Sigma}$ over $M$ with values in $\calS$ is a {\em $0$-cocycle}, if the following identities hold:
\[
\rho_{U_\alpha,U_{\alpha\beta}}(\Sigma_\alpha)=\rho_{U_\beta,U_{\alpha\beta}}(\Sigma_\beta),\quad U_{\alpha\beta}\neq\emptyset.
\]
\end{Def}
It is clear that, since $\calS$ is a sheaf, that a $0$-cocycle $\underline{\Sigma}$ over $M$ with values in $\calS$ is the same as a global section $\Sigma$ of $\calS$, i.e.\ an element $\Sigma$ of $\calS(M)$.
\begin{Rem}
To give a $0$-cochain $\underline{\Sigma}$ over $M$ with values in $\calS$ is equivalent to giving a $0$-cochain $\underline{s}$ over $M$ with values in the sheaf of arrows, $\calS$, and two $0$-cochains $\underline{\varepsilon^1}$ and $\underline{\varepsilon^2}$ over $M$ with values in the sheaf of objects $X_\calS$, such that the following condition holds:
\[
s_{U_\alpha}(\Sigma_\alpha)=\varepsilon^1_\alpha,\quad t_{U_\alpha}(\Sigma_\alpha)=\varepsilon^2_\alpha,\quad\forall \alpha.
\]
On the other hand, a $0$-cochain $\underline{\Sigma}$ over $M$ with values in $\calS$ gives rise to two $0$-cochains with values in the sheaf of objects, namely $t_{U_\alpha}(\Sigma_\alpha)$ and $s_{U_\alpha}(\Sigma_\alpha)$; these $0$-cochains are called {\em the target, resp.\ source, $0$-cochain of $\underline{\Sigma}$}, and are denoted by $\underline{t(\Sigma)}$, resp.\ $\underline{s(\Sigma)}$. 
\end{Rem}
\begin{Def}\label{def-1coch}
A {\em $1$-cochain $(\underline{\varepsilon},\underline{\Phi})$ over $M$ with values in $\calS$} (w.r.t.\ the open covering $\mathfrak{U}$) is given by $i)$ a $0$-cochain $\underline{\varepsilon}$ with values in the sheaf $X_{\calS}$ of objects of $\calS$ and $ii)$ a $1$-cochain over $M$ with values in the sheaf of arrows $\calS$, i.e.\ a collection of arrows $\Phi_{\alpha\beta}$ in $\calS(U_{\alpha\beta})$, such that the following additional condition holds:
\[
s_{U_{\alpha\beta}}(\Phi_{\alpha\beta})=r_{U_\beta,U_{\alpha\beta}}(\varepsilon_\beta),\quad t_{U_{\alpha\beta}}(\Phi_{\alpha\beta})=r_{U_\alpha,U_{\alpha\beta}}(\varepsilon_\alpha),\quad U_{\alpha\beta}\neq\emptyset.
\] 
A $1$-cochain $(\underline{\varepsilon},\underline{\Phi})$ over $M$ with values in $\calS$ is a {\em $1$-cocycle}, if it obeys the following identities:
\begin{align*}
\Phi_{\alpha\alpha}&=\iota_{U_\alpha}(\varepsilon_\alpha),\quad \forall \alpha;\\
\rho_{U_{\alpha\beta},U_{\alpha\beta\gamma}}(\Phi_{\alpha\beta})\rho_{U_{\beta\gamma},U_{\alpha\beta\gamma}}(\Phi_{\beta\gamma})&=\rho_{U_{\alpha\gamma},U_{\alpha\beta\gamma}}(\Phi_{\alpha\gamma}),\quad U_{\alpha\beta},U_{\beta\gamma},U_{\alpha\gamma},U_{\alpha\beta\gamma}\neq \emptyset.
\end{align*}
Notice that the left-hand side of the previous identity makes sense, because restriction maps are morphisms of Lie groupoids.
\end{Def}
\begin{Rem}
For the sake of simplicity, from now on I will drop from all formulae the restriction maps $\rho_{U,V}$ and $r_{U,V}$, if their presence is clear from the context.
Thus, the cocycle condition above will be written simply
\[
\Phi_{\alpha\beta}\Phi_{\beta\gamma}=\Phi_{\alpha\gamma}\in \calS(U_{\alpha\beta\gamma}).
\]
\end{Rem}
\begin{Exa}
Given a smooth manifold $M$, an open covering $\mathfrak{U}$ of $M$ and a Lie groupoid $\calG$, then a $1$-cocycle $\left(\underline{\varepsilon},\underline{\Phi}\right)$ over $M$ with values in $\calS_\calG$ is equivalent to local trivializing data over $M$ with values in $\calG$: in fact, the local momenta $\varepsilon_\alpha$ are the components of the $0$-cochain over $M$ with values in the sheaf of objects, and the transition functions $\Phi_{\alpha\beta}$ are the components of the $1$-cochain with values in the sheaf of arrows.
\end{Exa}

\begin{Rem}
Let me point out now the difference between Definition~\ref{def-1coch} of $1$-cochains and $1$-cocycles over $M$ with values in a sheaf of groupoids $\calS$ over $M$ and the one given by H{\"a}fliger in~\cite{H}.
What H{\"a}fliger called a $1$-cochain over $M$ with values in $\calS$ is simply a $1$-cochain over $M$ with values in the sheaf of arrows, which I denote by the same symbol as the sheaf of groupoids itself; analogously, a $1$-cocycle is a $1$-cochain satisfying the cocycle identities.
What is apparently missing is any information about the $0$-cochain with values in the sheaf of objects of $\calS$; I wrote ``apparently'', because the $0$-cochain is hidden in the ``diagonal sections'' $\Phi_{\alpha\alpha}$ over $U_\alpha$.
In Definition~\ref{def-1coch}, I assumed $\Phi_{\alpha\alpha}$ to be a unit of the groupoid $\calS(U_\alpha)$, i.e.\ to be the image w.r.t.\ $\iota_{U_\alpha}$ of a section $\varepsilon_\alpha$ of $X_{\calS(U_\alpha)}$, which can be obviously thought as a component of a $0$-cochain with values in the sheaf of objects of $\calS$, the so-called {\em unit $0$-cochain associated to the $1$-cochain $\underline{\Phi}$}; this is the idea that H{\"a}fliger had in mind, although he did not mention it explicitly. 
Nonetheless, I preferred to change the definition of H{\"a}fliger by mentioning explicitly the so-called unit $0$-cochain with values in $X_{\calG}$; this has the advantage of simplifying the notion of cohomologous cocycles, by the arguments that I will sketch in Remark~\ref{rem-cohomasleftspace}.
\end{Rem}

Assume now a $1$-cocycle $\left(\underline{\varepsilon},\underline{\Phi}\right)$ and a $0$-cochain $\underline{\Sigma}$ over $M$ with values in $\calS$ are given; assume additionally that the source $0$-cochain $\underline{s(\Sigma)}$ of $\underline{\Sigma}$ coincides with $\underline{\varepsilon}$.
\begin{Lem}\label{lem-cobound}
The $1$-cochain $\left(\underline{t(\Sigma)},\underline{\conj(\Sigma)\Phi}\right)$ over $M$ with values in $\calS$ is a $1$-cocycle, where
\[
\left(\conj(\Sigma)\Phi\right)_{\alpha\beta}\colon=\Sigma_\alpha \Phi_{\alpha\beta}\Sigma_\beta^{-1}\in \calS(U_{\alpha\beta}).
\]
\end{Lem}
\begin{proof}
First of all, one has
\begin{align*}
t_{U_{\alpha\beta}}\!\left(\left(\conj(\Sigma)\Phi\right)_{\alpha\beta}\right)&=t_{U_{\alpha\beta}}\!\left(\Sigma_\alpha \Phi_{\alpha\beta}\Sigma_\beta^{-1}\right)=\\
&=t_{U_{\alpha\beta}}(\Sigma_\alpha)=\\
&=r_{U_\alpha,U_{\alpha\beta}}(t_{U_\alpha}(\Sigma_\alpha)),
\end{align*}
where I used the fact that the restriction map is a morphism of Lie groupoids.
Similarly, one proves
\[
s_{U_{\alpha\beta}}\!\left(\left(\conj(\Sigma)\Phi\right)_{\alpha\beta}\right)=r_{U_\beta,U_{\alpha\beta}}(t_{U_\beta}(\Sigma_\beta)).
\]
Moreover, one gets
\begin{align*}
\left(\conj(\Sigma)\Phi\right)_{\alpha\alpha}&=\Sigma_\alpha\Phi_{\alpha\alpha}\Sigma_\alpha^{-1}=\\
&=\Sigma_\alpha \iota_{U_\alpha}(\varepsilon_\alpha) \Sigma_\alpha^{-1}=\\
&=\iota_{U_\alpha}(t_{U_\alpha}(\Sigma_\alpha)),
\end{align*}
using the properties of the unit map and the fact that $\varepsilon_\alpha=s_{U_\alpha}(\Sigma_\alpha)$.
Finally, the cocycle relation follows from the above computation:
\begin{align*}
\left(\conj(\Sigma)\Phi\right)_{\alpha\beta}\left(\conj(\Sigma)\Phi\right)_{\beta\gamma}=&\Sigma_\alpha \Phi_{\alpha\beta}\Sigma_\beta^{-1}\Sigma_\beta \Phi_{\beta\gamma}\Sigma_\gamma^{-1}=\\
&=\Sigma_\alpha \Phi_{\alpha\beta}\Phi_{\beta\gamma}\Sigma_\gamma^{-1}=\\
&=\Sigma_\alpha \Phi_{\alpha\gamma}\Sigma_\gamma^{-1}=\\
&=\left(\conj(\Sigma)\Phi\right)_{\alpha\gamma}.
\end{align*}
In the previous computations, I used implicitly the fact that restriction maps are morphisms of Lie groupoids, along with the properties of the unit map.
\end{proof}
\begin{Rem}\label{rem-cohomasleftspace}
Notice the particular notation I used for the $1$-cocycle $\underline{\conj(\Sigma)\Phi}$: in fact, the $0$-cochain $\underline{\Sigma}$ acts on the $1$-cocycle $\underline{\Phi}$ via generalized conjugation $\conj$.
On a more abstract level, the set $C^0\!\left(\mathfrak{U},\calS\right)$ of $0$-cochains over $M$ with values in the sheaf of groupoids $\calS$ (w.r.t.\ an open covering $\mathfrak{U}$) form an abstract groupoid, with the set of $0$-cochains over $M$ with values in the sheaf of objects of $\calS$ (w.r.t.\ the open covering $\mathfrak{U}$) as set of objects, the source and target map of this groupoid being
\[
s(\underline{\Sigma})\colon=\underline{s(\Sigma)},\quad t(\underline{\Sigma})\colon=\underline{t(\Sigma)},
\]
with unit map
\[
\iota(\underline{\varepsilon})\colon=\underline{\iota(\varepsilon)},\quad \iota(\varepsilon)_\alpha\colon=\iota_{U_\alpha}(\varepsilon_\alpha),\quad\forall \alpha,
\]
and inversion map
\[
j(\underline{\Sigma})\colon=\underline{j(\Sigma)},\quad j(\Sigma)_\alpha\colon=j_{U_\alpha}(\Sigma_\alpha),\quad\forall\alpha.
\]
Finally, the product of two composable $0$-cochains $\underline{\Sigma^1}$ and $\underline{\Sigma^2}$ is defined via
\[
\left(\Sigma^1\Sigma^2\right)_\alpha\colon=\Sigma^1_\alpha\Sigma_\alpha^2,\quad\forall \alpha.
\]
The set $C^1\!\left(\mathfrak{U},\calS\right)$ of $1$-cochains over $M$ with values in $\calS$ in the sense of Definition~\ref{def-1coch}, with momentum given by the projection onto $C^0\!\left(\mathfrak{U},X_{\calS}\right)$ is then a left $C^0\!\left(\mathfrak{U},\calS\right)$-space, with action given by
\[
\left(\underline{\Sigma},\left(\underline{\varepsilon},\underline{\Phi}\right)\right)\mapsto \left(t(\underline{\Sigma}),\underline{\conj(\Sigma)\Phi}\right),
\]
where $\underline{\varepsilon}=s(\underline{\Sigma})$, and 
\[
\left(\conj(\Sigma)\Phi\right)_{\alpha\beta}\colon=\Psi_{\conj}(\Sigma_\alpha,\Sigma_\beta;\Phi_{\alpha\beta}),
\]
where I used the notations of Subsection~\ref{ssec-defdivmap} for the (left) generalized conjugation; notice that generalized conjugation makes sense here, since $\left(\underline{\varepsilon},\underline{\Phi}\right)$ is a $1$-cocycle with values in $\calS$.
Notice also that I used an improper notation, hiding the restriction maps applied to the components of $\underline{\Sigma}$.
Lemma~\ref{lem-cobound} states that the previously described action of $C^0\!\left(\mathfrak{U},\calS\right)$ on $C^1\!\left(\mathfrak{U},\calS\right)$ descends to an action on $Z^1\!\left(\mathfrak{U},\calS\right)$ of $1$-cocycles over $M$ with values in $\calS$.
\end{Rem}

\begin{Def}\label{def-cohomology}
Two elements $\left(\underline{\varepsilon^1},\underline{\Phi^1}\right)$ and $\left(\underline{\varepsilon^2},\underline{\Phi^2}\right)$ of $Z^1\!\left(\mathfrak{U},\calS\right)$ are said to be {\em cohomologous}, if there exists an element $\underline{\Sigma}$ of $C^0\!\left(\mathfrak{U},\calS\right)$, such that
\[
\left(\underline{\varepsilon^2},\underline{\Phi^2}\right)=\underline{\Sigma}\left(\underline{\varepsilon^1},\underline{\Phi^1}\right).
\]
I will call the class of $\left(\underline{\varepsilon},\underline{\Phi}\right)$ in $Z^1\!\left(\mathfrak{U},\calS\right)$ by the previous equivalence relation the {\em cohomology class of $\left(\underline{\varepsilon},\underline{\Phi}\right)$}, and I will denote it by $\left[\underline{\varepsilon},\underline{\Phi}\right]$.
The set of all cohomology classes $\left(\underline{\varepsilon},\underline{\Phi}\right)$ of elements of $Z^1\!\left(\mathfrak{U},\calS\right)$ is called the {\em first nonabelian {\v C}ech cohomology group of $\calS$ w.r.t.\ the open covering $\mathfrak{U}$}, and is denoted by $\coh^1\!\left(\mathfrak{U},\calS\right)$.
\end{Def}
\begin{Rem}
It is clear that the first nonabelian {\v C}ech cohomology group of $\calS$ w.r.t.\ $\mathfrak{U}$ is the quotient space of $Z^1\!\left(\mathfrak{U},\calS\right)$ w.r.t.\ the left action of $C^0\!\left(\mathfrak{U},\calS\right)$ described in Remark~\ref{rem-cohomasleftspace}.
\end{Rem}

Consider a smooth manifold $M$ and a Lie groupoid $\calG$.
As I have proved previously, fixing an open covering $\mathfrak{U}$ of $M$, elements of $Z^1\!\left(\mathfrak{U},\calS_\calG\right)$ are in one-to-one correspondence with local trivializing data over $M$ with values in $\calG$; it is also clear that a local morphism $\Sigma$ between two such local trivializing data in the sense of Definition~\ref{def-nonabcobound} corresponds to the fact that the corresponding elements of $Z^1\!\left(\mathfrak{U},\calS_\calG\right)$ are cohomologous w.r.t.\ $\underline{\Sigma}$.
Hence, given an open covering $\mathfrak{U}$ of $M$, one gets the following result:

\fbox{\parbox{12cm}{\bf The first nonabelian {\v C}ech cohomology group $\coh^1\!\left(\mathfrak{U},\calS_\calG\right)$ over $M$ with values in the canonical sheaf of groupoids associated to $\calG$ is in one-to-one correspondence to isomorphism classes of principal bundles with structure groupoid $\calG$, trivialized w.r.t.\ the open covering $\mathfrak{U}$.}}

So far, I have discussed the first {\v C}ech nonabelian cohomology group over $M$ with values in a sheaf of Lie groupoids with a particular choice of an open covering $\mathfrak{U}$ of $M$.
But, as one can see directly from the above correspondence between this cohomology group and principal bundles with structure groupoid, there arises one problem, namely, one could consider different open coverings of $M$, leading to different local trivializations of such principal bundles, leading in turn to a priori different and unrelated first {\v C}ech nonabelian cohomology groups, which are but related to the same objects.
Solving this problem leads to the definition of the {\em first {\v C}ech nonabelian cohomology group $\coh^1\!\left(M,\calS\right)$}, for a sheaf of Lie groupoids $\calS$ over $M$.
\begin{Def}\label{def-refine}
Given two open coverings $\mathfrak{U}$ and $\mathfrak{V}$ of $M$, $\mathfrak{V}$ is said to be a {\em refinement} of $\mathfrak{U}$, if, denoting by $I$, resp.\ $J$, the set of indices of $\mathfrak{U}$, resp.\ $\mathfrak{V}$, there is a function $f$ from $J$ to $I$, such that
\[
V_j\subseteq U_{f(j)},\quad \forall j\in J.
\]
Equivalently, one says that the open covering $\mathfrak{V}$ is {\em finer} than $\mathfrak{U}$, and this property is denoted by $\mathfrak{U}<\mathfrak{V}$.
\end{Def}
In order to avoid  notational problems, let me switch momentarily from Greek indices for open coverings to Latin indices; in particular, given an open covering $\mathfrak{U}$ and a refinement $\mathfrak{V}$ thereof, the respective sets of indices will be denoted by $I$ and $J$, respectively.

Now, assume we are given an element $\underline{\Sigma}$ of $C^0\!\left(\mathfrak{U},\calS\right)$, for $\mathfrak{U}$ an open covering of $M$, and a refinement $\mathfrak{V}$ thereof, with associated function $f$.
The function defines a map $f^*$ as follows: given an element $\underline{\Sigma}$ of $C^0\!\left(\mathfrak{V},\calS\right)$, its image w.r.t.\ $f^*$, denoted by $f^*\underline{\Sigma}$, is defined as 
\[
(f^*\Sigma)_{j}\colon=\rho_{U_{f(j)},V_j}(\Sigma_{f(j)})\in\calS(V_j),\quad\forall j\in J.
\]
(The preceding operation of restriction makes sense, since $V_j\subseteq U_{f(j)}$, and $\Sigma_{f(j)}\in \calS(U_{f(j)})$.)
The map $f$ defines in a similar way a map $f^*$ on $C^0\!\left(\mathfrak{U},X_{\calS}\right)$ by the rule
\[
(f^*\varepsilon)_j\colon=r_{U_{f(j)},V_j}(\varepsilon_{f(j)})\in X_{\calS(V_j)},\quad\forall j\in J.
\]
Recall by Remark~\ref{rem-cohomasleftspace} that, for any open covering $\mathfrak{U}$ of $M$, $C^0\!\left(\mathfrak{U},\calS\right)$ has the structure of a groupoid over $C^0\!\left(\mathfrak{U},X_{\calS}\right)$.
Since restriction maps are morphisms of Lie groupoids, it follows easily that the source and target $0$-cochains of $\Sigma$ obey
\[
f^*(s(\underline{\Sigma}))=s(f^*(\underline{\Sigma})),\quad f^*(t(\underline{\Sigma}))=t(f^*(\underline{\Sigma}));
\]
moreover, given an element $\underline{\varepsilon}$ of $C^0\!\left(\mathfrak{U},X_{\calS}\right)$, it is easy to verify that
\[
f^*(\iota(\underline{\varepsilon}))=\iota(f^*(\underline{\varepsilon})),
\]
and $f^*$ commutes obviously with the inversion map.
Finally, it is clear that
\[
f^*\!\left(\underline{\Sigma^1}\underline{\Sigma^2}\right)=f^*(\underline{\Sigma^1})f^*(\underline{\Sigma^2}),
\] 
for any two composable $0$-cochains $\underline{\Sigma^1}$ and $\underline{\Sigma^2}$ in $C^0\!\left(\mathfrak{U},\calS\right)$.
All these computations can be summarized in the following
\begin{Lem}\label{lem-ref1}
Given an open covering $\mathfrak{U}$ of $M$ and a refinement $\mathfrak{V}$ thereof with associated function $f$, the map $f^*$ is a morphism of groupoids from $C^0\!\left(\mathfrak{U},\calS\right)$ to $C^0\!\left(\mathfrak{V},\calS\right)$.
\end{Lem}
The shape of this morphism depends explicitly on the choice of the map $f$; still, one could in principle consider the same refinement $\mathfrak{V}$ of $\mathfrak{U}$, but a different map, say $g:J\to I$, such that, for any $j\in J$, the inclusion holds 
\[
V_j\subseteq U_{g(j)}.
\]
Such a map $g$ determines a morphism $g^*$ from $C^0\!\left(\mathfrak{U},\calS\right)$ to $C^0\!\left(\mathfrak{V},\calS\right)$, which is a priori distinct from the one induced by $f$.
Still, considering the $0$-th nonabelian {\v C}ech cohomology group $\coh^0\!\left(\mathfrak{U},\calS\right)$ of $0$-cocycles over $M$ with values in $\calS$ (w.r.t.\ the open covering $M$), it is easy to prove that it is also a groupoid over the set of objects $Z^0\!\left(\mathfrak{U},X_\calS\right)$, consisting of all $0$-cocycles over $M$ with values in $X_{\calS}$, because restriction maps are morphisms of groupoids.
A choice of a map $f$ associated to $\mathfrak{U}<\mathfrak{V}$ defines a morphism $f^*$ of groupoids from $H^0\!\left(\mathfrak{U},\calS\right)$ to $H^0\!\left(\mathfrak{V},\calS\right)$. 
\begin{Lem}\label{lem-refin2}
Given an open covering $\mathfrak{U}$ of $M$ and a refinement $\mathfrak{V}$ thereof, and two maps $f,g\colon J\to I$, such that, for any $j\in J$, $V_j\subseteq U_{f(j)}$ and $V_j\subseteq U_{g(j)}$, then 
\[
f^*=g^*\quad\text{on $H^0\!\left(\mathfrak{U},\calS\right)$.}
\]
\end{Lem} 
\begin{proof}
Since, for any $j\in J$, $V_j$ is contained in $U_{f(j)}$ and $U_{g(j)}$, it follows $V_j\subseteq U_{f(j)g(j)}$.
Since any element $\underline{\Sigma}$ satisfies
\[
\rho_{U_{i_0},U_{i_0i_1}}(\Sigma_{i_0})=\rho_{U_{i_1},U_{i_0i_1}}(\Sigma_{i_1}),\quad U_{i_0i_1}\neq\emptyset,
\]
it follows, for any $j\in J$,
\[
\rho_{U_{f(j)},U_{f(j)g(j)}}(\Sigma_{f(j)})=\rho_{U_{g(j)},U_{f(j)g(j)}}(\Sigma_{g(j)}),
\]
whence the claim follows, by restricting further to $V_j$.
\end{proof}

Consider now an open covering $\mathfrak{U}$ of $M$ and a refinement $\mathfrak{V}$ thereof, with a map $f\colon J\to I$.
Then, $f$ determines a map $f^*$ from $C^1\!\left(\mathfrak{U},\calS\right)$ to $C^1\!\left(\mathfrak{V},\calS\right)$ as follows:
\[
f^*\!\left(\underline{\varepsilon},\underline{\Phi}\right)\colon=\left(f^*(\underline{\varepsilon}),f^*\!\left(\underline{\Phi}\right)\right),
\] 
where
\[
\left(f^*\Phi\right)_{j_0j_1}\colon=\rho_{U_{f(j_0)f(j_1)},V_{j_0j_1}}\!\left(\Phi_{f(j_0)f(j_1)}\right),\quad V_{j_0j_1}\neq\emptyset.
\]
(The previous definition makes sense, since $V_{j_k}\subseteq U_{f(j_k)}$, whence $V_{j_0j_1}\subseteq U_{f(j_0)f(j_1)}$.)
That $f^*$ is in fact a map from $C^1\!\left(\mathfrak{U},\calS\right)$ to $C^1\!\left(\mathfrak{V},\calS\right)$ follows again from the fact that restriction maps are morphisms of Lie groupoids.
By the very same reason, $f^*$ restricts to a map from $Z^1\!\left(\mathfrak{U},\calS\right)$ to $Z^1\!\left(\mathfrak{V},\calS\right)$.
I recall from Lemma~\ref{lem-cobound} that $Z^1\!\left(\mathfrak{U},\calS\right)$ is a left $C^0\!\left(\mathfrak{U},\calS\right)$-space.
\begin{Lem}\label{lem-refin3}
For an open covering $\mathfrak{U}$ of $M$ and a refinement $\mathfrak{V}$ thereof, with a map $f\colon J\to I$, $f^*$ defines a twisted equivariant map from the left $C^0\!\left(\mathfrak{U},\calS\right)$-space $Z^1\!\left(\mathfrak{U},\calS\right)$ to the left $C^0\!\left(\mathfrak{V},\calS\right)$-space $Z^1\!\left(\mathfrak{V},\calS\right)$, w.r.t.\ the morphisms of groupoids $f^*$.
\end{Lem}
\begin{proof}
It suffices to prove, for any element $\underline{\Sigma}$ of $C^0\!\left(\mathfrak{U},\calS\right)$ and $\left(\underline{\varepsilon},\underline{\Phi}\right)$ of $Z^1\!\left(\mathfrak{U},\calS\right)$, such that $s(\underline{\Sigma})=\underline{\varepsilon}$, that the following equality holds:
\[
f^*\!\left(\underline{\Sigma}\left(\underline{\varepsilon},\underline{\Phi}\right)\right)=f^*(\underline{\Sigma})f^*\left(\underline{\varepsilon},\underline{\Phi}\right).
\]
Writing down explicitly the left-hand side of the previous equation, one gets:
\begin{align*}
f^*\!\left(\underline{\Sigma}\left(\underline{\varepsilon},\underline{\Phi}\right)\right)&=f^*\!\left(\underline{t(\Sigma)},\underline{\conj(\Sigma)\Phi}\right)=\\
&=\left(f^*(t(\underline{\Sigma})),f^*\!\left(\underline{\conj(\Sigma)\Phi}\right)\right)
\end{align*}
By the very definition of $f^*$ acting on $Z^1\!\left(\mathfrak{U},\calS\right)$, and since $f^*$ is a morphism of groupoids, it follows immediately
\[
f^*(t(\underline{\Sigma}))=t(f^*(\underline{\Sigma})).
\]
On the other hand,  
\begin{align*}
f^*\!\left(\conj(\Sigma)\Phi\right)_{j_0j_1}&=\left(\conj(\Sigma)\Phi\right)_{f(j_0)f(j_1)}=\\
&=\Sigma_{f(j_0)}\Phi_{f(j_0)f(j_1)}\Sigma_{f(j_1)}^{-1}=\\
&=(f^*\Sigma)_{j_0}(f^*\Phi)_{j_0j_1}\left((f^*\Sigma)_{j_1}\right)^{-1}=\\
&=\left(\conj(f^*\Sigma)f^*\Phi\right)_{j_0j_1}.
\end{align*}
\end{proof}
Therefore, a refinement $\mathfrak{V}$ of any open covering $\mathfrak{U}$ of $M$, together with a map $f\colon J\to I$, defines a map between first nonabelian {\v C}ech cohomology groups with values in $\calS$ by
\[
f^*\!\left(\left[\underline{\varepsilon},\underline{\Phi}\right]\right)\colon=\left[f^*\left(\underline{\varepsilon},\underline{\Phi}\right)\right],
\] 
for any cohomology class $\left[\underline{\varepsilon},\underline{\Phi}\right]$ in $\coh^1\!\left(\mathfrak{U},\calS\right)$.

Let me come now to the main point of the construction.
\begin{Lem}\label{lem-refin4}
Consider, for an open covering $\mathfrak{U}$ of $M$ and a refinement $\mathfrak{V}$ thereof, two maps $f,g\colon J\to I$, such that $V_j\subseteq U_{f(j)}$ and $V_j\subseteq U_{g(j)}$ respectively, for any $j\in J$.
Then, the induced maps $f^*$ and $g^*$ on $\coh^1\!\left(\mathfrak{U},\calS\right)$ coincide.
\end{Lem}
\begin{proof}
Let me consider an element $\left(\underline{\varepsilon},\underline{\Phi}\right)$ of $Z^1\!\left(\mathfrak{U},\calS\right)$.
Then, I define the following element in $C^0\!\left(\mathfrak{V},\calS\right)$:
\[
\Sigma(f,g)_j\colon=\rho_{U_{f(j)g(j)},V_j}\!\left(\Phi_{f(j)g(j)}\right),\quad j\in J.
\]
The definition of $\underline{\Sigma(f,g)}$ makes sense, since $V_j$ is clearly contained in the intersection of $U_{f(j)}$ and $U_{g(j)}$.
I claim now that
\[
f^*\!\left(\underline{\varepsilon},\underline{\Phi}\right)=\underline{\Sigma(f,g)}\!\left(g^*\!\left(\underline{\varepsilon},\underline{\Phi}\right)\right).
\]
First of all, the right-hand side of the previous equation makes sense, because
\begin{align*}
s_{V_j}(\Sigma(f,g)_j)&=s_{V_j}\!\left(\rho_{U_{f(j)g(j)},V_j}\!\left(\Phi_{f(j)g(j)}\right)\right)=\\
&=r_{U_{f(j)g(j)},V_j}\!\left(s_{U_{f(j)g(j)}}\!\left(\Phi_{f(j)g(j)}\right)\right)=\\
&=r_{U_{f(j)g(j)},V_j}\!\left(r_{U_{g(j)},U_{f(j)g(j)}}(\varepsilon_{g(j)})\right)=\\
&=r_{U_{g(j)},V_j}(\varepsilon_{g(j)}).
\end{align*}
Then, let me compute $\underline{\conj(\Sigma(f,g))g^*\Phi}$:
\begin{align*}
\left(\conj(\Sigma(f,g))g^*\Phi\right)&=\Phi_{f(j_0)g(j_0)}\Phi_{g(j_0)g(j_1)}\Phi_{g(j_1)f(j_1)}=\\
&=\Phi_{f(j_0)g(j_1)}\Phi_{g(j_1)f(j_1)}=\\
&=\Phi_{f(j_0)f(j_1)}=\\
&=\left(f^*\Phi\right)_{j_0j_1}.
\end{align*}
Again, I used implicitly that restriction maps are morphisms of Lie groupoids, and that $\underline{\Phi}$ is a $1$-cocycle.
Thus, the image w.r.t.\ $f^*$ of the $1$-cocycle $\left(\underline{\varepsilon},\underline{\Phi}\right)$ is cohomologous to its image w.r.t.\ to $g^*$ in $Z^1\!\left(\mathfrak{V},\calS\right)$, thus their cohomology classes coincide.
\end{proof}

Lemma~\ref{lem-refin3} and~\ref{lem-refin4} imply together that, given an open covering $\mathfrak{U}$ of $M$, a refinement $\mathfrak{V}$ thereof, a map $f\colon J\to I$, such that $V_j\subseteq U_{f(j)}$, there is a map $f^*$ from $\coh^1\!\left(\mathfrak{U},\calS\right)$ to $\coh^1\!\left(\mathfrak{V},\calS\right)$, which does not depend on the choice of $f$.
Thus, the assignment $\mathfrak{U}\to \coh^1\!\left(\mathfrak{U},\calS\right)$, for any open covering $\mathfrak{U}$ of $M$, defines an {\em inductive system} over the partially ordered set of all open coverings of $M$, where the partial order is given by the refinement relation in Definition~\ref{def-refine}.
\begin{Def}\label{def-cohomol2}
Given a smooth manifold $M$ and a sheaf $\calS$ of Lie groupoids over $M$, the {\em first nonabelian {\v C}ech cohomology group of $M$ with values in $\calS$}, denoted by $\coh^1\!\left(M,\calS\right)$, is defined as the direct limit of the direct system $\mathfrak{U}\to \coh^1\!\left(\mathfrak{U},\calS\right)$ over the partially ordered set of all open coverings $\mathfrak{U}$ of $M$:
\[
\coh^1\!\left(M,\calS\right)\colon=\underset{\underset{\mathfrak{U}}\to}\lim\ \coh^1\!\left(\mathfrak{U},\calS\right).
\]
\end{Def}
As a consequence, consider a Lie groupoid $\calG$; then, the first non abelian cohomology group $\coh^1\!\left(M,\calS_\calG\right)$ of $M$ with values in $\calS_\calG$ parametrizes isomorphism classes of principal bundles over $M$ with structure groupoid $\calG$; it is customary to denote this cohomology group simply by $\coh^1\!\left(M,\calG\right)$.
E.g.\, if one considers a Lie group $G$ acting from the left on a manifold $X$, then
\[
\coh^1\!\left(M,G\ltimes X\right)\cong\coh^0\!\left(M,_G\! X\right)^{\coh^1\!\left(M,G\right)},
\]
where by $\coh^0\!\left(M,_G\! X\right)^{\coh^1\!\left(M,G\right)}$ I have denoted the map associating to any cohomology class in $\coh^1\!\left(M,G\right)$ the isomorphism class of a global section of the isomorphism class of a fiber bundle with typical fiber $X$ over $M$.

Given now a Lie groupoid $\calG$, a smooth manifold $M$ and a principal bundle $P$ over $M$ with structure groupoid $\calG$ as in Definition~\ref{def-princgroupoid}, consider the division map $\phi_P$ of $P$.
At this point, it should be clear why M{\oe}rdijk chose in~\cite{Moer1} for $\phi_P$ the name {\em cocycle of $P$ with values in $\calG$}: its invariance w.r.t.\ the action of morphisms of principal $\calG$-bundles over $M$, proved e.g.\ in~\cite{CR1} makes the division map an invariant of the isomorphism class of $P$ with values in $P$; moreover, composing it with the product of sections of $P$, gives the the natural $1$-cocycles over $M$ with values in $\calS_\calG$.
Thus, the division map can be thought of as a generator of the cohomology class of $P$ in $\coh^1\!\left(M,\calS_{\calG}\right)$, and, viceversa, a cohomology class in $\coh^1\!\left(M,\calS_{\calG}\right)$ gives rise to a unique division map; thus, division maps are a kind of global cocycles with values in $\calG$.

\section{Hilsum--Skandalis morphisms and a local version of generalized groupoids}\label{sec-hilskan}
A central notion in the theory of Lie groupoids is that of {\em generalized morphisms of Lie groupoids} or {\em Hilsum--Skandalis morphisms}; generalized morphisms, resp.\ {\em Morita equivalences} between Lie groupoids, mimic in the framework of Lie groupoids what ordinary generalized morphisms, resp.\ ordinary Morita equivalences, represent for algebra modules.
The former were first introduced by Connes~\cite{Con},~\cite{H} and were studied in great detail by Mr{\v c}un.
The main purpose here is to produce a different, ``local'' characterization of generalized morphisms, resp.\ Morita equivalence, in terms of cocycles as in Definition~\ref{def-trivdata} of Section~\ref{sec-localdata}.

Before entering into the details, let me briefly recall the notion of morphism of Lie groupoids.
\begin{Def}\label{def-groupmorph}
A morphism between two Lie groupoids $\calG$ and $\calH$ is a pair $\left(\Phi,\varphi\right)$, consisting of $i)$ a map $\Phi$ between the respective manifold of arrows $\calG$ and $\calH$ and $ii)$ a map between the respective manifolds of objects $X_{\calG}$ and $X_{\calH}$, such that the following conditions are satisfied:
\begin{itemize}
\item[a)] {\bf (Compatibility between the groupoid structures)} the following three diagrams must commute 
\begin{equation}\label{eq-commdiagmor}
\begin{CD} 
\calG  @>\Phi>> \calH\\
@Vs_\calG VV              @VV s_\calH V\\
X_\calG @>\varphi>> X_\calH
\end{CD}\quad,\quad 
\begin{CD} 
\calG  @> \Phi >> \calH\\
@V t_\calG VV              @VV t_\calH V\\
X_\calH          @>\varphi>> X_\calH
\end{CD}\quad\text{and}\quad
\begin{CD} 
X_\calG  @> \varphi >> X_\calH\\
@V \iota_\calG VV              @VV \iota_\calH V\\
\calG          @>\Phi>> \calH
\end{CD}\quad.
\end{equation}
\item[b)] {\bf (Homomorphism property)} For any composable pair $(g_1,g_2)$ of arrows, such that
\[
s_{\calG}(g_1)=t_{\calG}(g_2),
\]
the identity must hold
\begin{equation}\label{eq-homomor}
\Phi(g_1g_2)=\Phi(g_1)\Phi(g_2).
\end{equation} 
\end{itemize}
\end{Def}
In order to understand to significance of generalized morphisms, I need the following
\begin{Lem}\label{lem-genmorph1}
A morphism $\left(\Phi,\varphi\right)$ between the Lie groupoids $\calG$ and $\calH$ induces a principal bundle $\left(\varphi^*\calU_{\calH},\pr_1,s_{\calH}\circ\pr_2,X_\calG\right)$, such that there is a left $\calG$-action on this bundle w.r.t.\ the momentum $\pr_1$, which is compatible with the right $\calH$-action.
\end{Lem}
\begin{proof}
Consider the trivial bundle over the manifold of objects $X_{\calG}$ obtained by pulling the unit bundle of $\calH$ back w.r.t.\ the map $\varphi$; clearly, by previous arguments, the $4$-tuple $\left(\varphi^*\calU_{\calH},\pr_1,s_{\calH}\circ\pr_2,X_\calG\right)$ is a principal $\calH$-bundle over $X_{\calG}$.
Consider then the map $\pr_1$ from the trivial bundle to $X_{\calG}$ as a momentum for the following left action:
\begin{equation}\label{eq-morphaction}
\begin{aligned}
\calG\times_{\pr_1}\varphi^*\calU_{\calG}\ni\left(g,(x,h)\right),&\quad s_{\calG}(g)=x,\quad \varphi(x)=t_{\calH}(h)\mapsto\\
&\mapsto \left(t_{\calG}(g),\Phi(g)h\right). 
\end{aligned}
\end{equation}
First of all, notice that the action is well-defined, i.e.\ $Phi(g)$ and $h$ as above are composable, since
\[
s_{\calH}\!\left(\Phi(g)\right)=\varphi\!\left(s_{\calG}(g)\right)=\varphi(x)=t_{\calG}(h),
\]
and the result belongs to the trivial bundle $\varphi^*\calU_{\calH}$, since
\[
\varphi\!\left(t_{\calG}(g)\right)=t_{\calH}\!\left(\Phi(g)\right)=t_{\calH}\!\left(\Phi(h)g\right),
\]
by (\ref{eq-commdiagmor}).

One has now to show that (\ref{eq-morphaction}), together with the momentum $\pr_1$, defines a left $\calG$-action on $\varphi^*\calU_{\calH}$.
First of all, it is clear from (\ref{eq-morphaction}) that 
\[
\pr_1\!\left(g(x,h)\right)=t_{\calG}(g),\quad \forall (g,(x,h))\in\calG\times_{\pr_1}\varphi^*\calU_{\calG}. 
\]
Second, (\ref{eq-homomor}) implies directly that 
\[
(g_1g_2)(x,h)=\left(g_1\left(g_2(x,h)\right)\right),\quad \forall(g_2,(x,h))\in\calG\times_{\pr_1}\varphi^*\calU_{\calG},\quad s_{\calG}(g_1)=t_{\calG}(g_2), 
\]
and finally, by (\ref{eq-commdiagmor}) and (\ref{eq-homomor}), it follows immediately
\[
\iota_{\calG}(x)(x,h)=(x,h),\quad \forall (x,h)\in \varphi^*\calU_{\calH}.
\]
It remains to show that left $\calG$-action and the right $\calH$-action are compatible.
This is equivalent to showing $i)$ that the momentum for the right $\calH$-action is $\calG$-invariant (it is already known that the projection $\pr_1$ is $\calH$-invariant) and $ii)$ that the following identity holds
\begin{align*}
&\left(g(x,h_1)\right)h_2=g\left((x,h_1)h_2\right),\quad \forall g\in\calG,x\in X_{\calG},h_1,h_2\in\calH\ \text{such that}\\
&s_{\calG}(g)=x,\quad \varphi(x)=t_{\calH}(h_1),\quad s_{\calH}(h_1)=t_{\calH}(h_2).
\end{align*}
(Notice the $\calG$-invariance of the momentum for the right $\calH$-action makes the preceding expression well-defined.)
So, let me show $i)$:
\begin{align*}
\left(s_{\calH}\circ\pr_2\right)\!\left(g(x,h)\right)&=\left(s_{\calH}\circ\pr_2\right)\!(t_{\calG}(g),\Phi(g)h)=\\
&=s_{\calH}\!\left(\Phi(g)h\right)=\\
&=s_{\calH}(h)=\\
&=\left(s_{\calH}\circ\pr_2\right)\!(x,h),\quad \forall (g,(x,h))\in\calG\times_{\pr_1}\varphi^*\calU_{\calG}. 
\end{align*}
To show $ii)$, let me compute, by the associativity of the product structure in the groupoid $\calH$,
\begin{align*}
\left(g(x,h_1)\right)h_2&=\left(t_{\calG}(g),\Phi(g)h_1\right)h_2=\\
&=\left(t_{\calG}(g),\left(\Phi(g)h_1\right)h_2\right)=\\
&=\left(t_{\calG}(g),\Phi(g)\left(h_1h_2\right)\right)=\\
&=g(x,h_1h_2)=\\
&=g\left((x,h_1)h_2\right),\quad s_{\calG}(g)=x,\quad \varphi(x)=t_{\calG}(g).
\end{align*}
This proves the claim.
\end{proof}

\begin{Rem}
Consider the following map, which we denote by $\sigma_{\Phi}$
\[
X_{\calG}\ni x\overset{\sigma_{\Phi}}\mapsto(x,\iota_{\calH}(\varphi(x)))\in \varphi^*\calU_{\calH}. 
\]
The map is clearly well-defined, i.e.\ the image lies in fact in the trivial bundle $\varphi^*\calU_{\calH}$; it is also clearly a (global) section of the trivial bundle $\varphi^*\calU_{\calH}$.
Moreover, the composition of the section $\sigma_\Phi$ with the momentum for the right $\calH$-action equals simply $\varphi$.
A trivial, but remarkable property of the global section $\sigma_\Phi$ may be derived from the following computation:
\begin{align*}
g\sigma_{\Phi}\!\left(s_{\calG}(g)\right)&=g\left(s_{\calG}(g),\iota_{\calH}(\varphi(s_{\calG}(g)))\right)=\\
&=\left(t_{\calG}(g),\Phi(g)\iota_{\calH}(\varphi(s_{\calG}(g)))\right)=\\
&=\left(t_{\calG}(g),\Phi(g)\iota_{\calH}(s_{\calH}(\Phi(g)))\right)=\\
&=\left(t_{\calG}(g),\iota_{\calH}(t_{\calH}(\Phi(g)))\Phi(g)\right)=\\
&=\left(t_{\calG}(g),\iota_{\calH}(\varphi(t_{\calG}(g)))\Phi(g)\right)=\\
&=\sigma_{\Phi}(t_{\calG}(g))\Phi(g),\quad \forall g\in \calG,
\end{align*}
in other words, the left $\calG$-action on $\varphi^*\calU_{\calH}$ is intertwined by the global section $\sigma_{\Phi}$ with the right $\calH$-action via $\Phi$:
\begin{equation}\label{eq-intertw}
g\sigma_{\Phi}\!\left(s_{\calG}(g)\right)=\sigma_{\Phi}(t_{\calG}(g))\Phi(g),\quad \forall g\in \calG
\end{equation} 
A variation of Equation (\ref{eq-intertw}) will play a central r{\^o}le later in the discussion of local versions of Morita equivalences.
\end{Rem}

So, a morphism of two Lie groupoids $\calG$ and $\calH$ in the sense of Definition~\ref{def-groupmorph} defines automatically a right principal $\calH$-bundle on the manifold of objects of $\calG$, endowed with a left $\calG$-action, which is compatible with the principal bundle structure.
This motivates the following 
\begin{Def}\label{def-genmorph}[\cite{Mrcun},\cite{Moer1}]
Given two Lie groupoids $\calG$ and $\calH$, a {\em generalized morphism from $\calG$ to $\calH$ or a Hilsum--Skandalis morphism from $\calG$ to $\calH$} is a right principal $\calH$-bundle $\left(P,\pi,\varepsilon,X_{\calG}\right)$ over $X_{\calG}$, such that the following conditions are satisfied:
\begin{itemize}
\item[i)] the pair $\left(P,\pi\right)$ defines a left $\calG$-action on $P$ with momentum $\pi$;
\item[ii)] the momentum $\varepsilon$ for the right $\calH$-action is $\calG$-invariant, and moreover both actions are compatible in the sense that
\[
(gp)h=g(ph),\quad s_{\calG}(p)=\pi(p),\quad t_{\calH}(h)=\varepsilon(p).
\]
\end{itemize}
Notice that the $\calG$-invariance of the momentum $\varepsilon$ makes both sides of the compatibility condition between both actions well-defined.
\end{Def}

\begin{Rem}
By Lemma~\ref{lem-genmorph1}, it follows that every ``ordinary'' morphism between Lie groupoids gives rise to a generalized morphism between them in the sense of Definition~\ref{def-genmorph}.
\end{Rem}

Definition~\ref{def-genmorph} is an ``intrinsic'' definition of generalized morphisms: in fact, one does not see immediately the ``morphism'' hidden behind it.
However, in \cite{LTX}, the authors give at the beginning of the paper an equivalent characterization in terms of local versions of ``ordinary morphisms'', in the sense of Definition~\ref{def-groupmorph}.
More precisely, the decompose a generalized morphism from the Lie groupoid $\calG$ to the Lie groupoid $\calH$ into a Morita equivalence (see later for the precise definition) and morphisms defined w.r.t.\ a chosen open cover of the manifold of objects $X_{\calG}$ of the groupoid $\calG$.
This decomposition is a formal way of exposing the definition of Hilsum--Skandalis morphisms proposed in~\cite{Con},~\cite{HS}; again via the division map, I will soon come to this different, still equivalent, definition, although with some slight modifications.

In the next subsection, I will give a slightly different equivalent characterization of generalized morphism between Lie groupoids, based on the properties of local trivializing data in the sense of Definition~\ref{def-trivdata} of Subsection~\ref{ssec-trivdata}.

\subsection{From generalized morphisms to local generalized morphisms}\label{ssec-locgenmor}
Consider a generalized morphism $\left(P,\pi,\varepsilon,X_{\calG}\right)$ between the Lie groupoids $\calG$ and $\calH$ in the sense of Definition~\ref{def-genmorph}.
Lemma~\ref{lem-loctriv} of Subsection~\ref{ssec-localsec} enables us to find an open cover $\mathfrak{U}$ of the base manifold $X_{\calG}$ of the generalized morphism and corresponding local sections $\sigma_\alpha$ of the projection $\pi$, such that the restriction of the total space $P$ to any open set $U_\alpha$ of the cover $\mathfrak{U}$ is diffeomorphic to the trivial bundle $\varepsilon_\alpha^*\calU_{\calH}$, where $\varepsilon_\alpha$ is the local momentum w.r.t.\ the open set $U_\alpha$, defined by 
\[
\varepsilon_\alpha=\varepsilon\circ\sigma_\alpha.
\]
As was proved before, the local trivializations $\varphi_\alpha$, associated to the local sections $\sigma_\alpha$, give rise to the transition maps
\[
\Phi_{\alpha\beta}(x)=\phi_P\!\left(\sigma_\alpha(x),\sigma_\beta(x)\right),\quad \forall x\in U_{\alpha\beta}.
\]
Let me fix at this point some conventions: given two open sets $U_\alpha$ and $U_\beta$ of $X_{\calG}$, belonging to the trivializing open cover $\mathfrak{U}$ (not necessarily intersecting in a nontrivial way), define
\[
\calG_{U_\alpha,U_\beta}\colon=\left\{g\in\calG\colon s_{\calG}(g)\in U_{\alpha},\quad t_{\calG}(g)\in U_{\beta}\right\};
\]
I call $\calG_{\alpha,\beta}$ the {\em local $(\alpha,\beta)$-component of $\calG$}. 
\begin{Rem}
Given an open covering $\mathfrak{U}$ of $X_\calG$, the disjoint union of all local components of $\calG$ give rise to a groupoid, which is equivalent to $\calG$ itself; this can be found in~\cite{HS} and also in~\cite{Moer3}, although in a bit different shape.
Hilsum and Skandalis use the following notation
\[
\calG_{\alpha,\beta}=\calG_{U_\alpha}^{U_\beta},\quad \forall \alpha,\beta.
\]
\end{Rem} 
When the open cover is clear from the context, I will use the shorthand notation:
\[
\calG_{\alpha,\beta}\colon=\calG_{U_\alpha,U_\beta},\quad \calG_\alpha\colon=\calG_{U_\alpha,U_\alpha}.
\]
It is not difficult to prove the following
\begin{Lem}\label{lem-localgroupoid}
The $5$-tuple $\left(\calG_\alpha,U_\alpha,s_{\calG},t_{\calG},\iota_{\calG}\right)$ is a Lie subgroupoid of $\calG$, for any choice of the index $\alpha$. 
\end{Lem}
\begin{proof}
The topological spaces $\calG_\alpha$ and $U_\alpha$ are clearly smooth manifolds: in fact, $U_\alpha$ is an open subset of the smooth manifold $X_{\calG}$, and also $\calG_\alpha$, since
\[
\calG_\alpha=s_{\calG}^{-1}\!(U_\alpha)\cap t_{\calG}^{-1}\!(U_\alpha).
\]
It is easy to verify that the restrictions of source and target map $s_{\calG}$ and $t_{\calG}$ to $\calG_\alpha$ are still surjective submersions; moreover, the axioms for source, target and unit map of the Lie groupoid $\calG$ are clearly still valid for their respective restrictions on $\calG_\alpha$ and $U_\alpha$.

It remains only to prove that there is a well-defined associative multiplication in $\calG_\alpha$.
Consider two composable elements $g_1$ and $g_2$ in $\calG_\alpha$; this means that
\[
U_{\alpha}\ni s_{\calG}(g_1)=t_{\calG}(g_2)\in U_{\alpha}.
\] 
Then, consider the usual multiplication in $\calG$ as multiplication in $\calG_\alpha$:
\[
g_1g_2\in \calG.
\]
It is clear that $g_1g_2$ still belongs to $G_\alpha$, since
\[
s_{\calG}\!(g_1g_2)=s_{\calG}(g_2)\in U_\alpha,\quad t_{\calG}\!(g_1g_2)=t_{\calG}(g_1)\in U_\alpha.
\]
Associativity of the multiplication in $\calG_\alpha$ follows from the associativity of the multiplication in $\calG$, and, along the same lines, also the other axioms for the multiplication.
\end{proof}
Consider now the following map, for any open set $U_\alpha$:
\begin{equation}\label{eq-locgenmorph1}
\begin{aligned}
\Theta_\alpha\colon G_\alpha&\to \calH,\\
g&\mapsto \phi_P\!\left(\sigma_\alpha(t_{\calG}(g)),g\sigma_\alpha(s_{\calG}(g))\right),
\end{aligned}
\end{equation}
where $\sigma_\alpha$ is the local section of $\pi$ on $U_\alpha$.
I denoted by
\[
g\sigma_\alpha(s_{\calG}(g))
\]
the left $\calG$-action of $g$ on the element $\sigma_\alpha(s_{\calG}(g))$.

First, I need the following technical
\begin{Lem}\label{lem-wdeflocgenmor}
The map $\Theta_\alpha$ is well-defined.
\end{Lem} 
\begin{proof}
Let me check the following statements:
\begin{itemize}
\item[i)] the elements $\sigma_\alpha(t_{\calG}(g))$ and $\sigma_\alpha(s_{\calG}(g))$ are well-defined. 
\item[ii)] The left action of $g$ on $\sigma_\alpha(s_{\calG}(g))$ is well-defined. 
\item[iii)] The elements $\sigma_\alpha(t_{\calG}(g))$ and $g\sigma_\alpha(s_{\calG}(g))$ belong to the same fiber of $\pi$ 
\end{itemize}
The proof of $i)$ is trivial, since $g\in \calG_\alpha$, whence it follows
\[
s_{\calG}(g)\in U_\alpha,\quad t_{\calG}(g)\in U_\alpha.
\]
Part $ii)$ follows from the fact that $\sigma_\alpha$ is a section of $\pi$, whence
\[
\pi\!\left(\sigma_\alpha(s_{\calG}(g))\right)=s_{\calG}(\alpha),
\]
and from the fact that the left $\calG$-action has momentum $\pi$.

Part $iii)$ follows again because $\pi$ is the momentum of the left $\calG$-action: in fact,
\[
\pi\!\left(g\sigma_\alpha(s_{\calG}(g))\right)=t_{\calG}(g)=\pi\!\left(\sigma_\alpha(t_{\calG}(g))\right).
\] 
\end{proof}

The map $\Theta_\alpha$ is the explicit ``local morphism'' of~\cite{LTX}, as it is shown in the following
\begin{Prop}[Local morphism of Laurent-Gengoux, Tu and Xu]\label{prop-locmor1}
The pair $\left(\Theta_\alpha,\varepsilon_\alpha\right)$ is a morphism from the Lie groupoid $\calG_\alpha$ to the Lie groupoid $\calH$.
\end{Prop}
\begin{proof}
I have to show $i)$ the commutativity of the three diagrams (\ref{eq-commdiagmor}) and the homomorphism property (\ref{eq-homomor}).

Let me begin by showing the commutativity of the first diagram in (\ref{eq-commdiagmor}): this follows from Point $i)$ of Proposition~\ref{prop-prodivmap}, as the following computation shows
\begin{align*}
t_{\calH}\!\left(\Theta_\alpha(g)\right)&=t_{\calH}\!\left(\phi_P\!\left(\sigma_\alpha(t_{\calG}(g)),g\sigma_\alpha(s_{\calG}(g))\right)\right)=\\
&=\varepsilon\!\left(\sigma_\alpha(t_{\calG}(g)\right)=\\
&=\varepsilon_\alpha\!\left(t_{\calG}(g)\right).
\end{align*}
The commutativity of the second diagram follows again from Point $i)$ of Proposition~\ref{prop-prodivmap} and from the fact that the momentum $\varepsilon$ for the right $\calH$-action is $\calG$-invariant, so
\begin{align*}
s_{\calH}\!\left(\Theta_\alpha(g)\right)&=s_{\calH}\!\left(\phi_P\!\left(\sigma_\alpha(t_{\calG}(g)),g\sigma_\alpha(s_{\calG}(g))\right)\right)=\\
&=\varepsilon\!\left(g\sigma_\alpha(s_{\calG}(g)\right)=\\
&=\varepsilon\!\left(\sigma_\alpha(s_{\calG}(g)\right)=\\
&=\varepsilon_\alpha\!\left(t_{\calG}(g)\right).
\end{align*}
On the other hand, the commutativity of the third diagram follows from Point $ii)$ of Proposition~\ref{prop-prodivmap} and from the axioms for a left $\calG$-action; namely,
\begin{align*}
\Theta_\alpha(\iota_{\calG}(x))&=\phi_P\!\left(\sigma_\alpha(t_{\calG}(\iota_{\calG}(x))),\iota_{\calG}(x)\sigma_\alpha(s_{\calG}(\iota_{\calG}(x)))\right)=\\
&=\phi_P(\sigma_\alpha(x),\iota_\calG(x)\sigma_\alpha(x))=\\
&=\phi_P\!\left(\sigma_\alpha(x),\sigma_\alpha(x)\right)=\\
&=\iota_\calH\!\left(\varepsilon_\alpha(x)\right).
\end{align*}
It remains to show the homomorphism property (\ref{eq-homomor}); for this purpose, consider a composable pair $(g_1,g_2)$ in $\calG_\alpha\times \calG_\alpha$, i.e.
\[
s_{\calG}(g_1)=t_{\calG}(g_2).
\]
Then, one gets
\begin{align*}
\Theta_\alpha(g_1g_2)&=\phi_P\!\left(\sigma_\alpha(t_{\calG}(g_1g_2)),(g_1g_2)\sigma_{\alpha}(s_{\calG}(g_1g_2))\right)=\\
&=\phi_P\!\left(\sigma_\alpha(t_{\calG}(g_1)),g_1(g_2\sigma_{\alpha}(s_{\calG}(g_2)))\right).
\end{align*}
I want to prove now that $g_2\sigma_{\alpha}(s_{\calG}(g_2))$ belongs to the same fiber w.r.t.\ $pi$ as $\sigma_\alpha(s_{\calG}(g_1))$: in fact, a direct computation using the axioms of left $\calG$-action with momentum $\pi$ gives
\begin{align*}
\pi\!\left(g_2\sigma_\alpha(s_{\calG}(g_2))\right)&=t_{\calG}(g_2)=\\
&=s_{\calG}(g_1)=\\
&=\pi\!\left(\sigma_\alpha(s_{\calG}(g_1))\right),
\end{align*}
whence it follows
\begin{equation}\label{eq-genmortwist}
\begin{aligned}
g_2\sigma_\alpha(s_{\calG}(g_2))&=\sigma_{\alpha}(s_{\calG}(g_1))\phi_P\!\left(\sigma_{\alpha}(s_{\calG}(g_1)),g_2\sigma_\alpha(s_{\calG}(g_2))\right)=\\
&=\sigma_\alpha(s_{\calG}(g_1))\Theta_\alpha(g_2).
\end{aligned}
\end{equation}
Equation (\ref{eq-genmortwist}), which is the local analogon of Equation (\ref{eq-intertw}), implies now, together with Point $iv)$ of Proposition~\ref{prop-prodivmap}, that
\begin{align*}
\Theta_\alpha(g_1g_2)&=\phi_P\!\left(\sigma_\alpha(t_{\calG}(g_1)),g_1(g_2\sigma_{\alpha}(s_{\calG}(g_2)))\right)=\\
&=\phi_P\!\left(\sigma_\alpha(t_{\calG}(g_1)),g_1\sigma_\alpha(s_{\calG}(g_1))\Theta_\alpha(g_2)\right)=\\
&=\phi_P\!\left(\sigma_\alpha(t_{\calG}(g_1)),g_1\sigma_\alpha(s_{\calG}(g_1))\right)\Theta_\alpha(g_2)=\\
&=\Theta_\alpha(g_1)\Theta_\alpha(g_2).
\end{align*}
Hence, the claim follows.
\end{proof}
Let me point out the following transformation behaviour of the local morphisms $\Theta_\alpha$: for this purpose, assume to have two nontrivially intersecting open sets of the cover $\mathfrak{U}$, say $U_\alpha$ and $U_\beta$, and denote by $\Theta_\alpha$ and $\Theta_\beta$ the respective local morphisms.
It is immediate to see that
\[
\calG_\alpha\cap\calG_\beta=\calG_{\alpha\beta},
\]
which is, in virtue of Lemma~\ref{lem-localgroupoid}, a groupoid over the intersection $U_{\alpha\beta}$.

\begin{Cor}\label{cor-transflocmor}
For any two open subsets $U_\alpha$ and $U_\beta$, intersecting nontrivially, the corresponding local morphisms $\Theta_\alpha$ and $\Theta_\beta$ of Proposition~\ref{prop-locmor1} are related as follows:
\begin{equation}\label{eq-translocmor}
\Theta_\alpha(g)=\Phi_{\alpha\beta}\!(t_{\calG}(g))\Theta_\beta(g)\Phi_{\beta\alpha}\!(s_{\calG}(g)),\quad \forall g\in \calG_{\alpha\beta},
\end{equation}
where $\Phi_{\alpha\beta}$ denotes the transition map of $P$ associated to the local sections $\sigma_\alpha$ and $\sigma_\beta$.
\end{Cor}
\begin{proof}
By Equation (\ref{eq-locgenmorph1}), the claim follows immediately, since
\[
\Phi_{\alpha\beta}(x)=\phi_P\!\left(\sigma_\alpha(x),\sigma_\beta(x)\right).
\]
Moreover, the left $\calG$-action and the right $\calH$-action on $P$ are compatible; the claim then follows by the $\calH\times \calH$-equivariance of the division map.
\end{proof}
All the preceding computations are present in a different shape in Proposition 2.3 of~\cite{LTX}; but, in fact, the local morphisms $\Theta_\alpha$, arising directly from the ``bibundle'' structure of $P$, are only one piece of a more general construction, which I am now going to illustrate.

Consider the local component $\calG_{\alpha,\beta}$; it is an open submanifold of $\calG$, since it can be written in the following way:
\[
\calG_{\alpha,\beta}=s_{\calG}^{-1}(U_\alpha)\cap t_{\calG}^{-1}(U_\beta).
\]
Let me associate to the local component $\calG_{\alpha,\beta}$ the following map:
\begin{equation}\label{eq-genmor2}
\begin{aligned}
\Theta_{\beta\alpha}\colon \calG_{\alpha,\beta}&\to \calH,\\
g&\mapsto \phi_P\!\left(\sigma_\beta(t_{\calG}(g)),g\sigma_\alpha(s_{\calG}(g)\right).
\end{aligned}
\end{equation}
\begin{Rem}
The maps $\Theta_{\alpha\beta}$ were already introduced in~\cite{Con} and~\cite{HS} for the holonomy groupoids $\calG\colon=\calH_{V_1,F_1}$ and $\calH\colon=\calH_{V_2,F_2}$ of two foliated manifolds $(M_1,V_1)$ and $(M_2,V_2)$; the family of such maps was called by the authors a {\em cocycle over $\calG$ with values in $\calH$}.
The properties stated in the next Lemma are a slight generalization of the properties of the cocycle over $\calG$ with values in $\calH$ of~\cite{Con} and~\cite{HS}, in the sense that I consider an additional equation, which I will derive in the subsequent corollary (which follows along the same lines of Corollary~\ref{cor-transflocmor}).
Later, in a more implicit way, cocycles over $\calG$ with values in $\calH$ were considered also by M{\oe}rdijk in~\cite{Moer3}, although the author used a different notation in a slightly different framework. 
\end{Rem}
\begin{Rem}
For later purposes, let me just point out another interpretation of the local components $\calG_{\alpha,\beta}$: the product set $U_\alpha\times U_\beta$ is open in $X_{\calG}\times X_{\calG}$ w.r.t.\ the product topology.
The family of all such open sets is obviously an open covering of $X_{\calG}\times X_{\calG}$; then, $\calG_{\alpha,\beta}$ is the preimage of $U_\alpha\times U_\beta$ w.r.t.\ the smooth map $s_\calG\times t_\calG$ from $\calG$ to $X_\calG\times X_\calG$.
Since the $U_\alpha\times U_\beta$ cover $X_\calG\times X_\calG$, it is clear that the $\calG_{\alpha,\beta}$ form an open covering of $\calG$.
\end{Rem}

The same arguments used to prove Lemma~\ref{lem-wdeflocgenmor} lead to the following
\begin{Lem}\label{lem-wdeflocgenmor2}
The map $\Theta_{\beta\alpha}$ of Equation (\ref{eq-genmor2}) is well-defined, for any choice of indices $\alpha$, $\beta$.
\end{Lem}
Then, one has the following 
\begin{Thm}\label{thm-locgenmor}
The maps $\Theta_{\beta\alpha}$ enjoy the following properties:
\begin{itemize}
\item[a)] The following diagrams commute:
\begin{equation}\label{eq-commdiaglocmor}
\begin{CD} 
\calG_{\alpha,\beta}  @>\Theta_{\beta\alpha}>> \calH\\
@Vs_\calG VV              @VV s_\calH V\\
U_\alpha @>\varepsilon_\alpha>> X_\calH
\end{CD}\quad,\quad 
\begin{CD} 
\calG_{\alpha,\beta}  @> \Theta_{\beta\alpha} >> \calH\\
@V t_\calG VV              @VV t_\calH V\\
U_\beta          @>\varepsilon_\beta>> X_\calH\quad .
\end{CD}
\end{equation} 
\item[b)] For any three indices $\alpha$, $\beta$ and $\gamma$, and $g_1$ in $\calG_{\beta,\gamma}$ and $g_2\in\calG_{\alpha,\beta}$, such that
\[
s_{\calG}(g_1)=t_{\calG}(g_2),
\]
the following identity holds ({\bf generalized homomorphism property}):
\begin{equation}\label{eq-genmorprop}
\Theta_{\gamma\alpha}(g_1g_2)=\Theta_{\gamma\beta}(g_1)\Theta_{\beta\alpha}(g_2).
\end{equation}
Notice that both sides of Identity (\ref{eq-genmorprop}) are well-defined, since clearly, by the axioms of a groupoid, $g_1g_2$ belongs to $\calG_{\alpha,\gamma}$, and by the commutativity of the diagrams (\ref{eq-commdiaglocmor}).
\end{itemize}
\end{Thm}
\begin{proof}
\begin{itemize}
\item[a)] The commutativity of the diagrams (\ref{eq-commdiaglocmor}) follows using the same arguments displayed in the proof of Proposition~\ref{prop-locmor1} for showing the commutativity of the diagrams (\ref{eq-commdiagmor}).
\item[b)] First of all, since
\[
s_{\calG}(g_1)=t_{\calG}(g_2),\quad t_{\calG}(g_1g_2)=t_{\calG}(g_1),\quad s_{\calG}(g_1g_2)=s_{\calG}(g_2),
\]
one gets    
\begin{align*}
\pi\!\left(g_2\sigma_\alpha(s_{\calG}(g_2))\right)&=t_{\calG}(g_2)=\\
&=\pi\!\left(\sigma_\beta(t_{\calG}(g_2))\right),
\end{align*}
hence, one can write $g_2\sigma_\alpha(s_{\calG}(g_2))$ as follows
\begin{align*}
g_2\sigma_\alpha(s_{\calG}(g_2))&=\sigma_\beta(t_{\calG}(g_2))\phi_P\!\left(\sigma_\beta(t_{\calG}(g_2),g_2\sigma_\alpha(s_{\calG}(g_2))\right)=\\
&=\sigma_\beta(t_{\calG}(g_2))\Theta_{\beta\alpha}(g_2)=\\
&=\sigma_\beta(s_{\calG}(g_1))\Theta_{\beta\alpha}(g_2).
\end{align*}
Then, Identity (\ref{eq-genmorprop}) follows from the equivariance of the division map.
\end{itemize}
\end{proof}
The following Corollary is an immediate consequence of Equation (\ref{eq-genmor2}), recalling the Definition of the transition maps of $P$, and of the compatibility of both left $\calG$- and right $\calH$-action on $P$ and of the $\calH\times\calH$-equivariance of the division map:
\begin{Cor}\label{cor-transflocmor2}
Consider four open subsets of $X_{\calG}$, say $U_\alpha$, $U_\beta$, $U_\gamma$, $U_\delta$, such that
\[
U_{\alpha\gamma}\neq \emptyset,\quad U_{\beta\delta}\neq\emptyset,
\]
then the following Identity holds:
\begin{equation}\label{eq-translocmor2}
\Theta_{\beta\alpha}(g)=\Phi_{\beta\delta}\!\left(t_{\calG}(g)\right)\Theta_{\delta\gamma}(g)\Phi_{\gamma\alpha}\!\left(s_{\calG}(g)\right).
\end{equation}
\end{Cor} 
Hence, the following fact has been proved:

\fbox{\parbox{12cm}{\bf Given a generalized morphism $\left(P,\pi,\varepsilon,X_{\calG}\right)$ from the groupoid $\calG$ to $\calH$ in the sense of Definition~\ref{def-genmorph}, given local sections $\sigma_\alpha$ of $\pi$ subordinate to a trivializing cover $\mathfrak{U}$, there are local generalized morphisms $\Theta_{\beta\alpha}$ from any local component $\calG_{\alpha,\beta}$ of $\calG$ to $\calH$, transforming according to Equation (\ref{eq-translocmor2}).}}

\subsection{From local generalized morphisms to generalized morphisms}\label{ssec-locgentogenmor}
Assume to have two Lie groupoids $\calG$ and $\calH$, with respective manifolds of objects $X_{\calG}$ and $X_{\calH}$.
Assume additionally to have a fixed open cover $\mathfrak{U}$ of $X_{\calG}$ and corresponding local trivializing data $\left(\mathfrak{U},\varepsilon_\alpha,\Phi_{\alpha\beta}\right)$ on $X_{\calG}$ with values in $\calH$, with local momenta $\varepsilon_\alpha\colon U_\alpha\to X_{\calH}$ and cocycles $\Phi_{\alpha\beta}$. 
\begin{Def}[\cite{Con},~\cite{HS}]\label{def-locgenmor}
A {\bf local generalized morphism from $\calG$ to $\calH$ subordinate to the local trivializing data $\left(\mathfrak{U},\varepsilon_\alpha,\Phi_{\alpha\beta}\right)$} consists of smooth maps $\Theta_{\beta\alpha}$ from any local component $\calG_{\alpha,\beta}$ of $\calG$ to $\calH$, such that the following conditions hold:
\begin{itemize}
\item[a)] the following diagrams commute
\begin{equation}\label{eq-diagmormomen}
\begin{CD} 
\calG_{\alpha,\beta}  @>\Theta_{\beta\alpha}>> \calH\\
@Vs_\calG VV              @VV s_\calH V\\
U_\alpha @>\varepsilon_\alpha>> X_\calH
\end{CD}\quad,\quad 
\begin{CD} 
\calG_{\alpha,\beta}  @> \Theta_{\beta\alpha} >> \calH\\
@V t_\calG VV              @VV t_\calH V\\
U_\beta        @>\varepsilon_\beta>> X_\calH
\end{CD}\quad\text{and}\quad
\begin{CD} 
U_\alpha  @> \varepsilon_\alpha >> X_\calH\\
@V \iota_\calG VV              @VV \iota_\calH V\\
\calG_\alpha          @>\Theta_\alpha>> \calH
\end{CD}\quad,
\end{equation}
where the notation was used 
\[
\Theta_\alpha\colon=\Theta_{\alpha\alpha},\quad \forall \alpha.
\] 
\item[b)] The following identity must hold, for any choice of indices $\alpha$, $\beta$, $\gamma$:
\begin{equation}\label{eq-genhomom}
\Theta_{\gamma\alpha}(g_1g_2)=\Theta_{\gamma\beta}(g_1)\Theta_{\beta\alpha}(g_2),\quad g_1\in\calG_{\beta,\gamma},\quad g_2\in\calG_{\alpha,\beta},\quad s_{\calG}(g_1)=t_{\calG}(g_2); 
\end{equation} 
moreover, the local morphisms $\Theta_{\beta\alpha}$ are related to each other by the following equation:
\begin{equation}\label{eq-relgentrans}
\Theta_{\beta\alpha}(g)=\Phi_{\beta\delta}\!\left(t_{\calG}(g)\right)\Theta_{\delta\gamma}(g)\Phi_{\gamma\alpha}\!\left(s_{\calG}(g)\right),\quad \forall g\in \calG_{\alpha,\beta}\cap\calG_{\gamma,\delta},
\end{equation}
for any four open subsets $U_\alpha$, $U_\beta$, $U_\gamma$ and $U_\delta$ of $X_{\calG}$, such that
\[
U_{\alpha\gamma}\neq\emptyset,\quad U_{\beta\delta}\neq\emptyset.
\]
\end{itemize}
\end{Def}

\begin{Rem}
I will discuss later a cohomological interpretation of Equation (\ref{eq-relgentrans}) along the same patterns of Subsection~\ref{ssec-cohomint}: this suggests the interpretation of generalized morphisms between Lie groupoids as generalizations of {\em descent data}.
Moreover, it permits to construct more interesting examples of generalized morphisms: e.g.\ one feature that I will address in a subsequent paper is that generalized morphisms from the homotopy path groupoid of a smooth manifold $M$ to the trivial Lie groupoid associated to a Lie group $G$ corresponds to a flat principal $G$-bundle over $M$, and this also suggests another way of constructing flat bundles with groupoid structure. 
\end{Rem}

Given a local generalized morphism $\Theta$ from $\calG$ to $\calH$, subordinate to local trivializing data $\left(\mathfrak{U},\varepsilon_\alpha,\Phi_{\alpha\beta}\right)$ over $X_{\calG}$ with values in $\calH$, I want to prove that we can produce out of it a generalized morphism from $\calG$ to $\calH$ in the sense of Definition~\ref{def-genmorph}.
To do this, I first need the following technical
\begin{Lem}\label{lem-localaction}
For any choice of indices $\alpha$, $\beta$, there is an ``action'' of the local component $\calG_{\alpha,\beta}$ of $\calG$ from the trivial bundle $\varepsilon_\alpha^*\calU_{\calH}$ to the trivial bundle $\varepsilon_\beta^*\calU_{\calH}$, i.e.\ there is a smooth map $\Psi_{\alpha,\beta}^L$ with the following properties:
\begin{itemize}
\item[a)] $\Psi_{\alpha,\beta}^L$ is a map from
\[
\calG_{\alpha,\beta}\times_{\pr_1}\varepsilon_\alpha^*\calU_{\calH}\colon=\left\{\left(g;(x,h)\right)\in\calG_{\alpha,\beta}\times\varepsilon_\alpha^*\calU_{\calH}\colon s_{\calG}(g)=x\right\}
\] 
to the trivial bundle $\varepsilon_\beta^*\calU_{\calH}$; usually, the ``action map'' will simply be denoted by
\[
\calG_{\alpha,\beta}\times_{\pr_1}\varepsilon_\alpha^*\calU_{\calH}\ni(g;(x,h))\overset{\Psi_{\alpha,\beta}^L}\mapsto g(x,h)\in \varepsilon_\beta^*\calU_\calH.
\]
Moreover, choosing $\beta=\alpha$, the following equation holds:
\[
\iota_{\calG}(x)(x,h)=(x,h),\quad \forall (x,h)\in \varepsilon_\alpha^*\calU_\calH.
\]
\item[b)] For any choice of three indices $\alpha$, $\beta$, $\gamma$, the actions $\Psi_{\alpha,\beta}^L$, $\Psi_{\beta,\gamma}^L$ and $\Psi_{\alpha,\gamma}^L$ are compatible in the following way:
\[
(g_1g_2)(x,h)=g_1\left(g_2(x,h)\right),\quad \forall g_1\in\calG_{\beta,\gamma},\quad g_2\in \calG_{\alpha,\beta},\quad s_{\calG}(g_1)=t_{\calG}(g_2).
\] 
\end{itemize}
\end{Lem} 
\begin{proof}
Define the action $\Psi_{\alpha,\beta}^L$ as follows:
\begin{equation}\label{eq-genaction}
\calG_{\alpha,\beta}\times_{\pr_1}\varepsilon_\alpha^*\calU_{\calH}\ni (g,(x,h))\mapsto \left(t_{\calG}(g),\Theta_{\beta\alpha}(g)h\right).
\end{equation}
Notice that the action is well-defined: in fact, 
\[
s_{\calH}\!\left(\Theta_{\beta\alpha}(g)\right)=\varepsilon_\alpha(s_{\calG}(g))=\varepsilon_\alpha(x)=t_{\calH}(h).
\]
First thing, one has to show that the result belongs really to the trivial bundle $\varepsilon_\beta^*\calU_{\calH}$.
For this purpose, I use the commutativity of the second diagram in (\ref{eq-diagmormomen}), and one sees immediately that
\begin{align*}
t_{\calH}\!\left(\Theta_{\beta\alpha}(g)h\right)&=t_{\calH}\!\left(\Theta_{\beta\alpha}(g)\right)=\\
&=\varepsilon_\beta\!\left(t_{\calG}(g)\right),
\end{align*}
whence the claim follows.

An easy computation shows:
\[
\pr_1\!\left(g(x,h)\right)=\pr_1\!\left(t_{\calG}(g),\Theta_{\beta\alpha}(g)h\right)=t_{\calG}(g),
\]
hence the relations between action and momentum is verified.
Choosing $\beta=\alpha$, the ``identity axiom'' for the action (\ref{eq-genaction}) is an easy consequence of the commutativity of the third diagram in (\ref{eq-diagmormomen}). 

It remains to verify the compatibility condition for any choice of indices $\alpha$, $\beta$ and $\gamma$; this is but an easy consequence of Identity (\ref{eq-genhomom}) and Equation (\ref{eq-genaction}), recalling that, if $g_1$ is composable with $g_2$, and $g_1\in \calG_{\beta,\gamma}$ and $g_2\in\calG_{\alpha,\beta}$, then their product $g_1g_2$ lies in $\calG_{\alpha,\gamma}$:
\begin{align*}
(g_1g_2)(x,h)&=\left(t_{\calG}(g_1g_2),\Theta_{\gamma\alpha}(g_1g_2)h\right)=\\
&=\left(t_{\calG}(g_1),\Theta_{\gamma\beta}(g_1)\Theta_{\beta\alpha}(g_2)h\right)=\\
&=g_1\!\left(s_{\calG}(g_1),\Theta_{\beta\alpha}(g_2)h\right)=\\
&=g_1\!\left(t_{\calG}(g_2),\Theta_{\beta\alpha}(g_2)h\right)=\\
&=g_1\!\left(g_2(x,h)\right).
\end{align*}
\end{proof}

Theorem~\ref{thm-equivtrivdata} of Section~\ref{sec-localdata} states the equivalence between local trivializing data over a manifold $X$ with values in a groupoid $\calG$, in the sense of Definition~\ref{def-trivdata}, and principal bundles over $X$ with structure groupoid $\calG$ in the sense of Definition~\ref{def-princgroupoid} of Section~\ref{sec-divmap}.
Hence, given a local generalized morphism as in Definition~\ref{def-locgenmor}, subordinate to local trivializing data $\left(\mathfrak{U},\varepsilon_\alpha,\Phi_{\alpha\beta}\right)$ over $X_{\calG}$ with values in $\calH$, one gets automatically a right $\calH$-principal bundle $P$ over $X_{\calG}$ by the ``gluing procedure''.

Here comes now into play Identity (\ref{eq-relgentrans}), relating explicitly the cocycle $\Phi_{\alpha\beta}$ to the local generalized morphism $\Theta$: the actions of the local components $\calG_{\alpha,\beta}$ from the trivial bundles $\varepsilon_\alpha^*\calU_{\calH}$ to the trivial bundles $\varepsilon_\beta^*\calU_{\calH}$ from Lemma~\ref{lem-localaction} ``glue'' together to give a well-definite left $\calG$-action on $P$ along the projection $\pi$ from $P$ to $X_{\calG}$, which is compatible with the right $\calH$-action.

\begin{Lem}\label{lem-locmortogenmor}
Let $\calG$ and $\calH$ two Lie groupoids, and $\Theta$ be a local generalized morphism from $\calG$ to $\calH$ subordinate to the local trivializing data $\left(\mathfrak{U},\varepsilon_\alpha,\Phi_{\alpha\beta}\right)$ over $X_{\calG}$ with values in $\calH$, and let $\left(P,\pi,\varepsilon,X_{\calG}\right)$ the right $\calH$-principal bundle over $X_{\calG}$ associated to the above local trivializing data.

Then, $\Theta$ endows $P$ with a left $\calG$-action with momentum $\pi$, which is compatible with the right $\calH$-action.
\end{Lem}
\begin{proof}
Recall that the total space of $P$ is defined as the quotient of the disjoint union 
\[
\coprod_{\alpha}\varepsilon_\alpha^*\calU_{\calH}
\]
by the following equivalence relation:
\[
(\alpha,x_1,h_1)\sim (\beta,x_2,h_2)\Leftrightarrow \begin{cases}
&U_{\alpha\beta}\neq \emptyset,\\
&x_1=x_2,\\
&h_1=\Phi_{\alpha\beta}(x_1)h_2.
\end{cases}
\]
The projection $\pi$, resp.\ the momentum for the right $\calH$-action, is defined via
\[
[\alpha,x,h]\overset{\pi}\mapsto x,\quad\text{resp.}\quad [\alpha,x,h]\overset{\varepsilon}\mapsto s_{\calH}(h).
\]
Let me now define the left $\calG$-action simply by ``quotienting'' the left actions of the local components $\calG_{\alpha,\beta}$, i.e.\
\begin{equation}\label{eq-genleftact}
\left(g,[\alpha,x,h]\right)\mapsto\left[\beta,t_{\calG}(g),\Theta_{\beta\alpha}(g)h\right],
\end{equation} 
where 
\[
s_{\calG}(g)=x,\quad t_{\calG}(g)\in U_{\beta},
\]
for some index $\beta$.
One has to show first that the action defined in Equation (\ref{eq-genleftact}) is well-defined, i.e.\ it does not depend $i)$ neither on the choice of the representative $[\alpha,x,h]$ nor $ii)$ on the choice of the index $\beta$, such that the target of $g$ is in $U_\beta$.

Let me first show independence of the choice of the representative.
In fact, any other representative would have the form
\[
[\alpha,x,h]=[\gamma,x,\Phi_{\gamma\alpha}(x)h],\quad U_{\gamma\alpha}\neq\emptyset.
\]
Then, 
\begin{align*}
\left(g,[\gamma,x,\Phi_{\gamma\alpha}(x)h]\right)&\mapsto \left[\beta,t_{\calG}(g),\Theta_{\beta\gamma}(g)\Phi_{\gamma\alpha}(x)h\right]=\\
&=\left[\beta,t_{\calG}(g),\Theta_{\beta\gamma}(g)\Phi_{\gamma\alpha}(s_{\calG}(g))h\right]=\\
&=\left[\beta,t_{\calG}(g),\Theta_{\beta\alpha}(g)h\right],
\end{align*}
using Equation (\ref{eq-relgentrans}) with $\beta=\delta$.
On the other, when the target of $g$ belongs to the intersection $U_{\delta\beta}$ of two open subsets in the cover $\mathfrak{U}$ of $X_{\calG}$, consider the action of $g$ to be 
\[
\left(g,[\alpha,x,h]\right)\mapsto \left[\delta,t_{\calG}(g),\Theta_{\delta\alpha}(g)h\right].
\]
But then, again by Equation (\ref{eq-relgentrans}), one gets
\begin{align*}
\left[\delta,t_{\calG}(g),\Theta_{\delta\alpha}(g)h\right]&=\left[\beta,t_{\calG}(g),\Phi_{\beta\delta}(t_{\calG}(g))\Theta_{\delta\alpha}(g)h\right]=\\
&=\left[\beta,t_{\calG}(g),\Theta_{\beta\alpha}(g)h\right],
\end{align*}
by taking $\gamma=\alpha$.
Hence, the left $\calG$-action is well-defined.

One has to show that its momentum is $\pi$; this follows immediately by its very definition.
The fact that Equation (\ref{eq-genleftact}) defines really an action is a consequence of Lemma~\ref{lem-localaction}.

Therefore, it remains to show that the left $\calG$-action is compatible with the right $\calH$-action, i.e. one has to show that $i)$ the momentum $\varepsilon$ of the right $\calH$-action is $\calG$-invariant and $ii)$ both actions commute with each other.

To show $i)$, notice simply that the momentum is defined via projection onto the second factor of any pair $(x,h)$ in a trivial bundle $\varepsilon_\alpha^*\calU_{\calH}$ with the source of $\calH$; since the source map of $\calH$ is invariant w.r.t.\ left multiplication in $\calH$, the claim follows.

To show $ii)$, compute both actions:
\begin{align*}
\left(g\left[\alpha,x,h_1\right]\right)h_2&=\left[\beta,t_{\calG}(g),\Theta_{\beta\alpha}(g)h_1\right]h_2=\\
&=\left[\beta,t_{\calG}(g),(\Theta_{\beta\alpha}(g)h_1)h_2\right]=\\
&=\left[\beta,t_{\calG}(g),\Theta_{\beta\alpha}(g)(h_1h_2)\right]=\\
&=g\left[\alpha,x,h_1h_2\right]=\\
&=g\left(\left[\alpha,x,h_1\right]h_2\right),
\end{align*}
where
\[
s_{\calG}(g)=x,\quad t_{\calG}(g)\in U_\beta,\quad s_{\calH}(h_1)=t_{\calH}(h_2).
\]
Finally, notice that the left $\calG$-action is smooth; this follows from the fact that its local form is smooth and by the usual definition of a smooth structure on $P$.
\end{proof}
Hence, Lemma~\ref{lem-locmortogenmor}, together with Theorem~\ref{thm-locgenmor}, implies the following fact:

\fbox{\parbox{12cm}{\bf There is a one-to-one correspondence between generalized morphisms in the sense of Definition~\ref{def-genmorph} and local generalized morphisms in the sense of Definition~\ref{def-locgenmor}, provided one chooses an open covering of the manifold of objects of the source groupoid}} 
\begin{Rem}
Let me just notice that an equivalent statement was already present, for the special case of {\em monodromy groupoids of foliated manifolds}, in~\cite{HS} and in a more general form e.g.\ in~\cite{Moer3}: in fact, a generalized morphism from $\calG$ to $\calH$ is equivalent to a choice of an open covering of $X_\calG$ and a strict morphism of groupoids from $\coprod_{\alpha,\beta}\calG_{\alpha,\beta}$ (which is equivalent to $\calG$) to $\calH$.
\end{Rem}

\subsection{Some examples of local generalized morphisms}\label{ssec-exlocgenmor}
In this subsection, I construct some examples of local generalized morphisms in the sense of Definition~\ref{def-locgenmor}, and to relate them to the ordinary version of generalized morphisms, that of Definition~\ref{def-genmorph}.

\subsubsection{Local generalized morphisms between action groupoids}\label{sssec-genmoractgroup}
Assume given two action groupoids, say $G\ltimes X$ and $H\ltimes Y$, where $X$ and $Y$, resp.\ $G$ and $H$, are smooth manifolds, resp.\ Lie groups acting from the left on $X$ and $Y$ respectively; for the definition of action groupoids, I refer to Subsubsection~\ref{sssec-actiongroupoid}.
Consider local trivializing data $\left(\mathfrak{U},\varepsilon_\alpha,\Phi_{\alpha\beta}\right)$ over the manifold of objects of $G\ltimes X$, which is $X$, with values in $H\ltimes Y$; as was proved in Subsubsection~\ref{sssec-actiongroupoid}, this is equivalent to a $Y$-pointed principal $H$-bundle over $X$, in the language adopted in Subsubsection~\ref{sssec-actiongroupoid}.

Consider now a local generalized morphism $\Theta$ subordinate to the aforementioned local trivializing data.
First of all, the local components $\left(G\ltimes X\right)_{\alpha,\beta}$ take the form:
\[
\left(G\ltimes X\right)_{\alpha,\beta}=\left\{(g,x)\in G\times X\colon x\in U_\alpha, gx\in U_\beta\right\}.
\] 
Then, there are smooth maps $\Theta_{\beta\alpha}$ from $\left(G\ltimes X\right)_{\alpha,\beta}$ to $H\ltimes Y$, satisfying a certain number of properties, listed in Definition~\ref{def-locgenmor}.
First of all, write 
\[
\Theta_{\beta\alpha}(g,x)\colon=\left(\Theta_{\beta\alpha}^H(g,x),\Theta_{\beta\alpha}^Y(g,x)\right), 
\] 
where
\[
\Theta_{\beta\alpha}^H(g,x)\in H,\quad \Theta_{\beta\alpha}^Y(g,x)\in Y.
\]
Since (by the commutativity of the diagrams (\ref{eq-diagmormomen})) the following equations hold
\[
t_{\calH}\circ\Theta_{\beta\alpha}=\varepsilon_\beta\circ t_{\calG}\quad\text{and}\quad s_{\calH}\circ\Theta_{\beta\alpha}=\varepsilon_\alpha\circ s_{\calG},
\]
it follows immediately
\[
\Theta_{\beta\alpha}^Y(g,x)=\varepsilon_\alpha(x),\quad \varepsilon_\beta(gx)=\Theta_{\beta\alpha}^H(g,x)\varepsilon_\alpha(x).
\]
Let me introduce the following notation
\[
\theta_{\beta\alpha}(g,x)\colon=\Theta_{\beta\alpha}^H(g,x)\in H,\quad (g,x)\in \left(G\ltimes X\right)_{\alpha,\beta}.
\]
Moreover, by the commutativity of the third diagram of (\ref{eq-diagmormomen}), one has
\[
\theta_\alpha(e,x)=e,\quad \forall x\in U_\alpha.
\]
Writing down explicitly Identity (\ref{eq-genhomom}), one gets the following identity:
\begin{align*}
\Theta_{\gamma\alpha}\!\left(g_1g_2,x)\right)&=\left(\theta_{\gamma\alpha}(g_1g_2,x),\varepsilon_\alpha(x)\right)\overset{!}=\\
&\overset{!}=\Theta_{\gamma b\eta}\!\left(g_1,g_2x)\right)\Theta_{\beta\alpha}\!\left(g_2,x)\right)=\\
&=\left(\theta_{\gamma\beta}(g_1,g_2x),\varepsilon_\beta(g_2x)\right)\left(\theta_{\beta\alpha}(g_2,x),\varepsilon_\alpha(x)\right)=\\
&=\left(\theta_{\gamma\beta}(g_1,g_2x)\theta_{\beta\alpha}(g_2,x),\varepsilon_\alpha(x)\right),
\end{align*}
whence it follows
\[
\theta_{\gamma\alpha}(g_1g_2,x)=\theta_{\gamma\beta}(g_1,g_2x)\theta_{\beta\alpha}(g_2,x),\quad \forall x\in U_{\alpha},\quad g_2x\in U_\beta,\quad g_1g_2x\in U_\gamma. 
\]
There is also Identity (\ref{eq-relgentrans}) to be considered, relating the local generalized morphisms to transition maps of the bundle: namely, it takes the form
\begin{align*}
\Theta_{\beta\alpha}(g,x)&=\left(\theta_{\beta\alpha}(g,x),\varepsilon_\alpha(x)\right)\overset{!}=\\
&\overset{!}=\Phi_{\beta\delta}(gx)\Theta_{\delta\gamma}(g,x)\Phi_{\gamma\alpha}(x)=\\
&=\left(\Phi^H_{\beta\delta}(gx),\varepsilon_\delta(gx)\right)\left(\theta_{\delta\gamma}(g,x),\varepsilon_\gamma(x)\right)\left(\Phi_{\gamma\alpha}^H(x),\varepsilon_\alpha(x)\right)=\\
&=\left(\Phi^H_{\beta\delta}(gx)\theta_{\delta\gamma}(g,x)\Phi_{\gamma\alpha}^H(x),\varepsilon_\alpha(x)\right),
\end{align*}
whence one easily deduces the following relation between the maps $\theta_{\beta\alpha}$ and the transition maps of principal $H$-bundle $P$:
\[
\theta_{\beta\alpha}(g,x)=\Phi^H_{\beta\delta}(gx)\theta_{\delta\gamma}(g,x)\Phi_{\gamma\alpha}^H(x),\quad x\in U_\alpha\cap U_\gamma,\quad gx\in U_\beta\cap U_\gamma.
\] 
The right principal $H$-bundle $P$ is given explicitly by
\[
P=\coprod_\alpha U_\alpha\times H/ \sim,
\]
where the equivalence relation $\sim$ was described explicitly in Subsubsection~\ref{sssec-actiongroupoid}.

The maps $\theta_{\beta\alpha}$ define a lift of the left action of $G$ on $X$ to $P$. 
In fact, define the left action of $G$ on $P$ as follows:
\begin{equation}\label{eq-liftaction}
G\times P\ni \left(g,[\alpha,x,h]\right)\mapsto [\beta,gx,\theta_{\beta\alpha}(g,x)h]\in P,\quad gx\in U_\beta.
\end{equation}                                         
The left action (\ref{eq-liftaction}) is well-defined, since, picking up other indices $\gamma$ and $\delta$, such that
\[
x\in U_\gamma,\quad gx\in U_\delta, 
\]
so that
\[
[\alpha,x,h]=[\gamma,x,\Phi_{\gamma\alpha}^H(x)h],\quad [\beta,gx,\theta_{\beta\alpha}(g,x)h]=\left[\delta,gx,\Phi_{\delta\beta}^H(gx)\theta_{\beta\alpha}(g,x)h\right],
\]
one gets, by the above version of (\ref{eq-relgentrans}),
\begin{align*}
G\times P\ni\left(g,[\gamma,x,\Phi_{\gamma\alpha}^H(x)h]\right)&=\left(g,[\alpha,x,h]\right)\mapsto\\
&\mapsto\left[\delta,gx,\theta_{\delta\gamma}(g,x)\Phi_{\gamma\alpha}^H(x)h\right]=\\
&=\left[\delta,gx,\Phi_{\delta\beta}(gx)^H\theta_{\beta\alpha}(g,x)h\right]=\\
&=[\beta,gx,\theta_{\beta\alpha}(g,x)h]. 
\end{align*} 
Since $\theta_\alpha(e,x)=e$, it follows that 
\[
\left(e,[\alpha,x,h]\right)\mapsto [\alpha,x,h],
\]
and since $\theta_{\gamma\alpha}(g_1g_2,x)=\theta_{\gamma\beta}(g_1,g_2x)\theta_{\beta\alpha}(g_2,x)$, it follows
\begin{align*}
(g_1g_2)[\alpha,x,h]&=[\gamma,(g_1g_2)x,\theta_{\gamma\alpha}(g_1g_2,x)h]=\\
&=[\gamma,g_1(g_2x),\theta_{\gamma\beta}(g_1,g_2x)\theta_{\beta\alpha}(g_2,x)h]=\\
&=g_1\left[\beta,g_2x,\theta_{\beta\alpha}(g_2,x)h\right]=\\
&=g_1\left(g_2[\alpha,x,h]\right),\quad g_1g_2x\in U_\gamma,\quad g_2x\in U_\beta,\quad x\in U_\alpha.
\end{align*}
Finally, since the projection from $P$ to $X$ is simply
\[
[\alpha,x,h]\mapsto x,
\]
it is clear that $\pi$ is $G$-equivariant, and it is obvious that the left $G$-action is compatible with the right $H$-action.
Moreover, notice that the associated bundle $P\times_H Y$ may be constructed via the gluing procedure:
\[
P\times_H Y=\coprod_\alpha U_\alpha\times Y/ \sim,
\]
where now the equivalence relation takes the form
\[
[\alpha,x_1,y_1]\sim[\beta,x_2,y_2]\Leftrightarrow \begin{cases}
&U_{\alpha\beta}\neq \emptyset,\\
&x_1=x_2,\\
&y_2=\Phi_{\beta\alpha}^H(x_1)y_1.
\end{cases} 
\]
The left $G$-action on $X$ may be then also lifted to the associated bundle $P\times_H Y$ via the formula:
\[
P\times_H Y\ni \left(g,[\alpha,x,y]\right)\mapsto [\beta,gx,\theta_{\beta\alpha}(g,x)y]\in P\times_H Y.
\]
Repeating almost verbatim the arguments and the computations used for showing that $P$ inherits a lift of the left $G$-action on $X$, one can show that the left $G$-action on $P\times_H Y$ is well-defined and lifts exactly the left $G$-action on $X$.

Moreover, the global section $\eta$ of $P\times_H Y$ is defined by the local momenta $\varepsilon_\alpha$.
The section $\eta$ is $G$-equivariant, i.e.
\[
\eta(gx)=g\eta(x),\quad \forall g\in G,\quad x\in X.
\]
Namely, the following identities are already known:
\[
\varepsilon_\beta(x)=\Phi_{\beta\alpha}^H(x)\varepsilon_\alpha(x),\quad \forall x\in U_{\alpha\beta};\quad \varepsilon_\beta(gx)=\theta_{\beta\alpha}(g,x)\varepsilon_\alpha(x),\quad x\in U_\alpha,\quad gx\in U_\beta.
\]
The first identity together with the identity relating the transition maps of $P$ and the  
Now we have:
\begin{align*}
\eta(gx)&=\left[\beta,gx,\varepsilon_\beta(gx)\right]=\\
&=\left[\beta,gx,\theta_{\beta\alpha}(g,x)\varepsilon_\alpha(x)\right]=\\
&=g\left[\alpha,x,\varepsilon_\alpha(x)\right]=\\
&=g\eta(x),\quad x\in U_\alpha,\quad gx\in U_\beta.
\end{align*}
(The choice of the indices $\alpha$ and $\beta$ is completely uninfluent, because of the above relationship between the various local momenta and transition maps and between transition maps and the maps $\theta_{\beta\alpha}$.)  

Hence, the following equivalence has been established:

\fbox{\parbox{12cm}{\bf A local generalized morphism $\Theta$ from the action groupoid $G\ltimes X$ to the action groupoid $H\ltimes Y$, hence, a generalized morphism between $G\ltimes X$ and $H\ltimes Y$, subordinate to local trivializing data over $X$ with values in $H\ltimes Y$, is equivalent to $i)$ a right principal $H$-bundle $P$ over $X$, $ii)$ a smooth lift of the left $G$-action on $X$ to $P$ and to $P\times_H Y$ and $iii)$ a $G$-equivariant global section $\eta$ of $P\times_H Y$; Such a datum I will call a $G$-equivariant, $Y$-pointed principal $H$-bundle over $X$, provided one chooses an open covering of $X$.}}

\begin{Rem}
Notice that in the special case $Y$ is a point, on which $H$ acts trivially (so, the action groupoid reduces to the trivial groupoid $H$, with trivial source, target and unit maps), there is no condition about $Y$-points, and one gets simply a $G$-equivariant principal $H$-bundle over $X$, which is, in the language of~\cite{LTX}, a particular example of a {\em principal $H$-bundle over the action groupoid $G\ltimes X$}.
\end{Rem}    

\subsubsection{A local generalized morphism from the product groupoid $X\times X$ to the gauge groupoid $\calG(P)$, for a right principal $G$-bundle $P$ over $X$}\label{sssec-prodgauge}
Consider an ordinary right principal $G$-bundle over the manifold $X$, where $G$ is a Lie group; to $X$ is associated the product groupoid $X\times X$ (I refer to Subsubsection~\ref{sssec-productgroupoid} for more details), whereas to $P$ and $X$ is associated the gauge groupoid $\calG(P)$ (I refer to Subsubsection~\ref{sssec-gaugegroup} for more details).
The aim now is to produce a local generalized morphism $\Theta$ from $X\times X$ to $\calG(P)$, subordinate to the local trivializing data over $X$ with values in $\calG(P)$ associated to an open cover $\mathfrak{U}$, local momenta given by the canonical inclusions of $U_\alpha$ into $X$ and associated to the $P$-values cocycle (\ref{eq-gaugecoc1}) of Subsubsection~\ref{sssec-gaugegroup}; such local trivializing data give rise to the right $\calG(P)$-bundle $\left(X\times P,\pr_1,\pi\circ\pr_2,X\right)$.

Recall the explicit form of the transition maps of the local trivializing data:
\[
\Phi_{\alpha\beta}(x)=\left[\sigma_\alpha(x),\sigma_\beta(x)\right],\quad x\in U_{\alpha\beta},
\]
where $\sigma_\alpha$ denotes a local section of $P$ over $U_\alpha$

The maps
\begin{equation}\label{eq-locgengauge}
\Theta_{\beta\alpha}(x_1,x_2)\colon=\left[\sigma_\beta(x_2),\sigma_\alpha(x_1)\right],\quad x_1\in U_\alpha,\quad x_2\in U_\beta,
\end{equation}
define a local generalized morphism from $X\times X$ to $\calG(P)$ subordinate to the above local trivializing data; moreover, the associated generalized morphism from $X\times X$ to $\calG(P)$ is the bundle $\left(X\times P,\pr_1,\pi\circ\pr_2,X\right)$, with left $X\times X$-action given by
\begin{equation}\label{eq-genmorgauge}
(x,y)(y,p)\colon=(x,p),\quad x,y\in X,\quad p\in P.
\end{equation}
(It is immediate to verify that the above right $\calG(P)$-bundle, endowed with the left $X\times X$-action of Equation (\ref{eq-genmorgauge}), satisfies all the properties stated in Definition~\ref{def-genmorph}.) 

By the very definition of source map and target map of $X\times X$ and of $\calG(P)$, it follows immediately that all diagrams in (\ref{eq-diagmormomen}) commute.
Identity (\ref{eq-genhomom})is an easy consequence of the definition of the product in $\calG(P)$ and in $\Pi(X)$, namely:
\begin{align*}
\Theta_{\gamma\alpha}\!\left((x_1,x_2)(x_2,x_3)\right)&=\Theta_{\gamma\alpha}\!\left(x_1,x_3\right)=\\
&=\left[\sigma_\gamma(x_1),\sigma_\alpha(x_3)\right]=\\
&=\left[\sigma_\gamma(x_1),\sigma_\beta(x_2)\right]\left[\sigma_\beta(x_2),\sigma_\alpha(x_3)\right]=\\
&=\Theta_{\gamma\beta}(x_1,x_2)\Theta_{\beta\alpha}(x_2,x_3).
\end{align*}
Moreover, Identity (\ref{eq-relgentrans}), relating transition maps of local trivializing data and local generalized morphisms, follows from the very same arguments.
To see that in fact the generalized morphism $\Theta$ as defined in Equation (\ref{eq-locgengauge}) corresponds to the generalized morphism $\left(X\times P,\pr_1,\pi\circ\pr_2,X\right)$, with left $\Pi(X)$-action as defined in Equation (\ref{eq-genmorgauge}), recall the construction of Subsection~\ref{ssec-locgenmor}, where, given a generalized morphism $P$ from $\calG$ to $\calH$ in the sense of Definition~\ref{def-genmorph}, I constructed a local generalized morphism in the sense of Definition~\ref{def-locgenmor}.
In this particular case, one has to compute the division map of the bundle $\left(X\times P,\pr_1,\pi\circ\pr_2,X\right)$.
It is in fact simply given by
\begin{equation}\label{eq-divmapgenmor}
\phi_{X\times P}\!\left((x_1,p_1),(x_1,p_2)\right)\colon=[p_1,p_2].
\end{equation}
The verification that Equation (\ref{eq-divmapgenmor}) is an in fact the division map of the above bundle is immediate.

Then, it is known that the corresponding local generalized morphism is given by Equation (\ref{eq-genmor2}): let me compute explicitly the result
\begin{align*}
\Theta_{\beta\alpha}(x_1,x_2)&\overset{!}=\phi_{X\times P}\!\left(\sigma_\beta(t_{\Pi(X)}(x_1,x_2)),(x_1,x_2)\sigma_\alpha(s_{\Pi(X)}(x_1,x_2))\right)=\\
&=\phi_{X\times P}\!\left(\sigma_\beta(x_1),(x_1,x_2)\sigma_\alpha(x_2)\right)=\\
&=\phi_{X\times P}\!\left((x_1,\sigma_\beta(x_1)),(x_1,x_2)(x_2,\sigma_\alpha(x_2))\right)=\\
&=\phi_{X\times P}\!\left((x_1,\sigma_\beta(x_1)),(x_1,\sigma_\alpha(x_2))\right)=\\
&=\left[\sigma_\beta(x_1),\sigma_\alpha(x_2)\right],\quad x_1\in U_\beta,\quad x_2\in U_\alpha.
\end{align*}
Notice that I have adopted the same notation for local sections of $X\times P$ w.r.t.\ $\pr_1$ and local sections of $P$ over $U_\alpha$: in fact, a local section $\sigma_\alpha$ of $P$ over $U_\alpha$ specifies a local section of $X\times P$ over $U_\alpha$ in the following way:
\[
\sigma_\alpha:U_\alpha\to P\leadsto \sigma_\alpha:U_\alpha\to X\times P,\quad x\mapsto (x,\sigma_\alpha(x)).
\]

\subsubsection{A local generalized morphism from the action groupoid $G\ltimes X$ to the gauge groupoid $\calG(P)$ for a right principal $H$-bundle $P$ over $X$}\label{sssec-actiongauge}
Consider in this subsubsection a manifold $X$, acted on from the left by the Lie group $G$, and a right principal $H$-bundle $P$ over $X$, for $H$ a Lie group; consider furthermore the action groupoid $G\ltimes X$ and the gauge groupoid $\calG(P)$ associated to $P$.

Given the local trivializing data consisting of an open trivializing cover $\mathfrak{U}$ of $X$ w.r.t.\ $P$, the local momenta given by the inclusions $U_\alpha\hookrightarrow X$ and the cocycle over $X$ with values in $\calG(P)$ given by Equation (\ref{eq-gaugecoc1}), I construct a local generalized morphism from $G\ltimes X$ to $\calG(P)$ subordinate to the above local trivializing data, which, like in the previous subsubsection, define the right $\calG(P)$-principal bundle $\left(X\times P,\pr_1,\pi\circ\pr_2,X\right)$.
The maps
\begin{equation}\label{eq-locactiongroup}
\Theta_{\beta\alpha}(g,x)\colon=\left[\sigma_\beta(gx),\sigma_\alpha(x)\right],\quad x\in U_\alpha,\quad gx\in U_\beta,
\end{equation}
and $\sigma_\alpha$ denotes a local section over $U_\alpha$ of the right principal $H$-bundle $P$.

First of all, the very definition of target map, source map and unit map of both groupoids $G\ltimes X$ and $\calG(P)$ show immediately that the three diagrams in (\ref{eq-diagmormomen}) do indeed commute.
Identity (\ref{eq-genhomom}) follow from the following computation, where I make use again of the product laws in both $G\ltimes X$ and $\calG(P)$:
\begin{align*}
\Theta_{\gamma\alpha}\!\left((g_1,g_2x),(g_2,x)\right)&\overset{!}=\Theta_{\gamma\alpha}\!\left(g_1g_2,x\right)=\\
&=\left[\sigma_\gamma(g_1g_2x),\sigma_\alpha(x)\right]=\\
&=\left[\sigma_\gamma(g_1g_2x),\sigma_\beta(g_2x)\right]\left[\sigma_\beta(g_2x),\sigma_\alpha(x)\right]=\\
&=\Theta_{\gamma\beta}(g_1,g_2x)\Theta_{\beta\alpha}(g_2,x),\quad x\in U_\alpha,\quad g_2x\in U_\beta,\quad g_1g_2x\in U_\gamma.
\end{align*}
Analogous computations, recalling the explicit shape of transition maps of the right $\calG(P)$-bundle $\left(X\times P,\pr_1,\pi\circ\pr_2,X\right)$ associated to the local sections $\sigma_\alpha$ of $P$, show that Identity (\ref{eq-relgentrans}) also hold.

Recalling the shape of the division map of the bundle $\left(X\times P,\pr_1,\pi\circ\pr_2,X\right)$, let me now prove that the local generalized morphism defined by Equation (\ref{eq-locactiongroup}) corresponds to the generalized morphism $\left(X\times P,\pr_1,\pi\circ\pr_2,X\right)$, where the left $G\ltimes X$-action is defined via
\begin{equation}\label{eq-genactiongroup}
(g,x)(x,p)\colon=(gx,p),\quad g\in G,\quad x\in X,\quad p\in P.
\end{equation}
(It is immediate to verify that $\left(X\times P,\pr_1,\pi\circ\pr_2,X\right)$, endowed with the left $G\ltimes X$-action, satisfies all the requirements of Definition~\ref{def-genmorph}.)
In fact, using the computations of Subsection~\ref{ssec-locgenmor}, one gets the following expression for the local generalized morphism associated to the above generalized morphism:
\begin{align*}
\Theta_{\beta\alpha}(g,x)&=\phi_{X\times P}\!\left(\sigma_\beta(t_{G\ltimes X}(g,x)),(g,x)\sigma_\alpha(t_{G\ltimes X}(g,x))\right)=\\
&=\phi_{X\times P}\!\left(\sigma_\beta(gx),(g,x)\sigma_\alpha(x)\right)=\\
&=\phi_{X\times P}\!\left((x,\sigma_\beta(gx)),(g,x)(x,\sigma_\alpha(x))\right)=\\
&=\phi_{X\times P}\!\left((x,\sigma_\beta(gx)),(gx,\sigma_\alpha(x))\right)=\\
&=\left[\sigma_\beta(gx),\sigma_\alpha(x)\right],\quad x\in U_\alpha,\quad gx\in U_\beta,
\end{align*}
which corresponds exactly to the local generalized morphism of Equation (\ref{eq-locactiongroup}).

\subsection{Equivalences between generalized morphisms}\label{ssec-morgenmor}
Finally, I want to discuss the local aspects of equivalences between generalized morphisms from a groupoid $\calG$ to another groupoid $\calH$; by these, I mean morphisms of right $\calH$-principal bundles in the sense specified in~\cite{Moer2},~\cite{CR1}, which are additionally left $\calG$-equivariant.
I discussed already global aspects of such morphisms in~\cite{CR1}; in particular, since they are morphisms between $\calH$-principal bundles, they are invertible (whence the choice of name ``equivalence'') and there is a one-to-one correspondence between these equivalences and $\calG$-invariant generalized gauge transformations with values in $\calH$.

Let me now consider two groupoids $\calG$ and $\calH$, two generalized morphisms $\calG\overset{P}\to\calH$ and $\calG\overset{P}\to\calH$; as usual, I denote by $\pi_1$ and $\varepsilon_1$, resp.\ $\pi_2$ and $\varepsilon_2$, the projection and the momentum respectively of $P$, resp.\ $Q$.
I consider additionally an equivalence $\Sigma$ between $P$ and $Q$, in the sense specified above; by $K_{\Sigma}$, like in Subsection~\ref{ssec-morprincbun}, I denote the $\calG$-invariant generalized gauge transformation associated uniquely to $\Sigma$.
Assume that there is an open covering $\mathfrak{U}$ of $X_{\calG}$, such that there are local sections $\sigma_{\alpha}^1$ of $P$, resp.\ $\sigma_{\alpha}^2$ of $Q$, over any open set $U_\alpha$; the corresponding local trivializations, transition maps and local momenta are denoted as in Subsection~\ref{ssec-morprincbun}.
In analogy to Subsection~\ref{ssec-morprincbun}, Define
\begin{align*}
\Sigma_\alpha(x)\colon=K_\Sigma\!\left(\sigma_{\alpha}^1(x),\sigma_{\alpha}^2(x)\right),\quad \forall x\in U_\alpha.
\end{align*} 
Now, recall from Subsection~\ref{ssec-locgenmor} the formula for the local generalized morphisms $\Theta^P$ and $\Theta^Q$, associated respectively to $P$ and $Q$:
\[
\Theta^P_{\beta\alpha}(g)=\phi_P\!\left(\sigma_{\beta}^1(t_{\calG}(g)),g\sigma_{\alpha}^1(s_{\calG}(g))\right),\quad \Theta^Q_{\beta\alpha}(g)=\phi_Q\!\left(\sigma_{\beta}^2(t_{\calG}(g)),g\sigma_{\alpha}^2(s_{\calG}(g))\right),
\]
for $g$ belonging to the local component $\calG_{\alpha,\beta}$, defined explicitly also in the same subsection.
\begin{Prop}\label{prop-loceqgenmor}
The local maps $\Sigma_\alpha$ enjoy the following properties:
\begin{itemize}
\item[i)] (Compatibility with local momenta)
\[
t_{\calH}\circ \Sigma_\alpha=\varepsilon_{\alpha}^2,\quad s_{\calH}\circ\Sigma_\alpha=\varepsilon_{\alpha}^1.
\] 
\item[ii)] (Coboundary relation I)
\[
\Sigma_\beta(x)=\Phi_{\beta\alpha}^2(x)\Sigma_\alpha(x)\Phi_{\alpha\beta}^1(x),\quad x\in U_{\alpha\beta},
\] 
provided $U_\alpha$ and $U_\beta$ intersect nontrivially.
\item[iii)] (Coboundary relation II)
\[
\Theta^P_{\beta\alpha}(g)=\left(\Sigma_\beta(t_{\calG}(g))\right)^{-1}\Theta^Q_{\beta\alpha}(g)\Sigma_\alpha(s_{\calG}(g)),\quad \forall g\in \calG_{\alpha,\beta}.
\]
\end{itemize}
\end{Prop}  
\begin{proof}
The proof of Proposition~\ref{prop-locisomor} of Subsection~\ref{ssec-morprincbun} implies immediately the compatibility with local momenta and coboundary relation I; it remains therefore to show coboundary relation II.

Using the explicit form of the local generalized morphisms $\Theta_{\beta\alpha}^P$ and $\Theta_{\beta\alpha}^Q$, the relation
\[
\Sigma(\sigma_{\alpha}^1(x))=\sigma_{\alpha}^2(x)\Sigma_\alpha(x),\quad\forall\alpha, 
\]
as well as Theorem 5.11 of~\cite{CR1}, one gets:
\begin{align*}
\Theta^P_{\beta\alpha}(g)&=\phi_P\!\left(\sigma_{\beta}^1(t_{\calG}(g)),g\sigma_{\alpha}^1(s_{\calG}(g))\right)=\\
&=\phi_Q\!\left(\Sigma(\sigma_{\beta}^1(t_{\calG}(g))),\Sigma(g\sigma_{\alpha}^1(s_{\calG}(g)))\right)=\\
&=\phi_Q\!\left(\Sigma(\sigma_{\beta}^1(t_{\calG}(g))),g\Sigma(\sigma_{\alpha}^1(s_{\calG}(g)))\right)=\\
&=\phi_Q\!\left(\sigma_{\beta}^2(t_{\calG}(g))\Sigma_\beta(t_{\calG}(g)),g\sigma_{\alpha}^2(s_{\calG}(g))\Sigma_\alpha(s_{\calG}(g))\right)=\\
&=\left(\Sigma_\beta(t_{\calG}(g))\right)^{-1}\phi_Q\!\left(\sigma_{\beta}^2(t_{\calG}(g)),g\sigma_{\alpha}^2(s_{\calG}(g))\right) \Sigma_\alpha(s_{\calG}(g))=\\
&=\left(\Sigma_\beta(t_{\calG}(g))\right)^{-1}\Theta_{\beta\alpha}^Q(g)\Sigma_\alpha(s_{\calG}(g)),\quad \forall g\in \calG_{\alpha,\beta}.
\end{align*}
\end{proof}
The previous proposition motivates therefore the following
\begin{Def}\label{def-loceqgenmor}
Let $\calG$ and $\calH$ two groupoids.
Let $\left(\mathfrak{U},\varepsilon_{\alpha}^1,\Phi_{\alpha\beta}^1\right)$ and $\left(\mathfrak{U},\varepsilon_{\alpha}^2,\Phi_{\alpha\beta}^2\right)$ two local trivializing data over $X_{\calG}$ with values in $\calH$, with the same open covering, and let $\Theta$, resp. $\Eta$, a local generalized morphism in the sense of Definition~\ref{def-locgenmor} subordinate to the local trivializing data $\left(\mathfrak{U},\varepsilon_{\alpha}^1,\Phi_{\alpha\beta}^1\right)$, resp.\ $\left(\mathfrak{U},\varepsilon_{\alpha}^2,\Phi_{\alpha\beta}^2\right)$.

A {\em local equivalence} $\Sigma$ between $\Theta$ and $\Eta$ consists of a family of maps $\Sigma_\alpha$ from $U_\alpha$ to $\calH$, such that the following requirements hold:
\begin{itemize}
\item[a)] $\Sigma_\alpha$ puts the local momenta $\varepsilon_{\alpha}^1$ and $\varepsilon_{\alpha}^2$ in relationship as follows:
\[
t_{\calH}\circ\Sigma_\alpha=\varepsilon_{\alpha}^2,\quad s_{\calG}\circ \Sigma_{\alpha}=\varepsilon_{\alpha}^1.
\] 
\item[b)] For any two nontrivially intersecting open sets $U_\alpha$ and $U_\beta$, the {\em nonabelian {\v C}ech cocycles} $\Phi_{\alpha\beta}^1$ and $\Phi_{\alpha\beta}^2$ are cohomologous w.r.t.\ $\Sigma$:
\[
\Sigma_\beta(x)=\Phi_{\beta\alpha}^2(x)\Sigma_\alpha(x)\Phi_{\alpha\beta}^1(x),\quad x\in U_{\alpha\beta}.
\]
\item[c)] For any two open sets $U_\alpha$ and $U_\beta$, consider the local component $\calG_{\alpha,\beta}$; then, the local generalized morphisms $\Theta_{\beta\alpha}$ and $\Eta_{\beta\alpha}$ are ``cohomologous'' w.r.t.\ $\Sigma$ as follows:
\[
\Theta_{\beta\alpha}(g)=\left(\Sigma_\beta(t_{\calG}(g))\right)^{-1}\Eta_{\beta\alpha}(g)\Sigma_\alpha(s_{\calG}(g)),\quad \forall g\in \calG_{\alpha,\beta}.
\]
\end{itemize}
\end{Def}
\begin{Rem}
Notice that Condition $b)$ makes Condition $c)$ compatible with Equation (\ref{eq-relgentrans}) in Definition~\ref{def-locgenmor}: namely,
\begin{align*}
\Theta_{\beta\alpha}(g)&=\Phi_{\beta\delta}^1(t_{\calG}(g))\Theta_{\delta\gamma}(g)\Phi_{\gamma\alpha}^1(s_{\calG}(g))=\\
&=\Phi_{\beta\delta}^1(t_{\calG}(g))\left(\Sigma_\delta(t_{\calG}(g))\right)^{-1}\Eta_{\delta\gamma}(g)\Sigma_\gamma(s_{\calG}(g))\Phi_{\gamma\alpha}^1(s_{\calG}(g))=\\
&=\left(\Sigma_\delta(t_{\calG}(g))\Phi_{\delta\beta}^1(t_{\calG}(g))\right)^{-1}\Eta_{\delta\gamma}(g)\left(\Sigma_\gamma(s_{\calG}(g))\Phi_{\gamma\alpha}^1(s_{\calG}(g))\right)=\\
&=\left(\Phi_{\delta\beta}^2(t_{\calG}(g))\Sigma_{\beta}(t_{\calG}(g))\right)^{-1}\Eta_{\delta\gamma}(g)\left(\Phi_{\gamma\alpha}^2(s_{\calG}(g))\Sigma_\alpha(s_{\calG}(g))\right)=\\
&=\left(\Sigma_{\beta}(t_{\calG}(g))\right)^{-1}\Phi_{\beta\delta}^2(t_{\calG}(g))\Eta_{\delta\gamma}(g)\Phi_{\gamma\alpha}^2(s_{\calG}(g))\Sigma_\alpha(s_{\calG}(g))=\\
&=\left(\Sigma_{\beta}(t_{\calG}(g))\right)^{-1}\Eta_{\beta\alpha}(g)\Sigma_\alpha(s_{\calG}(g)),
\end{align*}
for $U_\alpha$ and $U_\gamma$, resp.\ $U_\beta$ and $U_\delta$, intersecting nontrivially, and $g\in \calG_{\alpha,\beta}\cap\calG_{\gamma,\delta}$.
\end{Rem}
It is already known that the generalized morphisms $\Theta$ and $\Eta$, together with their respective local trivializing data, give rise, by Lemma~\ref{lem-locmortogenmor}, to generalized morphisms $P_\Theta$ and $P_{\Eta}$ respectively between $\calG$ and $\calH$; it is natural to figure out that a local equivalence $\Sigma$ between $\Theta$ and $\Eta$ should give rise to an equivalence between the associated generalized morphisms $P_\Theta$ and $P_{\Eta}$.
This is the content of the following
\begin{Thm}\label{thm-loceqtoeq}
Given two local trivializing data $\left(\mathfrak{U},\varepsilon_{\alpha}^1,\Phi_{\alpha\beta}^1\right)$ and $\left(\mathfrak{U},\varepsilon_{\alpha}^2,\Phi_{\alpha\beta}^2\right)$ over $X_{\calG}$ with values in $\calH$, with the same open covering $\mathfrak{U}$, two local generalized morphisms $\Theta$ and $\Eta$ subordinate to $\left(\mathfrak{U},\varepsilon_{\alpha}^1,\Phi_{\alpha\beta}^1\right)$ and $\left(\mathfrak{U},\varepsilon_{\alpha}^2,\Phi_{\alpha\beta}^2\right)$ respectively, a local equivalence $\Sigma$ between $\Theta$ and $\Eta$ in the sense of Definition~\ref{def-loceqgenmor}, there is an equivalence $\Sigma$ between the generalized morphisms $P_\Theta$ and $P_{\Eta}$. 
\end{Thm}
\begin{proof}
By Theorem~\ref{thm-loctoglobmor}, it is already known that there is a morphism $\Sigma$ from $P_\Theta$ to $P_{\Eta}$ (i.e.\ a generalized gauge transformation on $P_\Theta\odot P_{\Eta}$ with values in $\calH$); it remains therefore to prove that it is $\calG$-equivariant.
Recall the construction of $P_\Theta$ and $P_{\Eta}$ from Subsection~\ref{ssec-locgentogenmor}:
\[
P_\Theta=\coprod_\alpha \varepsilon_{\alpha}^{1*}\calU_\calH/ \sim,\quad P_{\Eta}=\coprod_\alpha\varepsilon_{\alpha}^{2*}\calU_\calH/ \sim, 
\]
and the equivalence relations $\sim$ are in both cases induced by the cocycles $\Phi_{\alpha\beta}^1$ and $\Phi_{\alpha\beta}^2$ respectively.
The left $\calG$-action on $P_\Theta$ and $P_{\Eta}$ is defined respectively by
\[
g[\alpha,x,h]\colon=\left[\beta,t_{\calG}(g),\Theta_{\beta\alpha}(g)h\right],\quad g[\alpha,x,h]\colon=\left[\beta,t_{\calG}(g),\Eta_{\beta\alpha}(g)h\right],
\]
where $s_{\calG}(g)=x\in U_\alpha$ and $U_\beta$ is chosen so, that $t_{\calG}(g)\in U_\beta$ (so, $g\in\calG_{\alpha,\beta}$); it was already shown that the definition of left $\calG$-action is independent of any choices.
Recall also the definition of $\Sigma$:
\[
\Sigma\left([\alpha,x,h]\right)=\left[\alpha,x,\Sigma_\alpha(x)h\right],\quad x\in U_\alpha,\quad t_{\calH}(h)=\varepsilon_{\alpha}^1(x).
\]
Hence,
\begin{align*}
\Sigma\!\left(g[\alpha,x,h]\right)&=\Sigma\!\left(\left[\beta,t_{\calG}(g),\Theta_{\beta\alpha}(g)h\right]\right)=\\
&=\left[\beta,t_{\calG}(g),\Sigma_\beta(t_{\calG}(g))\Theta_{\beta\alpha}(g)h\right]=\\
&=\left[\beta,t_{\calG}(g),\Eta_{\beta\alpha}(g)\Sigma_\alpha(s_{\calG}(g))h\right]=\\
&=g\left[\alpha,x,\Sigma_\alpha(x)h\right]=\\
&=g\Sigma\!\left([\alpha,x,h]\right),\quad s_{\calG}(g)=x\in U_\alpha,
\end{align*}
which shows the $\calG$-equivariance of $\Sigma$.
The compatibility between Condition $c)$ in Definition~\ref{def-loceqgenmor} and Equation (\ref{eq-relgentrans}) in Definition~\ref{def-locgenmor} makes the above computation independent of any choice.
\end{proof}
The results of Theorem~\ref{thm-loceqtoeq} and the previous computations can be resumed as follows:

\fbox{\parbox{12cm}{\bf There is a one-to-one correspondence between equivalences of generalized morphisms from a groupoid $\calG$ to $\calH$ in the sense of~\cite{Moer2},~\cite{CR1}, and local equivalences between local generalized morphisms from $\calG$ to $\calH$ in the sense of Definition~\ref{def-loceqgenmor}.}}
 
\subsubsection{An example of equivalences between generalized morphisms: the case of action groupoids}\label{sssec-exaeqgenmor}
Consider two Lie groups $G$ and $H$, and two manifolds $X$ and $Y$, such that $G$ operates from the left on $X$ and $H$ from the left on $Y$.
The results of Subsubsection~\ref{sssec-genmoractgroup} imply that any generalized morphism from the action groupoid $G\ltimes X$ to $H\ltimes Y$ is, in the language introduced in the very same subsubsection, a $G$-equivariant, $Y$-pointed principal $H$-bundle $P$ over $X$, and the $X$-point corresponds in this case to a $G$-equivariant global section of the associated bundle $P\times_H Y$.

Consider thus two local trivializing data $\left(\mathfrak{U},\varepsilon_\alpha^i,\Phi_{\alpha\beta}^i\right)$, $i=1,2$, and local generalized morphism $\Theta$, $\Eta$, subordinate to $\left(\mathfrak{U},\varepsilon_\alpha^1,\Phi_{\alpha\beta}^1\right)$ and $\left(\mathfrak{U},\varepsilon_\alpha^2,\Phi_{\alpha\beta}^2\right)$ respectively; finally, consider a local equivalence $\Sigma$ between them in the sense of Definition~\ref{def-loceqgenmor}.
As such, it is possible to decompose the maps $\Sigma_\alpha$ as follows:
\[
\Sigma_\alpha(x)\colon=\left(\Sigma_\alpha^H(x),\Sigma_\alpha^Y(x)\right),\quad \forall x\in U_\alpha,
\]
where
\[
\Sigma_\alpha^H\colon U_\alpha\to H,\quad \Sigma_\alpha^Y\colon U_\alpha\to Y.
\]
First of all, I am going to inspect the compatibility condition with the local momenta:
\begin{align*}
\left(s_{H\ltimes Y}\circ\Sigma_\alpha\right)\!(x)&=s_{H\ltimes Y}\!\left(\Sigma_\alpha^H(x),\Sigma_\alpha^Y(x)\right)=\\
&=\Sigma_\alpha^Y(x)=\\
&=\varepsilon_\alpha^1(x),\quad \forall x\in U_\alpha,
\end{align*}
and
\begin{align*}
\left(t_{H\ltimes Y}\circ \Sigma_\alpha\right)\!(x)&=t_{H\ltimes Y}\!\left(\Sigma_\alpha^H(x),\Sigma_\alpha^Y(x)\right)=\\
&=\Sigma_\alpha^H(x)\Sigma_\alpha^Y(x)=\\
&=\varepsilon_\alpha^2(x),\quad\forall x\in U_\alpha,
\end{align*}
whence it follows:
\[
\Sigma_\alpha^H(x)\varepsilon_\alpha^1(x)=\varepsilon_\alpha^2(x).
\]
By the same arguments used in Paragraph~\ref{par-actgroup}, the nonabelian cohomological condition I translates into
\[
\Sigma_\beta^H(x)=\Phi_{\beta\alpha}^{H,2}(x)\Sigma_\alpha^H(x)\Phi_{\alpha\beta}^{1,H}(x),\quad \forall x\in U_{\alpha\beta}
\]
(having decomposed the cocycles $\Phi_{\alpha\beta}^i$ as in Subsubsection~\ref{sssec-actiongroupoid}), which tells us that the nonabelian {\v C}ech cocycles $\Phi_{\alpha\beta}^{H,i}$ over $X$ with values in $H$ are cohomologous, hence they give rise, by classical results, to isomorphic principal $H$-bundles $P_i$, the isomorphism being realized locally by the nonabelian $0$-cocycle $\Sigma_\alpha^H$.

Let me now write down the nonabelian cohomological condition II explicitly, after having decomposed $\Theta_{\beta\alpha}$ and $\Eta_{\beta\alpha}$ as in Subsubsection~\ref{sssec-actiongroupoid} (writing down only the $H$-component):
\begin{align*}
\theta_{\beta\alpha}(g,x)=\left(\Sigma_\beta^H(gx)\right)^{-1}\eta_{\beta\alpha}(g,x)\Sigma_\alpha^H(x),\quad \forall g\in G,\quad x\in U_\alpha,\quad gx\in U_\beta.
\end{align*}
This cocycle condition translates immediately into the fact that the isomorphism $\Sigma$ preserves also the left $G$-action on both $P_i$; moreover, $\Sigma$ gives rise also to an isomorphism of associated bundles $P_i\times_H Y$, which moreover preserves the induced left $G$-action on both of them. 
Moreover, the identity $\Sigma_\alpha^H(x)\varepsilon_\alpha^1(x)=\varepsilon_\alpha^2(x)$ means simply that $\Sigma$ maps the global section $\eta_1$ to $\eta_2$, and, since $\Sigma$ preserves the left $G$-action on $P_i\times_H Y$, $\Sigma$ preserves also the $G$-equivariance property of $\eta_1$ and $\eta_2$.
Hence

\fbox{\parbox{12cm}{\bf An equivalence between two generalized morphisms from the action groupoid $G\ltimes X$ to $H\ltimes Y$ corresponds uniquely to a morphism of $G$-equivariant, $Y$-pointed principal $H$-bundles over $X$, provided a choice of an open covering of $X$.}}

\subsection{Cohomological interpretation of local generalized morphisms: generalized morphisms as generalizations of descent data}\label{ssec-eqsheaves}
In this subsection, I discuss a more formal cohomological setting for local generalized morphisms (equivalently, for generalized morphisms) between Lie groupoids.
To fix notations once and for all in this Subsection, let $\calG$ be the ``source'' Lie groupoid, whereas $\calH$ is the ``target'' Lie groupoid of any local generalized morphism.
Let me go back to Definition~\ref{def-locgenmor} and try to reformulate in more abstract terms Equations (\ref{eq-diagmormomen}), (\ref{eq-genhomom}) and (\ref{eq-relgentrans}).
The first thing one can notice is the presence of local trivializing data on $X_{\calG}$: using the language of Subsection~\ref{ssec-cohomint}, local trivializing data, which encode the choice of an open covering of $X_\calG$, are in one-to-one correspondence with $1$-cocycles over $X_\calG$ with values in $\calS_{X_\calG,\calH}$, the sheaf of Lie groupoids over $X_\calG$ naturally associated to $\calH$.
(Let me just point out that in this Subsection, the natural sheaf of groupoids associated to a Lie groupoid $\calH$ over a manifold $M$ will be denoted $\calS_{M,\calH}$; I will explicitly label the sheaf by the base manifold, since in the subsequent computations, the base manifold is not always clear from the context.) 
Therefore, in an even more formal setting, the first ingredient should naturally be a sheaf $\calS$ of groupoids over $X_\calG$, an open covering $\mathfrak{U}$ of $X_\calG$ and a $1$-cocycle $\left(\underline{\varepsilon},\underline{\Sigma}\right)$ over $X_\calG$ with values in $\calS$, in the sense of Definition~\ref{def-1coch}.

The first important fact is that an open covering $\mathfrak{U}$ of $X_\calG$ defines three different open coverings of $\calG$, which I call the {\em source covering}, the {\em target covering} and the {\em source-target covering} of $\calG$ w.r.t.\ $\mathfrak{U}$, and denote by $\mathfrak{U}^s$, $\mathfrak{U}^t$ and $\mathfrak{U}^{s,t}$ respectively:
\begin{align*}
U_\alpha^s&\colon=s_{\calG}^{-1}\!\left(U_\alpha\right),\quad\forall\alpha,\\
U_\alpha^t&\colon=t_{\calG}^{-1}\!\left(U_\alpha\right),\quad\forall\alpha,\\
U_{\alpha,\beta}^{s,t}&\colon=\left(s_{\calG}\times t_{\calG}\right)^{-1}\!\left(U_\alpha\times U_\beta\right),\quad\forall\alpha,\beta.
\end{align*} 
Notice that the source-target covering is labelled by {\em two} indices; moreover, it is clear that
\[
U_{\alpha,\beta}^{s,t}=U_\alpha^s\cap U_\beta^t=\calG_{\alpha,\beta},
\]
with the notations of Subsection~\ref{ssec-locgenmor}, as was already mentioned. 
Notice also that the index set of the source-target covering $\mathfrak{U}^{s,t}$ is the Cartesian product of the index set $A$ of the covering $\mathfrak{U}$ with itself; there are two natural maps from $A\times A$ to $A$, namely the two projections $p_i$, $i=1,2$.
The projection $p_1$, resp.\ $p_2$, can be viewed as a (natural) map associated to the common refinement $\mathfrak{U}^{s,t}$ of $\mathfrak{U}^s$, resp.\ $\mathfrak{U}^t$. 

Let me now consider the three commutative diagrams of (\ref{eq-diagmormomen}).
As mentioned before, a local generalized morphism from $\calG$ to $\calH$ consists of local trivializing data $\left(\mathfrak{U},\varepsilon_\alpha,\Phi_{\alpha\beta}\right)$ over $X_\calG$ with values in $\calH$; these data correspond uniquely to an element $\left(\underline{\varepsilon},\underline{\Phi}\right)$ of $Z^1\!\left(\mathfrak{U},\calS_{X_\calG,\calH}\right)$, by the arguments of Subsection~\ref{ssec-cohomint}.
Consider now the surjective submersions $s_\calG$ and $t_\calG$ from $\calG$ to $X_\calG$: the pull-backs of the local momenta $\varepsilon_\alpha$ w.r.t.\ the source and target maps of $\calG$ define smooth maps from $U_\alpha^s$ and $U_\alpha^t$ respectively to $X_\calH$, whereas the pull-backs of the components $\Phi_{\alpha\beta}$ w.r.t.\ the $s_\calG$ and $t_\calG$ define smooth maps on $U_\alpha^s\cap U_\beta^s$ and $U_\alpha^t\cap U_\beta^t$.
Since the pull-backs w.r.t.\ $s_\calG$ and $t_\calG$ obviously commute with source, target, unit and inversion maps and product of the canonical sheaf $\calS_{\calG,\calH}$ of groupoids over $\calG$ associated to $\calH$, one has
\begin{Lem}\label{lem-pullbcocyc}
Writing $s_\calG^*\!\left(\underline{\varepsilon},\underline{\Phi}\right)$, resp.\ $t_\calG^*\!\left(\underline{\varepsilon},\underline{\Phi}\right)$, for the $1$-cochains over $\calG$ with values in $\calS_{\calG,\calH}$ w.r.t.\ the covering $\mathfrak{U}^s$, resp.\ $\mathfrak{U}^t$, obtained by pulling back the $1$-cocycle $\left(\underline{\varepsilon},\underline{\Phi}\right)$ w.r.t.\ $s_\calG$, resp\ $t_\calG$, one gets:
\begin{align*}
s_\calG^*\!\left(\underline{\varepsilon},\underline{\Phi}\right)&\in Z^1\!\left(\mathfrak{U}^s,\calS_{\calG,\calH}\right),\\
t_\calG^*\!\left(\underline{\varepsilon},\underline{\Phi}\right)&\in Z^1\!\left(\mathfrak{U}^t,\calS_{\calG,\calH}\right).
\end{align*} 
\end{Lem}
Let me then compute the image of both $1$-cocycles $s_\calG^*\!\left(\underline{\varepsilon},\underline{\Phi}\right)$, resp.\ $t_\calG^*\!\left(\underline{\varepsilon},\underline{\Phi}\right)$ w.r.t.\ the map $p_1^*$, resp.\ $p_2^*$ (see Subsection~\ref{ssec-cohomint} for more details). 
I begin by computing the respective restrictions of the $0$-cocycles $s_\calG^*\underline{\varepsilon}$ and $t_\calG^*\underline{\varepsilon}$: in fact, by the very definition, one has
\begin{align*}
p_1^*\!\left(s_\calG^*\underline{\varepsilon}\right)_{\alpha\beta}&=s_\calG^*\!\left(\underline{\varepsilon}\right)_\alpha=\\
&=\left(\varepsilon_\alpha\circ s_\calG\right)\arrowvert_{U_{\alpha,\beta}^{s,t}},\\
p_2^*\!\left(t_\calG^*\underline{\varepsilon}\right)_{\alpha\beta}&=t_\calG^*\!\left(\underline{\varepsilon}\right)_\beta=\\
&=\left(\varepsilon_\beta\circ t_\calG\right)\arrowvert_{U_{\alpha,\beta}^{s,t}}.
\end{align*}
Choosing now two nontrivially intersecting open sets $U_{\alpha,\beta}^{s,t}$ and $U_{\gamma,\delta}^{s,t}$ in the source-target covering of $\calG$, let me compute the respective restrictions of $s_\calG^*\underline{\Phi}$ and $t_\calG^*\underline{\Phi}$ w.r.t.\ $p_1^*$ and $p_2^*$ respectively:
\begin{align*}
p_1^*\!\left(s_\calG^*\underline{\Phi}\right)_{\alpha\beta,\gamma\delta}&=s_\calG^*\!\left(\underline{\Phi}\right)_{\alpha\gamma}=\\
&=\left(\Phi_{\alpha\gamma}\circ s_\calG\right)\arrowvert_{U_{\alpha\gamma,\beta\delta}^{s,t}},\\
p_2^*\!\left(t_\calG^*\underline{\Phi}\right)_{\alpha\beta,\gamma\delta}&=t_\calG^*\!\left(\underline{\Phi}\right)_{\beta\delta}=\\
&=\left(\Phi_{\beta\delta}\circ t_\calG\right)\arrowvert_{U_{\alpha\gamma,\beta\delta}^{s,t}}.
\end{align*}
The components $\Theta_{\beta\alpha}$ of the local generalized morphism $\Theta$ may be thought also as the components of a $0$-cochain $\underline{\Theta}$ over $\calG$ with values in the sheaf of groupoids $\calS_{\calG,\calH}$ w.r.t.\ the covering $\mathfrak{U}^{s,t}$; moreover, the three commutative diagrams of Equation (\ref{eq-diagmormomen}) of Definition~\ref{def-locgenmor} obviously imply together that the restrictions to the covering $\mathfrak{U}^{s,t}$ of $\calG$ of the $0$-cochains $s_\calG^*\underline{\varepsilon}$ and $t_\calG^*\underline{\varepsilon}$ are respectively the source and the target $0$-cochains of $\underline{\Theta}$, using arguments of Subsection~\ref{ssec-cohomint}.

Let me then consider Equation (\ref{eq-relgentrans}): it can be rewritten as follows
\[
\Theta_{\beta\alpha}(g)\Phi_{\alpha\gamma}(s_{\calG}(g))\left(\Theta_{\delta\gamma}(g)\right)^{-1}=\Phi_{\beta\delta}(t_\calG(g)),\quad\forall g\in U_{\alpha,\beta}^{s,t}\cap U_{\gamma,\delta}^{s,t},
\]
which can be further rewritten as (according to the above computations and recalling the arguments of Remark~\ref{rem-cohomasleftspace})
\begin{equation}\label{eq-restrcohom}
\underline{\Theta} p_1^*\!\left(s_\calG^*\!\left(\underline{\varepsilon},\underline{\Phi}\right)\right)=p_2^*\!\left(t_\calG^*\!\left(\underline{\varepsilon},\underline{\Phi}\right)\right),
\end{equation}
borrowing notations also from Subsection~\ref{ssec-cohomint}.
Therefore, combining Equation (\ref{eq-restrcohom}) with Lemma~\ref{lem-pullbcocyc} gives the following
\begin{Lem}\label{lem-cobound1}
Given an open covering $\mathfrak{U}$ of $X_\calG$, the manifold of objects of a Lie groupoid $\calG$, a local generalized morphism $\Theta$ from $\calG$ to a Lie groupoid $\calH$ in the sense of Definition~\ref{def-locgenmor} is the same as a $1$-cocycle $\left(\underline{\varepsilon},\underline{\Phi}\right)$ in $Z^1\!\left(\mathfrak{U},\calS_{X_\calG,\calH}\right)$, such that the two $1$-cocycles $s_\calG^*\!\left(\underline{\varepsilon},\underline{\Phi}\right)$ and $t_\calG^*\!\left(\underline{\varepsilon},\underline{\Phi}\right)$, respectively in $Z^1\!\left(\mathfrak{U^s},\calS_{\calG,\calH}\right)$ and $Z^1\!\left(\mathfrak{U}^t,\calS_{\calG,\calH}\right)$, are cohomologous when restricted to $Z^1\!\left(\mathfrak{U}^{s,t},\calS_{\calG,\calH}\right)$ w.r.t.\ to a $0$-cochain $\underline{\Theta}$, whose components are given by the components of the local generalized morphism $\Theta$.

Thus, using the results of Subsection~\ref{ssec-cohomint}, at a global level, a local generalized morphism $\Theta$ from $\calG$ to $\calH$ is equivalent to a principal $\calH$-bundle $P$ over $X_\calG$, such that there is a morphism $\Theta$ of $\calH$-principal bundles from $s_\calG^*P$ to $t_\calG^*P$.
\end{Lem}
\begin{Rem}
Given a generalized morphism $P$ from $\calG$ to $\calH$, in the sense of Definition~\ref{def-genmorph}, notice that the pull-back of $P$ w.r.t.\ the source map is the manifold of arrows of the groupoid version of the action groupoid; this plays an important in the characterization of principal $G$-bundles over groupoids $\Gamma$ of~\cite{L-GTX}.
\end{Rem}
It remains to inspect the meaning in this setting of Equation (\ref{eq-genhomom}) of Definition~\ref{def-locgenmor}; clearly, such an equation gives further informations about the bundle morphism $\Theta$. 
The first thing to notice is that, using the language from the first chapter of~\cite{L-GTX}, the Lie groupoid $\calG$ gives rise to a {\em simplicial manifold} $\calG_\bullet$; presently, I am interested only in the piece of $\calG_\bullet$ of degree $2$ (therefore, I will not define all now components of $\calG_\bullet$, deserving to it some attention later), which is the set of {\em composable arrows of $\calG$}, where multiplication in a groupoid makes sense:
\[
\calG_2=\left\{(g_1,g_2)\in \calG^2\colon s_\calG(g_1)=t_\calG(g_2)\right\}.
\]
There are three natural {\em face maps} (using again the language of simplicial manifolds) from $\calG_2$ to $\calG_1\colon=\calG$, namely $\pr_i$, $i=1,2$, and $\mu$, where $\pr_i$ is the projection onto the $i$-component and $\mu$ is the multiplication map; notice the following obvious identities, which follow from the groupoid axioms:
\begin{equation}\label{eq-pullident}
t_\calG\circ \pr_1=\mu\circ \pr_1,\quad s_\calG\circ\pr_2=s_\calG\circ\mu,\quad t_\calG\circ \pr_2=s_\calG\circ pr_1.
\end{equation}
The second thing to observe is that Equation (\ref{eq-genhomom}) makes sense on any open subset of $\calG_2$ of the form $U_{\alpha,\beta}^{s,t}\times U_{\beta,\gamma}^{s,t}$, intersected with $\calG_2$.
In fact, the collection of all such open sets, which I denote by $\mathfrak{U}^{s,t,2}$, is an open covering of $\calG_2$: namely, consider a pair $(g_1,g_2)$ in $\calG_2$, then, since $\mathfrak{U}$ is an open covering of $X_\calG$, it follows that there are indices $\alpha$, $\beta$ and $\gamma$ such that $s_\calG(g_2)$ is in $U_\alpha$, $s_\calG(g_1)=t_\calG(g_2)$ is in $U_\beta$ and $t_\calG(g_1)$ is in $U_\gamma$, which implies the claim.
There are three maps from the the triple product $A^3$ of the index set of the open covering $\mathfrak{U}$ onto $A^2$, namely $p_{ij}$, $1\leq i<j\leq 3$; with this in mind, Equation (\ref{eq-genhomom}) can be rewritten as
\begin{equation}\label{eq-compcocyc} 
p_{13}^*\!\left(\mu^*\underline{\Theta}\right)=p_{23}^*\!\left(\pr_1^*\underline{\Theta}\right)p_{12}^*\!\left(\pr_2^*\underline{\Theta}\right)
\end{equation}
as $0$-cochains over $\calG_2$ with values in $\calS_{\calG_2,\calH}$ w.r.t.\ the covering $\mathfrak{U}^{s,t,2}$; to be precise, the pull-backs of the $0$-cocycle $\underline{\Theta}$ w.r.t.\ $\mu$ and $\pr_i$ are $0$-cocycles with values in the sheaf $\calS_{\calG_2,\calH}$, but w.r.t.\ to different coverings of $\calG_2$, which admit a common refinement, which is in fact $\mathfrak{U}^{s,t,2}$.

Furthermore, the $1$-cocycles $s_\calG^*\!\left(\underline{\varepsilon},\underline{\Phi}\right)$ and $t_\calG^*\!\left(\underline{\varepsilon},\underline{\Phi}\right)$ can be further pulled back to give rise to six cocycles in $Z^1\!\left(\mathfrak{U}^{s,t,2},\calS_{\calG_2,\calH}\right)$, namely:
\begin{align*}
&\pr_1^*\!\left(s_\calG^*\!\left(\underline{\varepsilon},\underline{\Phi}\right)\right),\quad \pr_1^*\!\left(t_\calG^*\!\left(\underline{\varepsilon},\underline{\Phi}\right)\right),\quad \pr_2^*\!\left(s_\calG^*\!\left(\underline{\varepsilon},\underline{\Phi}\right)\right),\\
&\pr_2^*\!\left(t_\calG^*\!\left(\underline{\varepsilon},\underline{\Phi}\right)\right),\quad \mu^*\!\left(s_\calG^*\!\left(\underline{\varepsilon},\underline{\Phi}\right)\right),\quad \mu^*\!\left(t_\calG^*\!\left(\underline{\varepsilon},\underline{\Phi}\right)\right),
\end{align*}
Notice that each of these six $1$-cocycles takes value in the same sheaf $\calS_{\calG_2,\calH}$; but any is defined w.r.t.\ a different covering of $\calG_2$.
E.g., the first $1$-cocycle is defined w.r.t.\ the covering of $\calG_2$ with elements 
\[
\pr_1^{-1}\!\left(U_\alpha^s\right),\quad \forall\alpha,
\]
i.e.\ open subsets of $\calG_2$ of all pairs $(g_1,g_2)$ such that $s_\calG(g_1)\in U_\alpha$.
It is clear that $\mathfrak{U}^{s,t,2}$ is a common refinement of $\mathfrak{U}^{s,\pr_1}$ and of all the open coverings w.r.t.\ which the remaining pulled-back $1$-cocycles are defined.
For the sake of simplicity, I will omit all restriction maps to the common refinement $\mathfrak{U}^{s,t,2}$ of $\calG_2$.

Notice also that the first $1$-cocycle is equal to the fourth one, the second one is equal to the sixth one and the third one is equal to the fifth one, by Equations (\ref{eq-pullident}); at a global level, it means that there are three identities between principal $\calH$-bundles over $\calG_2$:
\[
\pr_1^*(s_\calG^*P)= \pr_2^*(t_\calG^*P),\quad \pr_1^*(t_\calG^*P)= \mu^*(t_\calG^*P),\quad \pr_2^*(s_\calG^*P)= \mu^*(s_\calG^*P),
\] 
for $P$ being the principal $\calH$-bundle over $X_\calG$ corresponding to the $1$-cocycle $\left(\underline{\varepsilon},\underline{\Phi}\right)$.

\begin{Lem}\label{lem-cobound2}
The following coboundary relations hold in $Z^1\!\left(\mathfrak{U}^{s,t,2},\calS_{\calG_2,\calH}\right)$:
\begin{align*}
\pr_1^*\!\left(t_\calG^*\!\left(\underline{\varepsilon},\underline{\Phi}\right)\right)&=\left(\pr_1^*\underline{\Theta}\right)\pr_1^*\!\left(s_\calG^*\!\left(\underline{\varepsilon},\underline{\Phi}\right)\right),\\
\pr_2^*\!\left(t_\calG^*\!\left(\underline{\varepsilon},\underline{\Phi}\right)\right)&=\left(\pr_2^*\underline{\Theta}\right)\pr_2^*\!\left(s_\calG^*\!\left(\underline{\varepsilon},\underline{\Phi}\right)\right),\\
\mu^*\!\left(t_\calG^*\!\left(\underline{\varepsilon},\underline{\Phi}\right)\right)&=\left(\mu^*\underline{\Theta}\right)\mu^*\!\left(s_\calG^*\!\left(\underline{\varepsilon},\underline{\Phi}\right)\right).
\end{align*}
\end{Lem}
\begin{proof}
The proof is a trivial consequence of Equation (\ref{eq-restrcohom}), by taking pull-backs and restrictions.
\end{proof}
Notice also the following coboundary relations:
\begin{align*}
\pr_1^*\!\left(s_\calG^*\!\left(\underline{\varepsilon},\underline{\Phi}\right)\right)&=\left(\pr_2^*\underline{\Theta}\right)\pr_2^*\!\left(s_\calG^*\!\left(\underline{\varepsilon},\underline{\Phi}\right)\right),\\
\mu^*\!\left(t_\calG^*\!\left(\underline{\varepsilon},\underline{\Phi}\right)\right)&=\left(\pr_1^*\underline{\Theta}\right)\pr_1^*\!\left(s_\calG^*\!\left(\underline{\varepsilon},\underline{\Phi}\right)\right),\\
\mu^*\!\left(t_\calG^*\!\left(\underline{\varepsilon},\underline{\Phi}\right)\right)&=\left(\mu^*\underline{\Theta}\right)\pr_2^*\!\left(s_\calG^*\!\left(\underline{\varepsilon},\underline{\Phi}\right)\right).
\end{align*}
At a global level, it means that there are morphisms of principal $\calH$-bundles over $\calG_2$, namely:
\begin{align*}
\pr_2^*(s_\calG^*P)=\mu^*(s_\calG^*P)&\overset{\pr_2^*\Theta}\cong \pr_2^*(t_\calG^*P)=\pr_1^*(s_\calG^*P),\\
\pr_2^*(t_\calG^*P)=\pr_1^*(s_\calG^*P)&\overset{\pr_1^*\Theta}\cong \pr_1^*(t_\calG^*P)=\mu^*(t_\calG^*P),\\
\pr_2^*(s_\calG^*P)=\mu^*(s_\calG^*P)&\overset{\mu^*\Theta}\cong \mu^*(t_\calG^*P)=\pr_1^*(t_\calG^*P).
\end{align*}
Combining the above coboundary relations, provided one restricts the $1$-cocycles and $0$-cochains to the common refinement $\mathfrak{U}^{s,t,2}$, one gets the following
\begin{Lem}\label{lem-cobound3}
Equation (\ref{eq-genhomom}) satisfied by the components of the local generalized morphism $\Theta$ is equivalent to the following compatibility condition for $\calH$-bundle morphisms over $\calG_2$:
\begin{equation}\label{eq-globgenhom}
\mu^*\Theta=\pr_1^*\Theta\circ\pr_2^*\Theta. 
\end{equation}
\end{Lem}
Putting Lemmata~\ref{lem-cobound1},~\ref{lem-cobound2} and~\ref{lem-cobound3} together, one gets the following 
\begin{Thm}\label{thm-genhomdescdata}
Given two Lie groupoids $\calG$ and $\calH$, local generalized morphisms $\Theta$ from $\calG$ to $\calH$ are in one-to-one correspondence with principal bundles $P$ over $X_\calG$ with structure groupoid $\calH$, such that there is a morphism $\Theta$ of principal $\calH$-bundles between $s_\calG^*P$ and $t_\calG^*P$, inducing in turn compatible morphisms $\pr_2^*\Theta$ from $\pr_2^*(s_\calG^*P)$ to $\pr_1^*(s_\calG^*P)$, $\pr_1^*\Theta$ from $\pr_1^*(s_\calG^*P)$ to $\pr_1^*(t_\calG^*P)$ and $\mu^*\Theta$ from $\pr_2^*(s_\calG^*P)$ to $\pr_1^*(t_\calG^*P)$ in the sense of Equation (\ref{eq-globgenhom}).
\end{Thm}
Let me now consider, for a given open covering $\mathfrak{U}$, two local generalized morphisms $\Theta$ and $\Eta$ from $\calG$ to $\calH$, and a local equivalence $\Sigma$ between them, in the sense of Definition~\ref{def-loceqgenmor}.
The generalized morphisms $\Theta$ and $\Eta$ correspond to respective local trivializing data, which in turn give rise to $1$-cocycles $\left(\underline{\varepsilon^1},\underline{\Phi^1}\right)$ and $\left(\underline{\varepsilon^2},\underline{\Phi^2}\right)$, and it is clear, from the arguments of Subsection~\ref{ssec-cohomint}, that an equivalence $\Sigma$ between them induces a coboundary $\underline{\Sigma}$ in $C^0\!\left(\mathfrak{U},\calS_{X_\calG,\calH}\right)$ between $\left(\underline{\varepsilon^1},\underline{\Phi^1}\right)$ and $\left(\underline{\varepsilon^2},\underline{\Phi^2}\right)$.
Furthermore, one can pull the $0$-cochain $\underline{\Sigma}$ back w.r.t.\ $s_\calG$ and $t_\calG$, obtaining two $0$-cochains $s_\calG^*\underline{\Sigma}$ and $t_\calG^*\underline{\Sigma}$ in $C^0\!\left(\mathfrak{U}^s,\calS_{\calG,\calH}\right)$ and $C^0\!\left(\mathfrak{U}^t,\calS_{\calG,\calH}\right)$ respectively; it is clear that these $0$-cochains induce coboundaries between the respective pull-backs of the $1$-cocycles $\left(\underline{\varepsilon^1},\underline{\Phi^1}\right)$ and $\left(\underline{\varepsilon^2},\underline{\Phi^2}\right)$.
This means that the principal $\calH$-bundles $P_1$ and $P_2$, associated respectively to the $1$-cocycles $\left(\underline{\varepsilon^1},\underline{\Phi^1}\right)$ and $\left(\underline{\varepsilon^2},\underline{\Phi^2}\right)$, are isomorphic; obviously, their respective pull-backs $s_\calG^*P_1$ and $s_\calG^*P_2$, as well as $t_\calG^*P_1$ and $t_\calG^*P_2$, are isomorphic. 

Consider now the restrictions of the two pulled-back $0$-cochains to the common refinement of $\mathfrak{U}^s$ and $\mathfrak{U}^t$ given by the source-target covering of $\calG$; then, by recalling the definition of the restriction morphism, the third equation in Definition~\ref{def-loceqgenmor} can be rewritten as
\begin{equation}\label{eq-globgenhom1}
p_2^*\!\left(t_\calG^*\underline{\Sigma}\right)\underline{\Theta}=\underline{\Eta}p_1^*\!\left(s_\calG^*\underline{\Sigma}\right);
\end{equation}
the above equation makes sense in $C^0\!\left(\mathfrak{U}^{s,t},\calS_{\calG,\calH}\right)$, which is a groupoid, as proved in Subsection~\ref{ssec-cohomint}.
It is already known that the pull-backs $s_\calG^*P_1$ and $t_\calG^*P_1$, as well as $s_\calG^*P_2$ and $t_\calG^*P_2$, are isomorphic as $\calH$-bundles.
Equation (\ref{eq-globgenhom1}) means simply that the following square of $\calH$-bundle morphisms commute:
\begin{equation}
\begin{CD} 
s_\calG^*P_1  @>\Theta>> t_\calG^*P_1\\
@Vs_\calG^*\Sigma VV              @VV t_\calG^*\Sigma V\\
s_\calG^*P_2 @>\Eta>> t_\calG^*P_2.
\end{CD}
\end{equation} 
By the compatibility condition expressed in the remark after Definition~\ref{def-loceqgenmor} between the Conditions $b)$ and $c)$, it follows that the higher pull-backs of $\Sigma$ w.r.t.\ the three maps from $\calG_2$ to $\calG$ previously introduced induce morphisms of $\calH$-bundles over $\calG_2$, which respect the compatibility condition (\ref{eq-globgenhom}) for the pull-backs of $\Theta$ and $\Eta$.

So, the following result is a consequence of all the previous computations:
\begin{Thm}\label{thm-cohomgenmor}
Given two Lie groupoids $\calG$ and $\calH$, equivalence classes (in the sense of Subsection~\ref{ssec-morgenmor}) of generalized morphisms from $\calG$ to $\calH$ (in the sense of Definition~\ref{def-genmorph}) are represented by isomorphism classes of principal $\calH$-bundles over $X_\calG$, which induce equivalence classes of bundle morphisms from the pull-backs of such bundles w.r.t.\ $s_\calG$ and $t_\calG$; these equivalence classes of bundle morphisms obey the compatibility condition (\ref{eq-globgenhom}) as morphisms of bundles over the term of degree $2$ in the simplicial manifold $\calG_\bullet$ associated to the groupoid $\calG$. 
\end{Thm}
\begin{Exa}
Let me discuss an interesting consequence of Theorem~\ref{thm-cohomgenmor}.
Let me consider a smooth, surjective local homeomorphism $f$ from a manifold $X$ to a manifold $Y$, and consider the {\em fibred product $X\times_Y X$}, namely the submanifold of the Cartesian product of $X$ with itself consisting of all pairs $(x_1,x_2)$, such that $f(x_1)=f(x_2)$.
It is clear that $X\times_Y X$ can be viewed as a Lie subgroupoid of the product groupoid $X\times X$ of $X$, with the same source, target, unit, inversion maps and product.
Let me also consider a Lie group $G$, which can be thought of as a Lie groupoid.
Theorem~\ref{thm-cohomgenmor} allows to give an explicit description of generalized morphisms from $X\times_Y X$ to $G$.
Namely, consider the manifold of objects of $X\times_YX$, which is simply $X$ itself; then, a generalized morphism from $X\times_YX$ to $G$ consists of a principal bundle over $X$ with structure groupoid $G$, which, by the arguments of Subsubsection~\ref{sssec-liegroup}, is simply a principal $G$-bundle $Q$ over $X$.
Later, the source, resp.\ target, map of $X\times_YX$ is the projection $\pr_2$, resp.\ $\pr_1$, from $X\times_YX$ to $X$; then, Theorem~\ref{thm-cohomgenmor} implies the existence of a $G$-bundle morphism $\Theta$ from $\pr_2^*Q$ to $\pr_1^*Q$.
Furthermore, let me consider the term of degree $2$ of the simplicial manifold associated to $X\times_YX$: it is not difficult to prove that it can be identified with the triple fibred product $X\times_YX\times_YX$; similarly, the face maps from $X\times_YX\times_YX$ to $X\times_YX$ are simply given by
\[
\pr_1=\pr_{12}, \quad \pr_2=\pr_{23},\quad \mu=\pr_{13}.
\] 
The compatibility condition (\ref{eq-globgenhom}) satisfied by the pull-backs of $\theta$ w.r.t.\ $\pr_{ij}$, $1\leq i<j\leq 3$, can be written as
\begin{equation}\label{eq-descocyc}
\pr_{13}^*\Theta=\pr_{12}^*\Theta \pr_{23}^*\Theta,
\end{equation}
where $\pr_{23}^*\Theta$ is a $G$-bundle morphism from $\pr_3^*Q$ to $\pr_2^*Q$, $\pr_{12}^*\Theta$ is a $G$-bundle morphism from $\pr_2^*Q$ to $\pr_1^*Q$, and $\pr_{13}^*\Theta$ is a $G$-bundle morphism from $\pr_3^*Q$ to $\pr_1^*Q$.
Recalling e.g.\ from~\cite{Moer2} or~\cite{CR1} that any $G$-bundle morphism is an isomorphism, if I consider the inverse $\phi$ of $\Theta$ as a $G$-bundle isomorphism from $p_1^*Q$ to $p_2^*Q$, then a generalized morphism from $X\times_YX$ to $G$ is equivalent to classical {\bf $G$-descent data w.r.t.\ a smooth, surjective local homeomorphism from $X$ to $Y$}, see e.g.~\cite{Bry}, Chapter 5, Section 1 for more details on such descent data.
Equivalently, this construction leads (see also~\cite{Bry}) to a principal $G$-bundle $P$ over $Y$, such that the $G$-bundle $Q$ over $X$ is naturally identified with the pull-back of $P$ w.r.t.\ $f$.    
To any smooth, surjective local homeomorphism $f\colon X\to Y$ corresponds the category of $G$-descent data w.r.t.\ $f$, whose objects consist of $G$-bundles $Q$ over $X$, such that there is a {\bf descent morphism} from $\pr_1^*Q$ to $\pr_2^*Q$ as $G$-bundles over the fibre product $X\times_YX$, satisfying the compatibility condition (\ref{eq-descocyc}), and whose morphisms are $G$-bundle morphisms, which are compatible with the respective descent morphisms.  
By Proposition 5.1.3. of~\cite{Bry}, the pull-back $f^*$ w.r.t.\ a smooth, surjective local homeomorphism $f\colon X\to Y$ induces an equivalence of categories between the category of $G$-bundles over $Y$ and the category of $G$-descent data w.r.t.\ $f$, i.e.\
\begin{Thm}\label{thm-eqgenmordescdata}
Given a smooth, surjective local homeomorphism $f\colon X\to Y$, equivalence classes of generalized morphisms from $X\times_YX$ to $G$ are in one-to-one correspondence via pull-back w.r.t.\ $f$ to isomorphism classes of $G$-bundles over $Y$, also to isomorphism classes of $G$-descent data w.r.t.\ $f$.
\end{Thm}
Theorem~\ref{thm-eqgenmordescdata} is the main reason, why generalized morphisms can be looked at as generalizations of classical descent data for principal bundles.
\end{Exa}

\begin{Exa}
Let me just mention another consequence of Theorem~\ref{thm-cohomgenmor}; this result, although at the beginning I did not expect it, was the main motivation, coming from Topological Quantum Field Theory, to pursue the theory of principal bundles with structure groupoid.
Namely, consider the {\em fundamental groupoid $\Pi(M)$} of a smooth manifold $M$, whose manifold of objects is $M$, and, for any two objects $x$, $y$, the set of arrows from $x$ to $y$ consists of all endpoints-preserving homotopy classes of smooth curves from $x$ to $y$; multiplication is given by composition of homotopy classes.
Since only endpoints-preserving homotopy classes of curves in $M$ are considered, and since one may assume that all such curves are parametrized over the unit interval $\unint$, then the target, resp.\ source, map of $\Pi(M)$ is the well-known {\em evaluation map} at $1$, resp.\ at $0$; let me also notice that the isotropy group at $x$ corresponds to the fundamental group of $M$ based at $x$.
Let me also consider a Lie group $G$, viewed as a trivial Lie groupoid.
Then, using the results of~\cite{CR} and of Theorem~\ref{thm-cohomgenmor}, one can view isomorphism classes of flat principal $G$-bundles as {\bf equivalence classes of generalized morphisms from $\Pi(M)$ to $G$.}
I will pursue this topic extensively in a subsequent paper, where I will prove the previous claim in all the details; let me just mention the main idea, namely to formalize in a fancy way the properties of the parallel transport, and then to obtain a connection starting from such an abstract parallel transport, such that the corresponding parallel transport equals the abstract one.
The fancy version of the parallel transport is given in terms of generalized gauge transformations, following the computations in~\cite{CR}.
Let me also mention that one could consider, more generally, generalized morphism from $\Pi(M)$ to a Lie groupoid $\calG$, using also methods introduced in~\cite{CR1}: this is expected to give a meaningful concept of flat connections of principal bundles with structure groupoid $\calG$; I plan to pursue also this topic in a subsequent paper in terms of abstract parallel transport.
\end{Exa}

Let me now formulate the results expressed in Theorem~\ref{thm-cohomgenmor} in the context of nonabelian {\v C}ech cohomology, restricting here to the canonical sheaf $\calS_\calH$ of groupoids associated to a Lie groupoid $\calH$.
I plan to come to the subject in more detail in subsequent works; in particular, it would be interesting to formulate a generalization of the theory of descent in the framework of Lie groupoids for general sheaves of groupoids, motivated by the previous results.

An open covering $\mathfrak{U}$ of $X_\calG$ gives rise to an open covering of the simplicial manifold $\calG_\bullet$ associated to $\calG$, i.e.\ it is possible to construct an open covering of any $\calG_n$, and these coverings are compatible with each other by the face maps in the sense that I am going to explain below.
Consider the piece $\calG_n$ of degree $n$, which is by definition the set of neighbour-to-neighbour pairwise composable elements in $\calG^n$:
\[
\calG_n=\left\{\left(g_1,\dots,g_n\right)\in \calG^n\colon s_{\calG}(g_i)=t_{\calG}(g_{i+1}),\quad i=1,\dots,n-1 \right\}.
\]
An open covering of $\calG_n$ constructed from $\mathfrak{U}$ is given by the collection of all sets of the form
\begin{align*}
U_{\alpha_{1},\dots,\alpha_{n+1}}&\colon=\left\{\left(g_1,\dots,g_n\right)\in \calG^n\colon t_{\calG}(g_1)\in U_{\alpha_1}, s_\calG(g_{n+1})\in U_{\alpha_{n+1}},\right.\\
&\phantom{\colon=}\left.s_{\calG}(g_i)=t_\calG(g_{i+1})\in U_{\alpha_{i+1}}\right\}.
\end{align*}
It is clear that any set $U_{\alpha_1,\dots,\alpha_{n+1}}$ is open: in fact, it can be written as the intersection with $\calG_n$ of the product of open sets in $\calG^n$ of the form
\[
\prod_{i=1}^n U_{\alpha_{i},\alpha_{i+1}}^{s,t},\quad U_{\alpha_i,\alpha_{i+1}}^{s,t}=(t_\calG\times s_\calG)^{-1}\!\left(U_{\alpha_i}\times U_{\alpha_{i+1}}\right).
\]
(Notice that I have inverted the r{\^o}les of source and target maps w.r.t.\ the previous notations; I did so for consistency reasons with the notations I have chosen for the indices.)
Furthermore, it is not difficult to prove that the collection of such open sets is truly an open covering of $\calG_n$: namely, let $(g_1,\dots,g_n)$ be an element of $\calG_n$, then there exist $n+1$ indices $\alpha_1,\dots,\alpha_n$, such that $t_\calG(g_1)\in U_{\alpha_1}$, $s_{\calG}(g_n)\in U_{\alpha_{n+1}}$ and $s_{\calG}(g_i)=t_\calG(g_{i+1})\in U_{\alpha_{i+1}}$, since $\mathfrak{U}$ is an open covering of $X_\calG$.
I denote by $\mathfrak{U}^n$ the covering of $\calG_n$ previously constructed; similarly, by $\mathfrak{U}^\bullet$ I denote the collection of all open coverings $\mathfrak{U}^n$.
I now explain the compatibility conditions for $\mathfrak{U}^\bullet$. 
The simplicial structure of $\calG_\bullet$ is encoded in the existence of the so-called face maps.
There are exactly $n+1$ such face maps from $\calG_{n}$ to $\calG_{n-1}$, denoted by $\iota_{k,n}$, $k=0,\dots,n$, which are explicitly given by
\[
\iota_{k,n}\!\left(g_1,\dots,g_n\right)=\begin{cases}
\left(g_2,\dots,g_n\right),& k=0\\
\left(g_1,\dots,g_{i}g_{i+1},\dots,g_n\right),& 1\leq k\leq n-1\\
\left(g_1,\dots,g_{n-1}\right),& k=n.
\end{cases}
\]
Set $\calG_0=X_\calG$, and $\iota_{0,1}=t_\calG$ and $\iota_{1,1}=s_\calG$.
Observe that, to be more precise, a simplicial manifold possesses another set of face maps, which go in the opposite direction (the face maps $\iota_{k,n}$ encode a sort of ``homological data'', while the second set encode ``cohomological data'').
Since, for my purposes, I need only the ``homological'' face maps, I will simply neglect the ``cohomological'' face maps; nonetheless, it was worth mentioning their existence. 
All face maps must satisfy a set of compatibility conditions, for whose explicit definition I refer to~\cite{L-GTX}; as an example, all compatibility conditions for the face maps from $\calG_2$ to $\calG_1$ with the face maps from $\calG_1$ to $\calG_0=X_\calG$ are written in (\ref{eq-pullident}).
Given the open covering $\mathfrak{U}^{n-1}$ of $\calG_{n-1}$, one can construct in principle $n+1$ different open coverings of $\calG_n$ by pulling back $\mathfrak{U}^{n-1}$ w.r.t.\ to the $n+1$ face maps; I denote by $\mathfrak{U}^{k,n}$ the pull-back of $\mathfrak{U}^{n-1}$ w.r.t.\ the face map $\iota_{k,n}$, i.e.\ the collections of all open sets of the form
\[
\iota_{k,n}^{-1}\!\left(U_{\alpha_1,\dots,\alpha_n}\right),\quad (\alpha_1,\dots,\alpha_n)\in A^n,
\]  
$A$ being the index set of the open covering $\mathfrak{U}$.
There are exactly $n+1$ projections $p_{i_1\cdots i_n}$ from $A^{n+1}$ to $A^n$, for $1\leq i_1<\cdots<i_n\leq n+1$.
\begin{Lem}\label{lem-commrefsimpl}
Any projection $p_{1\cdots\widehat{k+1}\cdots n+1}$, which forgets the $k+1$-th factor, makes $\mathfrak{U}^n$, whose index set is by construction $A^{n+1}$ a refinement of $\mathfrak{U}^{k,n}$, for any $1\leq k\leq n$; hence, $\mathfrak{U}^n$ is a common refinement of all open coverings $\mathfrak{U}^{k,n}$.
\end{Lem}
\begin{proof}
Consider first, for $k=0$, the pull-back $\mathfrak{U}^{0,n}$; one has to show that 
\[
U_{\alpha_1,\dots,\alpha_{n+1}}\subseteq \iota_{0,n}\!\left(U_{\alpha_2,\dots,\alpha_{n+1}}\right),\quad \forall (\alpha_1,\dots,\alpha_{n+1})\in A^{n+1}.
\]
Namely, consider $\left(g_1,\dots,g_n\right)$ in $U_{\alpha_1,\dots,\alpha_{n+1}}$, then, by its very definition, $\iota_{0,n}(g_1,\dots,g_n)=(g_2,\dots,g_n)$, and by the construction of $\mathfrak{U}^n$, is follows immediately that $(g_2,\dots,g_n)$ belongs to $U_{\alpha_2,\dots,\alpha_{n+1}}$.
Similarly, one proves the claim for $k=n$.
Let me consider $0<k<n$; then, one has to show
\[
U_{\alpha_1,\dots,\alpha_{n+1}}\subseteq \iota_{k,n}\!\left(U_{\alpha_1,\dots,\widehat{\alpha_{k+1}},\dots,\alpha_{n+1}}\right),\quad \forall (\alpha_1,\dots,\alpha_{n+1})\in A^{n+1}.
\] 
In this case, one has
\[
\iota_{k,n}(g_1,\dots,g_n)=(g_1,\dots,g_kg_{k+1},\dots,g_n);
\]
since $t_{\calG}(g_kg_{k+1})=t_{\calG}(g_k)$ and $s_{\calG}(g_kg_{k+1})=s_{\calG}(g_{k+1})$, one has
\[
t_{\calG}(g_kg_{k+1})\in U_{\alpha_k},\quad s_{\calG}(g_kg_{k+1})\in U_{\alpha_{k+2}},
\]
which implies the claim.
\end{proof}
This is the compatibility condition for the open covering $\mathfrak{U}^{\bullet}$.

Now, consider the first nonabelian {\v C}ech cohomology group $\coh^1\!\left(\mathfrak{U},\calS_{\calG_0,\calH}\right)$ w.r.t.\ the open covering $\mathfrak{U}$, introduced in Subsection~\ref{ssec-cohomint}; its elements correspond to isomorphism classes of principal $\calH$-bundles trivialized w.r.t.\ $\mathfrak{U}$.
Then, the face map $\iota_{0,1}$, resp.\ $\iota_{1,1}$, induce by pull-back a map from $\coh^1\!\left(\mathfrak{U},\calS_{\calG_0,\calH}\right)$ to $\coh^1\!\left(\mathfrak{U}^{0,1},\calS_{\calG_1,\calH}\right)$, resp.\ $\coh^1\!\left(\mathfrak{U}^{1,1},\calS_{\calG_1,\calH}\right)$.
By Lemma~\ref{lem-commrefsimpl} and Lemma~\ref{lem-refin3} of Subsection~\ref{ssec-cohomint}, there is a map $p_1^*$, resp.\ $p_2^*$, from $\coh^1\!\left(\mathfrak{U}^{0,1},\calS_{\calG_1,\calH}\right)$, resp.\ $\coh^1\!\left(\mathfrak{U}^{1,1},\calS_{\calG_1,\calH}\right)$, to $\coh^1\!\left(\mathfrak{U}^1,\calS_{\calG_1,\calH}\right)$.
Thus, there are two different maps from $\coh^1\!\left(\mathfrak{U},\calS_{\calG_0,\calH}\right)$ to $\coh^1\!\left(\mathfrak{U}^1,\calS_{\calG_1,\calH}\right)$.
Furthermore, there are three face maps $\iota_{k,2}$ from $\calG_2$ to $\calG_1$; these induce, by pull-back, three different maps from $\coh^1\!\left(\mathfrak{U}^1,\calS_{\calG_1,\calH}\right)$ to $\coh^1\!\left(\mathfrak{U}^{k,2},\calS_{\calG_2,\calH}\right)$ respectively, for $0\leq k\leq 2$.
Also, again by Lemma~\ref{lem-commrefsimpl} and Lemma~\ref{lem-refin3} of Subsection~\ref{ssec-cohomint}, there are three maps $p_{ij}^*$ from $\coh^1\!\left(\mathfrak{U}^{k,2},\calS_{\calG_2,\calH}\right)$ to $\coh^1\!\left(\mathfrak{U}^2,\calS_{\calG_2,\calH}\right)$; the index $k$, according to the above notations, is such that $k+1$ is the missing index in $\left\{1,2,3\right\}$, when sorting out $\left\{i,j\right\}$. 
So, there are three different natural maps from $\coh^1\!\left(\mathfrak{U}^1,\calS_{\calG_1,\calH}\right)$ to $\coh^1\!\left(\mathfrak{U}^2,\calS_{\calG_2,\calH}\right)$.
This procedure can be iterated at every order $n$, thus, one gets
\begin{Prop}\label{prop-cohomsimpl}
Given two Lie groupoids $\calG$ and $\calH$, an open covering $\mathfrak{U}^0=\mathfrak{U}$ of $\calG_0$, there is a sequence of first nonabelian {\v C}ech cohomology groups $\coh^1\!\left(\mathfrak{U}^n,\calS_{\calG_n,\calH}\right)$, denoted collectively by $\coh^1\!\left(\mathfrak{U}^\bullet,\calS_{\calG_\bullet,\calH}\right)$, and, for every $n$, $n+1$ maps $\rho_{k,n}$, $0\leq k\leq n$, from $\coh^1\!\left(\mathfrak{U}^{n-1},\calS_{\calG_{n-1},\calH}\right)$ to $\coh^1\!\left(\mathfrak{U}^n,\calS_{\calG_n,\calH}\right)$; these maps depend only on the face maps of the simplicial manifold $\calG_\bullet$ and satisfy, by their very construction, compatibility conditions contravariant to the ones satisfied by the face maps. 
\end{Prop}
\begin{proof}
It is clear from the construction; that in fact the maps $\rho_{k,n}$ depend only on the corresponding face maps is a consequence of Lemma~\ref{lem-refin4} of Subsection~\ref{ssec-cohomint}.
Taking the pull-backs of the compatibility conditions between face maps gives the compatibility conditions between the maps $\rho_{k,n}$.
\end{proof}
Consider now an equivalence class of a local generalized morphism $\Theta$; by Theorem~\ref{thm-cohomgenmor}, the cohomology class $\left[\underline{\varepsilon},\underline{\Phi}\right]$ in $\coh^1\!\left(\mathfrak{U}^0,\calS_{\calG_0,\calH}\right)$ is such that its images w.r.t.\ $\rho_{0,1}$ and $\rho_{1,1}$ in $\coh^1\!\left(\mathfrak{U}^1,\calS_{\calG_1,\calH}\right)$ coincide.
In fact, the restriction of $\Theta$, viewed as a $0$-cochain, provide the desired identification.
Moreover, the further images w.r.t.\ $\rho_{k,2}$ of this cohomology class also coincide, again by Theorem~\ref{thm-cohomgenmor}.
Iterating this procedure, one verifies that the images w.r.t.\ $\rho_{k,n}$ of the higher images of the cohomology class corresponding to the equivalence class of the local generalized morphism $\Theta$ always coincide.

Consider now a refinement $\mathfrak{V}^0$ of $\calG_0$; by Lemma~\ref{lem-commrefsimpl} one has a corresponding open covering $\mathfrak{V}^\bullet$ of $\calG_\bullet$.
\begin{Lem}\label{lem-simplrefin}
Denoting by $f$ any refinement map associated to $\mathfrak{V}^0$, the Cartesian products of $f$ with itself, collectively denoted by $f^\bullet$, make $\mathfrak{V}^\bullet$ into a refinement of $\mathfrak{U}^\bullet$.
\end{Lem}
\begin{proof}
Namely, consider $\mathfrak{V}^n$, for a general index $n$; one has to show that $f^{n+1}$, the Cartesian product of $f$ with itself $n+1$ times, has the property that
\[
V_{\beta_1,\dots,\beta_{n+1}}\subseteq U_{f(\beta_1),\dots,f(\beta_{n+1})},\quad \forall(\beta_1,\dots,\beta_{n+1})\in B^{n+1},  
\]
$B$ denoting the index set of $\mathfrak{V}^0$.
In fact, by definition, if $(g_1,\dots,g_n)$ lies in $V_{\beta_1,\dots,\beta_{n+1}}$, it means that
\[
t_{\calG}(g_1)\in V_{\beta_1},\quad s_{\calG}(g_i)=t_{\calG}(g_{i+1})\in V_{\beta_i},\quad 2\leq i\leq n,\quad s_{\calG}(g_{n+1})\in V_{\beta_{n+1}}.
\]
But since $V_\beta\subseteq U_{f(\beta)}$, for any $\beta\in B$, it follows
\[
t_{\calG}(g_1)\in U_{f(\beta_1)},\quad s_{\calG}(g_i)=t_{\calG}(g_{i+1})\in U_{f(\beta_i)},\quad 2\leq i\leq n,\quad s_{\calG}(g_{n+1})\in U_{f(\beta_{n+1})},
\]
which is equivalent to the result.
\end{proof}
If we are given a refinement $\mathfrak{V}^0$ of $\mathfrak{U}^0$, then this induces, by Lemma~\ref{lem-refin4} of Subsection~\ref{ssec-cohomint}, a well-defined restriction map in nonabelian cohomology; furthermore, by Lemma~\ref{lem-simplrefin}, it induces a well-defined restriction map from the nonabelian cohomology group $\coh^1\!\left(\mathfrak{U}^\bullet,\calS_{\calG_\bullet,\calH}\right)$ to $\coh^1\!\left(\mathfrak{V}^\bullet,\calS_{\calG_\bullet,\calH}\right)$.
It remains to check that the restriction maps are compatible with the maps $\rho_{k,n}$, for every $n$.
\begin{Lem}\label{lem-comprefsimpl}
The following square of maps is commutative, for any index $n$ and for every $0\leq k\leq n$:
\begin{equation*}
\begin{CD} 
\coh^1\!\left(\mathfrak{U}^{n-1},\calS_{\calG_{n-1},\calH}\right)  @>\rho_{k,n}>>\coh^1\!\left(\mathfrak{U}^{n},\calS_{\calG_{n},\calH}\right)  \\
@V(f^{n-1})^* VV              @VV (f^n)^* V\\
\coh^1\!\left(\mathfrak{V}^{n-1},\calS_{\calG_{n-1},\calH}\right) @>\rho_{k,n}>>\coh^1\!\left(\mathfrak{V}^{n},\calS_{\calG_{n},\calH}\right),
\end{CD}
\end{equation*}
where $f$ is any refinement map associated to the refinement $\mathfrak{V}^0$ of $\mathfrak{U}^0$.
\end{Lem}
\begin{proof}
By Lemma~\ref{lem-refin3} of Subsection~\ref{ssec-cohomint}, it is sufficient to work at the level of representatives; then, the claim follows clearly by the constructions of the maps.
Notice that by Lemma~\ref{lem-refin4}, the vertical maps in the above commutative square do not depend on the choice of $f$.
\end{proof}
So, the ``cosimplicial'' structure of the nonabelian cohomology groups $\coh^1\!\left(\mathfrak{U}^\bullet,\calS_{\calG_\bullet,\calH}\right)$, given by the maps $\rho_{k,n}$, is preserved by the refinement procedure.
The main consequence is one can define the nonabelian cohomology sequence $\coh^1\!\left(\calG_\bullet,\calS_{\calG_\bullet,\calH}\right)$ as a direct limit of the cohomology sequence $\coh^1\!\left(\mathfrak{U}^\bullet,\calS_{\calG_\bullet,\calH}\right)$ w.r.t.\ the refinement relation for open coverings $\mathfrak{U}^0$ of $\calG_0$; moreover, the ``cosimplicial structure'' encoded in the maps $\rho_{k,n}$ descends to $\coh^1\!\left(\calG_\bullet,\calS_{\calG_\bullet,\calH}\right)$.
Thus, equivalence classes of local generalized morphisms can be characterized in nonabelian {\v C}ech cohomology also by the following
\begin{Thm}\label{thm-cechgenmor}
Given two Lie groupoids $\calG$ and $\calH$, equivalence classes of generalized morphisms from $\calG$ to $\calH$ give rise to a sequence of cohomology classes in $\coh^1\!\left(\calG\bullet,\calS_{\calG_\bullet,\calH}\right)$, which is compatible with all maps $\rho_{k,n}$.
\end{Thm}
\begin{Rem}
A natural question arises at this point, which I think should deserve some attention: equivalence classes of generalized morphisms give rise to a compatible sequence of cohomology classes; do then all such compatible sequences of cohomology classes correspond to generalized morphisms? Or are there other objects, and which properties do they have? Do they somehow ``generalize'' the properties of generalized morphisms?
\end{Rem}

 \section{The composition law for generalized morphisms}\label{sec-genmorcomp}
The aim of this Section is to construct a suitable composition law for local generalized morphisms in the sense of Definition~\ref{def-locgenmor} of Subsection~\ref{ssec-locgentogenmor}.
Due to their (tautological!) local nature of , it is not possible to compose cocycles in the obvious way; this will be particularly clear when analyzing the local structure of the composition of generalized morphisms in the sense of Definition~\ref{def-genmorph}.

Thus, the task can be divided into two pieces: in the next subsection, I discuss the composition law for generalized morphisms in the sense of Definition~\ref{def-genmorph}, following~\cite{Mrcun} and~\cite{Moer2}, and I will in particular discuss the division map of the composition of two generalized morphisms.
Many details, like the general construction of associated bundles, are skipped, deserving to them special attention later.
Using then arguments from Subsection~\ref{ssec-localsec} and Subsection~\ref{ssec-locgenmor}, I will display a formula for the composition of local generalized morphisms arising naturally from generalized morphisms.

In the subsequent Subsection, I will generalize the notion of composition of local generalized morphisms, pointing out in particular to the ``refinement trick'', where I will introduce what may be called a ``birefinement''.
Once this is known, it is then relatively easy to define the composition law of local generalized morphisms, which corresponds to the composition law for generalized morphisms.

\subsection{From generalized morphisms to local generalized morphisms: the division map of the composition of generalized morphisms}\label{ssec-gencomplaw}
From now on, to simplify notations, I will write $\calG\overset{P}\to \calH$ for a generalized morphism $\left(P,\pi,\varepsilon,X_\calG\right)$ with a left $\calG$-action, as in Definition~\ref{def-genmorph}.
Likewise, I will denote usually by $\calG\overset{\Theta^P}\to \calH$, or simply by $\Theta^P$, the equivalent local generalized morphism, discussed in Subsection~\ref{ssec-locgenmor}.
Given a local generalized morphism $\Theta$, denote by $\calG\overset{P_\Theta}\to \calH$ the equivalent generalized morphism.
It is already known that there is an equivalence
\[
\left(\calG\overset{\Theta}\to\calH\right)\leftrightarrow \left(\calG\overset{P_\Theta}\to\calH\right).
\]
Given three groupoids $\calF$, $\calG$ and $\calH$, and two generalized morphisms $\calF\overset{P}\to\calG$, $\calG\overset{Q}\to\calH$, is there a natural notion of composition? And how is this notion of composition related to the corresponding local generalized morphisms $\Theta^P$ and $\Theta^Q$?
The composition law for generalized morphisms has already been studied by many authors; I will follow here the approach and conventions of e.g.~\cite{Moer2} and~\cite{Mrcun}.
I will only sketch some details; I will but discuss in great detail the division map of the composite of the generalized morphisms $P$ and $Q$.

Given now two generalized morphisms $P$ and $Q$ as above, i.e.\ a right principal $\calG$-bundle $\left(P,\pi_1,\varepsilon_1,X_{\calF}\right)$ endowed with a compatible left $\calF$-action with momentum $\pi_1$ and a right principal $\calH$-bundle $\left(Q,\pi_2,\varepsilon_2,X_{\calG}\right)$ endowed with a left $\calG$-action with momentum $\pi_2$, define their composition $P\circ Q$ as the ``associated bundle'' $P\times_{\calG} Q$.
Let me point out that there is a natural notion of bundle associated to a principal bundle $\left(P,\pi,\varepsilon,X\right)$ with structure groupoid $\calG$, for some (left) ``representation'' of $\calG$; here, by a left representation of $\calG$ is meant a $3$-tuple $\left(\calE,J_{\calE},X_{\calG}\right)$, where $i)$ $\calE$ is a fiber bundle over $X_{\calG}$ in the usual sense, the manifold of objects of $\calG$, with natural projection $J_{\calE}$ and $ii)$ $\calE$ is a left $\calG$-space with momentum $J_{\calE}$.
I will discuss in detail the construction of associated bundle, from the global and local point of view in subsequent work.

In the present situation, given the $\calG$-bundle $\left(P,\pi_1,\varepsilon_1,X_{\calF}\right)$, let me consider the $3$-tuple $\left(Q,\pi_2,X_{\calG}\right)$ as a left $\calG$-representation; all requirements are satisfied, because $Q$ is a generalized morphism from $\calG$ to $\calH$.
It thus makes sense to consider the associated bundle $P\times_{\calG} Q$, which is a bundle over $X_{\calF}$ with typical fibre $Q$: an element thereof consists of an equivalence class of pairs of shape
\[
(p,q)\in P\times Q,\quad \varepsilon_1(p)=\pi_2(q),
\] 
and the equivalence relation is given by
\[
(p_1,q_1)\sim (p_2,q_2)\Leftrightarrow\exists g\in\calG\colon (p_1g,g^{-1}q_1)=(p_2,q_2),\quad t_{\calG}(g)=\varepsilon_1(p_1). 
\]
Since $\pi_1$ and $\varepsilon_2$ are both $\calG$-invariant, they descend to well-defined maps from $P\times_{\calG}Q$ to $X_{\calF}$ and $X_{\calH}$, which I denote by $\overline{\pi}$ and $\overline{\varepsilon}$ respectively:
\begin{equation}\label{eq-genprodmaps}
\overline{\pi}\!\left([p,q]\right)\colon=\pi_1(p),\quad \overline{\varepsilon}\!\left([p,q]\right)\colon=\varepsilon_2(q).
\end{equation}
There is a left $\calF$-action on $P\times_{\calG}Q$ with momentum $\overline{\pi}$
\begin{equation}\label{eq-genprodleft}
\calG\times_{\overline{\pi}}(P\times_{\calG}Q)\ni\left(f,[p,q]\right)\mapsto [gp,q],
\end{equation}
and a right $\calH$-action with momentum $\overline{\varepsilon}$
\begin{equation}\label{eq-genprodright}
(P\times_{\calG}Q)\times_{\overline{\varepsilon}}\calH\ni\left([p,q],h\right)\mapsto [p,qh].
\end{equation}
It is not difficult to verify that both actions (\ref{eq-genprodleft}) and (\ref{eq-genprodright}) are well-defined and that they are truly actions; moreover, it is clear that $\overline{\pi}$ is $\calH$-invariant, while $\overline{\varepsilon}$ is $\calF$-invariant.

There is a division map w.r.t.\ the right $\calH$-action, which implies immediately by Definition~\ref{def-princgroupoid} of Section~\ref{sec-divmap} that the $4$-tuple $\left(P\times_{\calG}Q,\overline{\pi},\overline{\varepsilon},X_{\calF}\right)$ is a right principal $\calH$-bundle over $X_{\calF}$: namely, the map $\phi_{Q\circ P}$ defined by
\begin{equation}\label{eq-genprodivmap}
\left([p_1,q_1],[p_2,q_2]\right)\mapsto \phi_Q\!\left(q_1,\phi_P(p_1,p_2)q_2\right)\in \calH,\quad \begin{cases}
&\varepsilon_1(p_1)=\pi_2(q_1),\\
&\varepsilon_1(p_2)=\pi_2(q_2),\\
&\pi_1(p_1)=\pi_1(p_2),
\end{cases}
\end{equation}
is the division map of the bundle $P\times_{\calG}Q$.
First of all, let me show that the map of Equation (\ref{eq-genprodivmap}) is well-defined.
I show first that the operations involved make sense:
\begin{itemize}
\item[i)] $\phi_P(p_1,p_2)$ makes sense, since $\pi_1(p_1)=\pi_1(p_2)$; 
\item[ii)] the product $\phi_P(p_1,p_2)q_2$ also makes sense, since
\[
\pi_2(q_2)=\varepsilon_1(p_2)=s_{\calG}\!\left(\phi_P(p_1,p_2)\right);
\] 
\item[iii)] the elements $q_1$ and $\phi_P(p_1,p_2)q_2$ lie in the same fiber of $\pi_2$, since
\[
\pi_2(q_1)=\varepsilon_1(p_1)=t_{\calG}\!\left(\phi_P(p_1,p_2)\right)=\pi_2\!\left(\phi_P(p_1,p_2)q_2\right).
\] 
\end{itemize}
Let me show now independence on the chosen representatives of the classes; for this purpose, I need the following general
\begin{Lem}\label{lem-invdivmap}
Given a generalized morphism $\calG\overset{Q}\to \calH$, the division map $\phi_Q$ is $\calG$-invariant, where we consider the diagonal action of $\calG$ on the fibred product $Q\odot Q$ w.r.t.\ the momentum given by the natural projection from $Q\odot Q$ to $X_{\calG}$.  
\end{Lem}
For the proof, see e.g.~\cite{Moer2} or~\cite{CR1}
\begin{Rem}
Lemma~\ref{lem-invdivmap} will play a fundamental r{\^o}le later in the local theory of Morita equivalences, in particular in the ``factorization formula'' for the division map.
\end{Rem}

It is now possible to show that the map (\ref{eq-genprodivmap}) does not depend on the choice of representatives of the two classes involved, namely, I take different representatives as follows:
\begin{multline*}
(p_1,q_1)\leadsto (p_1g_1,g_1^{-1}q_1),\quad (p_2,q_2)\leadsto (p_2g_2,g_2^{-1}q_2),\\
\varepsilon_1(p_1)=t_{\calG}(g_1),\quad \varepsilon_1(p_2)=t_{\calG}(g_2).
\end{multline*}
Then, the map (\ref{eq-genprodivmap}) behaves as follows:
\begin{align*}
\phi_{Q\circ P}\!\left([p_1g_1,g_1^{-1}q_1],[p_2g_2,g_2^{-1}q_2]\right)&=\phi_Q\!\left(g_1^{-1}q_1,\phi_P(p_1g_1,p_2g_2)g_2^{-1}q_2\right)=\\
&=\phi_Q\!\left(g_1^{-1}q_1,g_1^{-1}\phi_P(p_1,p_2)g_2g_2^{-1}q_2\right)=\\
&=\phi_Q\!\left(g_1^{-1}q_1,g_1^{-1}\phi_P(p_1,p_2)q_2\right)=\\
&=\phi_Q\!\left(q_1,\phi_P(p_1,p_2)q_2\right)=\\
&=\phi_{Q\circ P}\!\left([p_1,q_1],[p_2,q_2]\right),
\end{align*}
where the last identity is a consequence of Lemma~\ref{lem-invdivmap}, and I have used of the $\calG$-equivariance of the division map $\phi_P$ of $P$.

It remains to show that the map (\ref{eq-genprodivmap}) satisfies Equation (\ref{eq-divmap}), which shows definitely that it is the division map for $Q\circ P$.
Namely, consider any two equivalence classes $[p_1,q_1]$ and $[p_2,q_2]$ in $P\times_\calG Q$ lying in the same fiber of $\overline{\pi}$, i.e.
\[
\pi_1(p_1)=\pi_1(p_2).
\] 
Then, one gets
\begin{align*}
[p_1,q_1]\phi_{Q\circ P}\!\left([p_1,q_1],[p_2,q_2]\right)&=[p_1,q_1]\phi_Q\!\left(q_1,\phi_P(p_1,p_2)q_2\right)=\\
&=\left[p_1,q_1\phi_Q\!\left(q_1,\phi_P(p_1,p_2)q_2\right)\right]=\\
&=\left[p_1,\phi_P(p_1,p_2)q_2\right]=\\
&=\left[p_1\phi_P(p_1,p_2),q_2\right]=\\
&=[p_2,q_2],
\end{align*}
where I have used Equation (\ref{eq-divmap}) for $\phi_P$ and $\phi_Q$, as well as of the $\calG$-invariance of any class in $P\times_{\calG}Q$.

Hence, the $4$-tuple $\left(P\times_{\calG}Q,\overline{\pi},\overline{\varepsilon},X_{\calF}\right)$, endowed with the left $\calF$-action (\ref{eq-genprodleft}) and right $\calH$-action (\ref{eq-genprodright}), is a generalized morphism in the sense of Definition~\ref{def-genmorph}, which I denote by $\calF\overset{Q\circ P}\to \calH$.

The shape of the division map (\ref{eq-genprodivmap}) is very important: in fact, the local generalized morphism associated to a generalized morphism is explicitly constructed via the division map of the generalized morphism, viewed as a principal $\calH$-bundle. 
This is the aim of what will follow: i.e.\, I want to relate the local generalized morphism $\Theta_{Q\circ P}$ associated to the composition $\calF\overset{Q\circ P}\to\calH$ to the respective local generalized morphisms $\Theta_P$ and $\Theta_Q$.
This will give the answer to the natural question: since local generalized morphisms (viewed as generalized morphisms) can be composed, what is the shape of their composition as a local generalized morphism?

\subsection{From generalized morphisms to local generalized morphisms: the composition law in terms of the division map}\label{ssec-compdivmap}
Following the construction in Subsection~\ref{ssec-locgenmor}, what I need in order to construct the local generalized morphism $\Theta_{Q\circ P}$ is an open cover of $X_{\calF}$, as well as local sections w.r.t.\ to the chosen cover of the projection $\overline{\pi}$. 
This is a subtle point of the construction: in fact, one cannot simply take any cover of $X_{\calF}$, because of the definition of the associated bundle $P\times_{\calG}Q$.
In fact, any representative of a general equivalence class $[p_1,q_1]$ in $P\times_{\calG}Q$ obeys the equation
\[
\varepsilon_1(p_1)=\pi_2(q_1),
\]
which puts some constraints on the possible sections of $\overline{\pi}$.
To deal with this compatibility relation, I work as follows:
\begin{itemize}
\item[i)] choose first any open cover $\mathfrak{U}_{\calF}$ of $X_{\calF}$ and local sections $\sigma_{\alpha}^1$ of $\pi_1$; the associated local momenta are, as usual, simply $\varepsilon_{\alpha}^1=\varepsilon_1\circ\sigma_{\alpha}^1:U_{\alpha}\to X_{\calG}$. 
\item[ii)] consider then an open cover $\mathfrak{V}_{\calG}$ of $X_{\calG}$ and associated local sections $\sigma_{i}^2$ of $\pi_2$; notice that, for the cover $\mathfrak{V}_{\calG}$, Latin indices are used, instead of Greek indices as for the cover $\mathfrak{U}_{\calF}$ of $X_{\calF}$.
Since $\mathfrak{V}_{\calG}$ is an open cover of $X_{\calG}$, it follows that any set $\varepsilon_{\alpha}^1(U_\alpha)$ (which is not necessarily open; in fact, this may happen when the local momenta are open maps, i.e.\ when the global momentum $\varepsilon_1$ is an open map) is covered by open sets $V_i$; the intersection $\varepsilon_{\alpha}^1(U_\alpha)\cap V_i$ is open in the relative topology of $\varepsilon_{\alpha}^1(U_\alpha)$, and will be denoted by $V_{\alpha_i,\varepsilon_1}$.
Since $\varepsilon_{\alpha}^1$ is smooth, it follows that the preimage $U_{\alpha_i}\colon=(\varepsilon_{\alpha}^1)^{-1}\!\left(V_{\alpha_i,\varepsilon_1}\right)$ is relatively open in $U_\alpha$ (hence, it is open in $X_{\calF}$).
Moreover, since
\[
\varepsilon_{\alpha}^1(U_\alpha)=\bigcup_{i\colon\varepsilon_{\alpha}^1(U_\alpha)\cap V_i\neq\emptyset}\varepsilon_{\alpha}^1(U_\alpha)\cap V_i,
\]
it follows
\[
U_\alpha=\bigcup_{i\colon\varepsilon_{\alpha}^1(U_\alpha)\cap V_i\neq\emptyset}U_{\alpha,i,\varepsilon_1}, 
\]
i.e.\ the open sets $U_{\alpha_i}$ cover $U_\alpha$, with the index $i$ such that $\varepsilon_{\alpha}^1(U_\alpha)\cap V_i$ is nontrivial.
\item[iii)] I define a refinement $\overline{\mathfrak{U}}_{\calF}$ of the cover $\mathfrak{U}_{\calF}$ of $X_{\calF}$ as follows: any open set $U_\alpha$ of the cover $\mathfrak{U}_{\calF}$ is covered by the open sets $U_{\alpha_i}$, which I have constructed in Part $ii)$.
The local sections $\sigma_{\alpha_i}^1$ are simply set to be the restrictions of the section $\sigma_\alpha^1$ to $U_{\alpha_i}$, which are clearly again sections of $\pi_1$; but now, the local momenta $\varepsilon_{\alpha_i}^1$ map $U_{\alpha_i}$ to $V_i$, hence it makes sense to consider the composition $\sigma_{i}^2\circ \varepsilon_{\alpha_i}^1\colon U_{\alpha_i}\to Q$.
\end{itemize}
With the previous notations in mind, I consider the following local sections of $\overline{\pi}$ w.r.t.\ the open cover $\overline{\mathfrak{U}}_{\calF}$ of $X_{\calF}$:
\begin{equation}\label{eq-genprodsect}
\overline{\sigma}_{\alpha_i}\colon\begin{cases}
U_{\alpha_i}&\to P\times_{\calG}Q,\\
x&\mapsto \left[\sigma_{\alpha_i}^1(x),\sigma_{i}^2(\varepsilon_{\alpha_i}^1(x))\right].
\end{cases}
\end{equation}
It is easy to verify that the maps $\overline{\sigma}_{\alpha_i}$ are well-defined sections of $\overline{\pi}$; moreover, an immediate computation shows, recalling Equations (\ref{eq-genprodmaps}), that the associated local momenta are simply
\[
\overline{\varepsilon}_{\alpha_i}=\varepsilon_{i}^2\circ \varepsilon_{\alpha_i}^1:U_{\alpha_i}\to X_{\calH}.
\]
Now that a convenient set of local sections for the composition $Q\circ P$ has been found, I need for later purposes to compute the associated transition maps.
Recalling the constructions in Subsections~\ref{ssec-localsec}, let me perform the following computation:
\begin{align*}
\Phi_{\alpha_i\beta_j}^{Q\circ P}(x)&\colon=\phi_{Q\circ P}\!\left(\overline{\sigma}_{\alpha_i}(x),\overline{\sigma}_{\beta_j}(x)\right)=\\
&=\phi_{Q\circ P}\!\left(\left[\sigma_{\alpha_i}^1(x),\sigma_{i}^2(\varepsilon_{\alpha_i}^1(x))\right],\left[\sigma_{\beta_j}^1(x),\sigma_{j}^2(\varepsilon_{\beta_j}^1(x))\right]\right)=\\
&=\phi_Q\!\left(\sigma_{i}^2(\varepsilon_{\alpha_i}^1(x)),\phi_P\!\left(\sigma_{\alpha_i}^1(x),\sigma_{\beta_j}^1(x)\right)\sigma_{j}^2(\varepsilon_{\beta_j}^1(x))\right)=\\
&=\phi_Q\!\left(\sigma_{i}^2\!\left(t_{\calG}\!\left(\Phi_{\alpha_i\beta_j}^P(x)\right)\right),\Phi_{\alpha_i\beta_j}^P(x)\sigma_{j}^2\!\left(s_{\calG}\!\left(\Phi_{\alpha_i\beta_j}^P(x)\right)\right)\right)=\\
&=\Theta_{ij}^Q\!\left(\Phi_{\alpha_i\beta_j}^P(x)\right),\quad x\in U_{\alpha_i}\cap U_{\beta_j},
\end{align*}
where I have made use of the first set of identities in Definition~\ref{def-trivdata}. 
I have also used the following notation:
\[
\Phi_{\alpha_i\beta_j}^P\colon=\Phi_{\alpha\beta}^P\big\vert_{U_{\alpha_i}\cap U_{\beta_j}},
\]
where $\Phi_{\alpha\beta}^P$ denotes the transition map of $P$ associated to the local sections $\sigma_\alpha$ and $\sigma_\beta$.
Hence, one can summarize the result of the previous computation as follows: the generalized morphism $Q\circ P$ has transition maps of the following form
\begin{equation}\label{eq-comptrans}
\Phi_{\alpha_i\beta_j}^{Q\circ P}(x)=\Theta_{ij}^Q\!\left(\Phi_{\alpha_i\beta_j}^P(x)\right),
\end{equation}
where $\Theta^Q$ is the local generalized morphism subordinate to the local trivializing data for $Q$ with cover $\mathfrak{V}_{\calG}$.
With the same choice of local sections of $\overline{\pi}$, one gets the following expression for the associated generalized morphism $\Theta^{Q\circ P}$, where I recall Identity (\ref{eq-genmor2}):
\begin{align*}
\Theta_{\beta_j,\alpha_i}^{Q\circ P}(f)&\colon=\phi_{Q\circ P}\!\left(\overline{\sigma}_{\beta_j}(t_{\calF}(f)),f\overline{\sigma}_{\alpha_i}(s_{\calF}(f))\right)=\\
&=\phi_{Q\circ P}\!\left(\left[\sigma_{\beta_j}^1(t_{\calF}(f)),\sigma_{j}^2(\varepsilon_{\beta_j}^1(t_{\calF}(f)))\right],\right.\\
&\phantom{=}\left.f\left[\sigma_{\alpha_i}^1(s_{\calF}(f)),\sigma_{i}^2(\varepsilon_{\alpha_i}^1(s_{\calF}(f)))\right]\right)=\\
&=\phi_{Q\circ P}\!\left(\left[\sigma_{\beta_j}^1(t_{\calF}(f)),\sigma_{j}^2(\varepsilon_{\beta_j}^1(t_{\calF}(f)))\right],\right.\\
&\phantom{=}\left.\left[f \sigma_{\alpha_i}^1(s_{\calF}(f)),\sigma_{i}^2(\varepsilon_{\alpha_i}^1(s_{\calF}(f)))\right]\right)=\\
&=\phi_Q\!\left(\sigma_{j,2}(\varepsilon_{\beta_j}^1(t_{\calF}(f))),\phi_P\!\left(\sigma_{\beta_j}^1(t_{\calF}(f)),\right.\right.\\
&\phantom{=}\left.\left.f \sigma_{\alpha_i}^1s_{\calF}(f))\right)\sigma_{i}^2(\varepsilon_{\alpha_i}^2(s_{\calF}(f)))\right)=\\
&=\phi_Q\!\left(\sigma_{j}^2(\varepsilon_{\beta_j}^2(t_{\calF}(f))),\phi_P\!\left(\sigma_{\beta_j}^1(t_{\calF}(f)),\right.\right.\\
&\phantom{=}\left.\left.f \sigma_{\alpha_i}^1(s_{\calF}(f))\right)\sigma_{i}^2(\varepsilon_{\alpha_i}^1(s_{\calF}(f)))\right)=\\
&=\phi_Q\!\left(\sigma_{j}^2(\varepsilon_{\beta_j}^1(t_{\calF}(f))),\Theta_{\beta_j,\alpha_i}^P(f)\sigma_{i}^2(\varepsilon_{\alpha_i}^1(s_{\calF}(f)))\right)=\\
&=\phi_Q\!\left(\sigma_{j}^2\!\left(t_{\calG}\!\left(\Theta_{\beta_j,\alpha_i}^P(f)\right)\right),\Theta_{\beta_j,\alpha_i}^P(f)\sigma_{i}^2\!\left(s_{\calG}\!\left(\Theta_{\beta_j,\alpha_i}^P(f)\right)\right)\right)=\\
&=\Theta_{j,i}^Q\!\left(\Theta_{\beta_j,\alpha_i}^P(f)\right),\quad t_{\calF}(f)\in U_{\beta_j},\quad s_{\calF}(f)\in U_{\alpha_i}.
\end{align*}
In the above computations, I have used the commutativity of the diagrams (\ref{eq-diagmormomen}) of Definition~\ref{def-locgenmor}.
Finally, with all the above notations, I have the following formula:
\begin{equation}\label{eq-comploc}
\Theta_{\beta_j,\alpha_i}^{Q\circ P}(f)=\Theta_{j,i}^Q\!\left(\Theta_{\beta_j,\alpha_i}^P(f)\right),\quad t_{\calF}(f)\in U_{\beta_j},\quad s_{\calF}(f)\in U_{\alpha_i}.
\end{equation}
Hence, to the composition of two generalized morphisms $\calF\overset{P}\to\calG$ and $\calG\overset{Q}\to\calH$, whose associated local generalized morphisms are $\Theta^P$ and $\Theta^Q$ respectively, can be associated a local generalized morphism $\Theta^{Q\circ P}$.
Taking some care towards the local nature of $\Theta^P$ and $\Theta^Q$, $\Theta^{Q\circ P}$ can be viewed as the composition of the local generalized morphisms $\Theta^Q$ and $\Theta^P$.

\subsection{The composition law for local generalized morphisms: the refinement trick and the general arguments}\label{ssec-locgencomplaw}
Motivated by the result (\ref{eq-comploc}) at the end of the preceding subsection, I now generalize the previous construction to the case of general local generalized morphisms in the sense of Definition~\ref{def-locgenmor}.
Hence, let me consider three Lie groupoids $\calF$, $\calG$ and $\calH$, and two local generalized morphisms $\calF\overset{\Theta}\to \calG$ and $\calG\overset{\Eta}\to\calH$.
Let $\Theta$ be subordinate to the local trivializing data $\left(\mathfrak{U}_\calF,\varepsilon_{\alpha}^1,\Phi_{\alpha\beta}\right)$ over $X_{\calF}$ with values in $\calG$, resp.\ $\Eta$ be subordinate to the local trivializing data $\left(\mathfrak{V}_{\calG},\varepsilon_{i}^2,\Psi_{ij}\right)$ over $X_{\calG}$ with values in $\calH$.

First of all, before entering into the details of the construction of the composition of $\Eta$ with $\Theta$, I need a discussion on the {\em refinement trick} sketched in the previous subsection.
I proceed along the following steps:
{\bf \begin{itemize}
  \item[i)] I consider all subsets of $X_{\calG}$ of the shape 
\[
\varepsilon_{\alpha}^1(U_\alpha)\subset X_{\calG},\quad U_\alpha\in \mathfrak{U}_\calF.
\]
  \item[ii)] Since $\mathfrak{V}_\calG$ is a cover of $X_{\calG}$, it is possible to write 
\[
\varepsilon_{\alpha}^1(U_\alpha)=\bigcup_{i\colon V_i\cap\varepsilon_{\alpha}^1(U_\alpha)\neq\emptyset}\left(V_i\cap \varepsilon_{\alpha}^1(U_\alpha)\right).
\] 
Since the subsets $V_i$ are all open, the nontrivial subsets $V_i\cap \varepsilon_{\alpha}^1(U_\alpha)$ are relatively open in $\varepsilon_{\alpha}^1(U_\alpha)$ (even the trivial ones, but they are of no interest).
  \item[iii)] Since the local momenta $\varepsilon_{\alpha}^1$ are smooth, the preimages
\[
U_{\alpha_i}\colon=(\varepsilon_{\alpha}^1)^{-1}\!\left(V_i\cap \varepsilon_{\alpha}^1(U_\alpha)\right)
\] 
are relatively open in $U_\alpha$, and hence they are open in $X_{\calF}$; moreover, 
\[
\varepsilon_{\alpha}^1(U_\alpha)=\bigcup_{i\colon V_i\cap\varepsilon_{\alpha}^1(U_\alpha)\neq\emptyset}\left(V_i\cap \varepsilon_{\alpha}^1(U_\alpha)\right)\Rightarrow U_\alpha=\bigcup_{i\colon V_i\cap\varepsilon_{\alpha}^1(U_\alpha)\neq\emptyset}U_{\alpha_i}.
\]
Since the $U_{\alpha_i}$ cover $U_\alpha$, it follows that
\[
\overline{\mathfrak{U}}_{\calF}\colon=\left\{U_{\alpha_i}\subset X_{\calF}\right\}_{\alpha,i}
\]
is also a cover and, by its very definition, it is moreover a refinement of the cover $\mathfrak{U}$, where the refinement map is simply $\alpha_i\mapsto i$.
\item[iv)] The refinement $\overline{\mathfrak{U}}_{\calF}$ has the following important property: considering the new local momenta $\varepsilon_{\alpha_i}^1$ as the restrictions of the local momenta $\varepsilon_{\alpha}^1$ to $U_{\alpha_i}\subset U_{\alpha,1}$, it follows by the very construction:
\begin{equation}\label{eq-compat}
\varepsilon_{\alpha_i}^1\colon U_{\alpha_i}\subset U_{\alpha,1}\to V_i\cap\varepsilon_{\alpha,1}(U_\alpha)\subset V_i,
\end{equation}
which is the main property that I need to construct the composition of $\Eta$ with $\Theta$.
\end{itemize}}
Given now the generalized morphism $\Theta$ from $\calF$ to $\calG$ subordinate to the local trivializing data $\left(\mathfrak{U}_{\calF},\varepsilon_{\alpha}^1,\Phi_{\alpha\beta}\right)$, define its {\em refinement $\overline{\Theta}$ w.r.t.\ the cover $\overline{\mathfrak{U}}_{\calF}$} as follows:
\begin{itemize}
\item[i)] the corresponding local trivializing data are now $\left(\overline{\mathfrak{U}}_{\calF},\varepsilon_{\alpha_i}^1,\Phi_{\alpha_i\beta_j}\right)$, where
\begin{align}
\varepsilon_{\alpha_i}^1&\colon =\varepsilon_{\alpha}^1\big\vert_{U_{\alpha_i}},\label{eq-refinmoment}\\
\Phi_{\alpha_i\beta_j}&\colon=\Phi_{\alpha\beta}\big\vert_{U_{\alpha_i}\cap U_{\beta_j}},\quad U_{\alpha_i}\cap U_{\beta_j}\neq\emptyset.\label{eq-refintrans}
\end{align}
By standard arguments in ordinary nonabelian {\v C}ech cohomology, it is not difficult to verify that the $3$-tuple $\left(\overline{\mathfrak{U}}_\calF,\varepsilon_{\alpha_i}^1,\Phi_{\alpha_i\beta_j}\right)$ defines local trivializing data in the sense of Definition~\ref{def-trivdata}.
\item[ii)] The local morphism $\overline{\Theta}$ is defined as follows:
\begin{equation}\label{eq-refinlocmor}
\overline{\Theta}_{\beta_j\alpha_i}(f)\colon=\Theta_{\beta\alpha}(f),\quad f\in\calF_{\alpha_i,\beta_j}.
\end{equation}
\end{itemize}
Then, one gets the following 
\begin{Prop}\label{prop-refingenmor}
The $4$-tuple $\left(\overline{\mathfrak{U}}_{\calF},\varepsilon_{\alpha_i}^1,\Phi_{\alpha_i\beta_j},\overline{\Theta}_{\beta_j\alpha_i}\right)$ is a local generalized morphism in the sense of Definition~\ref{def-locgenmor}.
\end{Prop}
\begin{proof}
It is already known that the $3$-tuple $\left(\overline{\mathfrak{U}}_{\calF},\varepsilon_{\alpha_i}^1,\Phi_{\alpha_i\beta_j}\right)$ defines local trivializing data; hence, it remains to check that the three diagrams in (\ref{eq-diagmormomen}) commute, that (\ref{eq-genhomom}) and (\ref{eq-relgentrans}) hold.

Let me check the commutativity of the first diagram in (\ref{eq-diagmormomen}); the proof of the commutativity if the remaining two is immediate.
Namely, 
\begin{align*}
s_{\calG}\left(\overline{\Theta}_{\beta_j\alpha_i}(f)\right)&=s_{\calG}\left(\Theta_{\beta\alpha}(f)\right)=\\
&=\varepsilon_{\alpha}^1\!\left(s_{\calF}(f)\right)=\\
&=\varepsilon_{\alpha_i}^1\!\left(s_{\calF}(f)\right),\quad f\in\calF_{\alpha_i,\beta_j}\Rightarrow s_{\calF}(f)\in U_{\alpha_i},
\end{align*}
by (\ref{eq-refinmoment}).

The proof of Identity (\ref{eq-genhomom}) is straightforward by Equation (\ref{eq-refinlocmor}), since $\Theta$ is a local generalized morphism.

The proof of Identity (\ref{eq-relgentrans}) is also immediate, by (\ref{eq-refintrans}) and (\ref{eq-refinlocmor}) and since $\Theta$ is subordinate to $\left(\mathfrak{U},\varepsilon_{\alpha}^1,\Phi_{\alpha\beta}\right)$.
\end{proof}

I have now all the elements needed to define the composition of two local generalized morphisms, provided they are composable.
\begin{Def}\label{def-complocmor}
Given two local generalized morphisms $\calF\overset{\Theta}\to\calG$ and $\calG\overset{\Eta}\to\calH$ in the sense of Definition~\ref{def-locgenmor}, subordinate respectively to the local trivializing data $\left(\mathfrak{U}_{\calF},\varepsilon_{\alpha}^1,\Phi_{\alpha\beta}\right)$ and $\left(\mathfrak{V}_{\calG},\varepsilon_{i}^2,\Psi_{ij}\right)$, the composition $\Eta\circ\Theta$ is defined as follows: consider the open cover $\overline{\mathfrak{U}}_{\calF}$ of $X_{\calF}$, constructed at the beginning of the subsection, the refinement of $\mathfrak{U}_{\calF}$ ``matched'' with the local momenta $\varepsilon_{\alpha}^1$ and $\varepsilon_{i}^2$; consequently, consider the corresponding refinement $\overline{\Theta}$, which, by the previous Lemma, defines also a local generalized morphism.   

Then, consider also the $3$-tuple $\left(\overline{\mathfrak{U}}_{\calF},\overline{\varepsilon}_{\alpha_i},\overline{\Phi}_{\alpha_i\beta_j}\right)$, where
\begin{align}
\label{eq-loccompmom}\overline{\varepsilon}_{\alpha_i}&\colon=\varepsilon_{i}^2\circ\varepsilon_{\alpha_i}^1,\\
\label{eq-loccomptrans}\overline{\Phi}_{\alpha_i\beta_j}&\colon=\Eta_{ij}\circ \Phi_{\alpha_i\beta_j}.
\end{align} 
Moreover, define
\begin{equation}\label{eq-loccomp}
\left(\Eta\circ\Theta\right)_{\beta_j\alpha_i}\!(f)\colon=\Eta_{ji}\!\left(\overline{\Theta}_{\beta_j\alpha_i}(f)\right),\quad f\in\calF_{\alpha_i,\beta_j}.
\end{equation}
The $3$-tuple $\left(\overline{\mathfrak{U}}_{\calF},\overline{\varepsilon}_{\alpha_i},\overline{\Phi}_{\alpha_i\beta_j}\right)$, together with the maps $\Eta\circ\Theta$, defined in (\ref{eq-loccomp}), is said to be the composition of the local generalized morphisms $\Eta$ and $\Theta$.
\end{Def}
First of all, notice that the maps defined in (\ref{eq-loccompmom}), (\ref{eq-loccomptrans}) and (\ref{eq-loccomp}) all make sense, because of the construction of the refinement $\overline{\mathfrak{U}}_{\calF}$ and because of (\ref{eq-compat}).

Moreover, the following Theorem holds
\begin{Thm}
Given two local generalized morphisms $\calF\overset{\Theta}\to\calG$ and $\calG\overset{\Eta}\to\calH$ in the sense of Definition~\ref{def-locgenmor}, subordinate respectively to the local trivializing data $\left(\mathfrak{U}_{\calF},\varepsilon_{\alpha}^1,\Phi_{\alpha\beta}\right)$ and $\left(\mathfrak{V}_{\calG},\varepsilon_{i}^2,\Psi_{ij}\right)$, the map $\Eta\circ\Theta$, defined in (\ref{eq-loccomp}) is a local generalized morphism, subordinate to the local trivializing data $\left(\overline{\mathfrak{U}}_{\calF},\overline{\varepsilon}_{\alpha_i},\overline{\Phi}_{\alpha_i\beta_j}\right)$.
\end{Thm}
\begin{proof}
Proposition~\ref{prop-refingenmor} implies immediately that $\left(\overline{\mathfrak{U}}_{\calF},\overline{\varepsilon}_{\alpha_i},\overline{\Phi}_{\alpha_i\beta_j}\right)$ defines local trivializing data on $X_{\calF}$ with values in $\calH$, since $\Eta$ is a local generalized morphism subordinate to the local trivializing data $\left(\mathfrak{V}_{\calG},\varepsilon_{i}^2,\Psi_{ij}\right)$.

It remains therefore to check that the three diagrams in (\ref{eq-diagmormomen}) commute, and that both identities (\ref{eq-genhomom}) and (\ref{eq-relgentrans}) hold.
Let me check the commutativity of the first diagram in (\ref{eq-diagmormomen}):
\begin{align*}
s_{\calH}\!\left(\left(\Eta\circ\Theta\right)_{\beta_j\alpha_i}\!(f)\right)&=s_{\calH}\!\left(\Eta_{ji}\!\left(\overline{\Theta}_{\beta_j\alpha_i}(f)\right)\right)=\\
&=\varepsilon_{i}^2\!\left(s_{\calG}\!\left(\overline{\Theta}_{\beta_j\alpha_i}(f)\right)\right)=\\
&=\varepsilon_{i}^2\!\left(\varepsilon_{\alpha_i}^1\!(s_{\calF}(f))\right)=\\
&=\overline{\varepsilon}_{\alpha_i}\!(s_{\calF}(f)),\quad f\in\calF_{\alpha_i,\beta_j},
\end{align*}
and similarly for the remaining two diagrams.

To check Identity (\ref{eq-genhomom}), let me compute
\begin{align*}
\left(\Eta\circ\Theta\right)_{\gamma_k\alpha_i}\!(f_1f_2)&=\Eta_{ki}\!\left(\overline{\Theta}_{\gamma_k\alpha_i}(f_1f_2)\right)=\\
&=\Eta_{ki}\!\left(\overline{\Theta}_{\gamma_k\beta_j}(f_1)\overline{\Theta}_{\beta_j\alpha_i}(f_2)\right)=\\
&=\Eta_{kj}\!\left(\overline{\Theta}_{\gamma_k\beta_j}(f_1)\right)\Eta_{ji}\!\left(\overline{\Theta}_{\beta_j\alpha_i}(f_2)\right)=\\
&=\left(\Eta\circ\Theta\right)_{\gamma_k\beta_j}\!(f_1)\left(\Eta\circ\Theta\right)_{\beta_j\alpha_i}\!(f_2),\quad f_1\in \calF_{\beta_j,\gamma_k},\quad f_2\in\calF_{\alpha_i,\beta_j},
\end{align*} 
since both $\Eta$ and $\overline{\Theta}$ are both local generalized morphisms.

Identity (\ref{eq-relgentrans}0, expressing the compatibility condition between the local generalized morphism and and the transition maps of the corresponding local trivializing data, is a consequence of Proposition~\ref{prop-refingenmor} and of the fact that $\Eta$ is a local generalized morphism:
\begin{align*}
\left(\Eta\circ\Theta\right)_{\beta_j\alpha_i}\!(f)&=\Eta_{ji}\!\left(\overline{\Theta}_{\beta_j\alpha_i}(f)\right)=\\
&=\Eta_{ji}\!\left(\Phi_{\beta_j\delta_l}\!(t_{\calF}(f))\overline{\Theta}_{\delta_l\gamma_k}(f)\Phi_{\gamma_k\alpha_i}\!(s_{\calF}(f))\right)=\\
&=\Eta_{jl}\!\left(\Phi_{\beta_j\delta_l}\!(t_{\calF}(f))\right)\Eta_{lk}\!\left(\overline{\Theta}_{\delta_l\gamma_k}(f)\right)\Eta_{ki}\!\left(\Phi_{\gamma_k\alpha_i}\!(s_{\calF}(f))\right)=\\
&=\overline{\Phi}_{\beta_j\delta_l}\!(t_{\calF}(f))\left(\Eta\circ\Theta\right)_{\delta_l\gamma_k}\overline{\Phi}_{\gamma_k\alpha_i}\!(t_{\calF}(f)),\quad f\in\calF_{\alpha_i,\beta_j}\cap\calF_{\gamma_k,\delta_l}.
\end{align*}
\end{proof}

Thus, it is possible to define in a meaningful way the composition of two composable local generalized morphisms $\calF\overset{\Theta}\to\calH$ and $\calG\overset{\Eta}\to\calH$; moreover, the computations performed in Subsection~\ref{ssec-gencomplaw}, combined with the preceding computations and with Lemma~\ref{lem-locmortogenmor} of subsection~\ref{ssec-locgentogenmor}, shows that the local generalized morphism associated to the composition of two generalized morphism is exactly the composition of the associated local generalized morphisms and viceversa.
Moreover, the computations done in the previous subsection show immediately that the composition of local generalized morphisms corresponds to the composition of the associated generalized morphisms, and viceversa.

\begin{Rem}
The final result $\overline{\mathfrak{U}}_{\calF}$ of the ``refinement trick'' I sketched above I call a {\em birefinement} of $\mathfrak{U}_\calF$ w.r.t.\ $\mathfrak{V}_\calG$.
Why such a choice of words? 
Because of its very nature, i.e.\ there are two maps defining the refinement: the first one, $\alpha_i\mapsto\alpha$, is a refinement map in the ordinary sense, whereas the second one, $\alpha_i\mapsto i$, takes the image of the ``birefinement'' w.r.t.\ the local momenta $\varepsilon_{\alpha_i,1}$ to the open cover $\mathfrak{V}_{\calG}$.
\end{Rem}

\section{Morita equivalence: a local characterization}\label{sec-morita}
In this section, I discuss the local nature of Morita equivalences, a special kind of generalized morphisms between groupoids.
I will therefore first introduce the notion of Morita equivalence, following closely~\cite{Moer2} and~\cite{Mrcun}.
I will study carefully how Morita equivalences give rise, in a canonical way, to two distinct division maps, whose properties I will display in detail.
Afterwards, I will define the ``inverse'' of a Morita equivalence, in which sense I will explain later.
The existence of an inverse of a Morita equivalence, combined with the results of Subsection~\ref{ssec-gencomplaw}, gives rise to the {\em factorization property of division maps}. 

Coming subsequently to the local nature of Morita equivalences, I will combine the factorization property of division maps with local data, in the same spirit of Subsection~\ref{ssec-localsec} with the results of Section~\ref{sec-hilskan}: this will give rise of a local version of Morita equivalences.

\subsection{Left and right division maps for a Morita equivalence}\label{ssec-divmorita}
Before entering into the details, I need the following
\begin{Def}\label{def-moritaeq}
Given two Lie groupoids $\calG$ and $\calH$, they are said to be {\em Morita equivalent} if there is a generalized morphism $\calG\overset{P}\to\calH$ in the sense of Definition~\ref{def-genmorph}, such that the $4$-tuple $\left(P,\varepsilon,\pi,X_{\calH}\right)$ is a left principal $\calG$-bundle with canonical projection $\varepsilon$ and momentum $\pi$.
\end{Def}

Given a Lie groupoid $\calG$, the definition and main properties of left principal $\calG$-bundle are completely analogous to those of right $\calG$-bundle: namely, a left principal $\calG$-bundle over $X$ is a $4$-tuple $\left(P,\pi,\varepsilon,X\right)$, such that $i)$ $P$ and $X$ are smooth manifolds, $ii)$ $\pi$ is a surjective submersion from $P$ to $X$ and $\varepsilon$ is a smooth map from $P$ to $X_{\calG}$ (called the {\em left momentum of the bundle}).
Moreover, the following requirements must be satisfied:
\begin{itemize}
\item[a)] the triple $\left(P,\calG,\varepsilon\right)$ defines a smooth left $\calG$-action on $P$; moreover, $\pi$ is $\calG$-invariant. 
\item[b)] The map
\[
\calG\times_\varepsilon P\ni (g,p)\mapsto (pg,p)\in P\times_XP
\] 
is a diffeomorphism.
\end{itemize}
In complete analogy to what was done in Subsection~\ref{ssec-defdivmap}, one can define the left fibred product bundle of a left principal $\calG$-bundle $\left(P,\pi,\varepsilon,X\right)$ with itself, denoted by $P\odot_LP$, by setting
\[
P\odot_LP\colon=\left\{(p,q)\in P\times P\colon \pi(p)=\pi(q)\right\}.
\]
It can then be proved that the $4$-tuple $\left(P\odot_LP,\overline{\pi},\varepsilon\times\varepsilon,X_{\calH}\right)$, where
\[
\overline{\pi}(p,q)\colon=\pi(p)=\pi(q),\quad \varepsilon\times\varepsilon(p_1,p_2)\colon=\left(\varepsilon(p),\varepsilon(q)\right),
\]
is a left principal $\calG^2$-bundle over $X$.

Recall now the construction of the {\em generalized conjugation of a Lie groupoid $\calG$} from Subsection~\ref{ssec-defdivmap}: I consider, in this case, the triple $\left(\calG^2,J_{\conj},\Psi_{\conj}\right)$, which defines, by Proposition~\ref{prop-genconj}, a left $\calG^2$-action on $\calG$ itself.

\begin{Def}\label{def-leftdivmap}
Given a left principal $\calG$-bundle $\left(P,\pi,\varepsilon,X\right)$, the {\em left division map of $P$} is defined by the following equation:
\begin{equation}\label{eq-leftdivmap}
p=\phi_P^L(p,q)q,\quad \pi(p)=\pi(q).
\end{equation}
\end{Def}
It follows immediately from Definition~\ref{def-leftdivmap} that the left division map $\phi_P^L$ is a smooth map from the fibred product $P\odot_LP$ to $\calG$; in fact, it is the first component of the smooth inverse of the map
\[
\calG\times P\ni (g,p)\mapsto (gp,p)\in P\odot_LP.
\]

Let me state without proof the following 
\begin{Prop}\label{prop-propleftdiv}
The left division map $\phi_P^L$ from $P\odot_LP$ to $\calG$ has the following properties:
\begin{itemize}
\item[i)] for any point $(p,q)$ of $P\odot_LP$, one has
\[
\phi_P^L(p,q)\in \calG_{\varepsilon(q),\varepsilon(p)}.
\] 
\item[ii)] On the diagonal submanifold of the total space of $P\odot_LP$, one has
\[
\phi_P^L(p,p)=\iota_{\calG}(\varepsilon(p)),\quad \forall p\in P.
\] 
\item[iii)] for any pair $(p,q)\in P\odot_LP$, the following equation holds
\[
\phi_P^L(p,q)=\phi_P^L(q,p)^{-1};
\] 
notice that the previous equation makes sense, since $(p,q)\in P\odot_LP$ implies that $(q,p)\in P\odot_LP$ also.
\item[iv)] The triple $\left(\phi_P^L,\id_{\calG^2},\id_{X_{\calG}^2}\right)$ is an equivariant map from the right $\calG^2$-space $P\odot_LP$ to the right $\calG^2$-space $\calG$ endowed with the left generalized conjugation defined by $\left(J_{\conj},\Psi_{\conj}\right)$.  
\end{itemize}
\end{Prop}
Consider now a Morita equivalence between Lie groupoids $\calG$ and $\calH$, say $\calG\overset{P}\to\calH$, hence, there are two principal bundles: the left principal $\calG$-bundle $\left(P,\varepsilon,\pi,X_{\calH}\right)$ and the right principal $\calH$-bundle $\left(P,\pi,\varepsilon,X_{\calG}\right)$.
Thus, to a Morita equivalence $\calG\overset{P}\to\calH$ are associated two canonical maps, which I denote by $\phi_P^L$ and $\phi_P^R$ respectively, namely the left division map of the bundle $\left(P,\varepsilon,\pi,X_{\calH}\right)$ and the right division map of the bundle $\left(P,\pi,\varepsilon,X_{\calG}\right)$:
\begin{align*}
p&=\phi_P^L(p,q)q,\quad \varepsilon(p)=\varepsilon(q),\\
q&=p\phi_P^R(p,q),\quad \pi(p)=\pi(q).
\end{align*}
For later purposes, I need the following
\begin{Lem}\label{lem-leftrightinv}
Given a Morita equivalence $\calG\overset{P}\to\calH$, the left, resp.\ right, division map $\phi_P^L$, resp.\ $\phi_P^R$, satisfies the following $\calH$-, resp.\ $\calG$-invariance:
\begin{align}
\label{eq-invleftdiv}\phi_P^L\!\left(p_1h,p_2h\right)&=\phi_P^L(p_1,p_2),\quad \varepsilon(p_1)=\varepsilon(p_2),\quad t_{\calH}(h)=\varepsilon(p_1)=\varepsilon(p_2),\\ 
\label{eq-invrightdiv}\phi_P^R\!\left(gp_1,gp_2\right)&=\phi_P^R(p_1,p_2),\quad \pi(p_1)=\pi(p_2),\quad s_{\calG}(g)=\pi(p_1)=\pi(p_2).
\end{align}
\end{Lem}
\begin{proof}
The $\calG$-invariance of the right division map of $P$ was already proven in Lemma~\ref{lem-invdivmap} of Subsection~\ref{ssec-gencomplaw}; Identity (\ref{eq-invleftdiv}) can be proved in the same way. 
\end{proof}

Now, since, in particular, a Morita equivalence is a generalized morphism in the sense of Definition~\ref{def-genmorph}, there is a local generalized morphism $\Theta^P$ from $\calG$ to $\calH$, constructed as in Subsection~\ref{ssec-locgenmor}, subordinate to the local trivializing data $\left(\mathfrak{U}_{\calG},\varepsilon_\alpha,\Phi_{\alpha\beta}^R\right)$ w.r.t.\ a choice of an open cover $\mathfrak{U}_{\calG}$ of $X_{\calG}$ and corresponding local sections $\sigma_\alpha$ of $\pi$:
\begin{align*}
\varepsilon_\alpha&\colon=\varepsilon\circ \sigma_\alpha,\\
\Phi_{\alpha\beta}(x)&\colon=\phi_P^R\!\left(\sigma_\alpha(x),\sigma_\beta(x)\right),\quad x\in U_{\alpha\beta},\\
\Theta_{\beta\alpha}^P(g)&\colon=\phi_P^R\!\left(\sigma_\beta(t_{\calG}(g)),g\sigma_\alpha(s_{\calG}(g))\right),\quad g\in \calG_{\alpha,\beta}.
\end{align*} 
One could wonder if there is a compatibility between the left and right division map associated to the Morita equivalence $\calG\overset{P}\to\calH$; this is the content of the following 
\begin{Prop}\label{prop-twistedinj}
The local generalized morphism $\Theta^P$ associated to the Morita equivalence $\calG\overset{P}\to\calH$ enjoys the following ``twisted injectivity'':
\begin{equation}\label{eq-moritainj}
\begin{aligned}
\Theta_{\beta\alpha}^P(g_1)&=\Theta_{\beta\alpha}^P(g_2),\quad g_1,g_2\in\calG_{\alpha,\beta}\Longleftrightarrow\\
\Longleftrightarrow \phi_P^L\!\left(\sigma_\beta(t_{\calG}(g_2)),\sigma_\beta(t_{\calG}(g_1))\right)g_1&=g_2\phi_P^L\!\left(\sigma_\alpha(s_{\calG}(g_2)),\sigma_\alpha(s_{\calG}(g_1))\right).
\end{aligned}
\end{equation}
\end{Prop}
\begin{proof}
Consider any two elements $g_1$ and $g_2$ in $\calG_{\alpha,\beta}$, such that
\[
\Theta_{\beta\alpha}^P(g_1)=\Theta_{\beta\alpha}^P(g_2).
\]
By applying on both sides of the previous equation the target, resp.\ the source map, of the groupoid $\calH$, one gets, by the commutativity of the diagrams (\ref{eq-diagmormomen}),
\begin{align*}
t_{\calH}\!\left(\Theta_{\beta\alpha}^P(g_1)\right)&=\varepsilon_\beta\!\left(t_{\calG}(g_1)\right)\overset{!}=\\
&\overset{!}=t_{\calH}\!\left(\Theta_{\beta\alpha}^P(g_2)\right)=\\
&=\varepsilon_\beta\!\left(t_{\calG}(g_2)\right),
\end{align*}
and similarly
\[
\varepsilon_\alpha\!\left(s_{\calG}(g_1)\right)=\varepsilon_\alpha\!\left(s_{\calG}(g_1)\right).
\]
Since $\varepsilon_\alpha=\varepsilon\circ\sigma_\alpha$, it follows immediately:
\begin{align*}
\sigma_\beta(t_{\calG}(g_2))&=\phi_P^L\!\left(\sigma_\beta(t_{\calG}(g_2)),\sigma_\beta(t_{\calG}(g_1))\right)\sigma_\beta(t_{\calG}(g_1))\quad\text{and}\\
\sigma_\beta(s_{\calG}(g_2))&=\phi_P^L\!\left(\sigma_\alpha(s_{\calG}(g_2)),\sigma_\alpha(s_{\calG}(g_1))\right)\sigma_\alpha(s_{\calG}(g_1)),
\end{align*}
by Equation (\ref{eq-leftdivmap}).
Set now
\begin{align*}
\phi_{\beta,t}^L(g_2,g_1)&\colon=\phi_P^L\!\left(\sigma_\beta(t_{\calG}(g_2)),\sigma_\beta(t_{\calG}(g_1))\right),\\
\phi_{\alpha,s}^L(g_2,g_1)&\colon=\phi_P^L\!\left(\sigma_\alpha(s_{\calG}(g_2)),\sigma_\alpha(s_{\calG}(g_1))\right).
\end{align*}
The elements $\phi_{\beta,t}^L(g_2,g_1)$ and $\phi_{\alpha,s}^L(g_2,g_1)$ of $\calG$ enjoy the following properties, which follow immediately from $i)$ of Proposition~\ref{prop-propleftdiv} and since $\sigma_\alpha$ and $\sigma_\beta$ are local sections of $\pi$:
\begin{align*}
t_{\calG}\!\left(\phi_{\beta,t}^L(g_2,g_1)\right)&=t_{\calG}(g_2),\quad s_{\calG}\!\left(\phi_{\beta,t}^L(g_2,g_1)\right)=t_{\calG}(g_1),\\
t_{\calG}\!\left(\phi_{\alpha,s}^L(g_2,g_1)\right)&=s_{\calG}(g_2),\quad s_{\calG}\!\left(\phi_{\alpha,s}^L(g_2,g_1)\right)=s_{\calG}(g_1).
\end{align*}
Thus, it makes sense to consider the element of $\calG$
\[
g_2^{-1}\phi_{\beta,t}^L(g_2,g_1)g_1.
\]
By $iv)$ of Proposition~\ref{prop-propleftdiv}, the above element may be written as
\begin{align*}
g_2^{-1}\phi_{\beta,t}^L(g_2,g_1)g_1&=\phi_P^L\!\left(g_2^{-1}\sigma_\beta(t_{\calG}(g_2)),g_1^{-1}\sigma_\beta(t_{\calG}(g_1))\right)=\\
&=\phi_P^L\!\left(\sigma_\alpha(s_{\calG}(g_2))\Theta_{\alpha\beta}\!\left(g_2^{-1}\right),\sigma_\alpha(s_{\calG}(g_1))\Theta_{\alpha\beta}\!\left(g_1^{-1}\right)\right)=\\
&=\phi_P^L\!\left(\sigma_\alpha(s_{\calG}(g_2))\Theta_{\beta\alpha}\!\left(g_2\right)^{-1},\sigma_\alpha(s_{\calG}(g_1))\Theta_{\beta\alpha}\!\left(g_1\right)^{-1}\right)=\\
&=\phi_P^L\!\left(\sigma_\alpha(s_{\calG}(g_2))\Theta_{\beta\alpha}\!\left(g_2\right)^{-1},\sigma_\alpha(s_{\calG}(g_1))\Theta_{\beta\alpha}\!\left(g_2\right)^{-1}\right)=\\
&=\phi_P^L\!\left(\sigma_\alpha(s_{\calG}(g_2)),\sigma_\alpha(s_{\calG}(g_1))\right)=\\
&=\phi_{\alpha,s}^L(g_2,g_1),
\end{align*}
where I have used compatibility between left $\calG$- and right $\calH$-action, as well as Identity (\ref{eq-invleftdiv}) of Lemma~\ref{lem-leftrightinv}.

On the other hand, assume to have two elements $g_1$ and $g_2$ of $\calG_{\alpha,\beta}$, related as in Equation~\ref{eq-moritainj}.
In particular, it must hold:
\[
\varepsilon_\beta(t_{\calG}(g_1))=\varepsilon_\beta(t_{\calG}(g_2)),\quad \varepsilon_\alpha(s_{\calG}(g_1))=\varepsilon_\alpha(s_{\calG}(g_2)).
\]
Using the same notations as above, the elements $\phi_{\beta,t}^L(g_2,g_1)g_1$ and $g_2\phi_{\alpha,s}^L(g_2,g_1)$ lie both again in $\calG_{\alpha,\beta}$, whereas $\phi_{\beta,t}^L(g_2,g_1)$, resp.\ $\phi_{\alpha,s}^L(g_2,g_1)$, lies in $\calG_{\beta}$, resp.\ $\calG_{\alpha}$.
Applying $\Theta_{\beta\alpha}^P$ to both sides of the second equality of (\ref{eq-moritainj}) and using (\ref{eq-genhomom}), one gets:
\begin{align*}
\Theta_{\beta\alpha}^P\!\left(\phi_{\beta,t}^L(g_2,g_1)g_1\right)&=\Theta_\beta^P\!\left(\phi_{\beta,t}^L(g_2,g_1)\right)\Theta_{\beta\alpha}^P(g_1)=\\
&=\Theta_{\beta\alpha}^P(g_2)\Theta_{\alpha}^P\!\left(\phi_{\alpha,s}^L(g_2,g_1)\right).
\end{align*}
Let me compute separately both $\Theta_\beta^P\!\left(\phi_{\beta,t}^L(g_2,g_1)\right)$ and $\Theta_{\alpha}^P\!\left(\phi_{\alpha,s}^L(g_2,g_1)\right)$, beginning with the former one:
\begin{align*}
\Theta_\beta^P\!\left(\phi_{\beta,t}^L(g_2,g_1)\right)&=\phi_P^R\!\left(\sigma_\beta(t_{\calG}(\phi_{\beta,t}^L(g_2,g_1))),\phi_{\beta,t}^L(g_2,g_1)\sigma_\beta(s_{\calG}(\phi_{\beta,t}^L(g_2,g_1)))\right)=\\
&=\phi_P^R\!\left(\phi_{\beta,t}^L(g_2,g_1)^{-1}\sigma_\beta(t_{\calG}(g_2)),\phi_{\beta,t}^L(g_2,g_1)\sigma_\beta(t_{\calG}(g_1))\right)=\\
&=\phi_P^R\!\left(\sigma_\beta(t_{\calG}(g_1)),\sigma_\beta(t_{\calG}(g_1))\right)=\\
&=\iota_{\calG}(t_{\calG}(g_2)),
\end{align*}
where I have used $ii)$ of Proposition~\ref{prop-propleftdiv}, Identity (\ref{eq-invrightdiv}) and Equation (\ref{eq-leftdivmap}).
Similarly, it follows
\[
\Theta_{\alpha}^P\!\left(\phi_{\alpha,s}^L(g_2,g_1)\right)=\iota_{\calG}(s_{\calG}(g_2)).
\]
Hence, it must hold:
\[
\Theta_{\beta\alpha}^P(g_1)=\Theta_{\beta\alpha}^P(g_2).
\]
\end{proof}

\subsection{The inverse of a Morita equivalence and the factorization property for division maps}\label{ssec-factprop}
Let $\calG$, $\calH$ two Lie groupoids, and $\calG\overset{P}\to\calH$ a Morita equivalence between them, in the sense of Definition~\ref{def-moritaeq}.
I define on the space $P$, viewed here as the total space of the Morita equivalence $\calG\overset{P}\to\calH$, the following left $\calH$-action and right $\calG$-action:
\begin{align}
\label{eq-morleftH}(h,p)\mapsto hp\colon=ph^{-1},\quad s_{\calH}(h)=\varepsilon(p),\\
\label{eq-morrightG}(p,g)\mapsto pg\colon=g^{-1}p,\quad t_{\calG}(g)=\pi(p).
\end{align}
It is immediate to verify that Equation (\ref{eq-morleftH}), resp.\ (\ref{eq-morrightG}), defines a left $\calH$-action on $P$ with momentum $\varepsilon$, resp.\ $\pi$.
This motivates the following
\begin{Def}\label{def-invmoreq}
Given a Morita equivalence $\left(P,\pi,\varepsilon,X_{\calG}\right)$, labelled by $\calG\overset{P}\to \calH$, its inverse $\calH\overset{P^{-1}}\to\calG$ is defined by the $4$-tuple $\left(P,\varepsilon,\pi,X_{\calH}\right)$, with left $\calH$-action defined by (\ref{eq-morleftH}) and right $\calG$-action defined by (\ref{eq-morrightG}). 
\end{Def}
First of all, I need the following technical 
\begin{Lem}\label{lem-invmoreq}
The inverse $\calH\overset{P^{-1}}\to\calG$ of a Morita equivalence $\calG\overset{P}\to\calH$ is again a Morita equivalence.
\end{Lem}
\begin{proof}
It is immediate to verify all requirements by the definition of the left $\calH$-action and the right $\calG$-action in (\ref{eq-morleftH}) and (\ref{eq-morrightG}).

For later purposes, I compute explicitly the canonical division maps associated to $P^{-1}$, which are, not surprisingly, related to the division maps associated to the Morita equivalence $P$.
In fact, by the very definition of left division map, resp.\ right division map, it must hold:
\begin{align*}
&\phi_{P^{-1}}^L\!\left(p_1,p_2\right)p_2=p_1,\quad \pi(p_1)=\pi(p_2)\Leftrightarrow p_2\phi_{P^{-1}}^L\!(p_1,p_2)^{-1}=p_1\Leftrightarrow\\
&\Leftrightarrow\phi_{P^{-1}}^L\!(p_1,p_2)^{-1}=\phi_P^R(p_2,p_1)\Leftrightarrow\phi_{P^{-1}}^L\!(p_1,p_2)=\phi_P^R(p_1,p_2);\\
&p_1 \phi_{P^{-1}}^R\!\left(p_1,p_2\right)=p_2,\quad \varepsilon(p_1)=\varepsilon(p_2)\Leftrightarrow \phi_{P^{-1}}^R\!(p_1,p_2)^{-1}p_1=p_2\Leftrightarrow\\
&\Leftrightarrow\phi_{P^{-1}}^R\!(p_1,p_2)^{-1}=\phi_P^L(p_2,p_1)\Leftrightarrow\phi_{P^{-1}}^R\!(p_1,p_2)=\phi_P^L(p_1,p_2).\\  
\end{align*}
\end{proof}
In particular, Lemma~\ref{lem-invmoreq} implies that there is a canonical generalized morphism $\calH\overset{P^{-1}}\to \calG$, for any Morita equivalence $\calG\overset{P}\to\calH$.
It thus makes sense to compute the composite generalized morphisms $P\circ P^{-1}$ and $P^{-1}\circ P$; let me compute e.g.\ $P\circ P^{-1}$; additionally, both generalized morphisms are Morita equivalences.
Recalling the construction of Subsection~\ref{ssec-gencomplaw}, the generalized morphism $P\circ P^{-1}$ has the following properties, viewed as a right principal $\calH$-bundle over $X_{\calH}$:
\begin{itemize}
\item[i)] its total space is the quotient $P\odot_\pi P/ \calG$ w.r.t.\ the right $\calG$-action defined via
\[
\left((p_1,p_2),g\right)\mapsto \left(p_1g,g^{-1}p_2\right)=\left(g^{-1}p_1,g^{-1}p_2\right),\quad t_{\calG}(g)=\pi(p_1)=\pi(p_2);
\] 
the base space is obviously $X_{\calG}$.
\item[ii)] The bundle projection $\overline{\pi}$ is 
\[
\overline{\pi}\!\left([p_1,p_2]\right)\colon=\varepsilon(p_1),
\] 
whereas the momentum map $\overline{\varepsilon}$ is
\[
\overline{\varepsilon}\!\left([p_1,p_2]\right)\colon=\varepsilon(p_2).
\] 
\item[iii)] The left $\calH$-action, resp.\ the right $\calH$-action, is defined via
\begin{align*}
h\left([p_1,p_2]\right)&=\left[p_1h^{-1},p_2\right],\quad s_{\calH}(h)=\varepsilon(p_1),\quad\text{resp.}\quad\\
\left([p_1,p_2]\right)h&=\left[p_1,p_2h\right],\quad t_{\calH}(h)=\varepsilon(p_2).
\end{align*}
\item[iv)] The right division map $\phi_{P\circ P^{-1}}^R$ of $P\circ P^{-1}$ is
\begin{align*}
&\phi_{P\circ P^{-1}}^R\!\left([p_1,p_2],\left[\overline{p}_1,\overline{p}_2\right]\right)=\phi_P^R\!\left(p_2,\phi_P^L(p_1,\overline{p}_1)\overline{p}_2\right),\\
&\pi(p_1)=\pi(p_2),\quad \pi(\overline{p}_1)=\pi(\overline{p}_2),\quad \varepsilon(p_1)=\varepsilon(\overline{p}_1).
\end{align*}
\end{itemize} 
The generalized morphism $P\circ P^{-1}$ is in fact a Morita equivalence, as follows from
\begin{Lem}\label{lem-leftdivmorita}
The generalized morphism $P\circ P^{-1}$ possesses a left division map, namely
\begin{equation}\label{eq-leftdivmor}
\begin{aligned}
&\phi_{P\circ P^{-1}}^L\!\left([p_1,p_2],\left[\overline{p}_1,\overline{p}_2\right]\right)=\phi_P^R\!\left(p_1,\phi_P^L(p_2,\overline{p}_2)\overline{p}_1\right),\\
&\pi(p_1)=\pi(p_2),\quad \pi(\overline{p}_1)=\pi(\overline{p}_2),\quad \varepsilon(p_2)=\varepsilon(\overline{p}_2).
\end{aligned}
\end{equation}
\end{Lem}
\begin{proof}
A slight modification of the arguments used in Subsection~\ref{ssec-gencomplaw} imply immediately the claim. 
\end{proof}

\begin{Prop}\label{prop-morinvert}
The $4$-tuple $\left(P\odot_\pi P/ \calG,\overline{\pi},\overline{\varepsilon},X_{\calH}\right)$ is isomorphic to the unit bundle $\calU_\calH$, viewed as a Morita equivalence of $\calH$.
\end{Prop}
\begin{proof}
By a slight modification of the arguments used to prove Proposition 4.9 of~\cite{CR}, it suffices to show the existence of a fibre-preserving, $\calH$-equivariant map from $P\circ P^{-1}$ to $U_{\calH}$, which additionally has to be left $\calH$-equivariant.
The map is defined as follows (recall that the total space of $\calU_\calH$ is simply $\calH$):
\begin{equation}\label{eq-morcompisom}
P\odot_\pi P\ni [p_1,p_2]\overset{\overline{\varphi}_{P\circ P^{-1}}}\mapsto \phi_P^R(p_1,p_2)\in \calH.
\end{equation}
Notice first that the map (\ref{eq-morcompisom}) makes sense, since any representative of $[p_1,p_2]$ belongs to $P\odot_\pi P$.
By Lemma~\ref{lem-leftrightinv}, it follows that the map (\ref{eq-morcompisom})is well-defined, i.e.\ it does not depend on the choice of the representative of the class $[p_1,p_2]$; it is clearly smooth, because the right division map $\phi_P^R$ is smooth.

It remains to check that the map (\ref{eq-morcompisom}) is fibre-preserving and left and right $\calH$-equivariant.
Point $i)$ of Proposition~\ref{prop-prodivmap} of Subsection~\ref{ssec-defdivmap} implies immediately that the map (\ref{eq-morcompisom}) is fibre-preserving and momentum-preserving (the latter proving a part of the equivariance); let me only recall that the bundle projection of $\calU_{\calH}$ is the target map of $\calH$, while the momentum is the source map of $\calH$. 
The left and right $\calH$-equivariance follows immediately from Point $iv)$ of Proposition~\ref{prop-prodivmap} of Subsection~\ref{ssec-defdivmap}, recalling Formula (\ref{eq-morleftH}).
\end{proof}

\begin{Rem}
Notice that the unit bundle $\calU_\calH$ may be interpreted as a generalized morphism from $\calH$ to itself: in fact, applying Lemma~\ref{lem-genmorph1} of Section~\ref{sec-hilskan} to the identity morphism of $\calH$, the resulting generalized morphism can be straightforwardly identified with the unit bundle $\calU_\calH$.
Moreover, it follows easily that the left $\calH$-action on $\calU_\calH$ is free and transitive on each fiber of the source map, which implies that $\calU_\calH$ is a Morita equivalence, which, in spite of Proposition~\ref{prop-morinvert} and the next Proposition~\ref{prop-morinvert2}, may be viewed as the unit w.r.t.\ the composition of generalized morphisms, since, moreover, it is not difficult to verify that left and right composition of generalized morphisms with the corresponding unit bundles are isomorphic to the initial generalized morphisms. 
\end{Rem}

Similarly, consider the composite generalized morphism $P^{-1}\circ P$, from $\calG$ to itself; without going into the details, let me only write its main properties, which can be deduced with the correct substitution from Subsection~\ref{ssec-gencomplaw}:
\begin{itemize}
\item[i)] the total space of $P^{-1}\circ P$ is the quotient space $P\odot_\varepsilon P/ \calH$, where $\calH$ acts from the right on $P\odot_\varepsilon P$ via the rule
\[
\left((p_1,p_2),h\right)\mapsto (p_1 h,h^{-1}p_2)=(p_1h,p_2h),\quad t_{\calH}(h)=\varepsilon(p_1)=\varepsilon(p_2)
\]
(diagonal action); the base space is clearly $X_{\calG}$.
\item[ii)] The bundle projection $\widetilde{\pi}$ is 
\[
\widetilde{\pi}\!\left([p_1,p_2]\right)=\pi(p_1),
\] 
and the momentum $\widetilde{\varepsilon}$ is
\[
\widetilde{\varepsilon}\!\left([p_1,p_2]\right)=\pi(p_2).
\]
\item[iii)] The left $\calG$-action on $P\circ P^{-1}$ is 
\[
\left(g,[p_1,p_2]\right)\mapsto \left[g p_1,p_2\right],\quad s_{\calG}(g)=\pi(p_1),
\] 
while the right $\calG$-action on $P\circ P^{-1}$ is 
\[
\left([p_1,p_2],g\right)\mapsto\left[p_1,g^{-1}p_2\right],\quad t_{\calG}(g)=\pi(p_2).
\]
\item[iv)] The right division map $\phi_{P^{-1}\circ P}^R$ of the bundle $P^{-1}\circ P$ is 
\begin{align*}
&\phi_{P^{-1}\circ P}^R\!\left([p_1,p_2],\left[\overline{p}_1,\overline{p}_2\right]\right)=\phi_P^L\!\left(p_2\phi_P^R(p_1,\overline{p}_1),\overline{p}_2\right),\\
&\varepsilon(p_1)=\varepsilon(p_2),\quad \varepsilon(\overline{p}_1)=\varepsilon(\overline{p}_2),\quad \pi(p_1)=\pi(\overline{p}_1).
\end{align*}
\end{itemize}
Moreover, the left $\calG$-action is free and transitive on every fiber of $\widetilde{\varepsilon}$, as stated in the following
\begin{Lem}\label{lem-leftdivmorita2}
The generalized morphism $P^{-1}\circ P$ possesses a left division map, given explicitly by the formula
\begin{equation}\label{eq-leftdivmor2}
\begin{aligned}
&\phi_{P^{-1}\circ P}\!\left(\left[\overline{p}_1,\overline{p}_2\right],[p_1,p_2]\right)=\phi_P^L\!\left(p_1,\overline{p}_1\phi_P^R(\overline{p}_2,p_2)\right),\\
&\varepsilon(p_1)=\varepsilon(p_2),\quad \varepsilon(\overline{p}_1)=\varepsilon(\overline{p}_2),\quad \pi(p_2)=\pi(\overline{p}_2).
\end{aligned}
\end{equation}
\end{Lem}

Analogously to Proposition~\ref{prop-morinvert}, the following Proposition holds
\begin{Prop}\label{prop-morinvert2}
The $4$-tuple $\left(P\odot_\varepsilon P/ \calH,\widetilde{\pi},\widetilde{\varepsilon},X_{\calG}\right)$ is isomorphic to the unit bundle $\calU_\calG$, viewed as a Morita equivalence of $\calG$ (again, the unit bundle is the canonical generalized morphism from $\calH$ to itself associated to the identity morphism of $\calG$).
\end{Prop}
\begin{proof}
The proof follows along the same lines of the proof of Proposition~\ref{prop-morinvert2}, once a fibre-preserving, left and right $\calG$-equivariant map from $P^{-1}\circ P$ to $\calG$ is defined; this turns out to be simply
\begin{equation}\label{eq-morcompisom2}
\left[p_1,p_2\right]\overset{\widetilde{\varphi}_{P^{-1}\circ P}}\mapsto\phi_P^L(p_1,p_2).
\end{equation}
\end{proof}
The following important formulae are the corollary of Proposition~\ref{prop-morinvert} and~\ref{prop-morinvert2} and Theorem 4.10 of Section 4 of~\cite{CR1}:
\begin{Thm}[Factorization formula for Morita equivalences]\label{thm-factmorita}
Let $\calG$ and $\calH$ be two Lie groupoids, and $\calG\overset{P}\to\calH$ a Morita equivalence between them; then the following formulae hold
\begin{align}
\label{eq-fact1}\phi_P^R\!\left(p_2,\phi_P^L(p_1,\overline{p}_1)\overline{p}_2\right)&=\phi_P^R(p_2,p_1)\phi_P^R(\overline{p}_1,\overline{p}_2),\\
\label{eq-fact2}\phi_P^L\!\left(p_2\phi_P^R(p_1,\overline{p}_1),\overline{p}_2\right)&=\phi_P^L(p_2,p_1)\phi_P^L(\overline{p}_1,\overline{p}_2).
\end{align}
In Identity (\ref{eq-fact1}), resp.\ (\ref{eq-fact2}), the elements $p_1$, $p_2$, $\overline{p}_1$ and $\overline{p}_2$ are chosen as follows:
\begin{align*}
&\pi(p_1)=\pi(p_2),\quad \pi(\overline{p}_1)=\pi(\overline{p}_2),\quad \varepsilon(p_1)=\varepsilon(\overline{p}_1),\quad\text{resp.}\\
&\varepsilon(p_1)=\varepsilon(\overline{p}_1),\quad \varepsilon(p_2)=\varepsilon(\overline{p}_2),\quad \pi(p_1)=\pi(\overline{p}_1).
\end{align*}
\end{Thm}
\begin{proof}
Let me prove Identity (\ref{eq-fact1}); Identity (\ref{eq-fact2}) follows by the very same arguments.

Observe first that the left-hand side of Identity (\ref{eq-fact2}) is simply the right division map of the generalized morphism $P\circ P^{-1}$, evaluated at a representative of the pair 
\[
\left([p_1,p_2],\left[\overline{p}_1,\overline{p}_2\right]\right)\in \left(P\odot_\pi P\right)\times_{\overline{\pi}}\left(P\odot_\pi P\right).
\]
On the other hand, Proposition~\ref{prop-morinvert} implies that the map (\ref{eq-morcompisom}) is a fibre-preserving, left and right $\calH$-equivariant from $P\circ P^{-1}$ to $\calU_{\calH}$, i.e.\ an isomorphism of right principal $\calH$-bundles.

Theorem 4.10 of Subsection 4.3 of~\cite{CR} implies therefore the following identity:
\[
\phi_{\calU_{\calH}}^R\circ\!\left(\overline{\varphi}_{P\circ P^{-1}}\times\overline{\varphi}_{P\circ P^{-1}}\right)=\phi_{P\circ P^{-1}}^R,
\]
whence Identity (\ref{eq-fact1}) follows immediately, evaluating both sides of the previous identity on a pair in $\left(P\odot_\pi P\right)\times_{\overline{\pi}}\left(P\odot_\pi P\right)$ and recalling the shape of the right division map of the unit bundle $\calU_\calH$.
\end{proof}

\subsection{A criterion for a generalized morphism to be a Morita equivalence}\label{ssec-critmor}
I discuss in this subsection a (global) criterion for guaranteeing that a generalized morphism in the sense of Definition~\ref{def-genmorph} is a Morita equivalence in the sense of Definition~\ref{def-moritaeq}; this criterion will play a central r{\^o}le in Subsection~\ref{ssec-loctoglobmor}.
Let me warn the reader that there is in principal nothing new in this subsection, as e.g.\ the same result may be found in~\cite{Moer4}; still, I thought it was worth writing down all the details of the proof step-by-step.
Moreover, I need this result for different purposes than the ones M{\oe}rdijk and Mr{\v c}un had in minds (and, to be really sincere, after I came to the result from my local point of view, I realized that surely a more clever mathematician than me should have proved it before, which, by the way, was exactly the case!).

Let therefore $\calG\overset{P}\to \calH$ be a generalized morphism; I denote by $\pi_1$, resp.\ $\varepsilon_1$, the projection from $P$ to $X_{\calG}$, resp.\ the momentum from $P$ to $X_{\calH}$.
\begin{Thm}\label{thm-critmor}
The generalized morphism $\calG\overset{P}\to \calH$ is a Morita equivalence between $\calG$ and $\calH$ if and only if there exists a generalized morphism $\calH\overset{Q}\to P$, such that the composite generalized morphism $Q\circ P$, resp.\ $P\circ Q$, is isomorphic to the unit bundle of $\calG$, resp.\ $\calH$, as a generalized morphism. 
\end{Thm}
The Factorization formula for Morita equivalences of Theorem~\ref{thm-factmorita} is a direct consequence of the fact that a Morita equivalence $\calG\overset{P}\to \calH$ satisfies the ``only if''-part of Theorem~\ref{thm-critmor}: namely, one can choose $Q\colon=P^{-1}$, the inverse generalized morphism of $P$ of Definition~\ref{def-invmoreq}, and the arguments of the preceding subsection show that $P^{-1}\circ P$, resp.\ $P\circ P^{-1}$, is isomorphic to the unit bundle of $\calG$, resp.\ $\calH$.
It remains therefore to prove the ``if''-part of Theorem~\ref{thm-critmor}; the proof follows from a series of technical Lemmata.
Let me first introduce and remind some notations: the projection, resp.\ momentum, of the generalized morphism $Q$ is denoted by $\pi_2$, resp.\ $\varepsilon_2$, the projection, resp.\ momentum, of $Q\circ P$ is denoted by $\overline{\pi}_1$, resp.\ $\overline{\varepsilon}_1$ (see Subsection~\ref{ssec-gencomplaw} for more details).
\begin{Lem}\label{lem-leftprincact}
Under the hypotheses of the ``if''-part of Theorem~\ref{thm-critmor}, the composite generalized morphism $Q\circ P$ is a Morita equivalence.
\end{Lem}  
\begin{proof}
The proof consists in two steps: $i)$ one has to show that the momentum of $Q\circ P$ is a surjective submersion and $ii)$ there is a right division map $\phi_{Q\circ P}^L$ on $Q\circ P$.
I postpone the proof of $i)$, which is a consequence of following technical Lemmata, showing first that there is a right division map for $Q\circ P$.

Let $\Phi_1$ be the isomorphism between $Q\circ P$ and $\calU_{\calG}$, the unit bundle of $\calG$.
Recall the construction of $Q\circ P$ from Subsection~\ref{ssec-gencomplaw}; consider two equivalence classes $[p_1,q_1]$ and $[p_2,q_2]$ in $Q\circ P$, such that
\[
\overline{\varepsilon}_1([p_1,q_1])=\varepsilon_2(q_1)=\varepsilon_2(q_2)=\overline{\varepsilon}_1([p_2,q_2]);\quad \varepsilon_1(p_1)=\pi_2(q_1),\quad \varepsilon_1(p_2)=\pi_2(q_2).
\]
I claim that the map
\[
\phi_{Q\circ P}^L\!\left([p_1,q_1],[p_2,q_2]\right)\colon=\Phi_1([p_1,q_1])\left(\Phi_1([p_2,q_2])\right)^{-1}
\]
is well-defined and that it is the left division map for $Q\circ P$.

First of all, $\phi_{Q\circ P}^L$ is well-defined: namely, recalling the properties of $\Phi_1$, one gets
\begin{align*}
s_{\calG}\!\left(\Phi_1([p_1,q_1])\right)&=\overline{\varepsilon}_1([p_1,q_1])=\\
&=\varepsilon_2(q_1)=\\
&=\varepsilon_2(q_2)=\\
&=s_{\calG}\!\left(\Phi_1([p_2,q_2]\right)=\\
&=t_{\calG}\!\left(\left(\Phi_1([p_2,q_2])\right)^{-1}\right),
\end{align*}
for any two pairs $[p_1,q_1]$ and $[p_2,q_2]$ in $Q\circ P$ as above.
It remains to show the following equation:
\[
\phi_{Q\circ P}^L\!\left([p_1,q_1],[p_2,q_2]\right)[p_2,q_2]=[p_1,q_1].
\]
This follows by applying $\Phi_1$ to the left-hand side of the previous equation; recall that $\Phi_1$ is (left) $\calG$-equivariant and is injective.
\end{proof}
The second lemma I need follows from the previous one
\begin{Lem}\label{lem-twisteqgenmor}
There exists a smooth map $\Psi_1$ from $P$ to $Q$, satisfying the following properties:
\begin{align*}
\pi_2\circ \Psi_1&=\varepsilon_1,\quad \varepsilon_2\circ\Psi_1=\pi_1;\\
\Psi_1(ph)&=h^{-1}\Psi_1(p),\quad \varepsilon_1(p)=t_{\calH}(h);\\
\Psi_1(gp)&=\Psi_1(p)g^{-1},\quad \pi_1(p)=s_{\calG}(g).
\end{align*}
\end{Lem}
\begin{proof}
Recall that the total space of the unit bundle $\calU_{\calG}$ is $\calG$.
Take $p\in P$; then there is a canonical element in $\calG$ associated to $p$, namely
\[
\iota_{\calG}(\pi_1(p)).
\]
Since $Q\circ P$ is isomorphic via $\Phi_1$ to $\calU_{\calG}$, it follows that there is a unique equivalence class $[p_1,q_1]$ in $Q\circ P$, such that
\[
\Phi_1([p_1,q_1])=\iota_{\calG}(\pi_1(p)).
\] 
A direct computation shows
\begin{align*}
t_{\calG}(\iota_{\calG}(\pi_1(p)))&=\pi_1(p)=\\
&=t_{\calG}(\Phi_1([p_1,q_1]))=\\
&=\overline{\pi}_1([p_1,q_1])=\\
&=\pi_1(p_1),
\end{align*}
i.e.\ $p_1$ lies in the same fiber as $p$, whence it follows that
\[
p_1=ph,
\]
for a unique element $h\in\calH$.
By the construction of $Q\circ P$, it follows
\[
[p_1,q_1]=[p,q],
\]
for a unique $q$ in $Q$, since the action of $\calH$ on $P$ is free.

Set $\Psi_1(p)\colon=q$, where $q$ is uniquely determined by
\[
\Phi_1([p,q])=\iota_{\calG}(\pi_1(p)).
\]
It follows immediately, by local arguments, that $\Psi_1$ is a smooth map, as
\begin{equation}\label{eq-gaugemorinv}
\left[p,\Psi_1(p)\right]=\Phi_1^{-1}\!\left(\iota_{\calG}(\pi_1(p))\right).
\end{equation}
It remains to check the properties of $\Psi_1$.
Since $\Psi_1(p)$ is uniquely determined by (\ref{eq-gaugemorinv}), it follows that
\[
\varepsilon_1(p)=\pi_2\!\left(\Psi_1(p)\right),\quad \forall p\in P.
\]
On the other hand, the computation above for showing the existence of a unique $q\in Q$, such that Equation (\ref{eq-gaugemorinv}) holds implies immediately that
\[
\pi_1(p)=\varepsilon_2\!\left(\Psi_1(p)\right),\quad \forall p\in P.
\]
Let me show now the two ``twisted'' equivariance properties; I begin with the $\calH$-twisted equivariance.
Let $p$ be an element of $P$ and $h\in \calH$, such that $\varepsilon_1(p)=t_{\calH}(h)$; then it follows
\begin{align*}
\Phi_1\!\left(\left[ph,\Psi_1(ph)\right]\right)&=\iota_{\calG}(\pi_1(ph))=\\
&=\iota_{\calG}(\pi_1(p))=\\
&=\Phi_1\!\left(\left[p,\Psi_1(p)\right]\right),
\end{align*}
by $\calH$-invariance of $\pi_1$.
The injectivity of $\Phi_1$ implies
\[
\left[ph,\Psi_1(ph)\right]=\left[p,h\Psi_1(ph)\right]=\left[p,\Psi_1(p)\right],
\]
whence the claim follows.

The $\calG$-twisted equivariance is a bit more complicated.
First of all, for $p\in P$ and $g\in \calG$ such that $s_{\calG}(g)=\pi_1(p)$, one gets
\begin{align*}
\Phi_1\!\left(\left[gp,\Psi_1(gp)\right]\right)&=\iota_{\calG}(\pi_1(gp))=\\
&=\iota_{\calG}(t_{\calG}(g))=\\
&=\Phi_1\!\left(g\left[p,\Psi_1(gp)\right]\right)=\\
&=g\Phi_1\!\left(\left[p,\Psi_1(gp)\right]\right),
\end{align*}
whence
\begin{align*}
\Phi_1\!\left(\left[p,\Psi_1(gp)\right]\right)&=g^{-1}=\\
&=\iota_{\calG}(s_{\calG}(g))g^{-1}=\\
&=\Phi_1\!\left(\left[p,\Psi_1(p)\right]\right)g^{-1}=\\
&=\Phi_1\!\left(\left[p,\Psi_1(p)g^{-1}\right]\right).
\end{align*}
Again, the injectivity of $\Phi_1$ implies
\[
\left[p,\Psi_1(gp)\right]=\left[p,\Psi_1(p)g^{-1}\right],
\]
and the claim follows immediately.
\end{proof}
Since the generalized morphism $P\circ Q$ from $\calH$ to itself is isomorphic to the unit bundle $\calU_{\calH}$, it follows that there exists a smooth map $\Psi_2$ from $Q$ to $P$, such that
\begin{align*}
\pi_1\circ\Psi_2&=\varepsilon_2,\quad \varepsilon_1\circ\Psi_2=\pi_2;\\
\Psi_2(hq)&=\Psi_2(q)h^{-1},\quad \pi_2(q)=s_{\calH}(h);\\
\Psi_2(qg)&=g^{-1}\Psi_2(q),\quad \varepsilon_2(q)=t_{\calG}(g).
\end{align*}
\begin{Lem}\label{lem-comptwist}
The composite map $\Psi_2\circ\Psi_1$, resp.\ $\Psi_1\circ\Psi_2$, is a gauge transformation of the generalized morphism $P$, resp.\ $Q$.
\end{Lem}
\begin{proof}
Consider e.g.\ the map $\Psi_2\circ\Psi_1$ from $P$ to $P$.
It follows
\begin{align*}
&\pi_1\circ\!\left(\Psi_2\circ\Psi_1\right)=\varepsilon_2\circ\Psi_1=\pi_1,\\
&\varepsilon_1\circ\!\left(\Psi_2\circ\Psi_1\right)=\pi_2\circ\Psi_1=\varepsilon_1.
\end{align*}
Moreover,
\begin{align*}
&\left(\Psi_2\circ\Psi_1\right)\!(ph)=\Psi_2\!\left(h^{-1}\Psi_1(p)\right)=\left(\Psi_2\circ\Psi_1\right)(p)h,\quad \varepsilon_1(p)=t_{\calH}(h),\\
&\left(\Psi_2\circ\Psi_1\right)\!(gp)=\Psi_2\!\left(\Psi_1(p)g^{-1}\right)=g\left(\Psi_2\circ\Psi_1\right)(p),\quad \pi_1(p)=s_{\calG}(g).
\end{align*}
\end{proof}
An important consequence of the existence of $\Psi_1$, resp.\ $\Psi_2$, is encoded in the following
\begin{Cor}\label{cor-comptwist}
The left $\calG$-, resp.\ left $\calH$-, action on $P$, resp.\ $Q$, is free.
\end{Cor}
\begin{proof}
Let me show that the left $\calH$-action on $Q$ is free.
Let $q$ be an element of $Q$ and $h\in\calH$, such that $s_{\calH}(h)=\pi_2(q)$ and
\[
hq=q.
\]
Then, applying $\Psi_2$ to both sides of the previous identity, one gets
\[
\Psi_2(hq)=\Psi_2(q)h^{-1}=\Psi_2(q).
\]
The right $\calH$-action on $P$ is free, whence
\[
h^{-1}=\iota_{\calH}\!\left(\varepsilon_1(\Psi_2(q))\right)=\iota_{\calH}(\pi_2(q)),
\]
and the claim follows.
\end{proof}
Let me now return to the proof of Theorem~\ref{thm-critmor}.
\begin{proof}[Proof of Theorem~\ref{thm-critmor}]
Also in this case, I show $i)$ that $\varepsilon_1$ is a surjective submersion and $ii)$ that there is a left division map $\phi_P^L$ for the left $\calG$-action.

That $\varepsilon_1$ is surjective, it follows from the identity $\varepsilon_1=\pi_2\circ\Psi_1$: since $\Psi_1$ is bijective by Lemma~\ref{lem-comptwist} and $\pi_2$ is surjective.
Moreover, since $\Psi_2\circ\Psi_1$, resp.\ $\Psi_1\circ\Psi_2$, is a gauge transformation of the generalized morphism $P$, resp.\ $Q$, by Lemma~\ref{lem-comptwist}, it follows from the chain rule that the tangent maps of $\Psi_1$ and $\Psi_2$ are linear isomorphisms on the corresponding tangent spaces; using again the chain rule on the identity $\varepsilon_1=\pi_2\circ\Psi_1$, together with the fact that $\pi_2$ is a submersion, it follows immediately that $\varepsilon_1$ is a submersion.

Let me now construct the left division map $\phi_P^L$ for the left $\calG$-action on $P$.
Consider thus two elements $p_1$, $p_2$ in $P$, such that $\varepsilon_1(p_1)=\varepsilon_1(p_2)$, and choose an element $q\in Q$, such that
\[
\pi_2(q)=\varepsilon_1(p_1)=\varepsilon_1(p_2).
\] 
Set
\begin{equation}\label{eq-ldivcrit}
\phi_P^L(p_1,p_2)\colon=\phi_{Q\circ P}^L\!\left([p_1,q],[p_2,q]\right),
\end{equation}
where $\phi_{Q\circ P}^L$ is the left division map for $Q\circ P$, constructed in Lemma~\ref{lem-leftprincact}.
First of all, I show that $\phi_P^L$ is well-defined, i.e.\ that it does not depend on the choice of $q$.
Namely, choosing $q_1$, such that
\[
\pi_2(q_1)=\varepsilon_1(p_1)=\varepsilon_1(p_2)=\pi_2(q),
\]
whence it follows
\[
q_1=qg,
\]
for a unique $g\in\calG$.
Then, one computes
\begin{align*}
\phi_{Q\circ P}^L\!\left([p_1,q_1],[p_2,q_1]\right)&=\Phi_1([p_1,q_1])\left(\Phi_1([p_2,q_1])\right)^{-1}=\\
&=\Phi_1([p_1,qg])\left(\Phi_1([p_2,qg])\right)^{-1}=\\
&=\Phi_1([p_1,q])g\left(\Phi_1([p_2,q]g)\right)^{-1}=\\
&=\Phi_1([p_1,q])\left(\Phi_1([p_2,q])\right)^{-1}=\\
&=\phi_{Q\circ P}^L\!\left([p_1,q],[p_2,q]\right),
\end{align*}
and this shows that $\phi_P^L$ is well-defined.
It remains to show that $\phi_P^L$ satisfies the following identity:
\[
\phi_P^L(p_1,p_2)p_2=p_1,\quad \varepsilon_1(p_1)=\varepsilon_1(p_2);
\]
Notice that $\phi_P^L(p_1,p_2)p_2$ makes sense, as
\begin{align*}
s_{\calG}\!\left(\phi_P^L(p_1,p_2)\right)&=s_{\calG}\!\left(\Phi_1([p_1,q])\left(\Phi_1([p_2,q])\right)^{-1}\right)=\\
&=t_{\calG}\!\left(\Phi_1([p_2,q])\right)=\\
&=\overline{\pi}_1([p_2,q])=\\
&=\pi_1(p_2).
\end{align*}
Now choose $q\in Q$, such that $\pi_2(q)=\varepsilon_1(p_1)=\varepsilon_1(p_2)$, then
\begin{align*}
\left[\phi_P^L(p_1,p_2)p_2,q\right]&=\left[\phi_{Q\circ P}^L\!\left([p_1,q],[p_2,q]\right)p_2,q\right]=\\
&=\phi_{Q\circ P}^L\!\left([p_1,q],[p_2,q]\right)[p_2,q]=\\
&=[p_1,q].
\end{align*}
Thus, there exist a unique $h\in \calH$, such that
\[
\phi_P^L(p_1,p_2)p_2=p_1h,\quad h^{-1}q=q.
\]
But Corollary~\ref{cor-comptwist} implies that the left $\calH$-action on $Q$ is free, and this shows
\[
\phi_P^L(p_1,p_2)p_2=p_1.
\]
\end{proof}
\begin{Rem}
Analogous arguments imply that the momentum $\varepsilon_2$ of $Q$ is also a surjective submersion; this implies, by the very construction of the generalized morphism $Q\circ P$, that its momentum is a surjective submersion, completing the proof of Lemma~\ref{lem-leftprincact}.
\end{Rem}

\subsection{Consequences of the Factorization Formula: from global to local Morita equivalences}\label{ssec-consfactform}
In this subsection, I want to deal explicitly with the local form of Morita equivalences: namely, I consider a Morita equivalence $\calG\overset{P}\to \calH$ in the sense of Definition~\ref{def-moritaeq} and the corresponding inverse Morita equivalence $\calH\overset{P^{-1}}\to \calG$ in the sense of Definition~\ref{def-invmoreq}; these two generalized morphisms give rise respectively to local generalized morphisms $\Theta_P$ and $\Theta_{P^{-1}}$, which are composable in the way explained in Subsection~\ref{ssec-compdivmap}.
I want to compute explicitly the composite local generalized morphisms $\Theta_{P^{-1}}\circ \Theta_P$ and $\Theta_P\circ \Theta_{P^{-1}}$, in order to find a local criterion for the characterization of Morita equivalences in local terms.
Theorem~\ref{thm-factmorita} plays the main r{\^o}le in these computations.

Let me compute the local generalized morphism associated to the generalized morphism $P^{-1}\circ P$, for a Morita equivalence $\calG\overset{P}\to\calH$, first w.r.t.\ to a particular choice of local sections of $P^{-1}\circ P$.
This does not correspond to the canonical choices of local sections of the composite $P^{-1}\circ P$ as specified in Subsection~\ref{ssec-compdivmap}; nonetheless, the following computations are helpful in order to understand completely the result of the computation of the composition $\Theta_{P^{-1}}\circ\Theta_P$ according to the rules of Subsection~\ref{ssec-compdivmap}.
Consider an open covering $\mathfrak{U}$ of $X_{\calG}$ and corresponding local sections $\sigma_\alpha$ of $P$ over $U_\alpha$; the associated local momenta are denoted, as usual, by $\varepsilon_\alpha$. 
Recall also that the total space of $P^{-1}\circ P$ is $P\odot_{\varepsilon}P/ \calH$, thus it makes sense to consider the following map:
\begin{equation}\label{eq-locsecmorita}
\begin{cases}
U_\alpha&\to P^{-1}\circ P\\
x&\overset{\widetilde{\sigma}_\alpha}\mapsto \left[\sigma_\alpha(x),\sigma_\alpha(x)\right].
\end{cases}
\end{equation}
It is immediate to verify that the map $\widetilde{\sigma}_\alpha$ defines a local section of $P^{-1}\circ P$; moreover, recalling the definition of the momentum of $P^{-1}\circ P$, it is immediate to see that the local momenta of $P^{-1}\circ P$ associated to the local sections defined in (\ref{eq-locsecmorita}) are the identity:
\begin{align*}
\widetilde{\varepsilon}_\alpha(x)&=\widetilde{\varepsilon}\!\left(\widetilde{\sigma}_\alpha(x)\right)=\\
&=\widetilde{\varepsilon}\!\left(\left[\sigma_\alpha(x),\sigma_\alpha(x)\right]\right)=\\
&=\pi(\sigma_\alpha(x))=\\
&=x.
\end{align*}
\begin{Rem}
Notice that the total space of the generalized morphism $P^{-1}\circ P$ is the analogon in the framework of principal bundles with groupoid structure of the {\em gauge groupoid $\calG(P)$}, introduced and discusses in Subsubsection~\ref{sssec-gaugegroup}; hence, the local sections $\widetilde{\sigma}_\alpha$ may be viewed as the unit map of the gauge groupoid of $P$.
\end{Rem}
Recalling now Equation (\ref{eq-genmor2} of Subsection~\ref{ssec-locgenmor} for local generalized morphisms obtained from generalized morphisms, one gets
\begin{align*}
\Theta^{P^{-1}\circ P}_{\beta\alpha}(g)&=\phi_{P^{-1}\circ P}\!\left(\widetilde{\sigma}_\beta\!\left(t_{\calG}(g)\right),g\widetilde{\sigma}_\alpha\!\left(s_{\calG}(g)\right)\right)=\\
&=\phi_{P^{-1}\circ P}\!\left(\left[\sigma_\beta(t_{\calG}(g)),\sigma_\beta(t_{\calG}(g))\right],\left[g\sigma_\alpha(s_{\calG}(g)),\sigma_\alpha(s_{\calG}(g))\right]\right)=\\
&=\phi_P^L\!\left(\sigma_\alpha(s_{\calG}(g)),\sigma_\beta(t_{\calG}(g))\phi_P^R\!\left(\sigma_\beta(t_{\calG}(g)),g\sigma_\alpha(s_{\calG}(g))\right)\right)=\\
&\overset{\text{Factorization property}}=\phi_P^L\!\left(\sigma_\beta(t_{\calG}(g)),\sigma_\beta(t_{\calG}(g))\right)^{-1}\phi_P^L\!\left(g\sigma_\alpha(s_{\calG}(g)),\sigma_\alpha(s_{\calG}(g))\right)=\\
&=\iota_\calG\!\left(t_{\calG}(g)\right)g\iota_{\calG}\!\left(s_{\calG}(g)\right)=\\
&=g,\quad \forall g\in \calG_{\beta,\alpha}.
\end{align*}
Hence, the local generalized morphism $\Theta_{P^{-1}\circ P}$ is the identity morphism of $\calG$, when one uses the specific sections $\widetilde{\sigma}_\alpha$ of $P^{-1}\circ P$ defined in Equation (\ref{eq-locsecmorita}).
Repeating the computations almost verbatim (with the due changes) leads to the result for the local generalized morphism $\Theta_{P\circ P^{-1}}$.

Let us now compute the local generalized morphism $\Theta_{P^{-1}\circ P}$ viewed as the composition of the local generalized morphisms $\Theta_P$ and $\Theta_{P^{-1}}$, using the results of Subsection~\ref{ssec-locgencomplaw}, from which I borrow all notations.
Let $\mathfrak{U}$, resp.\ $\mathfrak{V}$, be an open covering of $X_{\calG}$, resp.\ $X_\calH$; $\sigma_\alpha$, resp.\ $\tau_i$, denote smooth sections of $P$ (i.e.\ of $\pi$) on $U_\alpha$, resp.\ of $P^{-1}$ (i.e.\ of $\varepsilon$) on $V_i$.
$\overline{\mathfrak{U}}$ denotes the refined covering of $X_{\calG}$, constructed as in Subsection~\ref{ssec-locgencomplaw}, and, accordingly, $\sigma_{\alpha_i}$, resp.\ $\varepsilon_{\alpha_i}$, denote the restrictions of local sections of $P$, resp.\ local momenta of $P$, to $U_{\alpha_i}$.  
Consider the local sections $\widetilde{\sigma}_{\alpha_i}$, resp.\ $\widehat{\sigma}_{\alpha_i}$, of $P^{-1}\circ P$ over $U_{\alpha_i}$, where\[
\widehat{\sigma}_{\alpha_i}(x)\colon=\left[\sigma_{\alpha_i}(x),\tau_i(\varepsilon_{\alpha_i}(x))\right],
\]
and $\widetilde{\sigma}_{\alpha_i}$ are defined as in Equation (\ref{eq-locsecmorita}), only taking their restrictions to $U_{\alpha_i}$.
Both maps define local sections of $P^{-1}\circ P$, both subordinated to the refined open covering $\overline{\mathfrak{U}}$ and they are related to each other as follows (using the Factorization property):
\begin{align*}
\widehat{\sigma}_{\alpha_i}(x)&=\widetilde{\sigma}_{\alpha_i}(x)\phi_{P^{-1}\circ P}\!\left(\widetilde{\sigma}_{\alpha_i}(x),\widehat{\sigma}_{\alpha_i}(x)\right)=\\
&=\widetilde{\sigma}_{\alpha_i}(x)\phi_{P^{-1}\circ P}\!\left(\left[\sigma_{\alpha_i}(x),\sigma_{\alpha_i}(x)\right],\left[\sigma_{\alpha_i}(x),\tau_i(\varepsilon_{\alpha_i}(x))\right]\right)=\\
&=\widetilde{\sigma}_{\alpha_i}(x)\phi_P^L\!\left(\sigma_{\alpha_i}(x),\sigma_{\alpha_i}(x)\right)^{-1}\phi_P^L\!\left(\sigma_{\alpha_i}(x),\tau_i(\varepsilon_{\alpha_i}(x))\right)=\\
&=\widetilde{\sigma}_{\alpha_i}(x)\phi_P^L\!\left(\sigma_{\alpha_i}(x),\tau_i(\varepsilon_{\alpha_i}(x))\right);
\end{align*}
notice that the same result could have been found directly from the properties of $P^{-1}\circ P$.
The family of maps 
\[
U_{\alpha_i}\ni x\mapsto \phi_P^L\!\left(\sigma_{\alpha_i}(x),\tau_i(\varepsilon_{\alpha_i}(x))\right)\in \calG
\]
will be denoted by $\Phi_{\alpha_i}$ (notice that this is the common notation for transition functions); it is immediate to verify that the maps $\Phi_{\alpha_i}$ satisfy the properties
\begin{align*}
&t_{\calG}\circ \Phi_{\alpha_i}=\id,\quad s_{\calG}\circ \Phi_{\alpha_i}=\pi_i\circ \varepsilon_{\alpha_i};
&\Phi_{\beta_j}(x)=\Phi_{\alpha_i}(x)\Phi_{\alpha_i\beta_j}^{P^{-1}\circ P}(x),\quad \forall x\in U_{\alpha_i\beta_j},
\end{align*}
where by $\pi_i$ one denotes the composite function $\pi\circ\tau_i$ from $V_i$ to $P$, and provided $U_{\alpha_i}$ and $U_{\beta_j}$ intersect nontrivially; by $\Phi_{\alpha_i\beta_j}^{P^{-1}\circ P}$ I denote the transition function of $P^{-1}\circ P$.
The second identity follows immediately from the definition of $\Phi_{\alpha_i}$.

Now, I perform the explicit computation of $\Theta_{P^{-1}}\circ P$ as $\Theta_{P^{-1}}\circ \Theta_P$; notice that at some point I make use of the Factorization formula:
\begin{align*}
\Theta^{P^{-1}\circ P}_{\beta_j\alpha_i}(g)&=\Theta^{P^{-1}}_{ji}\!\left(\Theta^P_{\beta_j\alpha_i}(g)\right)=\\
&=\phi_{P^{-1}\circ P}\!\left(\widehat{\sigma}_{\beta_j}\!\left(t_{\calG}(g)\right),g\widehat{\sigma}_{\alpha_i}\!\left(s_{\calG}(g)\right)\right)=\\
&=\phi_{P^{-1}\circ P}\!\left(\left[\sigma_{\beta_j}(t_{\calG}(g)),\tau_j(\varepsilon_{\beta_j}(t_{\calG}(g)))\right],\left[g\sigma_{\alpha_i}(s_{\calG}(g)),\tau_i(\varepsilon_{\alpha_i}(s_{\calG}(g)))\right]\right)=\\
&=\phi_P^L\!\left(\sigma_{\beta_j}(t_{\calG}(g)),\tau_j(\varepsilon_{\beta_j}(t_{\calG}(g)))\right)^{-1}\phi_P^L\!\left(g\sigma_{\alpha_i}(s_{\calG}(g)),\tau_i(\varepsilon_{\alpha_i}(s_{\calG}(g)))\right)=\\
&=\phi_P^L\!\left(\sigma_{\beta_j}(t_{\calG}(g)),\tau_j(\varepsilon_{\beta_j}(t_{\calG}(g)))\right)^{-1} g \phi_P^L\!\left(\sigma_{\alpha_i}(s_{\calG}(g)),\tau_i(\varepsilon_{\alpha_i}(s_{\calG}(g)))\right)=\\
&=\Phi_{\beta_j}\!\left(t_{\calG}(g)\right)^{-1}\id_{\calG}(g)\ \Phi_{\alpha_i}\!\left(s_{\calG}(g)\right).
\end{align*}
In other words:

\fbox{\parbox{12cm}{\bf The composition of the local generalized morphisms $\Theta_{P^{-1}}$ and $\Theta_P$ in the sense explained in Subsection~\ref{ssec-locgencomplaw}, associated respectively to the generalized morphisms $P^{-1}$ and $P$, is the identity morphism of $\calG$ twisted by the transition functions $\Phi_{\alpha_i}$, in a way similar to the ``nonabelian {\v C}ech cohomological equation'' (\ref{eq-relgentrans}) of Definition~\ref{def-locgenmor} in Subsection~\ref{ssec-locgentogenmor}.}} 
A similar result can be proved, with due changes, also for the composition of $\Theta_P$ and $\Theta_{P^{-1}}$, i.e.\ $\Theta_P\circ\Theta_{P^{-1}}$ is ``cohomologous'' to the identity morphism of $\calH$ (clearly, this construction depends on the choice of open coverings of $X_{\calG}$ and $X_{\calH}$ and associated local sections and local momenta). 

\subsection{From local to global Morita equivalences}\label{ssec-loctoglobmor}
In this subsection, I give the definition of local Morita equivalences, motivated by the final results of the preceding subsection; subsequently, I show that there is a one-to-one correspondence between (global) Morita equivalences in the sense of Definition~\ref{def-moritaeq} and local Morita equivalences, as defined below.

Let first $\left(\mathfrak{U},\Phi_{\alpha\beta}^{\calG},\varepsilon_\alpha\right)$, resp.\ $\left(\mathfrak{V},\Phi_{ij}^{\calH},\pi_i\right)$, be local trivializing data over $X_{\calG}$ with values in $\calG$, resp.\ over $X_{\calH}$ with values in $\calH$; for the sake of simplicity, I use Greek indices, resp.\ Latin indices, for the labels of open sets of the open covering $\mathfrak{U}$, resp.\ $\mathfrak{V}$.
Let $\overline{\mathfrak{U}}$, resp. $\overline{\mathfrak{V}}$, be a refinement of the cover $\mathfrak{U}$, resp.\ $\mathfrak{V}$, with open sets $U_{\alpha_i}$, resp. $V_{i_\alpha}$, such that:
\begin{itemize}
\item[i)] $U_{\alpha_i}\subset U_\alpha,\quad \varepsilon_{\alpha_i}(U_{\alpha_i})\subset V_i$, where $\varepsilon_{\alpha_i}$ is the restriction of $\varepsilon_\alpha$ to $U_{\alpha_i}$; 
\item[ii)] $V_{i_\alpha}\subset V_i,\quad \pi_{i_\alpha}(V_{i_{\alpha}})\subset U_\alpha$, where $\pi_{i_\alpha}$ is the restriction of $\pi_i$ to $V_{i_\alpha}$. 
\end{itemize} 
Thus, it is possible to consider on the open covering $\overline{\mathfrak{U}}$, resp.\ $\overline{\mathfrak{V}}$, the composite local momenta $\pi_i\circ \varepsilon_{\alpha_i}$ from $U_{\alpha_i}\subset X_{\calG}$ to $X_\calG$, resp.\ $\varepsilon_\alpha\circ \pi_{\alpha_i}$ from $V_{i_\alpha}\subset X_{\calH}$ to $X_{\calH}$.
\begin{Def}\label{def-locmoreq}
A {\em local Morita equivalence $\Mu$} between the groupoids $\calG$ and $\calH$ consists of two pairs $\left(\Theta,\Phi^{\Theta}\right)$ and $\left(\Eta,\Phi^{\Eta}\right)$, where $i)$ $\Theta$, resp.\ $\Eta$, is a local generalized morphism from $\calG$ to $\calH$, resp.\ from $\calH$ to $\calG$, subordinate to the local trivializing data $\left(\mathfrak{U},\Phi_{\alpha\beta}^{\calG},\varepsilon_\alpha\right)$, resp.\ $\left(\mathfrak{V},\Phi_{ij}^{\calH},\pi_i\right)$ and $ii)$ $\Phi^\Theta=\left\{\Phi_{\alpha_i}^{\Theta}\right\}$, resp. $\Phi^{\Eta}=\left\{\Phi_{i_\alpha}^{\Eta}\right\}$, is a family of maps from $U_{\alpha_i}$ to $\calG$, resp.\ $V_{i_\alpha}$ to $\calH$, such that the following requirements are satisfied:
\begin{itemize}
\item[a)] 
\begin{equation}\label{eq-nonabcocycldef}
\begin{aligned}
t_{\calG}\circ\Phi_{\alpha_i}^{\Theta}&=\id_{X_{\calG}},\quad s_{\calG}\circ \Phi_{\alpha_i}^{\Theta}&=\pi_i\circ\varepsilon_{\alpha_i},\\ 
t_{\calH}\circ \Phi_{i_\alpha}^{\Eta}&=\id_{X_{\calH}},\quad s_{\calH}\circ \Phi_{i_\alpha}^{\Eta}&=\varepsilon_\alpha\circ \pi_{i_\alpha}.
\end{aligned}
\end{equation}
\item[b)] 
\begin{equation}\label{eq-nonabcobound}
\begin{aligned}
\Phi^{\Theta}_{\alpha_i}(x)&=\Phi_{\beta_j}^{\Theta}(x)\Phi_{\beta_j\alpha_i}^{\Eta\circ\Theta}(x),\quad \forall x\in U_{\alpha_i\beta_j},\\
\Phi^{\Eta}_{i_\alpha}(y)&=\Phi_{j_\beta}^{\Eta}(y)\Phi_{j_\beta i_\alpha}^{\Theta\circ \Eta}(y),\quad \forall y\in V_{i_\alpha j_\beta},
\end{aligned}
\end{equation}  
provided $U_{\alpha_i\beta_j}\subset X_{\calG}$, resp.\ $V_{i_\alpha j_\beta}\subset X_{\calH}$, is nontrivial.
By $\Phi_{\beta_j\alpha_i}^{\Eta\circ\Theta}$, resp.\ $\Phi_{j_\beta i_\alpha}^{\Theta\circ \Eta}$, I denoted the transition functions of $\Eta\circ\Theta$, resp.\ $\Theta\circ \Eta$, constructed by means of the ``refinement trick'' for $\Theta$ and $\Eta$, see Subsection~\ref{ssec-locgencomplaw} for more details.
\item[c)] 
\begin{equation}\label{eq-nonabcocycleq}
\begin{aligned}
(\Eta\circ\Theta)_{\beta_j\alpha_i}(g)&=\left(\Phi_{\beta_j}^{\Theta}(t_{\calG}(g))\right)^{-1}\id_{\calG}(g)\Phi_{\alpha_i}^{\Theta}(s_{\calG}(g)),\quad \forall g\in \calG_{\beta_j,\alpha_i},\\
\left(\Theta\circ\Eta\right)_{j_\beta i_\alpha}(h)&=\left(\Phi_{j_\beta}^{\Eta}(t_{\calH}(h))\right)^{-1}\id_{\calH}(h)\Phi_{i_\alpha}^{\Eta}(s_{\calH}(h)),\quad \forall h\in \calH_{j_\beta,i_\alpha},
\end{aligned}
\end{equation} 
i.e.\ the composite local generalized morphism $\Eta\circ\Theta$, resp.\ $\Theta\circ\Eta$, is ``cohomologous'' to the identity morphism of $\calG$, resp.\ $\calH$, via the ``coboundary'' $\Phi^\Theta$, resp.\ $\Phi^{\Eta}$.
\end{itemize}
\end{Def}
Let me discuss the contents of the previous definition.
First of all, one has two local generalized morphisms $\Theta$ and $\Eta$, from $\calG$ to $\calH$ and from $\calH$ to $\calG$ respectively, and subordinate to the local trivializing data $\left(\mathfrak{U},\Phi_{\alpha\beta}^{\calG},\varepsilon_\alpha\right)$, resp.\ $\left(\mathfrak{V},\Phi_{ij}^{\calH},\pi_i\right)$ respectively: by Lemma~\ref{lem-locmortogenmor}, $\Theta$, resp.\ $\Eta$, gives rise to a generalized morphism $\calG\overset{P_\Theta}\to\calH$, resp.\ $\calH\overset{P_{\Eta}}\to\calG$, in the sense of Definition~\ref{def-genmorph}.
On the other hand, the arguments of Subsection~\ref{ssec-locgencomplaw} imply that the composite local generalized morphisms $\Eta\circ\Theta$ and $\Theta\circ\Eta$ give rise to the composite generalized morphisms $\calG\overset{P_{\Eta}\circ P_\Theta}\to\calG$ and $\calH\overset{P_\Theta\circ P_{\Eta}}\to\calH$ respectively.
Recalling Definition~\ref{def-loceqgenmor}, it is not difficult to see that the family of maps $\Phi^{\Theta}$, resp.\ $\Phi^{\Eta}$, defines a local equivalence between $\Eta\circ\Theta$ and the identity morphism of $\calG$, resp.\ $\Theta\circ\Eta$ and the identity morphism of $\calH$.
Since the identity morphism of a groupoid $\calG$, viewed as a local generalized morphism subordinate to the ``trivial'' local data $\left(\mathfrak{U},\id,\iota_\calG\right)$, correspond to the unit bundle $\calU_{\calG}$, with obvious left and right $\calG$-action, Theorem~\ref{thm-loceqtoeq} implies that there is an equivalence between the generalized morphism $P_{\Eta}\circ P_\Theta$ and the unit bundle $\calU_{\calG}$, resp.\ between the generalized morphism $P_\Theta\circ P_{\Eta}$ and the unit bundle $\calU_{\calH}$.
Putting together these results with the arguments of Subsection~\ref{ssec-critmor}, in particular Theorem~\ref{thm-critmor}, the following Theorem holds
\begin{Thm}
Let $\left(\mathfrak{U},\Phi_{\alpha\beta}^{\calG},\varepsilon_\alpha\right)$, resp.\ $\left(\mathfrak{V},\Phi_{ij}^{\calH},\pi_i\right)$, be local trivializing data over $X_{\calG}$ with values in $\calG$, resp.\ over $X_{\calH}$ with values in $\calH$, let $\Theta$, resp.\ $\Eta$, a local generalized morphism from $\calG$ to $\calH$ subordinate to $\left(\mathfrak{U},\Phi_{\alpha\beta}^{\calG},\varepsilon_\alpha\right)$, resp.\ from $\calH$ to $\calG$ subordinate to $\left(\mathfrak{V},\Phi_{ij}^{\calH},\pi_i\right)$, and let $\Mu$ be a local Morita equivalence between $\Theta$ and $\Eta$ in the sense of Definition~\ref{def-locmoreq}.

Then, $\calG$ is Morita equivalent to $\calH$ w.r.t.\ the Morita equivalence $P_{\Theta}$, and $\calH$ is Morita equivalent to $\calG$ w.r.t.\ the Morita equivalence $P_{\Eta}$. 
\end{Thm}     

\subsubsection{A local Morita equivalence between the gauge groupoid $\calG(P)$ of a principal $G$-bundle $P$ and its structure group $G$}\label{sssec-locmorgaugegr}
In this subsubsection I want to compute an explicit example of a local Morita equivalence in the sense of Definition~\ref{def-locmoreq}.
For this purpose, I consider a principal $G$-bundle $P$ over a manifold $M$, and I consider the following groupoids: the gauge groupoid $\calG(P)$, for whose main properties I refer to Subsubsection~\ref{sssec-gaugegroup}, and the trivial groupoid associated to the Lie group $G$.
That these two groupoids are Morita-equivalent is already known, see e.g.~\cite{L-GTX}, Proposition 2.14; however, it is nice to prove it by a different perspective.

First of all, I construct a local generalized morphism from the gauge groupoid $\calG(P)$ to $G$.
As local trivializing data $\left(\mathfrak{U},\varepsilon_\alpha,\Phi_{\alpha\beta}\right)$ over $M$ (notice that $M$ is the manifold of objects of the gauge groupoid $\calG(P)$) with values in $G$, I choose an open covering of $M$ with local sections $\sigma_\alpha$ of $P$ as a $G$-bundle (in fact, by the arguments of Subsubsection~\ref{sssec-liegroup}, this exhausts all possible local trivializing data, up to isomorphism); the local sections, along with the division map of $P$, specify transition functions $\Phi_{\alpha\beta}$ on $U_{\alpha\beta}$ nontrivial with values in $G$.
The local component $\calG(P)_{\alpha,\beta}$, for any two open subsets $U_\alpha$, $U_\beta$ of $M$, takes the form:
\[
\calG(P)_{\alpha,\beta}\colon=\left\{[p_1,p_2]\in \calG(P)\colon \pi(p_1)\in U_\beta,\quad \pi(p_2)\in U_\alpha\right\}.
\]
Then, let me define local maps from $\calG(P)_{\alpha,\beta}$ to $G$ as follows:
\begin{equation}\label{eq-locgenmorgaugegr}
\begin{aligned}
\Theta_{\beta\alpha}\colon \calG(P)_{\alpha,\beta}&\to G\\
[p_1,p_2]&\mapsto \phi_P\!\left(\sigma_\beta(\pi(p_1)),p_1\right)\phi_P\!\left(p_2,\sigma_\alpha(\pi(p_2))\right),\quad \forall \alpha,\beta;
\end{aligned}
\end{equation}
as usual, $\phi_P$ is the division map of $P$.
As a consequence, one has the following
\begin{Lem}
The local maps (\ref{eq-locgenmorgaugegr}) are well-defined and thus give rise to a local generalized morphism $\Theta$ from $\calG(P)$ to $G$ subordinate to the trivializing data $\left(\mathfrak{U},\varepsilon_\alpha,\Phi_{\alpha\beta}\right)$.
\end{Lem}
\begin{proof}
First of all, let me show that $\Theta_{\beta\alpha}$ is well-defined, i.e.\ it does not depend on the choice of the representative of the class $[p_1,p_2]$.
In fact, choosing another representative $(\tildep_1,\tildep_2)$ of $[p_1,p_2]$, it follows that there exists a unique element $g\in G$, such that $\tildep_i=p_i g$, for $i=1,2$.
Then, it follows from the $G$-invariant of $\pi$ and of the $G\times G$-equivariance of the division map $\phi_P$ of $P$:
\begin{align*}
\Theta_{\beta\alpha}\!\left([\tildep_1,\tildep_2]\right)&=\phi_P\!\left(\sigma_\beta(\pi(\tildep_1)),\tildep_1\right)\phi_P\!\left(\tildep_2,\sigma_\alpha(\pi(\tildep_2))\right)=\\
&=\phi_P\!\left(\sigma_\beta(\pi(p_1g)),p_1g\right)\phi_P\!\left(p_2g,\sigma_\alpha(\pi(p_2g))\right)=\\
&=\phi_P\!\left(\sigma_\beta(\pi(p_1)),p_1\right)g g^{-1}\phi_P\!\left(p_2,\sigma_\alpha(\pi(p_2))\right)=\\
&=\Theta_{\beta\alpha}\!\left([p_1,p_2]\right),\quad \forall \alpha,\beta.
\end{align*} 
Let me now show that $\Theta$ is truly a local generalized morphism from $\calG(P)$ to $G$.
I have to verify the three conditions in Definition~\ref{def-locgenmor}.
Condition~\ref{eq-diagmormomen} is easily verified, since $G$ is endowed with a trivial groupoid structure.
Let me check Condition~\ref{eq-genhomom}: for this purpose, let me consider elements $[p_1,p_2]\in \calG(P)_{\alpha,\beta}$ and $[\tildep_1,\tildep_2]\in\calG_{\beta,\gamma}$, such that $\pi(p_1)=\pi(\tildep_2)$, and let me compute the following expression:
\begin{align*}
\Theta_{\gamma\alpha}\!\left([\tildep_1,\tildep_2][p_1,p_2]\right)&=\Theta_{\gamma\alpha}\!\left([\tildep_1\phi_P(\tildep_2,p_1),p_2]\right)=\\
&=\phi_P\!\left(\sigma_\gamma(\pi(\tildep_1\phi_P(\tildep_2,p_1))),\tildep_1\phi_P(\tildep_2,p_1)\right)\phi_P(p_2,\sigma_\alpha(\pi(p_2)))=\\
&=\phi_P\!\left(\sigma_\gamma(\pi(\tildep_1)),\tildep_1\right) \phi_P(\tildep_2,p_1) \phi_P(p_2,\sigma_\alpha(\pi(p_2)))=\\
&=\phi_P\!\left(\sigma_\gamma(\pi(\tildep_1)),\tildep_1\right)\phi_P(\tildep_2,\sigma_\beta(\pi(\tildep_2)))\\
&\phantom{==}\phi_P(\sigma_\beta(\pi(p_1)),p_1)\phi_P(p_2,\sigma_\alpha(\pi(p_2)))=\\
&=\Theta_{\gamma\beta}([\tildep_1,\tildep_2])\Theta_{\beta\alpha}([p_1,p_2]),\quad \forall\alpha,\beta,\gamma;
\end{align*}
notice that the equality
\[
\phi_P(\tildep_2,p_1)=\phi_P(\tildep_2,\sigma_\beta(\pi(\tildep_2)))\phi_P(\sigma_\beta(\pi(p_1)),p_1)
\]
follows immediately from the definition of the division map and from $\pi(\tildep_2)=\pi(p_1)$.
It remains to show Condition~\ref{eq-relgentrans}.
For this, notice first the identity
\[
\sigma_\beta(x)=\sigma_\alpha(x)\Phi_{\alpha\beta}(x),\quad \forall x\in U_{\alpha\beta},
\]
provided the intersection $U_{\alpha\beta}$ is nontrivial.
Then, the following holds, again using the $G\times G$-equivariance of the division map $\phi_P$:
\begin{align*}
\Theta_{\beta\alpha}([p_1,p_2])&=\phi_P(\sigma_\beta(\pi(p_1)),p_1)\phi_P(p_2,\sigma_\alpha(\pi(p_2)))=\\
&=\phi_P\!\left(\sigma_\delta(\pi(p_1))\Phi_{\delta\beta}(\pi(p_1)),p_1\right)\phi_P\!\left(p_2,\sigma_\gamma(\pi(p_2))\Phi_{\gamma\alpha}(\pi(p_2))\right)=\\
&=\Phi_{\beta\delta}(\pi(p_1))\phi_P(\sigma_\delta(\pi(p_1)),p_1)\phi_P(p_2,\sigma_\gamma(\pi(p_2)))\Phi_{\gamma\alpha}(\pi(p_2))=\\
&=\Phi_{\beta\delta}\!\left(t_{\calG(P)}([p_1,p_2])\right)\Theta_{\delta\gamma}([p_1,p_2])\Phi_{\gamma\alpha}\!\left(s_{\calG(P)}([p_1,p_2])\right),
\end{align*}
provided $U_{\beta\delta}$ and $U_{\alpha\gamma}$ are non empty.
\end{proof}
\begin{Rem}
Notice that the generalized morphism $P^\Theta$, associated to the local generalized morphism $\Theta$ of (\ref{eq-locgenmorgaugegr}), is $P$ itself; it remains to compute the explicit expression of the left $\calG(P)$-action.
Since $P$ is reconstructed from its local trivializing data, one has
\[
P=\coprod_\alpha U_\alpha\times G/ \sim,
\]
and for all notations I refer to Subsubsection~\ref{sssec-liegroup}; then, the left $\calG(P)$-action takes the form
\[
[p_1,p_2][x,g]=\left[\pi(p_1),\phi_P(\sigma_\beta(\pi(p_1)),p_1)\phi_P(p_2,\sigma_{\alpha}(\pi(p_2)))g\right],
\]
where $x\in U_{\alpha}$, $\pi(p_2)=x$ and $\beta$ is chosen so, that $\pi(p_1)$ lies in $U_\beta$.
Using the $G\times G$-equivariance of $\phi_P$ and the $G$-invariance of $\pi$, it follows
\[
[p_1,p_2][x,g]=\left[\pi(p_1),\phi_P(\sigma_\beta(\pi(p_1)),p_1)\phi_P(p_2,\sigma_{\alpha}(x)g)\right],
\]
and using the (well-defined, $G$-equivariant) isomorphism from $P$ to $\coprod_\alpha U_\alpha\times G/ \sim$
\[
[x,g]\mapsto \sigma_\alpha(x)g,\quad x\in U_\alpha,
\]
it follows immediately that the left $\calG(P)$-action takes the explicit form
\[
[p_1,p_2]p=p_1\phi_P(p_2,p),\quad \pi(p)=\pi(p_2).
\]
\end{Rem}

At this point, I need a local generalized morphism from $G$ to $\calG(P)$.
First of all, the manifold of objects of $G$ is a point $*$; therefore, an open covering for $*$ consists only of $*$ itself (with the trivial topology), and, as a consequence, one has to define only one local (in truth, global) momentum $\widetilde{\varepsilon}$ and one cocycle $\widetilde{\Phi}$:
\[
\widetilde{\varepsilon}(*)=x_0,\quad \widetilde{\Phi}(*)=\iota(x_0),
\]
where $x_0$ is a base point of $M$.
A local generalized morphism $\Eta$ from $G$ to $\calG(P)$ is then defined by
\[
g\mapsto [p_0,p_0g^{-1}],
\]
where $p_0$ is a fixed lift of $x_0$ in $P$; it is easy to verify that $\Eta$ is truly a local generalized morphism from $G$ to $\calG(P)$, by recalling the composition law in the gauge groupoid.

Now, consider any two open subsets $U_\alpha$, $U_\beta$ of $M$ in the covering $\mathfrak{U}$ of $M$, and compute the composition of $\Eta$ with $\Theta_{\beta\alpha}$ with the help of the arguments used in Subsection~\ref{ssec-locgencomplaw}:
\begin{align*}
\left(\Eta\circ\Theta_{\beta\alpha}\right)\!([p_1,p_2])&=\Eta\!\left(\phi_P(\sigma_\beta(\pi(p_1)),p_1)\phi_P(p_2,\sigma_\alpha(\pi(p_2)))\right)=\\
&=\left[p_0\phi_P(\sigma_\beta(\pi(p_1)),p_1)\phi_P(p_2,\sigma_\alpha(\pi(p_2))),p_0\right]=\\
&=\left[p_0,\sigma_\beta(\pi(p_1))\right][p_1,p_2]\left[\sigma_\alpha(\pi(p_2)),p_0\right],
\end{align*}
by the very definition of the composition law in the gauge groupoid.
On the other hand, let me fix an index $\alpha_0$, such that $x_0$ is in $U_{\alpha_0}$, and consider the corresponding local section $\sigma_{\alpha_0}$ of $P$; then, let me compute the composite morphism $\Theta_{\alpha_0}\!\circ\Eta$:
\begin{align*}
\left(\Theta_{\alpha_0}\!\circ\Eta\right)\!(g)&=\Theta_{\alpha_0}([p_0,p_0g^{-1}])=\\
&=\phi_P(\sigma_{\alpha_0}(\pi(p_0)),p_0)\phi_P(p_0g^{-1},\sigma_{\alpha_0}(\pi(p_0g^{-1})))=\\
&=\phi_P(\sigma_{\alpha_0}(\pi(p_0)),p_0)g\phi_P(p_0,\sigma_{\alpha_0}(\pi(p_0))).
\end{align*}
Therefore, let me define the following maps
\begin{align*}
\Phi^\Theta_\alpha\colon U_\alpha&\to \calG(P)\\
x&\mapsto [\sigma_\alpha(x),p_0]
\end{align*}
and
\begin{align*}
\Phi^{\Eta}\colon\left\{*\right\}&\to G\\
*&\to \phi_P(p_0,\sigma_{\alpha_0}(\pi(p_0))).
\end{align*}
it follows from the very definition of the gauge groupoid that
\[
\Phi^\Theta_\alpha(x)=\Phi^\Theta_\beta(x)[p_0,p_0\Phi_{\alpha\beta}(x)],\quad \forall x\in U_{\alpha\beta}\neq\emptyset,
\]
and it is immediate to verify that the maps
\[
U_{\alpha\beta}\ni x\mapsto [p_0,p_0\Phi_{\alpha\beta}(x)]\in\calG(P)
\]
are the composite cocycles $\Phi_{\beta\alpha}^{\Eta\circ\Theta}$ in the notations of Subsection~\ref{ssec-locgencomplaw}; on the other hand, it is obvious that the trivial cocycle
\[
*\mapsto e
\]
is the composite cocycle $\Phi^{\Theta\circ\Eta}$.
The above arguments altogether imply immediately that the pairs $\left(\Theta,\Phi^\Theta\right)$ and $\left(\Eta,\Phi^{\Eta}\right)$, constructed starting from a principal $G$-bundle $P$ over $M$, define a local Morita equivalence between the trivial groupoid $G$ and the gauge groupoid $\calG(P)$ of $P$.

\end{document}